%% file: brauer-soergel.tex
\documentclass[11pt,reqno]{amsart}
\input{macros}

\title[Cyclotomic nil-Brauer and Soergel bimodules]{Cyclotomic nil-Brauer and Singular Soergel bimodules of type D}

\author{Elijah Bodish}
\address[E.B.]{
  Department of Mathematics, Indiana University Bloomington, Bloomington, IN, USA}
\email{ebodish@iu.edu}

\author{Jonathan Brundan}
\address[J.B.]{
    Department of Mathematics,
    University of Oregon,
    Eugene, OR, USA
}
\urladdr{\href{https://pages.uoregon.edu/brundan}{pages.uoregon.edu/brundan}, \textrm{\textit{ORCiD}:} \href{https://orcid.org/0009-0009-2793-216X}{orcid.org/0009-0009-2793-216X}}
\email{brundan@uoregon.edu}

\author{Ben Elias}
\address[B.E.]{
  Department of Mathematics,
  University of Oregon,
  Eugene, OR, USA
}
\email{belias@uoregon.edu}

\thanks{E.B. is partially supported by NSF grant MSPRF-2202897.
J.B. is supported in part by NSF grant DMS-2348840. 
B.E. is supported in part by NSF grant DMS-2201387, and his research group is supported by DMS-2039316.}


\begin{document}

\begin{abstract}
We study a new family of strict monoidal categories, which are
cyclotomic quotients of the nil-Brauer category.
We construct a monoidal functor from the
cyclotomic nil-Brauer category of level $l$
to another monoidal category constructed from singular Soergel bimodules of type $\operatorname{A}_{l-1} \backslash \operatorname{D}_l / \operatorname{A}_{l-1}$.
We conjecture that our functor is an equivalence of categories.
Although we can prove neither 
fullness nor faithfulness at this point,
we are able to show that the functor induces an isomorphism
at the level of Grothendieck rings.
We compute these rings and their canonical bases, and give 
diagrammatic descriptions of the corresponding
primitive idempotents.
\end{abstract}

\maketitle
\tableofcontents
\input{s1-introduction}
\input{s2-nilbrauer}
\input{s3-ssbim}
\input{s4-typeD}
\input{s5-cyclotomic}
\bibliographystyle{alphaurl}
\bibliography{brauer-soergel}
\end{document}

%% file: macros.tex
\usepackage{
  mathabx,
  amsmath,
  amsfonts,
  amssymb,
  amsthm,
  amscd,
  graphicx,
  mathtools,
  etoolbox,
  enumitem,
  mathdots,
  txfonts,
  fancybox,
  dsfont,
  stackrel
}
\usepackage[all]{xy}
\usepackage[dvipsnames]{xcolor}
\usepackage[colorlinks=true, linkcolor=black, citecolor=black, urlcolor=blue, breaklinks=true]{hyperref}

\makeatletter
\@namedef{subjclassname@2020}{%
  \textup{2020} Mathematics Subject Classification}
\makeatother

\usepackage{tikz}
\usetikzlibrary{decorations.markings,decorations.pathmorphing,decorations.pathreplacing}
\usetikzlibrary{calc}
\usetikzlibrary{cd}
\usetikzlibrary{patterns}
\usetikzlibrary{shapes}

\tikzstyle{densely dashdotted}=[dash pattern=on 2pt off .5pt on \the\pgflinewidth off .5pt]
\tikzstyle{densely dashed}=[dash pattern=on 1.2pt off .5pt]
\tikzstyle{densely dotted}=[dash pattern=on 1pt off 1.2pt]

\tikzstyle{colp} = [ForestGreen,semithick]
\tikzstyle{colm} = [BrickRed,semithick,densely dashdotted]
\tikzstyle{twist} = [blue,densely dashed]

\colorlet{cyclocolor}{gray!20!white}
\colorlet{couponcolor}{gray!20!white}
\colorlet{pincolor}{gray!20!white}

\tikzset{anchorbase/.style={outer sep=auto,baseline={([yshift=-0.5ex]current bounding box.center)}}}
\tikzset{
    centerzero/.style={>=To,baseline={([yshift=-0.5ex](#1))}},
    centerzero/.default={0,0}
}

\def\bubble#1{\begin{tikzpicture}[anchorbase]
\node [draw,rounded corners,inner sep=2pt,black,fill=couponcolor] at (0,0) {$\scriptstyle #1$};
\end{tikzpicture}}
\def\Bubble#1{\begin{tikzpicture}[anchorbase]
\node [draw,semithick,rounded corners,inner sep=2pt,black,fill=couponcolor] at (0,0) {$\scriptstyle #1$};
\end{tikzpicture}}

\newcommand{\up}{\;\begin{tikzpicture}[anchorbase,scale=1]
\draw[thin,-to] (0,-.18) to (0,.18);
\end{tikzpicture}\;}

\newcommand{\down}{\;\begin{tikzpicture}[anchorbase,scale=1]
\draw[thin,to-] (0,-.18) to (0,.18);
\end{tikzpicture}\;}

\newcommand{\obj}{\;\begin{tikzpicture}[anchorbase,scale=1]
\draw[semithick] (0,-.18) to (0,.18);
\end{tikzpicture}\;}

\newcommand{\txtdot}{\begin{tikzpicture}[anchorbase,scale=1.2]
\draw[semithick,-] (0,-.2) to (0,.2);
\closeddot{0,0};
\end{tikzpicture}\!\!}
\newcommand{\txtcrossing}{\begin{tikzpicture}[anchorbase,scale=1.2]
\draw[semithick,-] (-.15,-.2) to (.15,.2);
\draw[semithick,-] (.15,-.2) to (-.15,.2);
\end{tikzpicture}}
\newcommand{\txtcap}{\begin{tikzpicture}[anchorbase,scale=1.2]
\draw[semithick,-] (-.15,-.15) to [out=90,in=180] (0,.15) to [out=0,in=90] (.15,-.15);
\end{tikzpicture}}
\newcommand{\txtcup}{\begin{tikzpicture}[anchorbase,scale=1.2]
\draw[semithick,-] (-.15,.15) to [out=-90,in=180] (0,-.15) to [out=0,in=-90] (.15,.15);
\end{tikzpicture}}

\newcommand{\stringlabel}[1]{\scriptstyle\color{blue}#1}

\newcommand{\cross}[2][0]{
\draw[line width=1.3pt,rotate=#1] ($#2+(-.077,-.077)$) -- ($#2+(.077,.077)$);
\draw[line width=1.3pt,rotate=#1] ($#2+(-.077,.077)$) -- ($#2+(.077,-.077)$);
}
\newcommand{\smallercross}[2][0]{
\draw[line width=1.3pt,rotate=#1] ($#2+(-.066,-.066)$) -- ($#2+(.066,.066)$);
\draw[line width=1.3pt,rotate=#1] ($#2+(-.066,.066)$) -- ($#2+(.066,-.066)$);
}
\newcommand{\pin}[3]{
\path #1 node[inner sep=0pt] (x) {$\bullet$};
\path #2 node[rectangle,rounded corners,draw,fill=pincolor,inner sep=2.5pt](y) {$\scriptstyle#3$};
\draw[Triangle Cap-,thick] (x)--(y);
}

\newcommand\closeddot[2][black]{\node at (#2) {$\color{#1}\scriptstyle\bullet$}}

\newcommand\regionlabel{\color{olive}\scriptstyle}
\def\labelp{{\regionlabel+}}
\def\labelm{{\regionlabel-}}
\def\labelpm{{\regionlabel\pm}}
\def\labeltwo{{\regionlabel\hat 2}}

\def\labeln{{\regionlabel\hat n}}

\def\LR#1#2#3{{#1\cap #2 #3 #2^{-1}}}
\def\RR#1#2#3{{#2^{-1} #1 #2 \cap #3}}
\newcommand\Z{\mathbb{Z}}
\newcommand\Q{\mathbb{Q}}
\newcommand\N{\mathbb{N}}
\def\ob{\operatorname{ob}}
\newcommand\kk{\Bbbk}
\newcommand{\qbinom}[2]{\genfrac{[}{]}{0pt}{}{#1}{#2}}
\def\op{{\operatorname{op}}}
\def\rev{{\operatorname{rev}}}
\newcommand\one{\mathds{1}}

\def\sl{\mathfrak{sl}}
\def\Par{\operatorname{Par}}
\newcommand\cC{\mathbf{C}}
\newcommand\cI{\mathbf{I}}
\newcommand\cJ{\mathbf{J}}
\def\pmod#1{{\:\left(\operatorname{mod} #1\right)}}
\def\ev{\operatorname{ev}}
\newcommand\tT{{\scriptstyle\mathrm{T}}}

\newcommand\tR{{\scriptstyle\mathrm{R}}}

\newcommand\tD{{\scriptstyle\mathrm{D}}}

\newcommand\tDl{\widetilde{\scriptstyle\mathrm{D}}}
\newcommand\gVec{\mathbf{gVec}}

\def\gMod{\operatorname{\!-gmod}}
\def\pgMod{\operatorname{\!-pgmod}}
\def\HOM{\mathbf{Hom}}
\def\END{\mathbf{End}}
\def\lround{(\!(}
\def\rround{)\!)}

\DeclareMathOperator{\End}{End}

\DeclareMathOperator{\Hom}{Hom}
\def\gKar{\operatorname{gKar}}

\DeclareMathOperator{\id}{id}

\def\dash{\operatorname{-}}

\DeclareMathOperator{\Res}{Res}
\DeclareMathOperator{\Ind}{Ind}
\DeclareMathOperator{\tr}{tr}
\def\ch{\operatorname{ch}}

\def\finsub{\subseteq_{\operatorname{f}}}
\def\finsup{\supseteq_{\operatorname{f}}}
\def\gBim{\mathfrak{gBim}}
\def\SBim{\mathfrak{SBim}}
\def\BSBim{\mathfrak{BSBim}}
\def\EBSBim{\mathfrak{eBSBim}}
\def\ESBim{\mathfrak{eSBim}}
\def\e{\mathbf{e}}
\def\r{\mathbf{r}}
\def\f{\mathbf{f}}

\def\u{\mathbf{u}}

\def\w{\mathbf{w}}

\def\u{\mathbf{u}}
\def\v{\mathbf{v}}

\newcommand{\U}{\mathrm{U}}
\newcommand{\Ui}{\mathrm{U}_t^{\imath}}
\newcommand{\Uizero}{\mathrm{U}_0^{\imath}}
\newcommand{\Uione}{\mathrm{U}_1^{\imath}}
\newcommand{\I}{\mathrm{I}}

\newcommand\cNB{\mathbf{NB}}
\newcommand\CNB{\mathbf{cNB}}
\newcommand\NB{\mathrm{NB}}
\newcommand{\V}{\mathrm{V}}

\newcommand{\Schur}{\mathrm{S}}

\newcommand{\Hecke}{\mathrm{H}}
\def\Sym{A}
\def\Center{C}
\def\B{\mathbb{B}}
\def\O{O}
\newcommand{\m}{*}

\newtheorem{theo}{Theorem}[section]
\newtheorem{lem}[theo]{Lemma}
\newtheorem{cor}[theo]{Corollary}

\theoremstyle{definition}
\newtheorem{defin}[theo]{Definition}

\newtheorem{rem}[theo]{Remark}

\numberwithin{equation}{section}

\usepackage{zref-clever}
\zcsetup{cap,nameinlink=true,
countertype = {enumi=item, enumii=item, enumiii=item, enumiv=item},
counterresetby = {enumii=enumi, enumiii=enumii, enumiv=enumiii}}
    
\zcRefTypeSetup{item}{
    Name-sg = {} ,
    name-sg = {} ,
    Name-pl = {} ,
    name-pl = {} ,
}
\zcRefTypeSetup{cor}{
    Name-sg = Corollary ,
    name-sg = corollary ,
    Name-pl = Corollaries ,
    name-pl = corollaries ,
}
\zcRefTypeSetup{defin}{
    Name-sg = Definition ,
    name-sg = definition ,
    Name-pl = Definitions ,
    name-pl = definitions ,
}
\zcRefTypeSetup{eg}{
    Name-sg = Example ,
    name-sg = example ,
    Name-pl = Examples ,
    name-pl = examples ,
}
\zcRefTypeSetup{equation}{
    Name-sg = {} ,
    name-sg = {} ,
    Name-pl = {} ,
    name-pl = {} ,
}
\zcRefTypeSetup{lem}{
    Name-sg = Lemma ,
    name-sg = lemma ,
    Name-pl = Lemmas ,
    name-pl = lemmas ,
}
\zcRefTypeSetup{rem}{
    Name-sg = Remark ,
    name-sg = remark ,
    Name-pl = Remarks ,
    name-pl = remarks ,
}
\zcRefTypeSetup{conj}{
    Name-sg = Conjecture ,
    name-sg = conjecture ,
    Name-pl = Conjectures ,
    name-pl = conjectures ,
}
\zcsetup{countertype={subsection=subsection}}
\zcRefTypeSetup{subsection}{
    Name-sg = Subsection ,
    name-sg = subsection ,
    Name-pl = Subsections ,
    name-pl = subsections ,
}
\zcRefTypeSetup{table}{
    Name-sg = Table ,
    name-sg = table ,
    Name-pl = Tables ,
    name-pl = tables ,
}
\zcRefTypeSetup{theo}{
    Name-sg = Theorem ,
    name-sg = theorem ,
    Name-pl = Theorems ,
    name-pl = theorems ,
}

\newcommand{\cref}[1]{\zcref{#1}}
\newcommand{\Cref}[1]{\zcref[S]{#1}}

\AddToHook{env/conj/begin}{\zcsetup{countertype={theo=conj}}}
\AddToHook{env/cor/begin}{\zcsetup{countertype={theo=cor}}}
\AddToHook{env/defin/begin}{\zcsetup{countertype={theo=defin}}}
\AddToHook{env/eg/begin}{\zcsetup{countertype={theo=eg}}}
\AddToHook{env/lem/begin}{\zcsetup{countertype={theo=lem}}}
\AddToHook{env/rem/begin}{\zcsetup{countertype={theo=rem}}}
\AddToHook{env/table/begin}{\zcsetup{countertype={theo=table}}}

\allowdisplaybreaks

%% file: s1-introduction.tex
\setcounter{section}{0}
\section{Introduction}\label{secintro}

Let $\kk$ be a 
field with $\operatorname{char} \kk\neq 2$.
The {\em nil-Brauer category} 
$\cNB_t$ is a graded $\kk$-linear 
monoidal category introduced in \cite{BWWnilBrauer}
(see \cref{NBdef}).
In \cite{BWWiquantum}, it was
shown to categorify the split iquantum group of rank 
one, which is the simplest 
of the coideal
subalgebras of quantized enveloping algebras
introduced by Letzter \cite{letzter1}
corresponding to symmetric pairs.
There are two admissible choices for the parameter $t$ in the definition of $\cNB_t$, $0$ or $1$, corresponding to the two possible $\Z[q,q^{-1}]$-forms, $\Uizero$ or $\Uione$, for this particular iquantum group according to the theory developed by Bao and Wang \cite{BW18QSP}. 
In this paper, we take the next step in the study of $\cNB_t$,
introducing monoidal categories $\CNB_l$ for $l \in \N$ with $l \equiv t \pmod{2}$, which are the {\em cyclotomic quotients} 
of $\cNB_t$ (see \cref{letters}).

Consider the irreducible $U_q(\sl_2)$ module of highest weight $l$, with its usual $\Z[q,q^{-1}]$-form $\V(l)$. 
Assuming that $l \equiv t \pmod{2}$,
$\V(l)$ is naturally a cyclic $\Ui$-module generated by its highest weight vector $\eta_l$, hence, $\V(l)$ is isomorphic to 
a quotient of $\Ui$. Since $\Ui$ is itself commutative,
this gives $\V(l)$ structure both as a $\Ui$-module and as a commutative $\Z[q,q^{-1}]$-algebra. The integral form
$\Ui$ also has a distinguished basis $\left\{b^{(n)}\:\big|\:n \geq 0\right\}$, called the
{\em icanonical basis}. The image of $b^{(n)}$ in $\V(l)$ is zero if $n > l$, and the non-zero images give the {\em icanonical basis} 
$\left\{b^{(n)} \eta_l\:\big|\:0 \leq n \leq l\right\}$
for $\V(l)$.

Explicit formulae expressing the icanonical basis of $\V(l)$ in terms of the standard basis 
were worked out by Berman and Wang in \cite{BeW18}.
The starting point for our work was the observation that these formulae match the Kazhdan-Lusztig combinatorics 
in a piece of the $q$-Schur algebra of type D$_l$,
associated to its maximal parabolics of type A$_{l-1}$. 
Categorifying this observation, we construct a monoidal functor from the cyclotomic nil-Brauer category $\CNB_l$ to a monoidal category 
consisting of
singular Soergel bimodules of type $\mathrm A_{l-1}\backslash \mathrm D_l /\mathrm A_{l-1}$.
This functor is a powerful tool enabling the computation of the Grothendieck ring of (the graded Karoubi envelope of) $\CNB_l$.
We use it to show that this Grothendieck ring is isomorphic to $\V(l)$ as a $\Z[q,q^{-1}]$-algebra, and match the basis arising from indecomposable objects with the icanonical basis.

Ordinary quantum groups 
are categorified by the Kac-Moody 2-categories of Khovanov, Lauda and Rouquier \cite{KL3, Rou}. In that setting,
cyclotomic quotients are certain universal 2-representations which categorify the 
integrable highest weight modules of the underlying quantum group. They
play a central role in the theory of categorical Kac-Moody actions developed in \cite{Rou}. 
The cyclotomic quotients 
$\CNB_l$
of $\cNB_t$ are expected to play a similar role in the study of categorical actions of 
the split iquantum group of rank one.
From this perspective, it may seem surprising that $\CNB_l$ is a monoidal category, rather than merely being a $\cNB_t$-module category, but this reflects the commutativity of the iquantum group in the split rank one case.

\vspace{2mm}

In order to formulate our results in more detail, we need to recall some elementary aspects 
of the theory of singular Soergel bimodules. To set conventions, here is the Dynkin diagram 
(with its diagram automorphism $\gamma$) for the Weyl group $W$ of type D$_l$, alongside the Dynkin diagram of the slightly larger Weyl group of type B$_l$:
$$
\begin{tikzpicture}[anchorbase]
\node[gray] at (-.4,0) {$\scriptstyle \gamma$};
\draw[densely dotted,gray,<->] (.08,.4) to [out=-150,in=150,looseness=1.4] (.08,-.4);
 \draw[-] (.8,0) to (0.2,-.45);
\draw[-] (.2,.45) to (0.8,0);
\draw[-] (.8,0) to (1.5,0);
\draw[-] (2.2,0) to (1.5,0);
\draw[-] (3.5,0) to (4.2,0);
\draw[dotted] (2.2,0) to (3.5,0);
\node at (.2,.65) {$\scriptstyle 1$};
\node at (.2,-.65) {$\scriptstyle -1$};
\node at (.85,.25) {$\scriptstyle 2$};
\node at (1.5,.25) {$\scriptstyle 3$};
\node at (2.2,.25) {$\scriptstyle 4$};

\node at (3.5,.25) {$\scriptstyle l-2$};
\node at (4.2,.25) {$\scriptstyle l-1$};
\closeddot{0.2,-.45};
\closeddot{0.2,.45};
\closeddot{.8,0};
\closeddot{1.5,0};
\closeddot{2.2,0};
\closeddot{3.5,0};
\closeddot{4.2,0};
\end{tikzpicture}
\hspace{35mm}
\begin{tikzpicture}[baseline=-.9mm]
 \draw[-] (.8,.03) to (.1,.03);
  \draw[-] (.8,-.03) to (.1,-.03);
\draw[-] (.8,0) to (1.5,0);
\draw[-] (2.2,0) to (1.5,0);
\draw[-] (3.5,0) to (4.2,0);
\draw[dotted] (2.2,0) to (3.5,0);
\node at (.45,0) {$\scriptstyle <$};
\node at (.1,.25) {$\scriptstyle 0$};
\node at (.8,.25) {$\scriptstyle 1$};
\node at (1.5,.25) {$\scriptstyle 2$};
\node at (2.2,.25) {$\scriptstyle 3$};
\node at (3.5,.25) {$\scriptstyle l-2$};
\node at (4.2,.25) {$\scriptstyle l-1$};
\closeddot{0.1,0};
\closeddot{.8,0};
\closeddot{1.5,0};
\closeddot{2.2,0};
\closeddot{3.5,0};
\closeddot{4.2,0};
\end{tikzpicture}
$$
We include the cases $l=2$ and $l=3$,
identifying $\operatorname{D}_2$ with $\operatorname{A}_1 \times \operatorname{A}_1$
and $\operatorname{D}_3$ with $\operatorname{A}_3$.
The Weyl groups of type D$_l$ and B$_l$ both act by automorphisms on the polynomial algebra $\kk[x_1,\dots,x_l]$,
which is the coordinate algebra
of the standard reflection representation. Explicitly, 
$s_i$ for $i=1,\dots,l-1$ switches $x_i$ and $x_{i+1}$ while fixing the other generators, $s_0$ maps $x_1 \mapsto -x_1$ while fixing the other generators, and $s_{-1} = s_0 s_1 s_0$.
We view $\kk[x_1,\dots,x_l]$ as a graded algebra so that each $x_i$ is in degree 2.

Let $\gBim$ be the graded $\kk$-linear bicategory with 
objects, 1-morphisms and 2-morphisms that are graded $\kk$-algebras, graded bimodules and graded bimodule homomorphisms,
respectively. There is a full sub-bicategory $\BSBim$ of $\gBim$
consisting of {\em singular Bott-Samelson bimodules} (see \cref{newssbimdef}). The objects of $\BSBim$ are the algebras of invariants $\kk[x_1,\dots,x_l]^{W_I}$ for all parabolic subgroups $W_I \leq W$, and its $1$-morphisms are tensor-generated by induction and restriction bimodules.
Then 
the graded $\kk$-linear
bicategory $\SBim$ 
of {\em singular Soergel bimodules} of type D
from the title
is the closure of $\BSBim$ under
direct sums, grading shifts, and summands.
By general results of Soergel and Williamson \cite{Will} recalled in \cref{homformula,class,catthm} below, the Grothendieck ring of $\SBim$ may be identified with the natural $\Z[q,q^{-1}]$-form $\Schur$ of the $q$-Schur algebra of 
type D$_l$.
In fact, we will work in this paper 
with a mildly extended version $\ESBim$ of $\SBim$
which
is obtained by incorporating additional bimodules which allow for twisting by the graph automorphism $\gamma$ (see \cref{anywayyouwant}).

Now we focus on the monoidal category 
$\END_{\ESBim}(\Sym)$ which is the endomorphism category 
of a particular object $\Sym$ 
of the bicategory $\ESBim$, namely, the algebra
$$
\Sym := \kk[x_1,\dots,x_l]^{S_l}
$$ 
of symmetric polynomials viewed as an object of $\ESBim$ 
by identifying $S_l$ with the maximal parabolic subgroup of 
$W$ obtained by deleting $\{-1\}$ from the Dynkin diagram.
Thus, $\END_{\ESBim}(\Sym)$ is a full monoidal subcategory of the graded $\kk$-linear monoidal category of 
graded $(A,A)$-bimodules.
Let $S_1 \times S_{l-1} < S_l$ denote the parabolic subgroup obtained by deleting $\{-1,1\}$ from the Dynkin diagram,
and 
$$
B := q^{1-l}
\kk[x_1,\dots,x_l]^{S_1 \times S_{l-1}}\in\ob \END_{\ESBim}(\Sym)
$$
be the graded $(\Sym,\Sym)$-bimodule
$\kk[x_1,\dots,x_l]^{S_1 \times S_{l-1}}$
with grading down-shifted so that the identity 
is in degree $(1-l)$. The right action of $\Sym$ on $B$ is the obvious one arising from the inclusion of algebras
$\kk[x_1,\dots,x_l]^{S_l} \subset \kk[x_1,\dots,x_l]^{S_1 \times S_{l-1}}$,
while the left action is defined by restriction along the {\em twisted} embedding $f \mapsto s_0(f)$ induced by the graph automorphism $\gamma$.
We note also that $\END_{\ESBim}(\Sym)$ is actually a 
$C$-linear graded monoidal category, where
$$
C := \kk[x_1^2,\dots,x_l^2]^{S_l} \subset A
$$ 
is the full algebra of invariants of the Weyl group of type B$_l$,
which is freely generated by the elementary symmetric polynomials
$e_r(x_1^2,\dots,x_l^2)\:(1 \leq r \leq l)$.

Define $t \in \{0,1\}$ so that $l \equiv t \pmod{2}$.
Let $\Lambda^{[2]}$ be the 
subalgebra of the algebra $\Lambda$ of symmetric functions over $\kk$
generated by the elements 
$$
e_r^{[2]} := \sum_{s=0}^r (-1)^s e_{r-s} e_{s},
$$
where $e_r$ is the usual $r$th elementary symmetric function.
The {\em nil-Brauer category} $\cNB_t$ is a 
strict graded $\Lambda^{[2]}$-linear monoidal category with one
generating object, also denoted $B$, and generating 
morphisms that are represented string-diagrammatically
by the dot $\txtdot\ :B \rightarrow B$, the crossing $\ \txtcrossing:B \star B \rightarrow B \star B$, the cup
$\ \txtcup:\one \rightarrow B \star B$ and the cap 
$\ \txtcap:B \star B \rightarrow \one$, subject to relations 
recorded in \cref{NBdef}.
There is an explicit $\Lambda^{[2]}$-algebra 
isomorphism $\zeta:\Lambda\stackrel{\sim}{\rightarrow}\End_{\cNB_t}(\one)$ (see \cref{whereswaldo}).
Denoting $\zeta(e_r)$ by the bubble 
$\Bubble{e_r} \in \End_{\cNB_t}(\one)$,
the {\em cyclotomic nil-Brauer category}
$\CNB_l$ is the graded monoidal category 
obtained from $\cNB_t$ by quotienting by the tensor ideal
generated by 
$\ 
\Bubble{e_l}\ 
\begin{tikzpicture}[anchorbase]
	\draw[semithick,-] (0,-.1) to (0,.3);
\end{tikzpicture}
\:+\:
\begin{tikzpicture}[anchorbase]
	\draw[semithick,-] (0,-.1) to (0,.3);
\end{tikzpicture}\ \Bubble{e_l}\;$
and
$\ \Bubble{e_r}\ (r > l)$.
%In fact, the left or right tensor ideal generated by these elements is already a two-sided tensor ideal (see \cref{bensdumbideatocallitfoursided}).
%The cyclotomic nil-Brauer category
It is naturally $C$-linear with $e_r(x_1^2,\dots,x_l^2)$
acting as $e_r^{[2]}$, this being 0 on all morphisms of 
$\CNB_l$ for $r > l$.
We write $\gKar(\CNB_l)$ for 
the graded Karoubi envelope of $\CNB_l$.

The following summarizes the main results of the article, which are established in \cref{okinawa,doors,bigtheorem,asleftmodule} in the main body of the text:

\vspace{2mm}
\noindent
{\bf Main Theorem.}
{\em
There is a $C$-linear graded monoidal functor
$\Theta:\CNB_l \rightarrow \END_{\ESBim}(\Sym)$ mapping the generating object $B$ of $\CNB_l$ to the $(\Sym,\Sym)$-bimodule $B$.
Moreover:
\begin{enumerate}
\item
The functor $\Theta$ induces an isomorphism
between the Grothendieck rings
$K_0(\gKar(\CNB_l))$ and
$K_0(\END_{\ESBim}(\Sym))$.
\item Both Grothendieck rings are isomorphic to the $\Z[q,q^{-1}]$-form $\V(l)$ of the $(l+1)$-dimensional irreducible $U_q(\sl_2)$-module equipped with an explicit algebra structure. 
\item 
Up to grading shift, isomorphism classes both 
of indecomposable objects of 
$\gKar(\CNB_l)$ and of indecomposable bimodules in $\END_{\ESBim}(\Sym)$ 
recover the icanonical basis of $\V(l)$.
\end{enumerate}
}

\vspace{2mm}
The construction of the monoidal functor $\Theta$
is naive---we simply write down bimodule homomorphisms in $\END_{\ESBim}(\Sym)$ that are the images of the generating morphisms of $\CNB_l$ then check that they satisfy the defining relations. 
The endomorphism $\Theta\big(\txtdot\  \big):B \rightarrow B$ is defined by right multiplication by $x_1$, and $\Theta\big(\ \txtcup\ )$ and $\Theta\big(\ \txtcap\ \big)$ are homomorphisms which are also quite easy to describe explicitly:
\begin{align*}
\Theta\big(\ \txtcup\ )
&:\Sym \rightarrow B \otimes_\Sym B,&
1 &\mapsto \sum_{r=0}^{l-1} e_{r}(x_2,\dots,x_l)\otimes x_1^{l-1-r},\\
\Theta\big(\ \txtcap\ )
&:B \otimes_\Sym B \rightarrow \Sym,&
f \otimes g
&\mapsto
\partial_{l-1} \cdots \partial_2 \partial_1(s_0(f)g),
\end{align*}
where $e_r(x_2,\dots,x_l)$ is the $r$th elementary symmetric polynomial in these variables and 
$\partial_i$ is the Demazure operator
$f \mapsto \frac{f - s_i(f)}{x_i-x_{i+1}}$.
However, $\Theta\big(\ \txtcrossing\ \big)$,
which is an endomorphism of the bimodule $B \otimes_{\Sym} B$, is more difficult and we are only able to write it down
using an extended version 
the diagrammatic calculus for $\BSBim$ from \cite{ESW}.
Using conventions which will be explained fully in the main body of the text, the picture for the crossing is
$$
\Theta\big(\ \txtcrossing\ \big)=
\begin{tikzpicture}[anchorbase,scale=2.2]
\draw[colp,<-] (.3,-.8) to (.3,.8);
\draw[colp,->] (-.3,-.8) to (-.3,.8);
	\draw[colm,<-] (.8,.8) to[out=-90, in=0] (0,.2) to[out = 180, in = -90] (-0.8,.8);
	\draw[colm,->] (.8,-.8) to[out=90, in=0] (0,-.2) to[out = 180, in = 90] (-0.8,-.8);
\draw (0,0) ellipse (.85 and .6);
\draw[->] (.85,0)
arc[start angle=-0.5, end angle=.5,radius=.9];
\node at (1.05,0) {$\labelp$};
\node at (.55,0) {$\labeltwo\labelp$};
\node at (-.55,0) {$\labeltwo\labelp$};
\node at (-.04,0) {$\labeltwo\labelpm$};    
\node at (-1,0) {$\labelp$};
\node at (0,-0.7) {$\labelp$};
\node at (0,0.7) {$\labelp$};
\node at (.38,0.4) {$\labeltwo$};
\node at (-.38,0.4) {$\labeltwo$};
\node at (.38,-0.4) {$\labeltwo$};
\node at (-.38,-0.4) {$\labeltwo$};
\node at (-.04,0.4) {$\labeltwo\labelp$};
\node at (-.04,-0.4) {$\labeltwo\labelp$};
\end{tikzpicture}
=\begin{tikzpicture}[anchorbase,scale=2]
\draw[colm,<-] (.3,-.8) to (.3,.8);
\draw[colm,->] (-.3,-.8) to (-.3,.8);
	\draw[colp,<-] (.8,.8) to[out=-90, in=0] (0,.2) to[out = 180, in = -90] (-0.8,.8);
	\draw[colp,->] (.8,-.8) to[out=90, in=0] (0,-.2) to[out = 180, in = 90] (-0.8,-.8);
\draw (0,0) ellipse (.85 and .6);
\draw[->] (.85,0) arc[start angle=-0.5, end angle=.5,radius=.9];
\node at (1.05,0) {$\labelp$};
\node at (.55,0) {$\labeltwo\labelp$};
\node at (-.63,0) {$\labeltwo\labelp$};
\node at (-.04,0) {$\labeltwo\labelpm$};    
\node at (-1,0) {$\labelp$};
\node at (0,-0.7) {$\labelp$};
\node at (0,0.7) {$\labelp$};
\node at (.38,0.4) {$\labeltwo$};
\node at (-.38,0.4) {$\labeltwo$};
\node at (.38,-0.4) {$\labeltwo$};
\node at (-.38,-0.4) {$\labeltwo$};
\node at (-.04,0.4) {$\labeltwo\labelp$};
\node at (-.04,-0.4) {$\labeltwo\labelp$};
\end{tikzpicture}\ .
$$
When $l=2$, the circle and all labels $\hat 2$ should be omitted from this picture.
In the statement of \cref{okinawa}, we also draw 
the much simpler pictures representing the images of the other generators.

Our Main Theorem implies that
the bimodule $B$ generates $\END_{\ESBim}(\Sym)$ in the sense that any bimodule in $\END_{\ESBim}(\Sym)$ is isomorphic to a finite direct sum of grading shifts of summands of the tensor powers 
$$
B^{\otimes n} = B \otimes_\Sym \cdots \otimes_\Sym B\qquad(0 \leq n \leq l).
$$
Up to isomorphism and grading shift, the indecomposable bimodules in $\END_{\ESBim}(\Sym)$ are parametrized by 
the set $\{0,1,\dots,l\}$. The $n$th one, denoted $B^{[n]}$ for $0 \leq n\leq l$, is uniquely determined by the property that 
it appears as a summand of
$B^{\otimes n}$ with graded multiplicity $[n]^!_q$,
and all other indecomposable summands of $B^{\otimes n}$
are isomorphic to degree shifts of $B^{[m]}$ for 
$m < n$ with $m \equiv n\pmod{2}$. 
For a formula giving the 
complete decomposition of $B^{\otimes n}$ into indecomposables, see \cref{step} below.
Since this decomposition 
is independent of the characteristic of the ground field, so too are the formal
characters of the indecomposable bimodules $B^{[n]}\:(0 \leq l \leq n)$.

We also give an explicit diagrammatic description of the 
homogeneous primitive 
idempotent $\f_n$ in $\End_{\Sym\dash\Sym}(B^{\otimes n})$ that is the projection of $B^{\otimes n}$ 
onto the unique summand that is 
$B^{[n]}$ shifted up in degree by $\binom{n}{2}$ (see \cref{bingley} and the definitions \cref{bplus,rn}). 
For example, assuming that $n$ is even with $2 \leq n \leq l-1$, the string diagram for $\f_n$ is 
$$
\begin{tikzpicture}[anchorbase,scale=4]
\node at (1.5,.7) {$\labelp$};
\node at (1.1,.7) {$\labelp$};
\node at (-.9,.7) {$\labelp$};
\node at (-.5,.7) {$\labelp$};
\node at (-0.1,.7) {$\labelp$};
\node at (1.05,.4) {$\labeln$};\node at (1.1,.4) {$\labelp$};
\node at (1.45,.4) {$\labeln$};\node at (1.5,.4) {$\labelp$};
\node at (-.95,.4) {$\labeln$};\node at (-.9,.4) {$\labelp$};
\node at (-.55,.4) {$\labeln$};\node at (-.5,.4) {$\labelp$};
\node at (-.15,.4) {$\labeln$};\node at (-0.1,.4) {$\labelp$};
\node at (1.3,.4) {$\labeln$};
\node at (1.68,.4) {$\labeln$};
\node at (-1.08,.4) {$\labeln$};
\node at (-.7,.4) {$\labeln$};
\node at (-.3,.4) {$\labeln$};
\node at (1.75,0) {$\labeln$};\node at (1.8,0) {$\labelp$};
\node at (.5,0) {$\labeln$};\node at (.55,0) {$\labelpm$};
\node at (2.1,0) {$\labelp$};
\node at (-1.25,0) {$\labeln$};\node at (-1.2,0) {$\labelp$};
\node at (-1.5,0) {$\labelp$};
	\draw[colm,->] (-1.2,.8) to[out=-90, in=180] (-.9,.2) to[out = 0, in = -90] (-0.6,.8);
	\draw[colp,->] (-.8,.8) to[out=-90, in=180] (-.5,.2) to[out = 0, in = -90] (-0.2,.8);
	\draw[colm,->] (-.4,.8) to[out=-90, in=180] (-.1,.2) to[out = 0, in = -90] (0.2,.8);
 \node at (.5,.7) {$\cdots$};
\node at (1.1,-.7) {$\labelp$};
\node at (1.5,-.7) {$\labelp$};
\node at (-.9,-.7) {$\labelp$};
\node at (-.5,-.7) {$\labelp$};
\node at (-0.1,-.7) {$\labelp$};
\node at (1.05,-.4) {$\labeln$};\node at (1.1,-.4) {$\labelp$};
\node at (1.45,-.4) {$\labeln$};\node at (1.5,-.4) {$\labelp$};
\node at (-.95,-.4) {$\labeln$};\node at (-.9,-.4) {$\labelp$};
\node at (-.55,-.4) {$\labeln$};\node at (-.5,-.4) {$\labelp$};
\node at (-.15,-.4) {$\labeln$};\node at (-0.1,-.4) {$\labelp$};
\node at (1.3,-.4) {$\labeln$};
\node at (1.68,-.4) {$\labeln$};
\node at (-1.08,-.4) {$\labeln$};
\node at (-.7,-.4) {$\labeln$};
\node at (-.3,-.4) {$\labeln$};
 \draw[colp,->] (.8,.8) to[out=-90, in=180] (1.1,.2) to[out = 0, in = -90] (1.4,.8);
	\draw[colm,->] (1.2,.8) to[out=-90, in=180]
 (1.5,.2) to[out = 0, in = -90] (1.8,.8);
	\draw[colm,<-] (-1.2,-.8) to[out=90, in=180] (-.9,-.2) to[out = 0, in = 90] (-0.6,-.8);
	\draw[colp,<-] (-.8,-.8) to[out=90, in=180] (-.5,-.2) to[out = 0, in = 90] (-0.2,-.8);
	\draw[colm,<-] (-.4,-.8) to[out=90, in=180] (-.1,-.2) to[out = 0, in = 90] (0.2,-.8);
 \node at (.5,-.7) {$\cdots$};
	\draw[colp,<-] (.8,-.8) to[out=90, in=180] (1.1,-.2) to[out = 0, in = 90] (1.4,-.8);
	\draw[colm,<-] (1.2,-.8) to[out=90, in=180] (1.5,-.2) to[out = 0, in = 90] (1.8,-.8);
 \draw[colp,<-] (1.6,-.8) to (1.6,.8);
 \draw[colp,->] (-1,-.8) to (-1,.8);
 \draw[-] (-1.4,0) to[out=90,looseness=1,in=180] (0.3,.6);
 \draw[dotted] (0.3,.6) to (.7,.6);
 \draw[-] (.7,.6) to[out=0,in=90,looseness=1] (2,0);
\draw[-] (-1.4,0) to[out=-90,in=180,looseness=1] (0.3,-.6);
 \draw[dotted] (0.3,-.6) to (.7,-.6);
\draw[->] (.7,-.6) to[out=0,in=-90,looseness=1] (2,0);
\node [draw,rounded corners,inner sep=2pt,black,fill=gray!20!white] at (1.3,0.7) {$\scriptstyle x_1$};
\node [draw,rounded corners,inner sep=1.5pt,black,fill=gray!20!white] at (-.3,0.7) {$\scriptstyle x_1^{n-3}$};
\node [draw,rounded corners,inner sep=1.5pt,black,fill=gray!20!white] at (-.7,0.7) {$\scriptstyle x_1^{n-2}$};
\node [draw,rounded corners,inner sep=1.5pt,black,fill=gray!20!white] at (-1.1,0.7) {$\scriptstyle x_1^{n-1}$};
\end{tikzpicture}
$$
In this picture, 
there are $(n-1)$ cups at the top and $(n-1)$ caps at the bottom.
There is a similar picture when $n$ is odd.
When $n=l-1$, the big bubble can be ``popped''---the penultimate $\f_n$ is also equal to
$$
\begin{tikzpicture}[anchorbase,scale=4]
\node [draw,rounded corners,inner sep=2pt,black,fill=gray!20!white] at (1.3,0.52) {$\scriptstyle x_1$};
\node [draw,rounded corners,inner sep=1.5pt,black,fill=gray!20!white] at (-.3,0.52) {$\scriptstyle x_1^{n-3}$};
\node [draw,rounded corners,inner sep=1.5pt,black,fill=gray!20!white] at (-.7,0.52) {$\scriptstyle x_1^{n-2}$};
\node [draw,rounded corners,inner sep=1.5pt,black,fill=gray!20!white] at (-1.1,0.52) {$\scriptstyle x_1^{n-1}$};
\node at (1.5,.4) {$\labelp$};
\node at (1.1,.4) {$\labelp$};
\node at (-.9,.4) {$\labelp$};
\node at (-.5,.4) {$\labelp$};
\node at (-0.1,.4) {$\labelp$};
\node at (1.9,0) {$\labelp$};
\node at (.5,0) {$\labelpm$};
\node at (-1.3,0) {$\labelp$};
	\draw[colm,->] (-1.2,.6) to[out=-90, in=180] (-.9,.1) to[out = 0, in = -90] (-0.6,.6);
	\draw[colp,->] (-.8,.6) to[out=-90, in=180] (-.5,.1) to[out = 0, in = -90] (-0.2,.6);
	\draw[colm,->] (-.4,.6) to[out=-90, in=180] (-.1,.1) to[out = 0, in = -90] (0.2,.6);
 \node at (.5,.5) {$\cdots$};
\node at (1.1,-.4) {$\labelp$};
\node at (1.5,-.4) {$\labelp$};
\node at (-.9,-.4) {$\labelp$};
\node at (-.5,-.4) {$\labelp$};
\node at (-0.1,-.4) {$\labelp$};
 \draw[colp,->] (.8,.6) to[out=-90, in=180] (1.1,.1) to[out = 0, in = -90] (1.4,.6);
	\draw[colm,->] (1.2,.6) to[out=-90, in=180]
 (1.5,.1) to[out = 0, in = -90] (1.8,.6);
	\draw[colm,<-] (-1.2,-.6) to[out=90, in=180] (-.9,-.1) to[out = 0, in = 90] (-0.6,-.6);
	\draw[colp,<-] (-.8,-.6) to[out=90, in=180] (-.5,-.1) to[out = 0, in = 90] (-0.2,-.6);
	\draw[colm,<-] (-.4,-.6) to[out=90, in=180] (-.1,-.1) to[out = 0, in = 90] (0.2,-.6);
 \node at (.5,-.5) {$\cdots$};
	\draw[colp,<-] (.8,-.6) to[out=90, in=180] (1.1,-.1) to[out = 0, in = 90] (1.4,-.6);
	\draw[colm,<-] (1.2,-.6) to[out=90, in=180] (1.5,-.1) to[out = 0, in = 90] (1.8,-.6);
 \draw[colp,<-] (1.6,-.6) to (1.6,.6);
 \draw[colp,->] (-1,-.6) to (-1,.6);
\end{tikzpicture}
$$
which has been drawn again assuming that $n$ is even (see \cref{popping}).
A similar bubble-free diagram also represents the primitive idempotent $\f_n$
in the ultimate case $n=l$.
These statements are proven using the functor $\Theta$ at the same time as establishing a parallel diagrammatic description of the corresponding
primitive idempotents $\e_n\:(0 \leq n \leq l)$ in $\CNB_l$. An essential ingredient is the description of primitive idempotents in $\cNB_t$ discovered in \cite{BWWiquantum} (see \cref{aboutidempotents}).

In $\CNB_l$ one obtains the idempotent $\e_2$ by placing a polynomial above the crossing. Similarly, one obtains $\f_2$ by placing a polynomial above $\Theta\big(\ \txtcrossing\ \big)$. Perhaps this helps explain how we determined the image of the crossing under $\Theta$ in the first place. For more clues about the idempotents $\f_n$, see \cref{rmk:factoredidempotentfn}.

% The idempotents $\f_n$ are complex, but with enough savvy one could have guessed most of their features. Slicing the diagram in half horizontally, we view $\f_n$ as a composition which factors through a particular object in $\SBim$. This object also has $B^{[n]}$ as a summand. The bottom half of $\f_n$ has degree $-\binom{n}{2}$, which is hard to achieve with any other diagram. For more on this perspective, see \cref{rmk:factoredidempotentfn}.

\vspace{2mm}
\noindent
{\bf Conjecture.}
{\em The monoidal functor $\Theta$ from the Main Theorem 
is fully faithful, hence, it induces a graded monoidal equivalence $\gKar(\CNB_l)\rightarrow\END_{\ESBim}(\Sym)$.}

\vspace{2mm}
Despite the elementary nature of these categories, we do not have a proof of this conjecture at this point.
Via the Soergel-Williamson categorification theorem, graded dimensions of morphism spaces in $\END_{\ESBim}(\Sym)$ can be computed using Lusztig's symmetric bilinear form on $\V(l)$. 
However, we are not able to compute these dimensions on the $\CNB_l$ side (although of course this would follow from the truth of our conjecture). One problem is that $\V(l)$ is not irreducible as a $\Ui$-module, so
bilinear forms on $\V(l)$ with the appropriate adjunction properties are not uniquely determined.
It seems likely that any proof of the conjecture will also produce explicit bases for the morphism
spaces $\Hom_{\Sym\dash\Sym}(\Sym,B^{[2n]})\:(0 \leq n \leq \lfloor l/2\rfloor)$ (see \cref{handy1} for the graded dimension of this space). We expect such bases will be essential for future aspirations related to categorification of ibraid group actions
of the split iquantum group of rank one.

\iffalse
The rest of the article is organized as follows.
In \cref{secnilbrauer}, we review results from \cite{BWWnilBrauer,BWWiquantum}
about the nil-Brauer category $\cNB_t$, then introduce the cyclotomic quotient $\CNB_l$.
A detailed explanation of the diagrammatic calculus for singular Soergel bimodules being used here will be given in \cref{secssbim}, which is a purely expository effort.
We introduce the extended category $\ESBim$ in type D$_l$ in \cref{secD}, and make some computations in it using an extended version of the diagrammatic calculus.
Then the main results are proved in \cref{seccyclotomic};
see Theorems \ref{okinawa}, \ref{bigtheorem}, \ref{windows}, \cref{doors} and \cref{bingley}.
The Main Theorem as formulated in this introduction is easily deduced from them.
\fi

\vspace{2mm}
\noindent
{\em General conventions.} 
By a {\em graded category},
we mean a category that is enriched in $\gVec$, the symmetric 
monoidal category of $\Z$-graded vector spaces over the ground field $\kk$.
In a graded category, we use the symbol $\cong$ to denote isomorphism of degree 0.
We use $\otimes$ to denote tensor product over the field $\kk$.
For graded vector spaces $V$ and $W$, we have by definition that
$\Hom_\kk(V,W) = \bigoplus_{n \in \Z} \Hom_\kk(V,W)_n$, with
$f \in \Hom(V,W)_n$ mapping $V_i$ into $W_{i+n}$.
The {\em upward} grading shift functor on $\gVec$ is denoted $q$, i.e., $(qV)_n = V_{n-1}$.
Assuming all $V_n$ are finite-dimensional,
its {\em graded dimension} is
\begin{equation*}
\dim_q V := \sum_{n \in \Z} q^{n} \dim V_n.
\end{equation*}
Graded rank of a free graded module over a graded algebra is defined similarly.
For a series $f = \sum_{n \in \Z} a_n q^n$ with each $a_n \in \N$,
$f V$ denotes $\bigoplus_{n \in \Z} q^n V^{\oplus a_n}$. 
An additive map between $\Z[q,q^{-1}]$-modules is
said to be {\em anti-linear}
if it twists scalars by the {\em bar 
involution}
$-:\Z[q,q^{-1}]\rightarrow \Z[q,q^{-1}], f(q) \mapsto f(q^{-1})$.
Let $[n]_q$ be the quantum integer
$\frac{q^n-q^{-n}}{q-q^{-1}}$, and $[n]^!_q$
and $\qbinom{n}{r}_q$ be the corresponding quantum factorial and quantum binomial coefficient.

\vspace{2mm}
\noindent
{\em Acknowledgements.}
The second author thanks Weiqiang Wang for helpful discussions about the underlying iSchur-Weyl duality.
All authors thank the Okinawa Institute of Science and Technology for its hospitality in June 2023, when this project was initiated.

%% file: s2-nilbrauer.tex
\setcounter{section}{1}
\section{The cyclotomic nil-Brauer category}\label{secnilbrauer}

Throughout the section, $t \in \{0,1\}$ is a fixed parameter.
We begin by defining a modified
version $\cNB_t$ of the nil-Brauer category from \cite{BWWnilBrauer}, and explaining its relationship to a $\Z[q,q^{-1}]$-form
$\Ui$ for the split iquantum group of rank one, following \cite{BWWiquantum}. Then, we pass to the quotient of $\cNB_t$ by a certain two-sided tensor ideal $\cI_l$ for $l \in \N$ with $l \equiv t \pmod{2}$,
thereby defining the cyclotomic nil-Brauer category of level $l$, denoted 
$\CNB_l$.
Finally, we discuss the combinatorics of the
icanonical basis for a $\Z[q,q^{-1}]$-form $\V(l)$ for the $(l+1)$-dimensional irreducible 
$U_q(\sl_2)$-module, which was worked out originally in \cite{BeW18}.

\subsection{Symmetric functions}\label{ssec2-1}
Let $\Lambda$ be the algebra of symmetric functions over $\kk$.
It is freely generated by either the
elementary symmetric polynomials $e_r\:(r \geq 1)$ or the complete symmetric polynomials
$h_r\:(r \geq 1)$, and $e_0 = h_0 := 1$.
We work with the grading on $\Lambda$ in which $e_r$ and $h_r$ are of degree $2r$.
As usual when working with symmetric functions, we use generating functions, setting
\begin{align}
e(u) &:= \sum_{r \geq 0} e_r u^{-r},&
h(u) &:= \sum_{r \geq 0} h_r u^{-r},
\end{align}
viewed as elements of $1+u^{-1}\Lambda \llbracket u^{-1} \rrbracket$
for a formal variable $u$.
Then we have that $e(-u) h(u) = 1$.
Let 
$\Gamma$ be the subalgebra of 
$\Lambda$ generated
by Schur's $q$-functions
\begin{equation}\label{qfunc}
q_r := \sum_{s=0}^r h_s e_{r-s}
\end{equation}
for $r \geq 0$.
We have that
\begin{align}
q(u) &:= \sum_{r \geq 0} q_r u^{-r} = e(u) h(u).
\end{align}
Hence, $q(u) q(-u) = 1$.
Recall also that
$2(-1)^{r-1} q_{2r} = q_r^2 
+ 2\sum_{s=1}^{r-1} (-1)^{r-s} q_s q_{2r-s}$
for $r \geq 1$;
cf. \cite[(III.8.2$'$)]{Mac}.
Using this, it follows that $\Gamma$
is freely generated by
$q_{2r-1}\:(r \geq 1)$.

Let $e^{[2]}(u) := e(u) e(-u)$.
Since $e^{[2]}(u) = e^{[2]}(-u)$, its expansion as a formal power series only involves even powers of $u$, so we have that
\begin{equation}\label{clunky}
e^{[2]}(u) = \sum_{r \geq 0} (-1)^r e_{r}^{[2]} u^{-2r}
\end{equation}
for some $e_r^{[2]} \in \Lambda$.
Equating coefficients gives that that $e_{r}^{[2]} = 2 (-1)^re_{2r}+ e_r^2
+ 2 \sum_{s=1}^{r-1} (-1)^{r-s} e_s e_{2r-s}$.
Let $\Lambda^{[2]}$ be the subalgebra of $\Lambda$ 
generated by $e_r^{[2]}\:(r \geq 1)$.

\begin{lem}\label{clanky}
The symmetric functions $e_{r}^{[2]}\:(r \geq 1)$ are algebraically independent.
Moreover, 
multiplication defines an algebra isomorphism
$\Gamma \otimes \Lambda^{[2]} \stackrel{\sim}{\rightarrow} \Lambda.$
\end{lem}

\begin{proof}
The first statement follows because
$e_{r}^{[2]} = 2(-1)^r e_{2r}+$(a linear combination of monomials in $e_s$ for $s < 2r$).
Also, because $h_r = (-1)^{r-1} e_r+$(a linear combination of monomials in $e_s$ for $s < r$), the definition \cref{qfunc} implies that
$q_{2r-1} = 2 e_{2r-1}+$(a linear combination of monomials in $e_s$ for $s < 2r-1$).
Everything else now follows because $\Lambda$ is
freely generated by $e_r\:(r \geq 1)$.
\end{proof}

\begin{rem}\label{simonriche}
The symmetric functions $e_{r}^{[2]}$ are natural to consider from the point of view of symmetric polynomials.
To explain what we mean, let 
$\Sym$ denote the algebra of symmetric polynomials $\kk[x_1,\dots,x_l]^{S_l}$ for some fixed $l \geq 0$, and
\begin{equation}\label{maggiesmith}
\ev_l:\Lambda \twoheadrightarrow \Sym
\end{equation}
be the homomorphism sending $e_r$ to the $r$th
elementary symmetric polynomial $e_r(x_1,\dots,x_l)$ (which is 0 for $r > l$).
Then we have that 
\begin{equation}
\ev_l(e_{r}^{[2]}) = e_{r}(x_1^2,\dots,x_l^2).
\end{equation}
To prove this, note that
$\ev_l(e(u)) = u^{-l}(u+x_1)\cdots(u+x_l)$
and $\ev_l(e(-u)) = u^{-l}(u-x_1)\cdots (u-x_1)$.
Hence,
$\ev_l\big(e^{[2]}(u)\big)=\ev_l(e(u) e(-u)) = u^{-2l} (u^2-x_1^2) \cdots (u^2 - x_l^2).$
Now equate $u^{-2r}$-coefficients.
\end{rem}

\subsection{The nil-Brauer category}
In \cite[Def.~2.1]{BWWnilBrauer}, the nil-Brauer category was defined over the ground field $\kk$. In this article, we will 
work with a modified version which is 
the $\Lambda^{[2]}$-linear monoidal category obtained by extending scalars to the ground ring $\Lambda^{[2]}$.
To explain the definition, we use the string calculus for monoidal categories, adopting 
the standard conventions from most of the recent categorification literature, as in \cite{BWWnilBrauer}.

\begin{defin}\label{NBdef}
Let $\Lambda^{[2]}$ be the subalgebra of $\Lambda$ defined in \cref{clanky}. Let $t \in \{0,1\}$. 
The 
{\em nil-Brauer category} is the strict graded $\Lambda^{[2]}$-linear monoidal category $(\cNB_t,-\star-,\one)$
with one generating object $B$, whose identity endomorphism $\id_B$ is denoted by the string $\obj$, and four generating morphisms
\begin{align*}
\begin{tikzpicture}[anchorbase]
	\draw[semithick] (0.08,-.3) to (0.08,.4);
    \closeddot{0.08,.05};
\end{tikzpicture}
&:B\rightarrow B,&
\begin{tikzpicture}[anchorbase]
\draw[semithick] (0.28,-.3) to (-0.28,.4);
	\draw[semithick] (-0.28,-.3) to (0.28,.4);
\end{tikzpicture}
&:B\star B \rightarrow B \star B,
&
\begin{tikzpicture}[anchorbase]
	\draw[semithick] (0.4,0) to[out=90, in=0] (0.1,0.4);
	\draw[semithick] (0.1,0.4) to[out = 180, in = 90] (-0.2,0);
\end{tikzpicture}
&:B \star B\rightarrow\one,&
\begin{tikzpicture}[anchorbase]
	\draw[semithick] (0.4,0.4) to[out=-90, in=0] (0.1,0);
	\draw[semithick] (0.1,0) to[out = 180, in = -90] (-0.2,0.4);
\end{tikzpicture}
&:\one\rightarrow B\star B,
\\\notag
&\text{(degree $2$)}&
&\text{(degree $-2$)}&
&\text{(degree $0$)}&
&\text{(degree $0$)}
\end{align*}
subject to the following relations:
\begin{align}
\label{rels2}
\begin{tikzpicture}[anchorbase]
  \draw[semithick] (0.3,0) to (0.3,.4);
	\draw[semithick] (0.3,0) to[out=-90, in=0] (0.1,-0.2);
	\draw[semithick] (0.1,-0.2) to[out = 180, in = -90] (-0.1,0);
	\draw[semithick] (-0.1,0) to[out=90, in=0] (-0.3,0.2);
	\draw[semithick] (-0.3,0.2) to[out = 180, in =90] (-0.5,0);
  \draw[semithick] (-0.5,0) to (-0.5,-.4);
\end{tikzpicture}
&=
\begin{tikzpicture}[anchorbase]
  \draw[semithick] (0,-0.4) to (0,.4);
\end{tikzpicture}
=\begin{tikzpicture}[anchorbase]
  \draw[semithick] (0.3,0) to (0.3,-.4);
	\draw[semithick] (0.3,0) to[out=90, in=0] (0.1,0.2);
	\draw[semithick] (0.1,0.2) to[out = 180, in = 90] (-0.1,0);
	\draw[semithick] (-0.1,0) to[out=-90, in=0] (-0.3,-0.2);
	\draw[semithick] (-0.3,-0.2) to[out = 180, in =-90] (-0.5,0);
  \draw[semithick] (-0.5,0) to (-0.5,.4);
\end{tikzpicture}
\:,&
\begin{tikzpicture}[baseline=-2.5mm]
\draw[semithick] (0,-.15) circle (.3);
\end{tikzpicture}
&= t\id_\one,
\\\label{rels3}
\begin{tikzpicture}[anchorbase,scale=1.1]
	\draw[semithick] (0.5,0) to[out=90, in=0] (0.1,0.5);
	\draw[semithick] (0.1,0.5) to[out = 180, in = 90] (-0.3,0);
 \draw[semithick] (0.1,0) to[out=90,in=-90] (-.3,.7);
\end{tikzpicture}
&=
\begin{tikzpicture}[anchorbase,scale=1.1]
	\draw[semithick] (0.5,0) to[out=90, in=0] (0.1,0.5);
	\draw[semithick] (0.1,0.5) to[out = 180, in = 90] (-0.3,0);
 \draw[semithick] (0.1,0) to[out=90,in=-90] (.5,.7);
\end{tikzpicture}\:,&
\begin{tikzpicture}[anchorbase,scale=1.1]
	\draw[semithick] (0.35,.3)  to [out=90,in=-90] (-.1,.9) to [out=90,in=180] (.1,1.1);
 \draw[semithick] (-.15,.3)  to [out=90,in=-90] (.3,.9) to [out=90,in=0] (.1,1.1);
\end{tikzpicture}
&=0
\:,\\
\label{rels4}
\begin{tikzpicture}[anchorbase,scale=1.1]
	\draw[semithick] (0.4,0) to[out=90, in=0] (0.1,0.5);
	\draw[semithick] (0.1,0.5) to[out = 180, in = 90] (-0.2,0);
 \closeddot{.38,.2};
\end{tikzpicture}
&=
-\begin{tikzpicture}[anchorbase,scale=1.1]
	\draw[semithick] (0.4,0) to[out=90, in=0] (0.1,0.5);
	\draw[semithick] (0.1,0.5) to[out = 180, in = 90] (-0.2,0);
 \closeddot{-.18,.2};
\end{tikzpicture}
\:,&
\begin{tikzpicture}[anchorbase,scale=1.4]
	\draw[semithick] (0.25,.3) to (-0.25,-.3);
	\draw[semithick] (0.25,-.3) to (-0.25,.3);
 \closeddot{-0.12,0.145};
\end{tikzpicture}
-
\begin{tikzpicture}[anchorbase,scale=1.4]
	\draw[semithick] (0.25,.3) to (-0.25,-.3);
	\draw[semithick] (0.25,-.3) to (-0.25,.3);
     \closeddot{.12,-0.145};
\end{tikzpicture}
&=
\begin{tikzpicture}[anchorbase,scale=1.4]
 	\draw[semithick] (0,-.3) to (0,.3);
	\draw[semithick] (-.3,-.3) to (-0.3,.3);
\end{tikzpicture}
-
\begin{tikzpicture}[anchorbase,scale=1.4]
 	\draw[semithick] (-0.15,-.3) to[out=90,in=180] (0,-.05) to[out=0,in=90] (0.15,-.3);
 	\draw[semithick] (-0.15,.3) to[out=-90,in=180] (0,.05) to[out=0,in=-90] (0.15,.3);
\end{tikzpicture}
\:,
\\
\label{rels1}
\begin{tikzpicture}[anchorbase]
	\draw[semithick] (0.28,0) to[out=90,in=-90] (-0.28,.6);
	\draw[semithick] (-0.28,0) to[out=90,in=-90] (0.28,.6);
	\draw[semithick] (0.28,-.6) to[out=90,in=-90] (-0.28,0);
	\draw[semithick] (-0.28,-.6) to[out=90,in=-90] (0.28,0);
\end{tikzpicture}
&=0,
&\begin{tikzpicture}[anchorbase]
	\draw[semithick] (0.45,.6) to (-0.45,-.6);
	\draw[semithick] (0.45,-.6) to (-0.45,.6);
        \draw[semithick] (0,-.6) to[out=90,in=-90] (-.45,0);
        \draw[semithick] (-0.45,0) to[out=90,in=-90] (0,0.6);
\end{tikzpicture}
&=
\begin{tikzpicture}[anchorbase]
	\draw[semithick] (0.45,.6) to (-0.45,-.6);
	\draw[semithick] (0.45,-.6) to (-0.45,.6);
        \draw[semithick] (0,-.6) to[out=90,in=-90] (.45,0);
        \draw[semithick] (0.45,0) to[out=90,in=-90] (0,0.6);
\end{tikzpicture}\:.
\end{align}
\end{defin}

We will denote the $r$th power of $\txtdot$ by labeling the dot with $r$.
Using the zig-zag relations from \cref{rels2}, the following are easily deduced from \cref{rels3,rels4}:
\begin{align}\label{rels3b}
\begin{tikzpicture}[anchorbase,scale=1.1]
	\draw[-,semithick] (0.5,0) to[out=-90, in=0] (0.1,-0.5);
	\draw[-,semithick] (0.1,-0.5) to[out = 180, in = -90] (-0.3,0);
 \draw[-,semithick] (0.1,0) to[out=-90,in=90] (-.3,-.7);
\end{tikzpicture}
&=
\begin{tikzpicture}[anchorbase,scale=1.1]
	\draw[-,semithick] (0.5,0) to[out=-90, in=0] (0.1,-0.5);
	\draw[-,semithick] (0.1,-0.5) to[out = 180, in = -90] (-0.3,0);
 \draw[-,semithick] (0.1,0) to[out=-90,in=90] (.5,-.7);
\end{tikzpicture}\:,&
\begin{tikzpicture}[anchorbase,scale=1.1]
	\draw[-,semithick] (0.35,-.3)  to [out=-90,in=90] (-.1,-.9) to [out=-90,in=180] (.1,-1.1);
 \draw[-,semithick] (-.15,-.3)  to [out=-90,in=90] (.3,-.9) to [out=-90,in=0] (.1,-1.1);
\end{tikzpicture}
&=0,\\
\label{rels4b}
\begin{tikzpicture}[anchorbase,scale=1.1]
	\draw[-,semithick] (0.4,0) to[out=-90, in=0] (0.1,-0.5);
	\draw[-,semithick] (0.1,-0.5) to[out = 180, in = -90] (-0.2,0);
 \closeddot{.38,-.2};
\end{tikzpicture}
&=
-\begin{tikzpicture}[anchorbase,scale=1.1]
	\draw[-,semithick] (0.4,0) to[out=-90, in=0] (0.1,-0.5);
	\draw[-,semithick] (0.1,-0.5) to[out = 180, in = -90] (-0.2,0);
 \closeddot{-.18,-.2};
\end{tikzpicture}\:,&
\begin{tikzpicture}[anchorbase,scale=1.4]
	\draw[-,semithick] (0.25,.3) to (-0.25,-.3);
	\draw[-,semithick] (0.25,-.3) to (-0.25,.3);
 \closeddot{-0.12,-0.145};
\end{tikzpicture}
-
\begin{tikzpicture}[anchorbase,scale=1.4]
	\draw[-,semithick] (0.25,.3) to (-0.25,-.3);
	\draw[-,semithick] (0.25,-.3) to (-0.25,.3);
     \closeddot{.12,0.145};
\end{tikzpicture}
&=
\begin{tikzpicture}[anchorbase,scale=1.4]
 	\draw[-,semithick] (0,-.3) to (0,.3);
	\draw[-,semithick] (-.3,-.3) to (-0.3,.3);
\end{tikzpicture}
-
\begin{tikzpicture}[anchorbase,scale=1.4]
 	\draw[-,semithick] (-0.15,-.3) to[out=90,in=180] (0,-.05) to[out=0,in=90] (0.15,-.3);
 	\draw[-,semithick] (-0.15,.3) to[out=-90,in=180] (0,.05) to[out=0,in=-90] (0.15,.3);
\end{tikzpicture}
\:.
\end{align}
Now it follows that there are strict graded $\Lambda^{[2]}$-linear monoidal functors
\begin{align}
\tR:\cNB_t &\rightarrow \cNB_t^{\rev}, & s &\mapsto (-1)^{\bullet(s)} s^{\leftrightarrow},
\label{R}\\\label{T}
\tT:\cNB_t &\rightarrow \cNB_t^{\op}, &
s &\mapsto
s^{\updownarrow}.
\end{align}
Here, the $\op$ denotes the opposite category with the same monoidal product, and $\rev$ denotes the same category with the reverse monoidal product.
Also, for a string diagram $s$ 
we are using $s^{\updownarrow}$ 
and $s^{\leftrightarrow}$ to denote
its reflection in a horizontal or vertical axis,
and $\bullet(s)$ denotes the total number of dots
in the diagram.
It follows that the category $\cNB_t$ has a strict pivotal structure. The underlying duality functor is $\tR \circ \tT = \tT \circ\tR$, which rotates a string diagram $s$ through $180^\circ$ then scales by $(-1)^{\bullet(s)}$. 
For more details about all of this, see \cite[Sec.~2]{BWWnilBrauer}.

\begin{lem} \label{AhaJonYouForgotToComeUpWithARidiculousNameForThisLemma}
Using the defining relations \cref{rels2,rels3,rels4} but neither of the relations \cref{rels1}, it follows that
\begin{align}\label{puppy1}
\begin{tikzpicture}[scale=.8,centerzero={0,0}]
	\draw[semithick] (0.28,0) to[out=90,in=-90] (-0.28,.6);
	\draw[semithick] (-0.28,0) to[out=90,in=-90] (0.28,.6);
	\draw[semithick] (0.28,-.6) to[out=90,in=-90] (-0.28,0);
	\draw[semithick] (-0.28,-.6) to[out=90,in=-90] (0.28,0);
 \closeddot{-.22,-.45};
\end{tikzpicture}
&=
\begin{tikzpicture}[scale=.8,centerzero={0,0}]
 \closeddot{-.22,.45};
	\draw[semithick] (0.28,0) to[out=90,in=-90] (-0.28,.6);
	\draw[semithick] (-0.28,0) to[out=90,in=-90] (0.28,.6);
	\draw[semithick] (0.28,-.6) to[out=90,in=-90] (-0.28,0);
	\draw[semithick] (-0.28,-.6) to[out=90,in=-90] (0.28,0);
\end{tikzpicture}\:,
&
\begin{tikzpicture}[centerzero={0,0},scale=.8]
	\draw[semithick] (0.28,0) to[out=90,in=-90] (-0.28,.6);
	\draw[semithick] (-0.28,0) to[out=90,in=-90] (0.28,.6);
	\draw[semithick] (0.28,-.6) to[out=90,in=-90] (-0.28,0);
	\draw[semithick] (-0.28,-.6) to[out=90,in=-90] (0.28,0);
 \closeddot{.22,-.45};
\end{tikzpicture}
&=
\begin{tikzpicture}[centerzero={0,0},scale=.8]
 \closeddot{.23,.45};
	\draw[semithick] (0.28,0) to[out=90,in=-90] (-0.28,.6);
	\draw[semithick] (-0.28,0) to[out=90,in=-90] (0.28,.6);
	\draw[semithick] (0.28,-.6) to[out=90,in=-90] (-0.28,0);
	\draw[semithick] (-0.28,-.6) to[out=90,in=-90] (0.28,0);
\end{tikzpicture}\:.
\end{align}
Using \cref{rels2,rels3,rels4} and the first relation from \cref{rels1} but not the second (the braid relation), 
it follows that
\begin{align}\label{puppy2}
\begin{tikzpicture}[centerzero={0,0},scale=.7]
 \closeddot{-.33,-.45};
	\draw[semithick] (0.45,.6) to (-0.45,-.6);
	\draw[semithick] (0.45,-.6) to (-0.45,.6);
        \draw[semithick] (0,-.6) to[out=90,in=-90] (-.45,0);
        \draw[semithick] (-0.45,0) to[out=90,in=-90] (0,0.6);
\end{tikzpicture}
-
\begin{tikzpicture}[centerzero={0,0},scale=.7]
 \closeddot{-.33,-.45};
	\draw[semithick] (0.45,.6) to (-0.45,-.6);
	\draw[semithick] (0.45,-.6) to (-0.45,.6);
        \draw[semithick] (0,-.6) to[out=90,in=-90] (.45,0);
        \draw[semithick] (0.45,0) to[out=90,in=-90] (0,0.6);
\end{tikzpicture}&=
\begin{tikzpicture}[centerzero={0,0},scale=.7]
 \closeddot{.33,.45};
	\draw[semithick] (0.45,.6) to (-0.45,-.6);
	\draw[semithick] (0.45,-.6) to (-0.45,.6);
        \draw[semithick] (0,-.6) to[out=90,in=-90] (-.45,0);
        \draw[semithick] (-0.45,0) to[out=90,in=-90] (0,0.6);
\end{tikzpicture}
-
\begin{tikzpicture}[centerzero={0,0},scale=.7]
 \closeddot{.33,.45};
	\draw[semithick] (0.45,.6) to (-0.45,-.6);
	\draw[semithick] (0.45,-.6) to (-0.45,.6);
        \draw[semithick] (0,-.6) to[out=90,in=-90] (.45,0);
        \draw[semithick] (0.45,0) to[out=90,in=-90] (0,0.6);
\end{tikzpicture}\:,&
\begin{tikzpicture}[centerzero={0,0},scale=.7]
 \closeddot{-0.04,-.45};
	\draw[semithick] (0.45,.6) to (-0.45,-.6);
	\draw[semithick] (0.45,-.6) to (-0.45,.6);
        \draw[semithick] (0,-.6) to[out=90,in=-90] (-.45,0);
        \draw[semithick] (-0.45,0) to[out=90,in=-90] (0,0.6);
\end{tikzpicture}
-
\begin{tikzpicture}[centerzero={0,0},scale=.7]
 \closeddot{0.04,-.45};
	\draw[semithick] (0.45,.6) to (-0.45,-.6);
	\draw[semithick] (0.45,-.6) to (-0.45,.6);
        \draw[semithick] (0,-.6) to[out=90,in=-90] (.45,0);
        \draw[semithick] (0.45,0) to[out=90,in=-90] (0,0.6);
\end{tikzpicture}&=
\begin{tikzpicture}[centerzero={0,0},scale=.7]
 \closeddot{-0.04,.45};
	\draw[semithick] (0.45,.6) to (-0.45,-.6);
	\draw[semithick] (0.45,-.6) to (-0.45,.6);
        \draw[semithick] (0,-.6) to[out=90,in=-90] (-.45,0);
        \draw[semithick] (-0.45,0) to[out=90,in=-90] (0,0.6);
\end{tikzpicture}
-
\begin{tikzpicture}[centerzero={0,0},scale=.7]
 \closeddot{0.04,.45};
	\draw[semithick] (0.45,.6) to (-0.45,-.6);
	\draw[semithick] (0.45,-.6) to (-0.45,.6);
        \draw[semithick] (0,-.6) to[out=90,in=-90] (.45,0);
        \draw[semithick] (0.45,0) to[out=90,in=-90] (0,0.6);
\end{tikzpicture}
\:,&\begin{tikzpicture}[centerzero={0,0},scale=.7]
 \closeddot{.33,-.45};
	\draw[semithick] (0.45,.6) to (-0.45,-.6);
	\draw[semithick] (0.45,-.6) to (-0.45,.6);
        \draw[semithick] (0,-.6) to[out=90,in=-90] (-.45,0);
        \draw[semithick] (-0.45,0) to[out=90,in=-90] (0,0.6);
\end{tikzpicture}
-
\begin{tikzpicture}[centerzero={0,0},scale=.7]
 \closeddot{.33,-.45};
	\draw[semithick] (0.45,.6) to (-0.45,-.6);
	\draw[semithick] (0.45,-.6) to (-0.45,.6);
        \draw[semithick] (0,-.6) to[out=90,in=-90] (.45,0);
        \draw[semithick] (0.45,0) to[out=90,in=-90] (0,0.6);
\end{tikzpicture}&=
\begin{tikzpicture}[centerzero={0,0},scale=.7]
 \closeddot{-.33,.45};
	\draw[semithick] (0.45,.6) to (-0.45,-.6);
	\draw[semithick] (0.45,-.6) to (-0.45,.6);
        \draw[semithick] (0,-.6) to[out=90,in=-90] (-.45,0);
        \draw[semithick] (-0.45,0) to[out=90,in=-90] (0,0.6);
\end{tikzpicture}
-
\begin{tikzpicture}[centerzero={0,0},scale=.7]
 \closeddot{-.33,.45};
	\draw[semithick] (0.45,.6) to (-0.45,-.6);
	\draw[semithick] (0.45,-.6) to (-0.45,.6);
        \draw[semithick] (0,-.6) to[out=90,in=-90] (.45,0);
        \draw[semithick] (0.45,0) to[out=90,in=-90] (0,0.6);
\end{tikzpicture}\:.
\end{align}
\end{lem}

\begin{proof}
In view of the symmetry $\tR$, it suffices to check the first relation from \cref{puppy1} and the first two relations from \cref{puppy2}.
They are all quite easy.
For example, here are the details for the second relation from \cref{puppy2}. 
We may use \cref{rels3b,rels4b}
since they are consequences of \cref{rels2,rels3,rels4}. We have that
\begin{align*}
\begin{tikzpicture}[centerzero={0,0},scale=.8]
 \closeddot{-0.04,-.45};
	\draw[semithick] (0.45,.6) to (-0.45,-.6);
	\draw[semithick] (0.45,-.6) to (-0.45,.6);
        \draw[semithick] (0,-.6) to[out=90,in=-90] (-.45,0);
        \draw[semithick] (-0.45,0) to[out=90,in=-90] (0,0.6);
\end{tikzpicture}&=\begin{tikzpicture}[centerzero={0,0},scale=.8]
 \closeddot{-0.45,0};
	\draw[semithick] (0.45,.6) to (-0.45,-.6);
	\draw[semithick] (0.45,-.6) to (-0.45,.6);
        \draw[semithick] (0,-.6) to[out=90,in=-90] (-.45,0);
        \draw[semithick] (-0.45,0) to[out=90,in=-90] (0,0.6);
\end{tikzpicture}-
\begin{tikzpicture}[centerzero={0,0},scale=.8]
 	\draw[semithick] (0,.6) to (-0.45,-.6);
	\draw[semithick] (.45,.6) to (0,-.6);
        \draw[semithick] (-0.45,.6)to (.45,-0.6);
\end{tikzpicture}=\begin{tikzpicture}[centerzero={0,0},scale=.8]
 \closeddot{-0.04,.5};
	\draw[semithick] (0.45,.6) to (-0.45,-.6);
	\draw[semithick] (0.45,-.6) to (-0.45,.6);
        \draw[semithick] (0,-.6) to[out=90,in=-90] (-.45,0);
        \draw[semithick] (-0.45,0) to[out=90,in=-90] (0,0.6);
\end{tikzpicture}-
\begin{tikzpicture}[centerzero={0,0},scale=.8]
 	\draw[semithick] (0,.6) to (-0.45,-.6);
	\draw[semithick] (.45,.6) to (0,-.6);
        \draw[semithick] (-0.45,.6)to (.45,-0.6);
\end{tikzpicture}+
\begin{tikzpicture}[centerzero={0,0},scale=.8]
 	\draw[semithick] (0,-.6) to (-0.45,.6);
	\draw[semithick] (.45,-.6) to (0,.6);
        \draw[semithick] (-0.45,-.6)to (.45,0.6);
\end{tikzpicture}\:,\\
\begin{tikzpicture}[centerzero={0,0},scale=.8]
 \closeddot{0.04,-.45};
	\draw[semithick] (0.45,.6) to (-0.45,-.6);
	\draw[semithick] (0.45,-.6) to (-0.45,.6);
        \draw[semithick] (0,-.6) to[out=90,in=-90] (.45,0);
        \draw[semithick] (0.45,0) to[out=90,in=-90] (0,0.6);
\end{tikzpicture}
&=
\begin{tikzpicture}[centerzero={0,0},scale=.8]
 \closeddot{0.45,0};
	\draw[semithick] (0.45,.6) to (-0.45,-.6);
	\draw[semithick] (0.45,-.6) to (-0.45,.6);
        \draw[semithick] (0,-.6) to[out=90,in=-90] (.45,0);
        \draw[semithick] (0.45,0) to[out=90,in=-90] (0,0.6);
\end{tikzpicture}+
\begin{tikzpicture}[centerzero={0,0},scale=.8]
 	\draw[semithick] (0,-.6) to (-0.45,.6);
	\draw[semithick] (.45,-.6) to (0,.6);
        \draw[semithick] (-0.45,-.6)to (.45,0.6);
\end{tikzpicture}=
\begin{tikzpicture}[centerzero={0,0},scale=.8]
 \closeddot{0.04,.5};
	\draw[semithick] (0.45,.6) to (-0.45,-.6);
	\draw[semithick] (0.45,-.6) to (-0.45,.6);
        \draw[semithick] (0,-.6) to[out=90,in=-90] (.45,0);
        \draw[semithick] (0.45,0) to[out=90,in=-90] (0,0.6);
\end{tikzpicture}+
\begin{tikzpicture}[centerzero={0,0},scale=.8]
 	\draw[semithick] (0,-.6) to (-0.45,.6);
	\draw[semithick] (.45,-.6) to (0,.6);
        \draw[semithick] (-0.45,-.6)to (.45,0.6);
\end{tikzpicture}-
\begin{tikzpicture}[centerzero={0,0},scale=.8]
 	\draw[semithick] (0,.6) to (-0.45,-.6);
	\draw[semithick] (.45,.6) to (0,-.6);
        \draw[semithick] (-0.45,.6)to (.45,-0.6);
\end{tikzpicture}
\:.
\end{align*}
Subtracting gives the result.  In this calculation we have omitted several diagrams which are 0 due to the first relation from \cref{rels1}. 
\end{proof}

Next we exploit the isomorphism
$\Lambda \cong \Gamma \otimes \Lambda^{[2]}$
from \cref{clanky} and freeness of
$\Gamma = \kk[q_1,q_3,\dots]$ to see that
there is a well-defined
$\Lambda^{[2]}$-algebra homomorphism
\begin{align}\label{gammat}
\zeta:\Lambda &\rightarrow \End_{\cNB_t}(\one)
\end{align}
mapping $q_{2r-1} \in \Gamma$ to
$-2(-1)^t \begin{tikzpicture}[anchorbase]
\draw[semithick] (0,0) circle (.2);
\closeddot{.2,0};
\node at (.6,0) {$\scriptstyle 2r-1$};
\end{tikzpicture}$
for each $r \geq 1$.
By \cite[Cor.~2.6]{BWWnilBrauer},
it follows that $\zeta$
maps $q_{r}$ to
$-2(-1)^t \begin{tikzpicture}[anchorbase]
\draw[semithick] (0,0) circle (.2);
\closeddot{.2,0};
\node at (.35,0) {$\scriptstyle r$};
\end{tikzpicture}$
for every $r \geq 1$. However, we warn the reader that $q_0 = 1$ does not agree with $-2(-1)^t \begin{tikzpicture}[anchorbase]
\draw[semithick] (0,0) circle (.2);
\end{tikzpicture}$\ .

We denote the image of any $a \in \Lambda$ under $\zeta$ 
simply by the labelled bubble $\Bubble{a}$\,.
The morphism space
$\Hom_{\cNB_t}(B^{\star n}, B^{\star m})$
is naturally a $(\Lambda,\Lambda)$-bimodule
so that the left and right actions $a \in \Lambda$
are by horizontal composition on the left or right with $\Bubble{a}$, respectively. The left and right actions by elements $\Lambda^{[2]}$ coincide since we are considering a $\Lambda^{[2]}$-linear monoidal category, but this is seldom the case for elements of $\Gamma$.
Note also that the symmetries $\tR$ and $\tT$ from \cref{R,T} fix all of the 
bubbles $\Bubble{a}$ for any $a \in \Lambda$.

\begin{theo}\label{whereswaldo}
For $m,n\geq 0$, the morphism space $\Hom_{\cNB_t}(B^{\star n}, B^{\star m})$
is free as a right $\Lambda$-module with basis given by a
certain combinatorially-defined set $\operatorname{D}(m,n)$ of {\em dotted reduced $m \times n$ string diagrams}.
In particular, the homomorphism $\zeta$ from \cref{gammat} is an isomorphism.
\end{theo}

\begin{proof}
In \cite[Th.~5.3]{BWWnilBrauer},
it is proved that morphism spaces in the nil-Brauer category considered there are
free as right $\Gamma$-modules with basis
$\operatorname{D}(m,n)$.
In view of \cref{clanky}, the result here
follows since our nil-Brauer category is the one from \cite{BWWnilBrauer} base-changed from $\kk$ to $\Lambda^{[2]}$.
\end{proof}

\subsection{Identification of the Grothendieck ring
of \texorpdfstring{$\cNB_t$}{}}
A thick string $\,\begin{tikzpicture}[anchorbase]
\draw[ultra thick] (0,.15) to (0,-.15);
\node at (0,-.25) {$\stringlabel{n}$};
\end{tikzpicture}\,$
labelled by $n$ indicates $n$ parallel thin strings, i.e., it is the identity endomorphism of $B^{\star n}$.
The crossing of two thick strings denotes the minimal length composition of crossings of thin strings:
\begin{align*}
\begin{tikzpicture}[baseline=-1.5pt,scale=.9]
\draw [ultra thick] (-.3,-.4) to (.3,.4);
\draw [ultra thick] (-.3,0.4) to (.3,-.4);
\node at (-.3,-.55) {$\stringlabel{a}$};
\node at (.3,-.55) {$\stringlabel{b}$};
\end{tikzpicture}
&=
\begin{tikzpicture}[baseline=-1.5pt,scale=.9]
\draw [semithick,-] (-.3,-.4) to (.45,.4);
\draw [semithick,-] (-.45,-.4) to (.3,.4);
\draw [semithick,-] (-.6,-.4) to (0.15,.4);
\draw [semithick,-] (-.3,0.4) to (.3,-.4);
\draw [semithick,-] (-.6,0.4) to (0,-.4);
\draw [semithick,-] (-.45,0.4) to (.15,-.4);
\draw [semithick,-] (-.15,0.4) to (.45,-.4);
\node at (-.45,-.55) {$\stringlabel{a}$};
\node at (.3,-.55) {$\stringlabel{b}$};
\end{tikzpicture}\:.
\end{align*}
We use a cross on a string of thickness $n$
to indicate the composition of thin crossings
according to a reduced expression for the longest permutation in $S_n$. On  a string of thickness one, this is just the identity.
In general, by the braid relation, we have that
\begin{align*}
\begin{tikzpicture}[baseline=-1.5pt,scale=.9]
\draw[ultra thick] (0,.4) to (0,-.4);
\cross{(0,0)}
\node at (0,-.6) {$\stringlabel{a+b}$};
\end{tikzpicture}
&=
\begin{tikzpicture}[baseline=-1.5pt,scale=.9]
\draw[ultra thick] (.4,.4) to (-.4,-.4);
\draw[ultra thick] (-.4,.4) to (.4,-.4);
\cross[45]{(0,-.3)}
\cross[-45]{(0,-.3)}
\node at (-.4,-.6) {$\stringlabel{a}$};
\node at (.4,-.6) {$\stringlabel{b}$};
\end{tikzpicture}
\end{align*}
for any $a+b=n$.
Given $\alpha = (\alpha_1,\dots,\alpha_n) \in \N^n$, let
$$
\begin{tikzpicture}[anchorbase,scale=.9]
\draw[ultra thick] (0,.4) to (0,-.4);
\closeddot{0,0};
\node at (0,-.55) {$\stringlabel{n}$};
\node at (0.25,0) {$\scriptstyle{\alpha}$};
\node at (0,.55){$\stringlabel{\phantom{n}}$};
\end{tikzpicture}
:=
\begin{tikzpicture}[anchorbase,scale=.9]
\draw[semithick,-] (-0.2,.4) to (-0.2,-.4);
\draw[semithick,-] (0.4,.4) to (0.4,-.4);
\draw[semithick,-] (1.2,.4) to (1.2,-.4);
\node at (-0.47,0) {$\scriptstyle{\alpha_1}$};
\node at (0.14,0) {$\scriptstyle{\alpha_2}$};
\node at (1.48,0) {$\scriptstyle{\alpha_n}$};
\node at (.65,0) {$\scriptstyle .$};
\node at (.8,0) {$\scriptstyle .$};
\node at (.95,0) {$\scriptstyle .$};
\closeddot{-0.2,0};
\closeddot{0.4,0};
\closeddot{1.2,0};
\end{tikzpicture}
$$
Using this notation, we define
\begin{equation}\label{endef}
\e_n :=
\begin{tikzpicture}[baseline=-1.5pt]
\draw[ultra thick] (0,.4) to (0,-.4);
    \closeddot{0,0.2};
    \node at (.25,.2) {$\scriptstyle{\rho_n}$};
    \cross{(0,-.1)};
\node at (0,-.55) {$\stringlabel{n}$};
\end{tikzpicture}\:\in \End_{\cNB_t}\big(B^{\star n}\big),
\end{equation}
where $\rho_n :=(n-1,\cdots,2,1,0) \in \N^n$.
For example:
\begin{align*}
\e_0 &= \id_\one,&
\e_1
&= 
\begin{tikzpicture}[anchorbase,scale=1.4]
\draw[semithick] (0,-.4) to (0,.4);
\end{tikzpicture}\:,&
\e_2 &=
\begin{tikzpicture}[anchorbase,scale=1.4]
\draw[semithick] (0.3,-.4) to (-.3,.4);
\draw[semithick] (-.3,-.4) to (0.3,.4);
\closeddot{-.15,.2};
\end{tikzpicture}\:,&
\e_3 &=
\begin{tikzpicture}[anchorbase,scale=1.4]
\draw[semithick] (0.3,-.4) to (-.3,.4);
\draw[semithick] (-.3,-.4) to (0.3,.4);
\draw[semithick] (0,.4) to[out=-45,in=90] (.3,0) to[out=-90,in=45] (0,-.4);
\closeddot{.08,.33};
\closeddot{-.25,.33};
\closeddot{-.14,.19};
\end{tikzpicture}\:\;,&
\e_4 &=
\begin{tikzpicture}[anchorbase,scale=1.4]
\draw[semithick] (0.3,-.4) to (-.3,.4);
\draw[semithick] (-.3,-.4) to (0.3,.4);
\draw[semithick] (.1,.4) to[out=-50,in=90,looseness=.75] (.3,.15) to[out=-90,in=35,looseness=.75] (-.1,-.4);
\draw[semithick] (-.1,.4) to[out=-35,in=90,looseness=.75] (.3,-.15) to[out=-90,in=50,looseness=.75] (.1,-.4);
\closeddot{-.26,.345};
\closeddot{-0.03,.345};
\closeddot{.15,.345};
\closeddot{.07,.25};
\closeddot{-.19,.25};
\closeddot{-.12,.16};
\end{tikzpicture}\:\;.
\end{align*}
The significance of these endomorphisms is 
explained by the following result:

\begin{theo}\label{aboutidempotents}
For $n \geq 0$,
$\e_n$ is a primitive homogeneous idempotent, and any primitive homogeneous idempotent in $\cNB_t$ is equivalent\footnote{Homogeneous 
idempotents $\e:X \rightarrow X$ and $\f:Y \rightarrow Y$ in a graded 
category $\cC$ are {\em equivalent} if there exist
homogeneous morphisms $\u:Y \rightarrow X$ and $\v:X \rightarrow Y$
such that $\e = \u \circ \v$ and $\f = \v \circ \u$.}
to $\e_n$ for a unique $n\geq 0$.
\end{theo}

\begin{proof}
This is proved in 
\cite[Cor.~4.24]{BWWiquantum}, but working over the ground field $\kk$ rather than the ground ring $\Lambda^{[2]}$ here. 
This change has no effect on homogeneous idempotents,
since they necessarily have degree 0 and 
$\Lambda^{[2]}$ is a connected positively graded algebra.
\end{proof}

For a graded category $\cC$, let $\gKar(\cC)$
be its graded Karoubi envelope.
This is obtained by enlarging the category by formally adjoining an invertible grading shift functor $q$
in such a way that 
$$
\Hom(X,Y)_n = \left(q^{-n} \Hom(X,Y)\right)_0 = \Hom(X, q^{-n} Y)_0 = \Hom(q^{n} X, Y)_0,
$$
then passing to the usual additive Karoubi envelope. Thus, objects of $\gKar(\cC)$ are pairs $[X, \e]$
consisting of a formal finite direct sum $X$ of grading shifts of objects of $\cC$ and a matrix $\e$ of endomorphisms defining a homogeneous idempotent $\e:X \rightarrow X$.
We write $K_0(\gKar(\cC))$
for the split Grothendieck group consisting of degree 0 isomorphism classes of objects in $\gKar(\cC)$. It is a $\Z[q,q^{-1}]$-module with action of $q$ induced by the grading shift functor.
If $\cC$ is monoidal then $\gKar(\cC)$
is too, so we get an induced multiplication making
$K_0(\gKar(\cC))$ into a $\Z[q,q^{-1}]$-algebra.

In \cite{BWWiquantum}, it is shown that
the $\Z[q,q^{-1}]$-algebra
$K_0(\gKar(\cNB_t))$ is isomorphic to 
a certain $\Z[q,q^{-1}]$-form
$\Ui$ for the split iquantum group of rank one. We refer to \cite[Sec.~2]{BWWiquantum} for the full definition of this, just noting for now that $\Q(q) \otimes_{\Z[q,q^{-1}]} \Ui$
is the polynomial algebra $\Q(q)[b]$, and
the $\Z[q,q^{-1}]$-form
$\Ui$ is free as a $\Z[q,q^{-1}]$-module with basis $b^{(n)}\:(n \geq 0)$ defined from
\begin{equation}\label{idp}
b^{(n)} :=
\begin{dcases}
\frac{1}{[n]^!}
\prod_{\substack{k=0\\k\equiv t\pmod{2}}}^{n-1}
\!\!\!\left(b^2 - [k]^2\right)
&\text{if $n$ is even}\\
\frac{b}{[n]^!}
\prod_{\substack{k=1\\k\equiv t\pmod{2}}}^{n-1}
\!\!\!\left(b^2 - [k]^2\right)
&\text{if $n$ is odd.}
\end{dcases}
\end{equation}
The basis $b^{(n)}\:(n \geq 0)$ is
the {\em icanonical basis} of $\Ui$ in the general sense
of \cite{BW18QSP, BW18KL} associated to the parameter $t$. These elements are also known as {idivided powers}.
Instead of the closed
formula \cref{idp}, which was worked out originally in \cite{BeW18},
$b^{(n)}$ can be also defined recursively:
we have that
$b^{(0)} := 1$ and 
\begin{equation}\label{ryder}
b \cdot b^{(n)} = 
\begin{cases}
{[n+1]_q b^{(n+1)}+[n]_q b^{(n-1)}}&\text{if $n \equiv t\pmod{2}$}\\
{[n+1]_q b^{(n+1)}}&\text{if $n \not\equiv t\pmod{2}$.}
\end{cases}
\end{equation}
Let 
$\Par(r \times c)$ be the set of partitions whose Young diagram fits into an $r \times c$ rectangle, i.e., partitions with at most $r$ non-zero parts, all of which are $\leq c$.
Let $\Par_{\not\equiv t}(r\times c)$ be the subset consisting of those partitions whose non-zero parts
are $\not\equiv t\pmod{2}$.
Another notable formula \cite[Cor.~2.13]{BWWiquantum} gives that
\begin{equation}\label{ohio}
b^n= \sum_{i=0}^{\lfloor \frac{n}{2}\rfloor}
[n-2i]_q^!
\left(\sum_{\lambda \in \Par_{\not\equiv t}(i\times (n-2i))}
[\lambda_1+1]_q^2\cdots [\lambda_i+1]_q^2\right) 
b^{(n-2i)}
\end{equation}
for any $n \geq 0$.

\begin{theo}\label{eatyourselffitter}
There is an isomorphism of 
$\Z[q,q^{-1}]$-algebras
$K_0(\gKar(\cNB_t))
\stackrel{\sim}{\rightarrow}\Ui$
taking the isomorphism class of 
\begin{equation}\label{poorpolly}
B^{(n)} := q^{-\binom{n}{2}}[B^{\star n}, \e_n] \in \gKar(\cNB_t)
\end{equation}
to the icanonical basis vector $b^{(n)}$.
\end{theo}

\begin{proof}
This follows from
\cite[Th.~B]{BWWiquantum}, but a couple of comments are in order.
One is that
we have extended the ground ring from $\kk$ to $\Lambda^{[2]}$, but this causes no problem
since $\Lambda^{[2]}$ is a connected positively graded algebra.
More likely to cause confusion, in \cite{BWWiquantum}, results were not explained in terms of graded Karoubi envelopes, rather, they were formulated in terms of finitely generated graded
projective {\em left} modules over the path algebra $\NB$ of $\cNB_t$.
The graded category $\NB\pgMod$ of finitely generated projective graded left $\NB$-modules
is {\em contravariantly} equivalent to $\gKar(\cNB_t)$ via the Yoneda equivalence.
Being contravariant, the canonical isomorphism
$$
K_0(\gKar(\cNB_t)) \stackrel{\sim}{\rightarrow} K_0(\NB\pgMod)
$$
induced by the Yoneda equivalence is {\em anti-linear}. 
This accounts for the fact that the grading shift in \cref{poorpolly} is the negation of the grading shift in the definition of the corresponding indecomposable projective module $P(n)$ defined in \cite[(4.33)]{BWWiquantum}.
\end{proof}

From \cref{eatyourselffitter,ohio}, we get that
\begin{equation}\label{theworld}
B^{\star n} \cong 
\bigoplus_{i=0}^{\lfloor\frac{n}{2}\rfloor}
[n-2i]_q^!
\left(\sum_{\lambda \in \Par_{\not\equiv t}(i\times (n-2i))}
[\lambda_1+1]_q^2\cdots [\lambda_i+1]_q^2\right) B^{(n-2i)}
\end{equation}
for any $n \geq 0$.
There is a unique summand
in this direct sum decomposition equal
to $q^{\binom{n}{2}} B^{(n)}$. This is the ``image" of the idempotent $\e_n$.

\subsection{Bubble slides}
Given $f(x) = \sum_{r\geq 0} c_r x^r \in \kk[x]$
and a dot in some string diagram $s$,
we denote
$$
\sum_{r \geq 0} c_r \times 
(\text{the morphism obtained from $s$ by labeling the dot by }r)
$$
by attaching what we call a {\em pin} to the dot, labeling 
the node at the head of the pin by $f(x)$:
\begin{equation}\label{labelledpin}
\begin{tikzpicture}[anchorbase]
	\draw[semithick,-] (0,-.4) to (0,.4);
    \pin{(0,0)}{(.8,0)}{f(x)};
\end{tikzpicture}\;
:=
\sum_{r \geq 0} c_r 
\:\begin{tikzpicture}[anchorbase]
	\draw[semithick,-] (0,-.4) to (0,.4);
    \closeddot{0,0};
    \node at (.2,0) {$\scriptstyle{r}$};
\end{tikzpicture}
\in \End_{\cNB_t}(B).
\end{equation}
More generally, $f(x)$ here could be a polynomial
with coefficients in the algebra $\kk\lround u^{-1} \rround$
of formal Laurent series in an indeterminate $u^{-1}$;
then the string $s$ decorated with a pin labelled $f(x)$ defines a morphism in $\cNB_t\lround u^{-1}\rround$. We think of this as being a generating function for a family of morphisms.

Now we can discuss ``bubble slides".
We obviously have that
\begin{equation}\label{victor0}
\Bubble{e^{[2]}(u)}\:\:
\begin{tikzpicture}[anchorbase]
	\draw[semithick,-] (0,-.3) to (0,.5);
\end{tikzpicture}=
\begin{tikzpicture}[anchorbase]
	\draw[semithick,-] (0,-.3) to (0,.5);
\end{tikzpicture}\:\:\Bubble{e^{[2]}(u)}
\end{equation}
since $\cNB_t$ is a $\Lambda^{[2]}$-linear
monoidal category.
The following is \cite[Th.~2.5(5)]{BWWnilBrauer}:
\begin{equation}\label{victor1}
\Bubble{q(u)}\:\:
\begin{tikzpicture}[anchorbase]
	\draw[semithick,-] (0,-.3) to (0,.5);
\end{tikzpicture}=
\begin{tikzpicture}[anchorbase]
	\draw[semithick,-] (0,-.3) to (0,.5);
	\pin{(0,.1)}{(1,.1)}{\left(\frac{u-x}{u+x}\right)^2};
\end{tikzpicture}\:\:\Bubble{q(u)}\:.
\end{equation}
Writing $\sqrt{q(u)}$ for the unique square root of $q(u)$ in $1+u^{-1}\Gamma\llbracket u^{-1}\rrbracket$, 
\cref{victor1} is equivalent to
\begin{align}\label{victor8}
\Bubble{\sqrt{q(u)}}\:\:
\begin{tikzpicture}[anchorbase]
	\draw[semithick,-] (0,-.3) to (0,.5);
\end{tikzpicture}&=
\begin{tikzpicture}[anchorbase]
	\draw[semithick,-] (0,-.3) to (0,.5);
	\pin{(0,.1)}{(.8,.1)}{\frac{u-x}{u+x}};
\end{tikzpicture}\:\:\Bubble{\sqrt{q(u)}}\:.
\end{align}
The following lemma describes the bubble slides
for elementary and complete symmetric functions;
the formulae
are the same as the one for $\sqrt{q(u)}$:

\begin{lem}\label{victor}
We have that
\begin{align}\label{victor2}
\Bubble{e(u)}\:\:
\begin{tikzpicture}[anchorbase]
	\draw[semithick,-] (0,-.3) to (0,.5);
\end{tikzpicture}&=
\begin{tikzpicture}[anchorbase]
	\draw[semithick,-] (0,-.3) to (0,.5);
	\pin{(0,.1)}{(.8,.1)}{\frac{u-x}{u+x}};
\end{tikzpicture}\:\:\Bubble{e(u)}\:,&
\Bubble{h(u)}\:\:
\begin{tikzpicture}[anchorbase]
	\draw[semithick,-] (0,-.3) to (0,.5);
\end{tikzpicture}&=
\begin{tikzpicture}[anchorbase]
	\draw[semithick,-] (0,-.3) to (0,.5);
	\pin{(0,.1)}{(.8,.1)}{\frac{u-x}{u+x}};
\end{tikzpicture}\:\:\Bubble{h(u)}\:.
\end{align}
\end{lem}

\begin{proof}
Recall that $e^{[2]}(u) = e(u) e(-u)$
and $q(u) = e(u) e(-u)^{-1}$.
Hence, $e(u)^2 = e^{[2]}(u) q(u)$
and
$e(u) = \sqrt{e^{[2]}(u)}\sqrt{q(u)}$ (again, for $e^{[2]}(u)$, we take the square root with constant term 1).
The first identity in \cref{victor2} is now clear from \cref{victor0,victor8}.
To deduce the second one, replace $u$ by $-u$ and then take inverses on both sides.
\end{proof}

\begin{cor}\label{cyclotomicrelationexplained}
We have that
$\displaystyle\Bubble{e_r}\:
\begin{tikzpicture}[anchorbase]
	\draw[semithick,-] (0,-.3) to (0,.5);
\end{tikzpicture}
\:+\:
\begin{tikzpicture}[anchorbase]
	\draw[semithick,-] (0,-.3) to (0,.5);
\end{tikzpicture}\:\Bubble{e_r}
=
2\sum_{s=0}^r
\Bubble{e_{r-s}}
\begin{tikzpicture}[anchorbase]
\node at (.18,.1) {$\scriptstyle s$};
\closeddot{0,.1};
	\draw[semithick,-] (0,-.3) to (0,.5);
\end{tikzpicture}=2\sum_{s=0}^r (-1)^s
\begin{tikzpicture}[anchorbase]
\node at (.18,.1) {$\scriptstyle s$};
\closeddot{0,.1};
	\draw[semithick,-] (0,-.3) to (0,.5);
\end{tikzpicture}\Bubble{e_{r-s}}\:.
$
\end{cor}

\begin{proof}
For the first equality, 
equate coefficients of $u^{-r}$ on both sides
of the first relation from \cref{victor2}
using the identity
$\frac{u-x}{u+x}=2(1+xu^{-1})^{-1} - 1
= 2 \sum_{s \geq 0} (-1)^s x^s u^{-s} -1$.
The second equality follows from the first 
by applying $\tR$.
\end{proof}

\subsection{The cyclotomic nil-Brauer category}

The endomorphism of $B$ appearing in \cref{cyclotomicrelationexplained} will play an important role, so we introduce some special notation for it: let
\begin{align}\label{purplesquare}
\begin{tikzpicture}[anchorbase]
	\draw[semithick,-] (0,-.3) to (0,.5);
\node[rectangle,draw=black,inner sep=2.5,fill=cyclocolor] at (0,0.1) {$\scriptstyle r$};
\end{tikzpicture}\ 
&:=\ \Bubble{e_r}\:
\begin{tikzpicture}[anchorbase]
	\draw[semithick,-] (0,-.3) to (0,.5);
\end{tikzpicture}
\:+\:
\begin{tikzpicture}[anchorbase]
	\draw[semithick,-] (0,-.3) to (0,.5);
\end{tikzpicture}\:\Bubble{e_r}\ .
\end{align}
For $l \geq 0$, let $\cI_l$ be the four-sided ideal of
$\cNB_t$ generated by
$\ \begin{tikzpicture}[anchorbase,scale=.8]
	\draw[semithick,-] (0,-.3) to (0,.5);
\node[rectangle,draw=black,inner sep=2.5,fill=cyclocolor] at (0,0.1) {$\scriptstyle l$};
\end{tikzpicture}\ $ and the bubbles
$\ \Bubble{e_{2r}}\ $ for all $r$ with $2r > l$.
By ``four-sided ideal'' here, we mean that it is a {\em two-sided ideal}
of the $\Lambda$-linear category
$\cNB_t$, i.e. a
family of $\Lambda$-submodules
$\cI_l(B^{\star n}, B^{\star m})
\subseteq \Hom_{\cNB_t}(B^{\star n}, B^{\star m})$ for all $m,n \geq 0$
closed under vertical composition on top and bottom with any morphism, 
which is also a {\em two-sided tensor ideal},
i.e. it is closed under horizontal composition on left and right by any morphism. 

\begin{theo}\label{bensdumbideatocallitfoursided}
The four-sided ideal $\cI_l$ is equal both to the right tensor ideal and to the left tensor ideal
generated by
$\ \begin{tikzpicture}[anchorbase,scale=.8]
	\draw[semithick,-] (0,-.3) to (0,.5);
\node[rectangle,draw=black,inner sep=2.5,fill=cyclocolor] at (0,0.1) {$\scriptstyle l$};
\end{tikzpicture}\ $ and $\ \Bubble{e_{2r}}\:(2r > l)$.
Moreover, $\cI_l(\one,\one)$ contains
$\ \Bubble{e_{r}}\ $ for all $r > l$.
\end{theo}

\begin{proof}
Let $\cJ_l$ be the
right tensor ideal generated by 
$\ \begin{tikzpicture}[anchorbase,scale=.8]
	\draw[semithick,-] (0,-.3) to (0,.5);
\node[rectangle,draw=black,inner sep=2.5,fill=cyclocolor] at (0,0.1) {$\scriptstyle l$};
\end{tikzpicture}\ $ and $\ \Bubble{e_{2r}}\:(2r > l)$.
In a series of claims below, we will show that $\cJ_l$ is a two-sided tensor ideal and that $\cJ_l(\one,\one)$ contains $\ \Bubble{e_r}\ $ for all $r > l$.
This is sufficient to complete the proof of the theorem. Indeed,
since $\cJ_l$ has the same generators as $\cI_l$,
it follows that $\cI_l = \cJ_l$, so 
 $\cI_l$ is the right tensor ideal with these generators. Applying $\tR$, it also follows that $\cI_l$ is the left tensor ideal with these generators.

\vspace{2mm}
\noindent
\underline{Claim 0}: {\em $\cJ_l(\one,\one)$ contains
$\ \Bubble{e_r}\ $ for all $r > l$.}

\vspace{1mm}
\noindent
\underline{Proof}.
To prove this, we will use the following notation: for a formal series $f(u) = \sum_{r \geq 0} f_r u^{-r}$, we use  $[f(u)]_{> l}$ to denote
$\sum_{r > l} f_r u^{-r}$ and $[f(u)]_{\leq l}$
to denote $\sum_{r=0}^{l} f_r u^{-r}$.
By \cref{cyclotomicrelationexplained},
$\cJ_l(B,B)$ contains $\sum_{s=0}^l \Bubble{e_{l-s}}
\begin{tikzpicture}[anchorbase]
\node at (.18,.1) {$\scriptstyle s$};
\closeddot{0,.1};
	\draw[semithick,-] (0,-.1) to (0,.3);
\end{tikzpicture}$.
If we add $r-l$ dots to the string, then close on the right, we deduce that $\cJ_l(\one,\one)$ contains
$\sum_{s=0}^l \Bubble{e_{l-s}}
\begin{tikzpicture}[anchorbase,scale=1]
\draw[semithick] (0.4,0) circle (.2);
\closeddot{0.2,0};
\node at (-.35,0) {$\scriptstyle r+s-l$};
\end{tikzpicture}\ $ for all $r > l$.
Hence, $\cJ_l(\one,\one)\llbracket u^{-1}\rrbracket$ contains 
the image under $\zeta$ of $\left[[e(u)]_{\leq l} q(-u)\right]_{> l}$. Now we observe that
\begin{align*}
\left[[e(u)]_{\leq l} q(-u)\right]_{>l} &= 
\left[e(u)q(-u) - [e(u)]_{>l}q(-u)\right]_{>l}\\
&=
[e(u)h(-u) e(-u)]_{>l} - [e(u)]_{>l} q(-u)
\\&= [e(-u)]_{>l}  - [e(u)]_{>l} q(-u).
\end{align*}
We have that
$\zeta\left([e(-u)]_{> l}\right)\equiv \zeta\left(-[e(u)]_{> l}\right)\pmod{\cJ_l(\one,\one)\llbracket u^{-1}\rrbracket}$ 
as $\zeta(e_{2r}) \in \cJ_l(\one,\one)$ for $2r > l$.
We deduce that
$\zeta\left([e(u)]_{> l} (1+q(-u))\right)$ is in $\cJ_l(\one,\one)\llbracket u^{-1}\rrbracket$.
Hence, $\zeta\left([e(u)]_{>l}\right)$ is in $\cJ_l(\one,\one)\llbracket u^{-1}\rrbracket$.

\vspace{2mm}
\noindent
\underline{Claim 1}: {\em We have that
$\ \begin{tikzpicture}[anchorbase,scale=1.1]
	\draw[semithick] (0.4,0) to[out=90, in=0] (0.1,0.5);
	\draw[semithick] (0.1,0.5) to[out = 180, in = 90] (-0.2,0);
\node[rectangle,draw=black,inner sep=2.5,fill=cyclocolor] at (.35,0.3) {$\scriptstyle l$};
\end{tikzpicture} \in \cJ_l(B\star B, \one)$ and $\ \begin{tikzpicture}[anchorbase,scale=1.1]
	\draw[semithick] (0.4,0) to[out=-90, in=0] (0.1,-0.5);
	\draw[semithick] (0.1,-0.5) to[out = 180, in = -90] (-0.2,0);
\node[rectangle,draw=black,inner sep=2.5,fill=cyclocolor] at (.35,-0.3) {$\scriptstyle l$};
\end{tikzpicture} \in \cJ_l(\one,B\star B)$.}

\vspace{1mm}
\noindent
\underline{Proof}.
The first assertion follows because
$\begin{tikzpicture}[anchorbase,scale=1.1]
	\draw[semithick] (0.4,0) to[out=90, in=0] (0.1,0.5);
	\draw[semithick] (0.1,0.5) to[out = 180, in = 90] (-0.2,0);
\node[rectangle,draw=black,inner sep=2.5,fill=cyclocolor] at (.35,0.3) {$\scriptstyle l$};
\end{tikzpicture} = \begin{tikzpicture}[anchorbase,scale=1.1]
	\draw[semithick] (0.4,0) to[out=90, in=0] (0.1,0.5);
	\draw[semithick] (0.1,0.5) to[out = 180, in = 90] (-0.2,0);
\node[rectangle,draw=black,inner sep=2.5,fill=cyclocolor] at (-.15,0.3) {$\scriptstyle l$};
\end{tikzpicture}$\, , as is clear from the definition \cref{purplesquare}.
The second one is similar (or one can apply $\tT$).

\vspace{2mm}
\noindent
\underline{Claim 2}: {\em We have that
$\begin{tikzpicture}[anchorbase,scale=.8]
	\draw[semithick,-] (0,-.3) to (0,.5);
\node[rectangle,draw=black,inner sep=2.5,fill=cyclocolor] at (0,0.1) {$\scriptstyle r$};
\end{tikzpicture} \in \cJ_l(B,B)$ for all $r \geq l$.}

\vspace{1mm}
\noindent
\underline{Proof}.
This is proved by induction on $r$, the base case $r=l$ being immediate from the definition of $\cJ_l$. Now suppose $r > l$. By \cref{cyclotomicrelationexplained} we have that
\begin{equation}
\begin{tikzpicture}[anchorbase,scale=1.1]
	\draw[semithick,-] (0,-.3) to (0,.5);
\node[rectangle,draw=black,inner sep=2.5,fill=cyclocolor] at (0,0.1) {$\scriptstyle r$};
\end{tikzpicture}
=
2\sum_{s=0}^r
\Bubble{e_{r-s}}
\begin{tikzpicture}[anchorbase,scale=1.1]
\node at (.18,.1) {$\scriptstyle s$};
\closeddot{0,.1};
	\draw[semithick,-] (0,-.3) to (0,.5);
\end{tikzpicture}
=
2 \ 
\Bubble{e_{r}}\ 
\begin{tikzpicture}[anchorbase,scale=1.1]
	\draw[semithick,-] (0,-.3) to (0,.5);
\end{tikzpicture}
+2\sum_{s=1}^r \Bubble{e_{r-s}}
\begin{tikzpicture}[anchorbase,scale=1.1]
\node at (.35,.25) {$\scriptstyle s-1$};
\closeddot{0,.25};
\closeddot{0,-.05};
	\draw[semithick,-] (0,-.3) to (0,.5);
\end{tikzpicture}= 
2 \ 
\Bubble{e_{r}}\ 
\begin{tikzpicture}[anchorbase,scale=1.1]
	\draw[semithick,-] (0,-.3) to (0,.5);
\end{tikzpicture}
+2\sum_{s=0}^{r-1} \Bubble{e_{r-1-s}}
\begin{tikzpicture}[anchorbase,scale=1.1]
\node at (.185,.25) {$\scriptstyle s$};
\closeddot{0,.25};
\closeddot{0,-.05};
	\draw[semithick,-] (0,-.3) to (0,.5);
\end{tikzpicture}
=2 \ 
\Bubble{e_{r}}\ 
\begin{tikzpicture}[anchorbase,scale=1.1]
	\draw[semithick,-] (0,-.3) to (0,.5);
\end{tikzpicture}
+
\begin{tikzpicture}[anchorbase,scale=1.1]
	\draw[semithick,-] (0,-.3) to (0,.5);
\node[rectangle,draw=black,inner sep=2.5,fill=cyclocolor] at (0,0.2) {$\scriptstyle r-1$};
\closeddot{0,-.15};
\end{tikzpicture}\ .
\label{soundedlike}\end{equation}
This lies in $\cJ_l$ by induction (and because $\displaystyle\Bubble{e_r} \in \cJ_l$).

\vspace{2mm}
\noindent
\underline{Claim 3}: {\em We have that $\begin{tikzpicture}[anchorbase,scale=.8]
	\draw[semithick,-] (0,-.3) to (0,.5);
\end{tikzpicture}\ 
\Bubble{e_r} \in \cJ_l(B,B)$ for any $r > l$.}

\vspace{1mm}
\noindent
\underline{Proof}.
This follows from Claim 2 because (by \cref{purplesquare})
$$\begin{tikzpicture}[anchorbase,scale=1]
	\draw[semithick,-] (0,-.3) to (0,.5);
\end{tikzpicture}\:\Bubble{e_r}
=\begin{tikzpicture}[anchorbase]
	\draw[semithick,-] (0,-.3) to (0,.5);
\node[rectangle,draw=black,inner sep=2.5,fill=cyclocolor] at (0,0.1) {$\scriptstyle r$};
\end{tikzpicture}-
\Bubble{e_r}\:
\begin{tikzpicture}[anchorbase]
	\draw[semithick,-] (0,-.3) to (0,.5);
\end{tikzpicture}\; .
$$ 

\vspace{2mm}
\noindent
\underline{Claim 4}: {\em We have that $\begin{tikzpicture}[anchorbase,scale=.9]
	\draw[semithick] (0.28,0) to[out=90,in=-90] (-0.28,.6);
	\draw[semithick] (-0.28,0) to[out=90,in=-90] (0.28,.6);
	\draw[semithick] (0.28,-.6) to[out=90,in=-90] (-0.28,0);
	\draw[semithick] (-0.28,-.6) to[out=90,in=-90] (0.28,0);
 \node[rectangle,draw=black,inner sep=2.5,fill=cyclocolor] at (.28,0) {$\scriptstyle r$};
\end{tikzpicture}
\in \cJ_l(B\star B,B\star B)$ for all $r \geq l$.}

\vspace{1mm}
\noindent
\underline{Proof}.
By \cref{purplesquare}, we have that
$$
\begin{tikzpicture}[anchorbase]
	\draw[semithick] (0.28,0) to[out=90,in=-90] (-0.28,.6);
	\draw[semithick] (-0.28,0) to[out=90,in=-90] (0.28,.6);
	\draw[semithick] (0.28,-.6) to[out=90,in=-90] (-0.28,0);
	\draw[semithick] (-0.28,-.6) to[out=90,in=-90] (0.28,0);
 \node[rectangle,draw=black,inner sep=2.5,fill=cyclocolor] at (.28,0) {$\scriptstyle r$};
\end{tikzpicture}
=
\begin{tikzpicture}[anchorbase,scale=.9]
	\draw[semithick] (0.28,0) to[out=90,in=-90] (-0.28,.6);
	\draw[semithick] (-0.28,0) to[out=90,in=-90] (0.28,.6);
	\draw[semithick] (0.28,-.6) to[out=90,in=-90] (-0.28,0);
	\draw[semithick] (-0.28,-.6) to[out=90,in=-90] (0.28,0);
 \node[rectangle,draw=black,inner sep=2.5,fill=cyclocolor] at (-.28,0) {$\scriptstyle r$};
\end{tikzpicture}
\ +\ 
\begin{tikzpicture}[anchorbase,scale=.9]
	\draw[semithick] (0.28,0) to[out=90,in=-90] (-0.28,.6);
	\draw[semithick] (-0.28,0) to[out=90,in=-90] (0.28,.6);
	\draw[semithick] (0.28,-.6) to[out=90,in=-90] (-0.28,0);
	\draw[semithick] (-0.28,-.6) to[out=90,in=-90] (0.28,0);
\end{tikzpicture}\ \Bubble{e_r}\ 
-\ 
\Bubble{e_r}\ 
\begin{tikzpicture}[anchorbase,scale=.9]
	\draw[semithick] (0.28,0) to[out=90,in=-90] (-0.28,.6);
	\draw[semithick] (-0.28,0) to[out=90,in=-90] (0.28,.6);
	\draw[semithick] (0.28,-.6) to[out=90,in=-90] (-0.28,0);
	\draw[semithick] (-0.28,-.6) to[out=90,in=-90] (0.28,0);
\end{tikzpicture}\ .
$$
The latter two terms on the right-hand side are zero by \cref{rels1}, thus the claim follows.

\vspace{2mm}
\noindent
\underline{Claim 5}: {\em We have that 
$\begin{tikzpicture}[anchorbase,scale=.9]
	\draw[semithick] (-0.4,-.5) to (0.4,.5);
	\draw[semithick] (0.4,-.5) to (-.4,.5);
 \node[rectangle,draw=black,inner sep=2.5,fill=cyclocolor] at (.2,-.24) {$\scriptstyle l$};
\end{tikzpicture} \in \cJ_l(B\star B, B\star B)$.}

\vspace{1mm}
\noindent
\underline{Proof}.
By the second relation from \cref{rels4b}, we have that
\begin{align*}
\begin{tikzpicture}[anchorbase,scale=1.2]
	\draw[semithick] (-0.4,-.6) to (0.4,.6);
	\draw[semithick] (0.4,-.6) to (-.4,.6);
 \node[rectangle,draw=black,inner sep=2.5,fill=cyclocolor] at (.2,-.3) {$\scriptstyle l$};
\end{tikzpicture}
&=
\begin{tikzpicture}[anchorbase,scale=1.2]
	\draw[semithick] (0.28,0) to[out=90,in=-90] (-0.28,.6);
	\draw[semithick] (-0.28,0) to[out=90,in=-90] (0.28,.6);
	\draw[semithick] (0.28,-.6) to[out=90,in=-90] (-0.28,0);
	\draw[semithick] (-0.28,-.6) to[out=90,in=-90] (0.28,0);
 \closeddot{-.16,-.4};
 \node[rectangle,draw=black,inner sep=2.5,fill=cyclocolor] at (.28,0) {$\scriptstyle l$};
\end{tikzpicture}
-
\begin{tikzpicture}[anchorbase,scale=1.2]
	\draw[semithick] (0.28,0) to[out=90,in=-90] (-0.28,.6);
	\draw[semithick] (-0.28,0) to[out=90,in=-90] (0.28,.6);
	\draw[semithick] (0.28,-.6) to[out=90,in=-90] (-0.28,0);
	\draw[semithick] (-0.28,-.6) to[out=90,in=-90] (0.28,0);
 \closeddot{.2,-.19};
 \node[rectangle,draw=black,inner sep=2.5,fill=cyclocolor] at (.22,0.1) {$\scriptstyle l$};
\end{tikzpicture}+
\begin{tikzpicture}[anchorbase,scale=1.2]
	\draw[semithick] (0.28,0) to[out=90,in=-90] (-0.28,.6);
	\draw[semithick] (-0.28,0) to[out=90,in=-90] (0.28,.6);
	\draw[semithick] (0,-.24) to[out=180,in=-90] (-0.28,0);
	\draw[semithick] (0,-.24) to[out=00,in=-90] (0.28,0);
 \node[rectangle,draw=black,inner sep=2.5,fill=cyclocolor] at (.28,0.03) {$\scriptstyle l$};
 \draw[semithick] (.28,-.6) to [out=90,in=90,looseness=1.6] (-.28,-.6);
\end{tikzpicture}\ .
\end{align*}
The first and third terms on the right hand side lie in $\cJ_l(B\star B, B\star B)$ by Claims 1 and 4. It remains to show that the second term on the right hand side belongs to $\cJ_l(B\star B, B\star B)$ too. We have
$$
\begin{tikzpicture}[anchorbase,scale=1.3]
	\draw[semithick] (0.28,0) to[out=90,in=-90] (-0.28,.6);
	\draw[semithick] (-0.28,0) to[out=90,in=-90] (0.28,.6);
	\draw[semithick] (0.28,-.6) to[out=90,in=-90] (-0.28,0);
	\draw[semithick] (-0.28,-.6) to[out=90,in=-90] (0.28,0);
 \closeddot{.2,-.19};
 \node[rectangle,draw=black,inner sep=2.5,fill=cyclocolor] at (.22,0.1) {$\scriptstyle l$};
\end{tikzpicture}=
\begin{tikzpicture}[anchorbase,scale=1.3]
	\draw[semithick] (0.28,0) to[out=90,in=-90] (-0.28,.6);
	\draw[semithick] (-0.28,0) to[out=90,in=-90] (0.28,.6);
	\draw[semithick] (0.28,-.6) to[out=90,in=-90] (-0.28,0);
	\draw[semithick] (-0.28,-.6) to[out=90,in=-90] (0.28,0);
 \node[rectangle,draw=black,inner sep=2.5,fill=cyclocolor] at (.24,0) {$\scriptstyle l+1$};
\end{tikzpicture}
-2\ 
\begin{tikzpicture}[anchorbase,scale=1.3]
	\draw[semithick] (0.28,0) to[out=90,in=-90] (-0.28,.6);
	\draw[semithick] (-0.28,0) to[out=90,in=-90] (0.28,.6);
	\draw[semithick] (0.28,-.6) to[out=90,in=-90] (-0.28,0);
	\draw[semithick] (-0.28,-.6) to[out=90,in=-90] (0.28,0);
 \node at (0,0) {$\Bubble{\!e_{l+1}\!}$};
\end{tikzpicture}
$$
by \cref{soundedlike}. Both of the terms on the right hand side here lie in $\cJ_l(B\star B,B\star B)$ by Claims 4 and 3.

\vspace{2mm}
\noindent
\underline{Claim 6}: {\em We have that $\begin{tikzpicture}[anchorbase,scale=.8]
	\draw[semithick] (0.2,-.5) to (0.2,.5);
	\draw[semithick] (-0.2,-.5) to (-.2,.5);
 \node[rectangle,draw=black,inner sep=2.5,fill=cyclocolor] at (.2,0) {$\scriptstyle l$};
\end{tikzpicture}
\in \cJ_l(B\star B,B\star B)$.}

\vspace{1mm}
\noindent
\underline{Proof}.
By the second relation from \cref{rels4}, 
we have that
\begin{align*}
\begin{tikzpicture}[anchorbase,scale=1.2]
	\draw[semithick] (0.2,-.5) to (0.2,.5);
	\draw[semithick] (-0.2,-.5) to (-.2,.5);
 \node[rectangle,draw=black,inner sep=2.5,fill=cyclocolor] at (.2,0) {$\scriptstyle l$};
\end{tikzpicture}
&=
\begin{tikzpicture}[anchorbase,scale=1.2]
	\draw[semithick] (-0.4,-.5) to (0.4,.5);
	\draw[semithick] (0.4,-.5) to (-.4,.5);
 \node[rectangle,draw=black,inner sep=2.5,fill=cyclocolor] at (.15,-.19) {$\scriptstyle l$};
 \closeddot{-.2,.26};
\end{tikzpicture}-\begin{tikzpicture}[anchorbase,scale=1.2]
	\draw[semithick] (-0.4,-.5) to (0.4,.5);
	\draw[semithick] (0.4,-.5) to (-.4,.5);
 \node[rectangle,draw=black,inner sep=2.5,fill=cyclocolor] at (.15,-.19) {$\scriptstyle l$};
 \closeddot{.34,-.42};
\end{tikzpicture}+
\begin{tikzpicture}[anchorbase,scale=1.2]
	\draw[semithick] (0.4,-.5) to[out=90, in=0] (0.1,-0.1);
	\draw[semithick] (0.1,-0.1) to[out = 180, in = 90] (-0.2,-.5);
 	\draw[semithick] (0.4,.5) to[out=-90, in=0] (0.1,0.1);
	\draw[semithick] (0.1,0.1) to[out = 180, in = -90] (-0.2,.5);
\node[rectangle,draw=black,inner sep=2.5,fill=cyclocolor] at (.32,-0.26) {$\scriptstyle l$};
\end{tikzpicture}\ .
\end{align*}
This lies in $\cJ_l(B\star B, B\star B)$ by Claims 5 and 1.

\vspace{2mm}
Finally, we can complete the proof by showing that $\cJ_l$ is a two-sided tensor ideal. 
Since it is a right tensor ideal by definition, it suffices to show that
$B \star f \in \cJ_l(B^{\star (n+1)}, 
B^{\star (m+1)})$ for any $f \in \cJ_l(B^{\star n}, B^{\star m})$.
The proof of this reduces to the case that $f$ is one of the generators of $\cJ_l$, in which case it follows 
by Claims 6 and 3.
\end{proof}

\begin{defin}\label{letters}
Assume that $l \equiv t \pmod{2}$.
The {\em cyclotomic nil-Brauer category} of {\em level $l$} is the $\Lambda$-linear strict graded monoidal category
that is the quotient $\cNB_t / \cI_l$ of the nil-Brauer category
by the four-sided tensor ideal $\cI_l$. 
We denote it by $\CNB_l$.
\end{defin}

\begin{lem}\label{lemmabelow}
For $n > l$, the morphism
$\begin{tikzpicture}[anchorbase]
\draw[ultra thick] (0,.25) to (0,-.25);
\cross{(0,0)};
\node at (0,-.35) {$\stringlabel{n}$};
\end{tikzpicture}$
is 0 in $\CNB_l$.
Hence, the image in $\CNB_l$ of the primitive homogeneous idempotent $\e_n$ from \cref{endef}
is 0 for $n > l$.
\end{lem}

\begin{proof}
It suffices to prove this in the case that $n=l+1$,
when it follows because we have
that
$$
\begin{tikzpicture}[anchorbase,scale=1.2]
\draw[ultra thick] (0,-.4) to (0,.4);
\cross{(0,0)};
\node at (0,-.5) {$\stringlabel{l+1}$};
\node at (0,.5) {$\stringlabel{\phantom{l+1}}$};
\end{tikzpicture}=
\begin{tikzpicture}[anchorbase,scale=1.2]
\draw[ultra thick] (0,-.4) to (0,.4);
\closeddot{0.01,0};
\node at (0.35,0) {$\scriptstyle \rho_{l+1}$};
\cross{(0,.25)};
\cross{(0,-.25)};
\node at (0,-.5) {$\stringlabel{l+1}$};
\node at (.2,.5) {$\stringlabel{\phantom{l+1}}$};
\end{tikzpicture}=\begin{tikzpicture}[anchorbase,scale=1.2]
\draw[ultra thick] (-.4,.4) to[looseness=1.88,out=-20,in=20] (-.4,-.4);
\draw (0,.4) to[looseness=1.88,out=-160,in=160] (0,-.4);
    \smallercross[-140]{(-.07,-.2)};
    \smallercross[140]{(-.07,.2)};
\closeddot{-.42,0};
\node at (-.55,0) {$\scriptstyle{l}$};
\closeddot{.01,0};
\node at (.22,0) {$\scriptstyle \rho_l$};
\node at (-.45,-.5) {$\stringlabel{l}$};
\node at (0,.5){$\stringlabel{\phantom{l}}$};
\end{tikzpicture}
=\begin{tikzpicture}[anchorbase,scale=1.2]
\draw[ultra thick] (-.4,.4) to[looseness=1.88,out=-20,in=20] (-.4,-.4);
\draw (0,.4) to[looseness=1.88,out=-160,in=160] (0,-.4);
    \smallercross{(.01,0)};
\closeddot{-.42,0};
\node at (-.55,0) {$\scriptstyle{l}$};
\node at (-.45,-.5) {$\stringlabel{l}$};
\node at (0,.5){$\stringlabel{\phantom{l}}$};
\end{tikzpicture}
= 
\sum_{r=0}^{l}
\begin{tikzpicture}[anchorbase,scale=1.2]
\draw[ultra thick] (-.4,.4) to[looseness=1.88,out=-20,in=20] (-.4,-.4);
\draw (0,.4) to[looseness=1.88,out=-160,in=160] (0,-.4);
    \smallercross{(.01,0)};
\closeddot{-.42,0};
\node at (-.55,0) {$\scriptstyle{r}$};
\node at (-.45,-.5) {$\stringlabel{l}$};
\node at (0,.5){$\stringlabel{\phantom{l}}$};
\node [draw,rounded corners,inner sep=2pt,black,fill=gray!20!white] at (-1,0) {$\scriptstyle e_{l-r}$};
\end{tikzpicture}= {\textstyle\frac{1}{2}}\ 
\begin{tikzpicture}[anchorbase,scale=1.2]
\draw[ultra thick] (-.4,.4) to[looseness=1.88,out=-20,in=20] (-.4,-.4);
\draw (0,.4) to[looseness=1.88,out=-160,in=160] (0,-.4);
    \smallercross{(.01,0)};
\node at (-.45,-.5) {$\stringlabel{l}$};
\node at (0,.5){$\stringlabel{\phantom{l}}$};
\node[rectangle,draw=black,inner sep=2.5,fill=cyclocolor] at (-.42,0) {$\scriptstyle l$};
\end{tikzpicture},
$$
which is 0 in $\CNB_l$.
The first and third equalities here follow by \cite[Cor.~4.4]{BWWiquantum}, the second follows from definitions and the braid relation, the fourth equality follows as each of the terms added to the summation is zero
by \cite[Lem.~4.2]{BWWiquantum}, and the final equality is \cref{cyclotomicrelationexplained}.
\end{proof}

\begin{rem}
The assumption that $l \equiv t \pmod{2}$ in \cref{letters} is natural from the point of view of categorification since, according to the construction in \cite{BW18QSP},
the iquantum group $\Ui$ is an inverse limit of the $\Z[q,q^{-1}]$-forms $\V(l)$ of the irreducible
$U_q(\sl_2)$-modules of highest weights $l \equiv t \pmod{2}$ (see the next subsection).
Another way to justify this assumption can be seen from the following
relation in $\cNB_t$ from \cite[Cor.~3.5]{BWWiquantum}:
 \begin{align}\label{dialog}
\begin{tikzpicture}[anchorbase]
	\draw[-] (0,-.3)  to (0,-.1) to [out=90,in=0] (-.3,.3) to [out=180,in=90] (-.5,.1) to[out=-90,in=180] (-.3,-.1) to [out=0,in=-90] (0,.3) to (0,.5);
  \closeddot{-.5,.1};\node at (-.9,0.1) {$\scriptstyle r+1$};
\end{tikzpicture}&=
\begin{dcases}
\begin{tikzpicture}[anchorbase]
	\draw[-] (0,-.3) to (0,.3);
 \closeddot{0,0};\node at (.2,0) {$\scriptstyle r$};
\end{tikzpicture}+{\textstyle\frac{1}{2}}\sum_{s=1}^{r}
\begin{tikzpicture}[anchorbase]
	\draw[-] (.5,-.3) to (.5,.3);
\closeddot{-0.2,0};
\closeddot{0.5,0};
\node at (0.7,0) {$\scriptstyle s$};
\node at (-0.1,0) {\Bubble{q_{r-s}}};
\end{tikzpicture}
&\text{if $r \equiv t \pmod{2}$}\\
-{\textstyle\frac{1}{2}}
\sum_{s=1}^{r}\begin{tikzpicture}[anchorbase]
	\draw[-] (.5,-.3) to (.5,.3);
\closeddot{-0.2,0};
\closeddot{0.5,0};
\node at (0.7,0) {$\scriptstyle s$};
\node at (-0.1,0) {\Bubble{q_{r-s}}};
\end{tikzpicture}&\text{if $r \not\equiv t \pmod{2}$.}
\end{dcases}
\end{align}
If $l \not \equiv t \pmod{2}$, we can attach a curl to the right of
$\sum_{r=0}^l
\Bubble{e_{l-r}}\ 
\begin{tikzpicture}[anchorbase]
\node at (.2,.1) {$\scriptstyle r$};
\closeddot{0,.1};
	\draw[semithick,-] (0,-.1) to (0,.3);
\end{tikzpicture} \in \cI_l(B,B)$ then
expand using \cref{dialog} to deduce that
$\sum_{r=0}^{l-1}
\Bubble{f_{l-1-r}}\ 
\begin{tikzpicture}[anchorbase]
\node at (.2,.1) {$\scriptstyle r$};
\closeddot{0,.1};
	\draw[semithick,-] (0,-.1) to (0,.3);
\end{tikzpicture} \in \cI_l(B,B)$
for some $f_0,\dots,f_{l-1} \in \Lambda$ with $f_0 = 1$.
Then the argument from the proof of \cref{lemmabelow}
can be mimicked to show that the images of the primitive idempotents
$\e_n$ in the quotient $\cNB_t / \cI_l$ are zero not only for $n > l$ but also for $n=l$.
\end{rem}

\subsection{The module \texorpdfstring{$\V(l)$}{} and its \texorpdfstring{$\imath$}{}-canonical basis}\label{ssec2-6}
We continue working with a fixed
$l \geq 0$ such that $l \equiv t\pmod{2}$.
Let $U_q(\sl_2)$ be the quantized enveloping algebra of $\sl_2$
over $\Q(q)$ with its standard generators $e, f, k^{\pm 1}$.
Let $\U$ be the $\Z[q,q^{-1}]$-form
generated by the divided powers
$e^{(n)} := e^n / [n]^!_q$ and $f^{(n)} := f^n / [n]^!_q$. Let $\eta_l$ be a highest weight vector of weight $l$, that is,
a vector such that $e \eta_l = 0$ and 
$k \eta_l = q^l \eta_l$.
Then $\V(l) := \U \cdot \eta_l$ is a $\Z[q,q^{-1}]$-form for the
$(l+1)$-dimensional
irreducible $U_q(\sl_2)$-module. We have that $f^{(n)} \eta_l = 0$ for $n > l$, and $\V(l)$ 
is free as a $\Z[q,q^{-1}]$-module with 
{\em standard basis} (also the {\em canonical basis})
$f^{(n)} \eta_l\:(0 \leq n \leq l)$.

There is a $\Z[q,q^{-1}]$-algebra anti-involution
$\rho:\U \rightarrow \U$ induced by the 
$\Q(q)$-algebra anti-involution of $U_q(\sl_2)$
which maps
$f \mapsto q k^{-1} e, e \mapsto q^{-1} f k,
k \mapsto k$.
It is easy to check that
\begin{align}\label{nnn}
f \cdot f^{(n)} \eta_l &= [n+1]_q f^{(n+1)}\eta_l,&
\rho(f) \cdot f^{(n)} \eta_l &= q^{2n-l-1}[l+1-n]_q
f^{(n-1)} \eta_l,
\end{align}
interpreting $f^{(-1)} \eta_l$ as $0$.
We find it helpful to visualize this as follows: 
\begin{equation}\label{mmm}
\begin{tikzcd}
\arrow[dddd,dashed, "{f}" left]\phantom{axyzw}& \arrow[d,bend right,"{[1]_{q}}" left]f^{(0)} \eta_l&\phantom{a}
\\
\phantom{a}
&  \arrow[u,bend right,"{q^{1-l}[\ell]_{q}}" right] f^{(1)} \eta_l \arrow[bend right,d,"{[2]_{q}}" left]&\phantom{a}
\\
\phantom{a}
&  \arrow[u,bend right,"{q^{3-l}[\ell-1]_{q}}" right]\vdots\arrow[d,bend right,"{[\ell-1]_{q}}" left]&\phantom{a}
\\
\phantom{a}
&\arrow[u,bend right,"{q^{l-3}[2]_{q}}" right]  f^{(l-1)} \eta_l\arrow[d,bend right,"{[\ell]_{q}}" left]&\phantom{a}
\\
\phantom{a}
&\arrow[u,bend right,"{q^{l-1}[1]_{q}}" right]   f^{(l)}\eta_l&\phantom{axywesdrzw}\arrow[uuuu,dashed,"\rho(f)" right]
\end{tikzcd}
\end{equation}
Recall finally that there is a unique symmetric bilinear form
$(\cdot,\cdot)_l:\V(l)\times \V(l)$
such that $(\eta_l,\eta_l)_l = 1$
and $(u v_1, v_2)_l = (v_1, \rho(u) v_2)_l$
for all $v_1,v_2 \in \V(l)$
and $u \in \U$.
The standard basis is orthogonal with
\begin{equation}\label{bilform}
\overline{\big(f^{(n)} \eta_l, f^{(n)} \eta_l\big)_l}
= q^{n(l-n)}\qbinom{l}{n}_q \in 1 + q^2\N[q^2],
\end{equation}
which is the Poincar\'e polynomial for the cohomology of the Grassmannian $\operatorname{Gr}_{n,l}$.

The iquantum group $\Ui$
is naturally a $\Z[q,q^{-1}]$-subalgebra of $\U$,
with $b \in \Ui$ being the element
$f+\rho(f) \in \U$.
Of course, the $\U$-module $\V(l)$ can also be viewed as a $\Ui$-module by restriction.
By \cref{nnn}, we have that
$b \eta_l = f \eta_l$ and
\begin{equation}\label{rakes}
b \cdot f^{(n)} \eta_l = 
[n+1]_q f^{(n+1)} \eta_l
+ q^{2n-l-1} [l-n+1]_q
f^{(n-1)} \eta_l
\end{equation}
for $n \geq 1$.
The icanonical basis of $\Ui$ 
descends to the {\em 
icanonical basis}
$b^{(n)}\eta_l\:(0 \leq n \leq l)$ 
of $\V(l)$, and $b^{(n)} \eta_l = 0$ for $n > l$. These statements all follow from 
\cite{BW18QSP, BW18KL}.

\begin{rem}\label{oldpeculiar}
From the preceding discussion,
it follows that the kernel $\I_l$ of
the $\Ui$-module
homomorphism $\Ui \twoheadrightarrow \V(l),
u \mapsto u \eta_l$
is generated by
the icanonical basis elements 
$b^{(n)}$ for $n > l$.
%In fact, by \cref{ryder}, this ideal is generated already just by $b^{(l+1)}$.
The quotient $\Ui / \I_l$ is isomorphic to 
$\V(l)$ as a $\Ui$-module. In view of the commutativity of $\Ui$, this quotient has the additional structure of a commutative $\Z[q,q^{-1}]$-algebra.
\end{rem}

The icanonical basis of $\V(l)$
has an alternative characterization; in fact, this is the way the icanonical basis gets introduced in the first place in \cite{BW18QSP, BW18KL}. To explain this, we
need the {\em ibar involution}
$\V(l) \rightarrow \V(l), v \mapsto \overline{v}$. This is the unique anti-linear 
map fixing $\eta_l$ 
such that $\overline{bv} = b \overline{v}$ for all $v \in \V(l)$.
This is {\em not} the usual bar involution on $\V(l)$ as it is definitely not true that $\overline{fv} = f \overline{v}$.
Then, for $0 \leq n \leq l$,
$b^{(n)} \eta_l \in \V(l)$ is the unique bar-invariant vector
such that \begin{equation}\label{characterize}
b^{(n)} \eta_l = f^{(n)} \eta_l
+\text{(a $q^{-1} \Z[q^{-1}]$-linear combination of
the vectors $f^{(i)} \eta_l$ for $i < n$)}.
\end{equation}
There is also an explicit formula expressing the icanonical basis in terms of the standard bases:

\begin{theo}[{\cite[(2.16),(2.17),(3.8),(3.9)]{BeW18}}]\label{crazy}
We have that
$$
b^{(n)} \eta_l = 
\begin{cases}
\displaystyle\sum_{i=0}^{\lfloor\frac{n}{2}\rfloor}
q^{-i(l+2i-n-1)}\qbinom{
(l+2i-n)/2
}{i}_{q^2} f^{(n-2i)} \eta_l
&\text{if $0 \leq n \leq l$ and $n \equiv t \pmod{2}$,}\\
\displaystyle\sum_{i=0}^{\lfloor\frac{n}{2}\rfloor}
q^{-i(l+2i-n)}\qbinom{
(l+2i-n-1)/2
}{i}_{q^2} f^{(n-2i)} \eta_l
&\text{if $0 \leq n \leq l$ and $n \not\equiv t \pmod{2}$}
\\
0&\text{if $n > l$}.
\end{cases}
$$
\end{theo}

\begin{proof}
This follows from the four formulae cited from
\cite{BeW18}, but it is some work to convert these into the more concise form recorded here. Alternatively, one can give an independent proof by using \cref{nnn} to verify that our reformulation
satisfies the recurrence relation \cref{ryder}.
\end{proof}

\begin{cor}\label{cute1}
For $0 \leq n \leq \lfloor\frac{l}{2}\rfloor$, we have that
$\big(\eta_l, b^{(2n)}\eta_l\big)_l = 
q^{-n(l+t-1)} \qbinom{(l-t)/2}{n}_{q^2}$.
\end{cor}

Here is one more observation about the bilinear form on $\V(l)$.

\begin{lem}\label{cute2}
The top degree term of the Laurent polynomial
$\big(b^n \eta_l, b^n \eta_l)_l$
is $q^{2\binom{n}{2}}$
for $0 \leq n \leq l$, and it is
$q^{(2n-l)(l-1)}$ for $n \geq l$.
\end{lem}

\begin{proof}
Since $\big(b^n \eta_l, b^n \eta_l\big)_l
= \big(\eta_l, b^{2n} \eta_l\big)_l$,
we need to compute the $\eta_l$-coefficient of $$
b^{2n} \eta_l = (f + \rho(f))^{2n} \eta_l.
$$
Such coefficients arise from the monomials
in the expansion of
$(f + \rho(f))^{2n}$
that have $f$ appearing $n$ times
and $\rho(f)$ appearing $n$ times.
To produce a non-zero coefficient, 
the monomials also need to be Dyck words in the sense that at each position there are not more letters 
$\rho(f)$ in this position or to the right than 
there are letters $f$.
 The coefficient of $\eta_l$ arising when each of these monomials is applied to $\eta_l$
 is easily computed by contemplating \cref{mmm}.
 The Laurent polynomials on the edges of this diagram are monic of degrees
 $0,1,\dots,l-1$ from top to bottom. It follows that the biggest degree arises from a unique monomial,
 namely, $\rho(f)^n f^n$ if $n \leq l$
 or $\rho(f)^l (f \rho(f))^{n-l} f^l$ if $n \geq l$.
 Since the Laurent polynomials are monic,
 this produces the top degree term in the statement of the lemma.
\end{proof}

\subsection{The special cases \texorpdfstring{$l=0$ and $l=1$}{}}
The definition of $\CNB_l$ makes sense when $l=0$ and $l=1$, but these are degenerate cases and we treat them separately here. In the case $l=0$,
$\CNB_0$ is the trivial graded monoidal category with one object $\one$ and $\End_{\CNB_0}(\one) = \kk$. Its graded Karoubi envelope is monoidally equivalent to $\gVec_{\operatorname{fd}}$, with $K_0(\gKar(\CNB_0)) \cong \V(0)$ so that $[\one]$ is identified with $\eta_0$. There is nothing more to be said about this case.

If $l=1$ then $\CNB_l$ is slightly more interesting, but still easy to understand directly. It can be presented as a strict graded monoidal category with one generating object $B$ and three generating morphisms
\begin{align*}
\Bubble{e_1}&:\one \rightarrow \one,&
\begin{tikzpicture}[anchorbase]
	\draw[semithick] (0.4,0) to[out=90, in=0] (0.1,0.4);
	\draw[semithick] (0.1,0.4) to[out = 180, in = 90] (-0.2,0);
\end{tikzpicture}
&:B \star B\rightarrow\one,&
\begin{tikzpicture}[anchorbase]
	\draw[semithick] (0.4,0.4) to[out=-90, in=0] (0.1,0);
	\draw[semithick] (0.1,0) to[out = 180, in = -90] (-0.2,0.4);
\end{tikzpicture}
&:\one\rightarrow B\star B,
\\\notag
&\text{(degree $2$)}&
&\text{(degree $0$)}&
&\text{(degree $0$)}
\end{align*}
subject to the following relations:
\begin{align}
\label{stupidrels2}
\begin{tikzpicture}[baseline=-2.5mm]
\draw[semithick] (0,-.15) circle (.3);
\end{tikzpicture}
&= \id_\one,&
\begin{tikzpicture}[anchorbase,scale=1.2]]
	\draw[semithick,-] (0,-.3) to (0,.5);
\end{tikzpicture}\:\:\Bubble{e_1}\ &=\ 
-\:\displaystyle\Bubble{e_1}\:\:
\begin{tikzpicture}[anchorbase,scale=1.2]
	\draw[semithick,-] (0,-.3) to (0,.5);
\end{tikzpicture}
\ ,&
\begin{tikzpicture}[anchorbase,scale=1.2]
\draw[semithick] (-.5,-.2) to (-.5,.4);
 	\draw[semithick] (-0.2,.4) to[out=-90,in=180,looseness=1] (0,0) to[out=0,in=-90,looseness=1] (0.2,.4);
\end{tikzpicture}\ 
&=\ 
\begin{tikzpicture}[anchorbase,scale=1.2]
\draw[semithick] (.5,-.2) to (.5,.4);
 	\draw[semithick] (-0.2,.4) to[out=-90,in=180,looseness=1] (0,0) to[out=0,in=-90,looseness=1] (0.2,.4);
\end{tikzpicture}\ ,&
\begin{tikzpicture}[anchorbase,scale=1.2]
\draw[semithick] (-.5,-.4) to (-.5,.2);
 	\draw[semithick] (-0.2,-.4) to[out=90,in=180,looseness=1] (0,0) to[out=0,in=90,looseness=1] (0.2,-.4);
\end{tikzpicture}\ 
&=\ 
\begin{tikzpicture}[anchorbase,scale=1.2]
\draw[semithick] (.5,-.4) to (.5,.2);
 	\draw[semithick] (-0.2,-.4) to[out=90,in=180,looseness=1] (0,0) to[out=0,in=90,looseness=1] (0.2,-.4);
\end{tikzpicture}\ ,&
\begin{tikzpicture}[anchorbase,scale=1.2]
 	\draw[semithick] (0,-.3) to (0,.3);
	\draw[semithick] (-.3,-.3) to (-0.3,.3);
\end{tikzpicture}
&=
\begin{tikzpicture}[anchorbase,scale=1.2]
 	\draw[semithick] (-0.18,-.3) to[out=90,in=90,looseness=2.3] (0.18,-.3);
 	\draw[semithick] (-0.18,.3) to[out=-90,in=-90,looseness=2.3] (0.18,.3);
\end{tikzpicture}\ .
\end{align}
This follows quite easily from the defining relations above; the dot is redundant as a generator 
since $\begin{tikzpicture}[anchorbase,scale=.8] \draw[semithick](0,.3) to (0,-.3);\closeddot{0,0};\end{tikzpicture}
=\begin{tikzpicture}[anchorbase,scale=.8] \draw[semithick] (0,.3) to (0,-.3);\end{tikzpicture}\ \ \Bubble{e_1}\ $
in $\End_{\CNB_1}(B)$ by \cref{cyclotomicrelationexplained}, 
$\ \Bubble{e_r} = 0\ $ in $\End_{\CNB_1}(\one)$ for $r > 1$
by \cref{bensdumbideatocallitfoursided},
and $\begin{tikzpicture}[anchorbase,scale=.6]
\draw[semithick] (0.28,-.3) to (-0.28,.4);
	\draw[semithick] (-0.28,-.3) to (0.28,.4);
\end{tikzpicture}= 0$ 
in $\End_{\CNB_1}(B\star B)$. To prove the final assertion here about the vanishing of the crossing,
attach a crossing to the bottom of the second relation from \cref{rels4}, then use other relations to check that all of the terms except for the crossing are 0 in $\CNB_1$, hence, the crossing must be 0 too.

\begin{theo}
View 
the morphism space
$\Hom_{\CNB_1}(B^{\star n}, B^{\star m})$
as a right $\kk[x]$-module so that $x$ acts by horizontal composition on the right by $\Bubble{e_1}\ $.
Then this morphism space
is $\{0\}$ unless $m \equiv n \pmod 2$, in which case it is free of rank one as a right $\kk[x]$-module with basis
given by the string diagram obtained by drawing
$\lfloor m/2\rfloor$ side-by-side caps at the top and $\lfloor n/2\rfloor$ side-by-side cups at the bottom,
plus a vertical string on the right hand side if $m$ and $n$ are odd.
\end{theo}

\begin{proof}
It is an exercise using the relations 
to see that $\Hom_{\CNB_1}(B^{\star n}, B^{\star m})$ is $\{0\}$ if $m \not \equiv n\pmod{2}$, and that it is spanned as a $\Z[x,x^{-1}]$-module by the morphism just described otherwise. 
To prove the freeness, we let $\cC$ be the graded category with
objects $\N$ and $\Hom_{\cC}(n,m) := \kk[x]$ if $m \equiv n \pmod{2}$ or $\{0\}$ otherwise, with the obvious composition law defined by multiplication of polynomials.
We make this into a strict graded monoidal category by defining $m \star n := m+n$, and defining the tensor product of morphisms
$f(x):m \rightarrow m', g(x):n \rightarrow n'$
by setting $f(x) \star g(x) := f((-1)^n x) g(x):
m+n \rightarrow m'+n'$.
It is trivial to check from the defining relations that there is a graded monoidal functor
$\CNB_1 \rightarrow \cC$ taking $B \mapsto 1$,
$\Bubble{e_1} \mapsto (x:0 \rightarrow 0)$, and the cap and the cup to the
morphisms $2 \rightarrow 0$ and $0 \rightarrow 2$ 
represented by $1 \in \kk[x]$.
The desired freeness is now obvious.
\end{proof}

From this, it follows that $B^{\star n} \cong \one$ if $n$ is even, $B^{\star n} \cong B$ if $n$ is odd, $\End_{\CNB_1}(\one) \cong \End_{\CNB_1}(B)
\cong \kk[x]$, and $\Hom_{\CNB_1}(B,\one)
= \Hom_{\CNB_1}(\one,B) = \{0\}$. Hence, there is a $\Z[q,q^{-1}]$-module isomorphism
\begin{align}
K_0(\gKar(\CNB_1)) &\stackrel{\sim}{\rightarrow}
\V(1),&
[\one] \mapsto \eta_1,
\quad
[B] \mapsto b \eta_1.
\end{align}
Since $B\star B \cong \one$, we have that $[B]^2 = [\one]$, so as an algebra 
$K_0(\gKar(\CNB_1))
\cong \Z[q,q^{-1}][b] / (b^2-1)$.

%% file: s3-ssbim.tex
\setcounter{section}{2}
\section{Reminders about singular Soergel bimodules}\label{secssbim}

In this expository section, we review some fundamental results about singular Soergel bimodules and the graphical calculus for them following \cite{Will,ESW}.
Our general setup is essentially 
the same as in \cite[$\S$3.1]{Will}
with additional assumptions imposed so that we can appeal to the diagrammatic calculus
of \cite{ESW}; see also \cite[Ch.~24]{EMTW}.
Later, we will appeal to these results only 
for the reflection realization in finite type D, for which it is well known that all of the assumptions hold.

\subsection{Coxeter group realization}
\label{ssec3-1}
Let $(W,S)$ be a Coxeter system
with length function $\ell:W \rightarrow \N$
and Bruhat order $\leq$. 
We
assume simple reflections are parametrized by another set\footnote{This is the second author's pickiness. It is perfectly reasonable to index everything by the set $S$ but that leads to expressions like $\alpha_{s_i}$ not $\alpha_i$ to denote the simple root associated to $s_i \in S$.} $N$ (``nodes of the Coxeter diagram") via a given bijection
$N \rightarrow S, i \mapsto s_i$. 
%In fact, we excepting one use in a bullet point below, we no longer use the set $S$, to avoid conflict with the $q$-Schur algebra $\Schur$ to be defined later.

We fix a realization in the sense of \cite[Def.~3.1]{soergelcalculus}. So 
we are given a finite-dimensional
vector space $\mathfrak{h}$ over the ground field $\kk$ ($\operatorname{char} \kk\neq 2$) and
subsets 
$\{\alpha_i\:|\:i \in N\}\subset\mathfrak{h}^*$ and
$\{\alpha_i^\vee\:|\:i \in N\}\subset\mathfrak{h}$
such that
\begin{itemize}
\item $\langle \alpha_i, \alpha_i^\vee \rangle = 2$ for all $i \in N$;
\item
there is a well-defined representation $W \rightarrow 
GL(\mathfrak{h}^*)$ taking $s_i$ to
the map 
$\lambda \mapsto \lambda - \langle \lambda,\alpha_i^\vee\rangle \alpha_i$.
\end{itemize}
We assume in addition that the realization is
\begin{itemize}
\item {\em balanced} in the sense of \cite[Def.~3.6]{soergelcalculus};
\item
{\em finitary Demazure surjective} as in \cite[Sec.~24.3.6]{EMTW} (called {\em generalized Demazure surjective} in \cite[Def.~3.5]{EKLPdemazure}), meaning that
the longest Demazure operator
$\partial_I:R \rightarrow R^I$
is surjective for each $I \finsub N$ (the notation here will be explained shortly);
\item {\em faithful} and {\em reflection faithful} as in \cite[Sec.~4.1]{Will}.
This means for any $w \in W$ that
the fixed point space $\{v \in \mathfrak{h}^*\:|\:w(v) = v\}$
is $\mathfrak{h}^*$ if and only if $w=1$, and it
is of codimension 1 if and only if $w$ is
conjugate to a simple reflection. 
\end{itemize}
Using further language yet to be introduced, these assumptions are required to ensure
\begin{itemize}
\item
the Demazure operators $\partial_i\:(i \in N)$
satisfy the braid relations, and
the sets $\Phi_I^+$ of positive roots introduced below are of the claimed size;
\item
the upgraded Chevalley theorem as formulated in
\cite[Th.~24.36, Th.~24.40]{EMTW} 
(and proved again in \cite[Sec.~4]{EKLPdemazure})
 holds, i.e., we have available the appropriate squares of Frobenius extensions (see \cref{ssec3-2,ssec3-3});
\item
the hypotheses of \cite{Will} are satisfied (see \cref{ssec3-4} below where the important results deduced from these hypotheses discussed further).
\end{itemize}
For further justification of the assumptions, see \cite[Sec.~24.3.6]{EMTW} and \cite[Sec.~3]{EKLPdemazure}.

Let $R := \kk[\mathfrak{h}]$ be  the symmetric algebra on the dual space $\mathfrak{h}^*$, viewed as a graded algebra with $\mathfrak{h}^*$ in degree 2. The Coxeter group acts naturally 
on $R$ by graded algebra automorphisms.
For $i \in N$, the {\em Demazure operator} $\partial_i:R
\rightarrow R$ is the degree $-2$ linear map defined by
$$
\partial_i(f) := \frac{f - s_i(f)}{\alpha_i}.
$$
The endomorphisms $\partial_i\:(i \in N)$ generate a copy of the
{\em nil-Coxeter algebra} associated to the Coxeter group $W$. 
They satisfy the same braid relations as the simple reflections in the Coxeter group, but the quadratic relation $s_i^2 = 1$
is replaced by $\partial_i^2 = 0$.

For $I \subseteq N$, let $W_I$
be the parabolic subgroup $\langle s_i\:|\:i \in I \rangle$ of $W$, and
$R^I$ be the subalgebra of $R$ consisting of its fixed points.
We say that $I$ is {\em finitary} if $W_I$ is a finite group, and use the notation $I \finsub N$
to indicate that $I$ is such a subset. 
The letters $I, J, K$ will be reserved for finitary subsets.
For $I \finsub N$,
let $w_I$ be the longest element of $W_I$,
and $\partial_I$ be the product of Demazure operators corresponding to a reduced expression for $w_I$.
Let
\begin{align}\label{usualpoincare}
\pi^+_I &:= \sum_{w \in W_I}
q^{2\ell(w)},&
\pi_I &:= \sum_{w \in W_I} q^{2\ell(w)-\ell(w_I)}=
q^{-\ell(w_I)} \pi^+_I.
\end{align}
The first of these is the Poincar\'e polynomial of $W_I$, and the second is its bar-invariant renormalization.
Let $\Phi_I$ be the finite set 
$\{w (\alpha_i)\:|\:w \in W_I, i \in I\}$,
and 
\begin{equation}
\Phi_I^+ := \{w(\alpha_i)\:|\:w \in W_I, i \in I\text{ such that }\ell(ws_i) = \ell(w)+1\}.
\end{equation}
The assumption that the realization is balanced implies\footnote{This is proven for dihedral groups in \cite[Appendix]{Etwocolor}. The general case is not in the literature to our knowledge. The result will be obvious in our applications, since we make use of realizations associated to typical root data.} that $|\Phi^+_I| = \ell(w_I)$.
Finally, let
\begin{equation}
\mu_I := \prod_{\alpha \in \Phi^+_I} \alpha \; \in R.
\end{equation}
This polynomial is $W_I$-anti-invariant 
(meaning $s_i \mu_I = -\mu_I$ for $i \in I$) and has degree $2 \ell(w_I)$.

\subsection{Chains of Frobenius extensions}\label{ssec3-2}
Suppose that we are given $J\subseteq I \finsub N$.
We have that $W_J \subseteq W_I$, hence,
$R^I$ is a subalgebra of $R^J$.
We denote the set of minimal length $W_I / W_J$-coset representatives by
$(W_I / W_J)_{\min}$.
Let 
\begin{equation}\label{notstandard}
\eta_I^J := \prod_{\substack{J \subseteq K \subseteq I\\|K| \equiv |I|\pmod{2}}}\!\!\!\! \mu_K
\:\Bigg/\:
\prod_{\substack{J \subseteq K \subseteq I\\ |K| \not\equiv |I|\pmod{2}}}\!\!\!\! \mu_K.
\end{equation}
The inclusion-exclusion principle implies that this is a polynomial in $R$
(rather than a rational function). In fact, since each $\mu_K$ is $W_J$-anti-invariant, 
$\eta_I^J$ lies in $R^J$.
For example:
\begin{itemize}
\item
writing $Ii$ 
for $I \cup \{i\}$ for $i \in N - I$ and assuming this set is finitary,
 we have that $\eta^I_{Ii} = \frac{\mu_{Ii}}{\mu_I} \in R^{I}$;
\item
writing $Iij$ for $I \cup \{i,j\}$
for distinct $i,j \in N - I$ and assuming it is finitary, we have that $\eta^I_{Iij} = \frac{\mu_{Iij}\, \mu_I}{\mu_{Ii}\, \mu_{Ij}}\in R^I$;
\item
writing $Iijk$ for $I \cup \{i,j,k\}$
for distinct $i,j,k \in N - I$ and assuming it is finitary, we have that $\eta^I_{Iijk} = \frac{\mu_{Iij}\, \mu_{Iik}\,\mu_{Ijk}\, \mu_I}{\mu_{Iijk} \,\mu_{Ii} \,\mu_{Ij} \,\mu_{Ik}}\in R^I$.
\end{itemize}

Continuing  with $J \subseteq I\finsub N$,
a classical result of Demazure \cite{Demazure}
(see also \cite[Th.~4.3]{EKLPdemazure})
shows that $R^J$ is a symmetric graded Frobenius extension of $R^I$.
The
(unique up to a scalar) 
Frobenius trace is the $(R^I, R^I)$-bimodule homomorphism
of degree $2(\ell(w_J)-\ell(w_I))$
\begin{align}\label{tracedef}
\tr_I^J:R^J &\rightarrow R^I, &f &\mapsto \partial_I^J(f),
\end{align}
where $\partial_I^J$ denotes the element of 
the nil-Coxeter algebra defined by composing Demazure operators in the same order as
a reduced expression for the longest element\footnote{Although $w_J^{-1}=w_J$, we write $w_Iw_J^{-1}$ as a reminder that $\ell(w_Iw_J^{-1}) = \ell(w_I)-\ell(w_J)$.
} 
$w_I w_J^{-1}$ of $(W_I / W_J)_{\min}$. 
The comultiplication
is the $(R^J, R^J)$-bimodule homomorphism
of degree $2(\ell(w_I)-\ell(w_J))$
\begin{align}\label{comultdef}
\Delta_I^J&:R^J \rightarrow R^J \otimes_{R^I} R^J,
&1&\mapsto
\sum_{b \in
\B_I^J} b \otimes b^\vee,
\end{align}
where 
$\B_I^J$ is a homogeneous basis for $R^J$ as a free $R^I$-module and $b^\vee$ is the dual basis element defined by
$\tr_I^J(a b^\vee) = \delta_{a,b}\:(a \in \B_I^J)$.

The following observation about chains of Frobenius extensions appears in \cite[Sec.~2.2]{ESW}:

\begin{lem}\label{L:transitive-trace-and-bases}
Let $K \subseteq J \subseteq I \finsub N$. Then $\tr^K_I=\tr^J_I\circ \tr^K_J$. Moreover, $\{ab \:|\: a\in \B^K_J, b\in \B^J_I\}$ is a basis for $R^{K}$ as a free graded $R^{I}$-module, with dual basis $\{a^{\vee}b^{\vee} \:|\: a\in \B^K_J, b\in \B^J_I\}$. 
\end{lem}

\begin{proof}
The equality of Demazure operators $\partial^K_I=\partial^J_I\circ \partial^K_J$ implies that $\tr^K_I=\tr^J_I\circ \tr^K_J$. Using the fact that $\tr^K_J$ is $R^J$-linear, we compute
\[
\tr^K_I((ab)\cdot (a')^{\vee}(b')^{\vee}) = \tr^J_I\circ \tr^K_J(ab(a')^{\vee}(b')^{\vee}) = \tr^J_I\left(b(b')^{\vee} \tr^K_J(a(a')^{\vee})\right) = \delta_{a,a'}\delta_{b,b'}.
\]
Since $\B^K_J$ spans $R^K$ over $R^J$ and $\B^J_I$ spans $R^J$ over $R^I$, it follows that $\{ab \ | \ a\in \B^K_J, y\in \B^J_I\}$ spans $R^K$ over $R^I$. Linear independence is a consequence of the existence of a dual basis.
\end{proof}

Using \cref{L:transitive-trace-and-bases}, the following useful identities are established in 
\cite[(2.7)]{ESW}: for $K \subseteq J \subseteq I \finsub N$ and $f \in R^K$, we have that
\begin{align}\label{deadhead1}
\sum_{b \in \B_I^K} b \otimes \tr^K_J(f b^\vee)
&=\sum_{a \in \B_I^J} af \otimes a^\vee \in R^K \otimes_{R^I} R^J,\\\label{deadhead2}
\sum_{b \in \B_I^K} \tr^K_J(bf) \otimes b^\vee
&=\sum_{a \in \B_I^J} a \otimes f a^\vee \in R^J \otimes_{R^I} R^K.
\end{align}
In particular, taking $K=J$, this gives that
\begin{equation}\label{deadhead3}
\sum_{b \in \B_I^J} b f \otimes b^\vee =
\sum_{b \in \B_I^J} b \otimes  f b^\vee \in R^J \otimes_{R^I} R^J
\end{equation}
for all $f \in R^J$.
Also \cite[(2.8)]{ESW} gives
\begin{align}\label{deadhead4}
\sum_{b \in \B^K_I} b \tr^K_J(b^\vee) &=\mu_I / \mu_J \in R^J.
\end{align}
Taking $K = J$, this implies that
$\sum_{b \in \B^J_I} b b^\vee = \mu_I / \mu_J$. Another special case recovers the counit axiom
\begin{equation}\label{counit}
\sum_{b \in \B_I^J} b \tr_I^J(b^\vee) = 1.
\end{equation}

\subsection{Diagrammatics for singular Bott-Samelson bimodules}\label{ssec3-3}
Suppose that $J \subseteq I \finsub N$,
so that $R^I \subseteq R^J$.
Working with the usual graded 
categories $R^I\gMod$ and $R^J\gMod$ 
of graded left modules,
we have 
the restriction and induction functors
\begin{align*}
\Res_{I}^{J}&:R^J\gMod \rightarrow R^I\gMod,
&\Ind_{I}^{J}&:R^I\gMod\rightarrow R^J\gMod,
\end{align*}
which are defined by tensoring on the left 
with $R^J$ viewed as an $(R^I, R^J)$-bimodule
or as an $(R^J, R^I)$-bimodule, respectively.
They form an adjoint pair
$(\Ind_{I}^{J}, \Res_{I}^{J})$
via the canonical (degree 0) adjunction between tensor and hom. 
Since 
$R^I \subseteq R^J$ is a graded Frobenius extension, 
there is also an adjunction the other way around, after a suitable grading shift.
To make the grading shifts on the various units and counits of adjunction more balanced, we incorporate a grading shift into
{\em restriction}\footnote{This convention is consistent with the choice of duality made in \cref{rightduality} below; see also \cref{drivingmecrazy}.},
working henceforth with the functors
$q^{\ell(w_J)-\ell(w_I)} \Res^J_I$ and
$\Ind_I^J$.
We have that
\begin{align}\label{appearedalready}
q^{\ell(w_J)-\ell(w_I)} \Res^J_I&=B_{I,J} \otimes_{R^J} -,
& \Ind_I^J&=B_{J,I} \otimes_{R^I} -,
\end{align}
where $B_{I,J}$ denotes the graded
$(R^I, R^J)$-bimodule $q^{\ell(w_J)-\ell(w_I)} R^J$,
and $B_{J,I}$ denotes the graded $(R^J, R^I)$-bimodule $R^J$.
The two adjunctions give natural 
degree-preserving isomorphisms
\begin{align}\label{adj1}
\Hom_{R^J\dash}(B_{J,I} \otimes_{R^I} M, M')
&\cong
q^{\ell(w_I)-\ell(w_J)}\Hom_{R^I\dash}(M, B_{I,J} \otimes_{R^J} M'),\\\label{adj2}
\Hom_{R^I\dash}(B_{I,J} \otimes_{R^J} M', M)&\cong
q^{\ell(w_J)-\ell(w_I)}\Hom_{R^J\dash} (M', B_{J,I} \otimes_{R^I} M)
\end{align}
for any graded left $R^I$-module $M$ and $R^J$-module $M'$.

More generally, for $I, J \finsub N$
and $K \subseteq I \cap J$,
we define the
graded $(R^I, R^J)$-bimodule
\begin{align}
\label{inductionrestrictionbimodule}
B_{I,J}^K &:= 
q^{\ell(w_K)-\ell(w_I)}
R^K.
\end{align}
This is naturally isomorphic to
$B_{I,K} \otimes_{R^K} B_{K,J}$.
We denote $B_{I,J}^{I \cap J}$ simply by $B_{I,J}$.
When $I \subseteq J$ or $I \supseteq J$, this notation is consistent with the notation for induction/restriction bimodules that appeared in \cref{appearedalready}.

\begin{rem}
The general results explained in the next subsection (specifically, \cref{bscorollary} in the light of \cref{homformula,class,catthm}) imply that $B_{I,J}$ is an indecomposable bimodule for any $I, J \finsub N$. 
Moreover, for $K \subseteq I \cap J$,
we have that
$B_{I,J}^K \cong \frac{\pi_{I \cap J}}{\pi_K} B_{I,J}$
as a graded $(R^I, R^J)$-bimodule.
\end{rem}

Recall that $\gBim$ is the graded bicategory of graded algebras, 
graded bimodules, and graded bimodule homomorphisms.

\begin{defin}[{\cite[Def.~7.1]{Will}}]\label{newssbimdef}
The category $\BSBim$ 
of {\em singular Bott-Samelson bimodules} is the
full\footnote{We mean that it has the same 2-morphism spaces.}
sub-bicategory of $\gBim$
with objects $R^I$ for $I \finsub N$, and
1-morphisms $\Hom_{\BSBim}(R^J,R^I)$
given by the graded $(R^I, R^J)$-bimodules
\begin{align}\label{abs}
B_{I_0, I_1}^{K_1} \otimes_{R^{I_1}}
B_{I_1, I_2}^{K_2} \otimes_{R^{I_2}}\cdots
\otimes_{R^{I_{n-1}}} B_{I_{n-1},I_n}^{K_n},
\end{align}
for all $n \geq 1$ and
$I=I_0 \finsup K_1 \finsub I_1 \finsup \cdots
\finsub I_{n-1} \finsup K_n \finsub I_n=J$.

The graded bicategory $\SBim$ 
of {\em singular Soergel bimodules} is the graded Karoubian closure of $\BSBim$ in $\gBim$, that is, it is the full sub-bicategory of $\gBim$ with the same objects as $\BSBim$, but the 1-morphisms 
in $\Hom_{\SBim}(R^J, R^I)$ 
are the graded $(R^I, R^J)$-bimodules which are 
isomorphic to summands of finite direct sums of grading shifts of the bimodules \cref{abs}.

We will often denote the object $R^I$ of $\BSBim$ (resp., $\SBim$) simply by $I$, so that the set of objects
is identified with the set of finitary subsets of $N$.
The identity 1-endomorphism $\one_I$ of the object $I$ is the 
regular bimodule $R^I$. We use the notation
$\one_I \BSBim \one_J$ (resp., $\one_I \SBim \one_J$)
to denote the morphism category
$\HOM_{\BSBim}(J,I)$ (resp., $\HOM_{\SBim}(J,I)$).
\end{defin}

Now we introduce the string calculus for 
%(the strictification of) 
$\BSBim$
following \cite{ESW}.
By transitivity of induction and restriction,
to generate $\BSBim$,
it suffices to consider the bimodules
$B_{I,Ii}$
and
$B_{Ii,I}$
for $I \subset N$ and $i \in N - I$
such that $Ii$ is finitary.
We use the string diagrams
${\scriptstyle I}\up\scriptstyle{Ii}$
and ${\scriptstyle Ii}\down\scriptstyle{I}$
 to denote the identity endomorphisms of these bimodules.
For extra clarity, the strings can also be colored by the color $i$ but this is not essential since the string color is determined by the labels of the adjacent 2-cells.
We use the following rightward cap and rightward cup to denote 
the bimodule homorphisms arising from the 
natural adjunction between induction and restriction:
\begin{align*}
\begin{tikzpicture}[anchorbase,scale=2]
\draw[thin,-to] (-.15,-.15) to [out=90,in=180] (0,.15) to [out=0
,in=90] (.15,-.15);
\node at (0.25,0) {$\regionlabel I$};
\node at (0,0) {$\regionlabel Ii$};
\end{tikzpicture}
\hspace{.5mm}
&:B_{I,Ii} \otimes_{R^{Ii}} B_{Ii,I} \rightarrow R^I,&
f \otimes g &\mapsto fg,\\
\begin{tikzpicture}[anchorbase,scale=2]
\draw[thin,-to] (-.15,.15) to [out=-90,in=180] (0,-.15) to [out=0,in=-90] (.15,.15);
\node at (0,0) {$\regionlabel I$};
\node at (0.25,0) {$\regionlabel Ii$};
\end{tikzpicture}
&:R^{Ii}\rightarrow B_{Ii,I} \otimes_{R^I} B_{I,Ii},&
1 &\mapsto 1 \otimes 1.\\\intertext{The leftward cap and leftward cup 
denote the bimodule homomorphisms arising from
the Frobenius adjunction, i.e.,
the trace and comultiplication maps \cref{tracedef,comultdef}:}
\begin{tikzpicture}[anchorbase,scale=2]
\draw[thin,to-] (-.15,-.15) to [out=90,in=180] (0,.15) to [out=0
,in=90] (.15,-.15);
\node at (0.25,0) {$\regionlabel Ii$};
\node at (0,0) {$\regionlabel I$};
\end{tikzpicture}
&:
B_{Ii,I}\otimes_{R^I} B_{I,Ii}
\rightarrow R^{Ii},&
f \otimes g &\mapsto \tr_{Ii}^I(fg),\\
\begin{tikzpicture}[anchorbase,scale=2]
\draw[thin,to-] (-.15,.15) to [out=-90,in=180] (0,-.15) to [out=0,in=-90] (.15,.15);
\node at (0,0) {$\regionlabel Ii$};
\node at (0.25,0) {$\regionlabel I$};
\end{tikzpicture}\:
&:R^I \rightarrow B_{I,Ii}\otimes_{R^{Ii}} B_{Ii,I},&
1 &\mapsto \sum_{b \in \B_{Ii}^I} b \otimes b^\vee.
\end{align*}
The (clockwise) rightward cap and leftward cup are of positive degree $\ell(w_{Ii})-\ell(w_I)$ and the (counterclockwise) rightward cup and leftward cap are of negative degree 
$\ell(w_I)-\ell(w_{Ii})$.
The adjunctions give the zig-zag relations
\begin{align}\label{zigzag}
\begin{tikzpicture}[anchorbase]
  \draw[->] (0.3,0) to (0.3,.4);
	\draw[-] (0.3,0) to[out=-90, in=0] (0.1,-0.2);
	\draw[-] (0.1,-0.2) to[out = 180, in = -90] (-0.1,0);
	\draw[-] (-0.1,0) to[out=90, in=0] (-0.3,0.2);
	\draw[-] (-0.3,0.2) to[out = 180, in =90] (-0.5,0);
  \draw[-] (-0.5,0) to (-0.5,-.4);
\node at (0.5,0) {$\regionlabel Ii$};
\node at (-.7,0) {$\regionlabel I$};
\end{tikzpicture}
&=
\begin{tikzpicture}[anchorbase]
  \draw[->] (0,-0.4) to (0,.4);
\node at (0.2,0) {$\regionlabel Ii$};
\node at (-.2,0) {$\regionlabel I$};
\end{tikzpicture}=
\begin{tikzpicture}[anchorbase]
  \draw[-] (0.3,0) to (0.3,-.4);
	\draw[-] (0.3,0) to[out=90, in=0] (0.1,0.2);
	\draw[-] (0.1,0.2) to[out = 180, in = 90] (-0.1,0);
	\draw[-] (-0.1,0) to[out=-90, in=0] (-0.3,-0.2);
	\draw[-] (-0.3,-0.2) to[out = 180, in =-90] (-0.5,0);
  \draw[->] (-0.5,0) to (-0.5,.4);
\node at (0.5,0) {$\regionlabel Ii$};
\node at (-.7,0) {$\regionlabel I$};
\end{tikzpicture},&
\begin{tikzpicture}[anchorbase]
  \draw[->] (0.3,0) to (0.3,-.4);
	\draw[-] (0.3,0) to[out=90, in=0] (0.1,0.2);
	\draw[-] (0.1,0.2) to[out = 180, in = 90] (-0.1,0);
	\draw[-] (-0.1,0) to[out=-90, in=0] (-0.3,-0.2);
	\draw[-] (-0.3,-0.2) to[out = 180, in =-90] (-0.5,0);
  \draw[-] (-0.5,0) to (-0.5,.4);
\node at (0.5,0) {$\regionlabel I$};
\node at (-.7,0) {$\regionlabel Ii$};
\end{tikzpicture}
&=
\begin{tikzpicture}[anchorbase]
  \draw[<-] (0,-0.4) to (0,.4);
\node at (0.2,0) {$\regionlabel I$};
\node at (-.2,0) {$\regionlabel Ii$};
\end{tikzpicture}=
\begin{tikzpicture}[anchorbase]
  \draw[-] (0.3,0) to (0.3,.4);
	\draw[-] (0.3,0) to[out=-90, in=0] (0.1,-0.2);
	\draw[-] (0.1,-0.2) to[out = 180, in = -90] (-0.1,0);
	\draw[-] (-0.1,0) to[out=90, in=0] (-0.3,0.2);
	\draw[-] (-0.3,0.2) to[out = 180, in =90] (-0.5,0);
  \draw[->] (-0.5,0) to (-0.5,-.4);
\node at (0.5,0) {$\regionlabel I$};
\node at (-.7,0) {$\regionlabel Ii$};
\end{tikzpicture}.
\end{align}
It makes sense to include an element $a$ of the algebra $R^I$ into a 2-cell labelled $I$;
we draw this by putting $a$ into a lightly shaded node in the region.
We obviously have that
\begin{align}\label{bubbleslides}
\begin{tikzpicture}[anchorbase]
  \draw[->] (0,-.4) to (0,.4);
\node at (-.2,0) {$\regionlabel I$};
\node at (.6,0) {$\regionlabel Ii$};
\node [draw,rounded corners,inner sep=2pt,black,fill=gray!20!white] at (.25,0) {$\scriptstyle f$};
\end{tikzpicture}&=
\begin{tikzpicture}[anchorbase]
  \draw[->] (0,-.4) to (0,.4);
\node at (-.6,0) {$\regionlabel I$};
\node at (.2,0) {$\regionlabel Ii$};
\node [draw,rounded corners,inner sep=2pt,black,fill=gray!20!white] at (-.25,0) {$\scriptstyle f$};
\end{tikzpicture},&
\begin{tikzpicture}[anchorbase]
  \draw[<-] (0,-.4) to (0,.4);
\node at (.2,0) {$\regionlabel I$};
\node at (-.6,0) {$\regionlabel Ii$};
\node [draw,rounded corners,inner sep=2pt,black,fill=gray!20!white] at (-.25,0) {$\scriptstyle f$};
\end{tikzpicture}&=
\begin{tikzpicture}[anchorbase]
  \draw[<-] (0,-.4) to (0,.4);
\node at (.6,0) {$\regionlabel I$};
\node at (-.25,0) {$\regionlabel Ii$};
\node [draw,rounded corners,inner sep=2pt,black,fill=gray!20!white] at (.25,0) {$\scriptstyle f$};
\end{tikzpicture}
\end{align}
for $f \in R^{Ii}\subset R^I$.
In particular, elements $f \in R^N$
slide across strings of all colors, i.e., $\BSBim$ is an $R^N$-linear graded bicategory.
Closed bubbles evaluate according to the relations
\begin{align}\label{evaluations}
\begin{tikzpicture}[anchorbase]
\draw[->] (.5,0) arc[start angle=0, end angle=360,radius=.5];
\node [draw,rounded corners,inner sep=2pt,black,fill=gray!20!white] at (0,.2) {$\scriptstyle f$};
\node at (.75,0) {$\regionlabel Ii$};
\node at (0,-.2) {$\regionlabel I$};
\end{tikzpicture}&=\begin{tikzpicture}[anchorbase]
\node [draw,rounded corners,inner sep=2pt,black,fill=gray!20!white] at (0,0) {$\scriptstyle \tr_{Ii}^{I}(f)$};
\node at (.75,0) {$\regionlabel Ii$};
\end{tikzpicture},&
\begin{tikzpicture}[anchorbase]
\draw[<-] (.5,0) arc[start angle=0, end angle=360,radius=.3];
\node at (.75,0) {$\regionlabel I$};
\node at (.2,0) {$\regionlabel Ii$};
\end{tikzpicture}&=
\begin{tikzpicture}[anchorbase]
\node [draw,rounded corners,inner sep=2pt,black,fill=gray!20!white] at (0,0) {$\scriptstyle \eta^I_{Ii}$};
\node at (.5,0) {$\regionlabel I$};
\end{tikzpicture},
\end{align}
for $f \in R^I$.
Also \cref{counit} can be formulated diagrammatically as
\begin{align}\label{cupcap}
\begin{tikzpicture}[anchorbase]
\draw[<-] (-.4,-1) to (-.4,1);
\draw[->] (.4,-1) to (.4,1);
\node at (0,0) {$\regionlabel I$};
\node at (.65,0) {$\regionlabel Ii$};
\node at (-.65,0) {$\regionlabel Ii$};
\end{tikzpicture}&=
\sum_{b \in \B_{Ii}^I}\begin{tikzpicture}[anchorbase]
\draw[<-] (-.6,-1) to[out=90,in=90,looseness=2.4] (.6,-1);
\draw[->] (-.6,1) to[out=-90,in=-90,looseness=2.4] (.6,1);
\node at (0,.4) {$\regionlabel I$};
\node at (0,-.4) {$\regionlabel I$};
\node at (0,0) {$\regionlabel Ii$};
\node [draw,rounded corners,inner sep=2pt,black,fill=gray!20!white] at (0,.8) {$\scriptstyle b$};
\node [draw,rounded corners,inner sep=2pt,black,fill=gray!20!white] at (0,-.8) {$\scriptstyle b^\vee$};
\end{tikzpicture}.
\end{align}

There are also
upward and downward crossings of strings of different colors $i$ and $j$, which
are the obvious degree 0 bimodule isomorphisms 
arising from transitivity of induction and restriction:
\begin{align*}
\begin{tikzpicture}[anchorbase,scale=1.8]
	\draw[<-] (0.25,.3) to (-0.25,-.3);
	\draw[->] (0.25,-.3) to (-0.25,.3);
 \node at (.2,0) {$\regionlabel Iij$};
 \node at (0,.2) {$\regionlabel Ii$};
 \node at (0,-.2) {$\regionlabel Ij$};
 \node at (-.2,0) {$\regionlabel I$};
\end{tikzpicture}
&:B_{I,Ij} \otimes_{R^{Ij}} B_{Ij,Iij}
\rightarrow B_{I,Ii} \otimes_{R^{Ii}} B_{Ii,Iij},&
f \otimes 1 &\mapsto f \otimes 1,\\
\begin{tikzpicture}[anchorbase,scale=1.8]
	\draw[->] (0.25,.3) to (-0.25,-.3);
	\draw[<-] (0.25,-.3) to (-0.25,.3);
 \node at (-.2,0) {$\regionlabel Iij$};
 \node at (0,.2) {$\regionlabel Ij$};
 \node at (0,-.2) {$\regionlabel Ii$};
 \node at (.2,0) {$\regionlabel I$};
\end{tikzpicture}
&:B_{Iij,Ii}\otimes_{R^{Ii}} B_{Ii,I} \rightarrow B_{Iij,Ij}\otimes_{R^{Ij}} B_{Ij,I},&
1 \otimes f &\mapsto 1 \otimes f.
\end{align*}
This notation implicitly assumes that $i, j$ 
are different elements of $N - I$
such that $Iij$ is finitary. The following 
{\em easy} Reidemeister II relations are clear:
\begin{align}\label{easyRII}
\begin{tikzpicture}[anchorbase,scale=1.8]
	\draw[->] (0.25,-.5) to[looseness=2,out=150,in=-150] (0.25,.5);
	\draw[->] (-.25,-.5) to[looseness=2,out=30,in=-30] (-0.25,.5);
 \node at (.4,0) {$\regionlabel Iij$};
 \node at (0,.5) {$\regionlabel Ii$};
 \node at (0,0) {$\regionlabel Ij$};
 \node at (0,-.5) {$\regionlabel Ii$};
 \node at (-.35,0) {$\regionlabel I$};
\end{tikzpicture}
&=
\begin{tikzpicture}[anchorbase,scale=1.8]
	\draw[->] (0.25,-.5) to (0.25,.5);
	\draw[->] (-.2,-.5) to (-0.2,.5);
 \node at (.4,0) {$\regionlabel Iij$};
 \node at (-.3,0) {$\regionlabel I$};
 \node at (0,0) {$\regionlabel Ii$};
 \end{tikzpicture},
& \begin{tikzpicture}[anchorbase,scale=1.8]
	\draw[<-] (0.25,-.5) to[looseness=2,out=150,in=-150] (0.25,.5);
	\draw[<-] (-.25,-.5) to[looseness=2,out=30,in=-30] (-0.25,.5);
 \node at (-.4,0) {$\regionlabel Iij$};
 \node at (0,.5) {$\regionlabel Ii$};
 \node at (0,0) {$\regionlabel Ij$};
 \node at (0,-.5) {$\regionlabel Ii$};
 \node at (.35,0) {$\regionlabel I$};
\end{tikzpicture}
&=
\begin{tikzpicture}[anchorbase,scale=1.8]
	\draw[<-] (0.25,-.5) to (0.25,.5);
	\draw[<-] (-.2,-.5) to (-0.2,.5);
 \node at (-.4,0) {$\regionlabel Iij$};
 \node at (.35,0) {$\regionlabel I$};
 \node at (0,0) {$\regionlabel Ii$};
 \end{tikzpicture}.
 \end{align}
 
There are sideways crossings which are defined by rotating the upward or downward ones:
\begin{align*}
\begin{tikzpicture}[anchorbase,scale=1.6]
	\draw[<-] (0.25,.3) to (-0.25,-.3);
	\draw[<-] (0.25,-.3) to (-0.25,.3);
 \node at (.2,0) {$\regionlabel Ii$};
 \node at (0,.2) {$\regionlabel I$};
 \node at (0,-.2) {$\regionlabel Iij$};
 \node at (-.2,0) {$\regionlabel Ij$};
\end{tikzpicture}
&:=\!
\begin{tikzpicture}[anchorbase,scale=1]
\draw[->] (.2,-.5) to (-.2,.5);
\draw[->] (-.6,.5) to[out=-90,in=-135,looseness=2] (0,0) to [out=45,in=90,looseness=2] (.6,-.5);
\node at (0.4,0.4) {$\regionlabel Ii$};
\node at (-.4,-.4) {$\regionlabel Ij$};
\node at (.3,-.1) {$\regionlabel Iij$};
\node at (-.3,.1) {$\regionlabel I$};
\end{tikzpicture}
\!=\!\begin{tikzpicture}[anchorbase,scale=1]
\draw[->] (.2,.5) to (-.2,-.5);
\draw[->] (-.6,-.5) to[out=90,in=135,looseness=2] (0,0) to [out=-45,in=-90,looseness=2] (.6,.5);
\node at (0.4,-.4) {$\regionlabel Ii$};
\node at (-.4,.4) {$\regionlabel Ij$};
\node at (.3,.1) {$\regionlabel I$};
\node at (-.3,-.1) {$\regionlabel Iij$};
\end{tikzpicture}
:B_{Ij,Iij}\otimes_{R^{Ij}}B_{Iij,Ii}
\rightarrow
B_{Ij,I}\otimes_{R^{I}}B_{I,Ii}
,&1\otimes 1 &\mapsto 1\otimes 1,\\ 
\begin{tikzpicture}[anchorbase,scale=1.6]
	\draw[->] (0.25,.3) to (-0.25,-.3);
	\draw[->] (0.25,-.3) to (-0.25,.3);
 \node at (.2,0) {$\regionlabel Ij$};
 \node at (0,.2) {$\regionlabel Iij$};
 \node at (0,-.2) {$\regionlabel I$};
 \node at (-.2,0) {$\regionlabel Ii$};
\end{tikzpicture}
&:=\!
\begin{tikzpicture}[anchorbase,scale=1]
\draw[->] (-.2,-.5) to (.2,.5);
\draw[<-] (-.6,-.5) to[out=90,in=135,looseness=2] (0,0) to [out=-45,in=-90,looseness=2] (.6,.5);
\node at (0.4,-.4) {$\regionlabel Ij$};
\node at (-.4,.4) {$\regionlabel Ii$};
\node at (.3,.1) {$\regionlabel Iij$};
\node at (-.3,-.1) {$\regionlabel I$};
\end{tikzpicture}
\!=\!\begin{tikzpicture}[anchorbase,scale=1]
\draw[->] (-.2,.5) to (.2,-.5);
\draw[<-] (-.6,.5) to[out=-90,in=-135,looseness=2] (0,0) to [out=45,in=90,looseness=2] (.6,-.5);
\node at (0.4,0.4) {$\regionlabel Ij$};
\node at (-.4,-.4) {$\regionlabel Ii$};
\node at (.3,-.1) {$\regionlabel I$};
\node at (-.3,.1) {$\regionlabel Iij$};
\end{tikzpicture}
:B_{Ii,I}\otimes_{R^{I}}B_{I,Ij}
\rightarrow
B_{Ii,Iij}\otimes_{R^{Iij}}B_{Iij,Ij}
,&f \otimes 1&\mapsto \sum_{b \in \B^{Ij}_{Iij}}
\tr^I_{Ii}(fb)\otimes b^\vee.
\end{align*}
The equalities here reflect the cyclicity of this bicategory. Both sideways crossings have degree
$\ell(w_{Iij})+\ell(w_I)-\ell(w_{Ii})-\ell(w_{Ij})$, which is half the degree of the element $\eta_{Iij}^I$. The mnemonic to remember this degree is: ``big plus small minus middle minus middle."
This degree is always non-negative, and it is 0 if and only if $i$ and $j$ live in different connected components of $Iij$.
The {\em hard}
Reidemeister II relations are
\begin{align}\label{hardRII}
\begin{tikzpicture}[anchorbase,scale=1.8]
	\draw[->] (0.25,-.5) to[looseness=2,out=150,in=-150] (0.25,.5);
	\draw[<-] (-.25,-.5) to[looseness=2,out=30,in=-30] (-0.25,.5);
 \node at (.4,0) {$\regionlabel Ij$};
 \node at (0,.5) {$\regionlabel I$};
 \node at (0,0) {$\regionlabel Iij$};
 \node at (0,-.5) {$\regionlabel I$};
 \node at (-.4,0) {$\regionlabel Ii$};
\end{tikzpicture}
&=
\begin{tikzpicture}[anchorbase,scale=1.8]
	\draw[->] (0.25,-.5) to (0.25,.5);
	\draw[<-] (-.2,-.5) to (-0.2,.5);
 \node at (.4,.2) {$\regionlabel Ij$};
 \node at (-.4,.2) {$\regionlabel Ii$};
 \node at (0,.2) {$\regionlabel I$};
\node [draw,rounded corners,inner sep=1.5pt,black,fill=gray!20!white] at (0.02,-.2) {$\scriptstyle \eta^I_{Iij}$};
 \end{tikzpicture},
 &\begin{tikzpicture}[anchorbase,scale=1.8]
	\draw[<-] (0.25,-.5) to[looseness=2,out=150,in=-150] (0.25,.5);
	\draw[->] (-.25,-.5) to[looseness=2,out=30,in=-30] (-0.25,.5);
 \node at (-.4,0) {$\regionlabel Ii$};
 \node at (0,.5) {$\regionlabel Iij$};
 \node at (0,.1) {$\regionlabel I$};
 \node at (0,-.5) {$\regionlabel Iij$};
 \node at (.4,0) {$\regionlabel Ij$};
\node [draw,rounded corners,inner sep=2pt,black,fill=gray!20!white] at (0,-.1) {$\scriptstyle f$};
\end{tikzpicture}
&=
 \sum_{b \in \B_{Iij}^{Ij}}\begin{tikzpicture}[anchorbase,scale=1.8]
	\draw[<-] (0.2,-.5) to (0.2,.5);
	\draw[->] (-.2,-.5) to (-0.2,.5);
 \node at (.4,.15) {$\regionlabel Ij$};
 \node at (-.4,.15) {$\regionlabel Ii$};
 \node at (0,.15) {$\regionlabel Iij$};
\node [draw,rounded corners,inner sep=2pt,black,fill=gray!20!white] at (-.6,-.15) {$\scriptstyle \tr_{Ii}^I(bf)$};
\node [draw,rounded corners,inner sep=1.5pt,black,fill=gray!20!white] at (.4,-.15) {$\scriptstyle b^\vee$};
 \end{tikzpicture}
 \end{align}
for $f \in R^I$.
The second of these can be rewritten in several equivalent ways using the following identity,
which follows from \cref{deadhead1,deadhead2,deadhead3}:
\begin{equation}\label{deadhead5}
\sum_{b \in \B_{Iij}^{Ij}}
\tr_{Ii}^I(bf) \otimes b^\vee
=\sum_{b \in \B_{Iij}^I} \tr^I_{Ii}(b) \otimes \tr^I_{Ij}(fb^\vee)
=\sum_{b \in \B_{Iij}^I} \tr^I_{Ii}(bf) \otimes \tr^I_{Ij}(b^\vee)
=\sum_{b \in \B_{Iij}^{Ii}} b\otimes \tr_{Ij}^I (fb^\vee)\end{equation}
in $R^{Ii} \otimes_{R^{Iij}} R^{Ij}$.
There is a special case of \cref{hardRII} we wish to call out explicitly: when $i$ and $j$ are in distinct connected components of $Iij$,
 we have the {\em distant} Reidemeister II relations
\begin{align}\label{distantRII}
\begin{tikzpicture}[anchorbase,scale=1.8]
	\draw[->] (0.25,-.5) to[looseness=2,out=150,in=-150] (0.25,.5);
	\draw[<-] (-.25,-.5) to[looseness=2,out=30,in=-30] (-0.25,.5);
 \node at (.4,0) {$\regionlabel Ij$};
 \node at (0,.5) {$\regionlabel I$};
 \node at (0,0) {$\regionlabel Iij$};
 \node at (0,-.5) {$\regionlabel I$};
 \node at (-.4,0) {$\regionlabel Ii$};
\end{tikzpicture}
&=
\begin{tikzpicture}[anchorbase,scale=1.8]
	\draw[->] (0.25,-.5) to (0.25,.5);
	\draw[<-] (-.2,-.5) to (-0.2,.5);
 \node at (.4,0) {$\regionlabel Ij$};
 \node at (-.4,0) {$\regionlabel Ii$};
 \node at (0,0) {$\regionlabel I$};
 \end{tikzpicture},
 &\begin{tikzpicture}[anchorbase,scale=1.8]
	\draw[<-] (0.25,-.5) to[looseness=2,out=150,in=-150] (0.25,.5);
	\draw[->] (-.25,-.5) to[looseness=2,out=30,in=-30] (-0.25,.5);
 \node at (-.4,0) {$\regionlabel Ii$};
 \node at (0,.5) {$\regionlabel Iij$};
 \node at (0,0) {$\regionlabel I$};
 \node at (0,-.5) {$\regionlabel Iij$};
 \node at (.4,0) {$\regionlabel Ij$};
\end{tikzpicture}
&=
\begin{tikzpicture}[anchorbase,scale=1.8]
	\draw[<-] (0.2,-.5) to (0.2,.5);
	\draw[->] (-.2,-.5) to (-0.2,.5);
 \node at (.4,0) {$\regionlabel Ij$};
 \node at (-.4,0) {$\regionlabel Ii$};
 \node at (0,0) {$\regionlabel Iij$};
 \end{tikzpicture}\;.
 \end{align}
Thus, in this case, the sideways crossings are mutually inverse isomorphisms.

The final general relations from \cite{ESW}
are the Reidemeister III relations, {\em easy} and {\em hard}:
\begin{align}\label{RIII}
\begin{tikzpicture}[anchorbase,scale=2.5]
\draw[->] (-.3,-.5) to (.3,.5);
\draw[->] (.3,-.5) to (-.3,.5);
	\draw[->] (0,-.5) to[looseness=1,out=30,in=-30] (0,.5);
 \node at (.4,0) {$\regionlabel Iijk$};
 \node at (.16,.48) {$\regionlabel Ijk$};
 \node at (0,.25) {$\regionlabel Ik$};
 \node at (0,-.25) {$\regionlabel Ii$};
 \node at (0.16,-.48) {$\regionlabel Iij$};
 \node at (-.2,0) {$\regionlabel I$};
 \node at (.14,0) {$\regionlabel Iik$};
\end{tikzpicture}&=\begin{tikzpicture}[anchorbase,scale=2.5]
\draw[->] (-.3,-.5) to (.3,.5);
\draw[->] (.3,-.5) to (-.3,.5);
	\draw[->] (0,-.5) to[looseness=1,out=150,in=-150] (0,.5);
 \node at (.2,0) {$\regionlabel Iijk$};
 \node at (-.18,.45) {$\regionlabel Ik$};
 \node at (-.18,-.45) {$\regionlabel Ii$};
 \node at (0,.25) {$\regionlabel Ijk$};
 \node at (0,-.25) {$\regionlabel Iij$};
 \node at (-.35,0) {$\regionlabel I$};
 \node at (-.15,0) {$\regionlabel Ij$};
\end{tikzpicture}
,&
\begin{tikzpicture}[anchorbase,scale=2.5]
\draw[->] (-.3,-.5) to (.3,.5);
\draw[->] (.3,-.5) to (-.3,.5);
	\draw[<-] (0,-.5) to[looseness=1,out=30,in=-30] (0,.5);
 \node at (.38,0) {$\regionlabel Iik$};
 \node at (.18,.45) {$\regionlabel Ik$};
 \node at (0,.25) {$\regionlabel Ijk$};
 \node at (0,-.25) {$\regionlabel Iij$};
 \node at (0.18,-.45) {$\regionlabel Ii$};
 \node at (-.2,0) {$\regionlabel Ij$};
 \node at (.13,0) {$\regionlabel Iijk$};
\end{tikzpicture}&=\begin{tikzpicture}[anchorbase,scale=2.5]
\draw[->] (-.3,-.5) to (.3,.5);
\draw[->] (.3,-.5) to (-.3,.5);
	\draw[<-] (0,-.5) to[looseness=1.8,out=150,in=-150] (0,.5);
 \node at (.2,0) {$\regionlabel Iik$};
 \node at (-.16,.5) {$\regionlabel Ijk$};
 \node at (-.16,-.5) {$\regionlabel Iij$};
 \node at (0,.25) {$\regionlabel Ik$};
 \node at (0,-.25) {$\regionlabel Ii$};
 \node at (-.58,0) {$\regionlabel Ij$};
 \node at (-.25,0.1) {$\regionlabel I$};
\node [draw,rounded corners,inner sep=1.5pt,black,fill=gray!20!white] at (-.25,-.1) {$\scriptstyle \eta^I_{Iijk}$};
\end{tikzpicture}.
\end{align}
Note also that if $\{i,j,k\}$ is not contained in one connected component of $Iijk$, then $\eta^I_{Iijk} = 1$, leading to the \emph{distant} Reidemeister III relation.
There are several more variations on Reidemeister III, which may be obtained from the ones above by rotating using the cyclic structure. 

Another way to obtain further relations 
is to apply the graded 2-functor
\begin{align}\label{tR}
\tR&:\BSBim\rightarrow \BSBim^{\rev}
\end{align}
defined as follows:
it fixes objects,
takes the graded $(R^I, R^J)$-bimodule $M$ to $q^{\ell(w_I)-\ell(w_J)} M$
viewed using commutativity as a graded $(R^J, R^I)$-bimodule,
and takes a bimodule homomorphism to the same underlying function, which is automatically also a homomorphism with respect to the new bimodule structure.
In particular, it maps the bimodule $B_{I,K}^J$ to $B_{K,I}^J$.
In terms of string diagrams, 
$\tR$ reflects diagrams in a vertical axis and reverses the orientation on all strings.

Finally, we record a useful consequence of the relations described so far.
This was originally discovered as \cite[Theorem 6.7]{EKLPSingLL}.
For $i \in I \finsub N$, we write $I \hat \imath$ as a shorthand for $I - \{i\}$.

\begin{lem}[Switchback relation]\label{switchback}
Suppose we are given $I$ and distinct $j,k \in N - I$ with $Ijk \finsub N$,
and also $i \in Ik$ such that
$\partial_{Ijk}^{Ij} = \partial_{I\hat\imath jk}^{I\hat\imath k}\circ\partial_{Ik}^I$.
Then the following relation holds:
\begin{align}\label{engage}
\begin{tikzpicture}[anchorbase,scale=1.2]
\draw[<-] (.3,-.6) to (.3,.6);
\draw[->] (-.3,-.6) to (-.3,.6);
\draw[-] (0,0) ellipse (.8 and .4);
\draw[->] (.8,0)
arc[start angle=-.5, end angle=.5,radius=.8];
\node at (-1.1,0) {$\regionlabel I\hat\imath jk$};
\node at (1,0) {$\regionlabel Ij$};
\node at (-.55,0) {$\regionlabel I\hat\imath k$};
\node at (.5,0) {$\regionlabel I$};
\node at (0,0) {$\regionlabel Ik$}; 
\node at (0,-0.54) {$\regionlabel Ijk$};
\node at (0,0.54) {$\regionlabel Ijk$};
\end{tikzpicture}
&=\begin{tikzpicture}[anchorbase,scale=1.2]
\draw[<-] (.25,-.6) to (.25,.6);
\draw[->] (-.25,-.6) to (-.25,.6);
\node at (0,0) {$\regionlabel Ijk$};
\node at (-.55,0) {$\regionlabel I\hat\imath jk$};
\node at (.5,0) {$\regionlabel Ij$};
\end{tikzpicture}\:.
\end{align}
\end{lem}

\begin{proof}
By \cref{hardRII} and easy relations, we have that
\begin{align*}
\begin{tikzpicture}[anchorbase,scale=1.2]
\draw[<-] (.3,-.8) to (.3,.8);
\draw[->] (-.3,-.8) to (-.3,.8);
\draw[-] (0,0) ellipse (.8 and .4);
\draw[->] (.8,0)
arc[start angle=-.5, end angle=.5,radius=.8];
\node at (-1.1,0) {$\regionlabel I\hat\imath jk$};
\node at (1,0) {$\regionlabel Ij$};
\node at (-.55,0) {$\regionlabel I\hat\imath k$};
\node at (.5,0) {$\regionlabel I$};
\node at (0,0) {$\regionlabel Ik$}; 
\node at (0,0.54) {$\regionlabel Ijk$};
\node at (0,-0.54) {$\regionlabel Ijk$};
\end{tikzpicture}\!
&=
\sum_{b \in \B_{Ijk}^{Ij}}
\begin{tikzpicture}[anchorbase]
\draw[<-] (1.25,-.8) to (1.25,.8);
\draw[->] (-.8,-.8) to (-.8,.8);
\draw[-] (-.65,0) ellipse (1.1 and .6);
\draw[->] (-1.75,0)
arc[start angle=179.5, end angle=180.5,radius=1.1];
\node at (2.25,0) {$\regionlabel Ij$};
\node at (-2.15,0) {$\regionlabel I\hat\imath jk$};
\node at (-1.3,0) {$\regionlabel I\hat\imath k$};
\node at (-.62,.35) {$\regionlabel Ik$}; 
\node at (-0.23,-.1) [draw,rounded corners,inner sep=2pt,black,fill=gray!20!white] {$\scriptstyle\tr_{Ik}^I(b)$};
\node at (1.68,0) [draw,rounded corners,inner sep=2pt,black,fill=gray!20!white] {$\scriptstyle b^\vee$};
\node at (.8,0) {$\regionlabel Ijk$}; 
\end{tikzpicture}
=
\sum_{b \in \B_{Ijk}^{Ij}}
\begin{tikzpicture}[anchorbase]
\draw[<-] (1.25,-.8) to (1.25,.8);
\draw[->] (.3,-.8) to (.3,.8);
\draw[-] (-.8,0) ellipse (.8 and .6);
\draw[->] (-1.6,0)
arc[start angle=179.5, end angle=180.5,radius=.8];
\node at (2.25,0) {$\regionlabel Ij$};
\node at (-2.15,0) {$\regionlabel I\hat\imath jk$};
\node at (-1.3,0) {$\regionlabel I\hat\imath k$};
\node at (-0.6,0) [draw,rounded corners,inner sep=2pt,black,fill=gray!20!white] {$\scriptstyle\tr_{Ik}^I(b)$};
\node at (1.68,0) [draw,rounded corners,inner sep=2pt,black,fill=gray!20!white] {$\scriptstyle b^\vee$};
\node at (.75,0) {$\regionlabel Ijk$}; 
\end{tikzpicture}
\\
&
=\sum_{b \in \B_{Ijk}^{Ij}}
\begin{tikzpicture}[anchorbase]
\draw[<-] (-.2,-.6) to (-.2,.6);
\draw[->] (-.7,-.6) to (-.7,.6);
\node at (-2.9,0) {$\regionlabel I\hat\imath jk$};
\node at (.9,0) {$\regionlabel Ij$};
\node at (-.45,0) {$\regionlabel Ijk$};
\node [draw,rounded corners,inner sep=2pt,black,fill=gray!20!white] at (.4,0) {$\scriptstyle b^\vee$};
\node [draw,rounded corners,inner sep=2pt,black,fill=gray!20!white] at (-1.8,0) {$\scriptstyle \tr^{I\hat\imath k}_{I\hat\imath jk}\tr_{Ik}^I(b)$};
\end{tikzpicture}
=\sum_{b \in \B_{Ijk}^{Ij}}
\begin{tikzpicture}[anchorbase]
\draw[<-] (-.2,-.6) to (-.2,.6);
\draw[->] (-.7,-.6) to (-.7,.6);
\node at (-2.6,0) {$\regionlabel I\hat\imath jk$};
\node at (.85,0) {$\regionlabel Ij$};
\node at (-.45,0) {$\regionlabel Ijk$};
\node [draw,rounded corners,inner sep=2pt,black,fill=gray!20!white] at (.4,0) {$\scriptstyle b^\vee$};
\node [draw,rounded corners,inner sep=2pt,black,fill=gray!20!white] at (-1.6,0) {$\scriptstyle \tr^{Ij}_{Ijk}(b)$};
\end{tikzpicture}\\
&=\sum_{b \in \B_{Ijk}^{Ij}}
\begin{tikzpicture}[anchorbase]
\draw[<-] (-.1,-.6) to (-.1,.6);
\draw[->] (-1.4,-.6) to (-1.4,.6);
\node at (-1.8,0) {$\regionlabel I\hat\imath jk$};
\node at (.9,0) {$\regionlabel Ij$};
\node at (-.7,0.3) {$\regionlabel Ijk$};
\node [draw,rounded corners,inner sep=2pt,black,fill=gray!20!white] at (.35,0) {$\scriptstyle b^\vee$};
\node [draw,rounded corners,inner sep=2pt,black,fill=gray!20!white] at (-.75,-.15) {$\scriptstyle \tr^{Ij}_{Ijk}(b)$};
\end{tikzpicture}=\begin{tikzpicture}[anchorbase,scale=1.2]
\draw[<-] (.25,-.6) to (.25,.6);
\draw[->] (-.25,-.6) to (-.25,.6);
\node at (0,0) {$\regionlabel Ijk$};
\node at (-.55,0) {$\regionlabel I\hat\imath jk$};
\node at (.5,0) {$\regionlabel Ij$};
\end{tikzpicture}.
\end{align*}
The fourth equality is where we used the hypothesis
$\partial_{Ijk}^{Ij} = \partial_{I\hat\imath jk}^{I\hat\imath k}\circ\partial_{Ik}^I$. The last equality uses a bubble slide \cref{bubbleslides}, followed by the counit axiom \cref{counit}.
\end{proof}

 The hypothesis $\partial_{Ijk}^{Ij} = \partial_{I\hat\imath jk}^{I\hat\imath k}\circ\partial_{Ik}^I$ appearing in \cref{switchback} is
an instance of the {\em switchback relation} introduced in \cite{EK}. A complete list of all situations in which it holds can be found in \cite[Sec.~6]{EK}.

\begin{rem}
The relations recorded in this subsection are not intended to be a complete presentation of $\BSBim$, and indeed no such presentation exists currently in the literature.
\end{rem}

\subsection{The Soergel-Williamson categorification theorem}\label{ssec3-4}
Next we give an account of
the main results of \cite{Will}. 
When it comes to matters related to the Hecke algebra, our $q$ is equal to $v$ in \cite{Will}, but on bimodules Williamson identifies multiplication by $v$ with the downward degree shift, which is our $q^{-1}$. As well as updating definitions to take this into account, 
we have made a few other expository changes.
Our setup is biased towards writing translation functors on the left, and consequently many of the statements below have been obtained from the ones in \cite{Will}
by twisting with the anti-involution $\rho$ on the 
$q$-Schur algebra from \cref{tau} or the symmetry $\tR$ on singular Bott-Samelson bimodules from
\cref{tR}. See \cref{drivingmecrazy} for further discussion of differences compared to \cite{Will}.

Let $\Hecke$ be the Hecke algebra associated to the Coxeter group $W$ over the ground ring $\Z[q,q^{-1}]$.
This is the free $\Z[q,q^{-1}]$-module with basis
$\{h_w\:|\:w \in W\}$ viewed as an $\Z[q,q^{-1}]$-algebra with multiplication satisfying
\begin{align}
h_i h_w
 &= \begin{cases}
h_{s_i w}&\text{if $\ell(s_i w) > \ell(w)$}\\
h_{s_i w} - (q-q^{-1}) h_w&\text{if $\ell(s_i w) < \ell(w)$}
\end{cases}
\end{align}
for each $w \in W$ and $i \in I$, where
$h_i$ denotes $h_{s_i}$ for short.
In particular, $(h_i+q)(h_i-q^{-1})=1$
and $h_i^{-1} = h_i + (q-q^{-1})$.
The {\em bar involution}
is the anti-linear algebra involution 
$\Hecke \rightarrow \Hecke, h \mapsto \overline{h}$
defined so that $\overline{h_w} := h_{w^{-1}}^{-1}$ for each $w \in W$.
Also useful is the linear algebra anti-involution 
$\rho: \Hecke \rightarrow \Hecke, h_w \mapsto h_{w^{-1}}$.
It commutes with the bar involution, so 
$\omega:\Hecke \rightarrow \Hecke, h \mapsto 
\rho(\overline{h})$ is an anti-linear
algebra anti-involution
fixing each basis vector $h_w\:(w \in W)$.
We denote the {\em Kazhdan-Lusztig basis}
of $\Hecke$ (renormalized as in \cite{soergelKL})
by $b_w\:(w \in W)$.
So $b_w$ is the unique bar-invariant element of
$h_w + \sum_{w' < w} 
q\Z[q] h_{w'}$.
We have that $\rho(b_w) = b_{w^{-1}}$.

For $I \finsub N$, let $\Hecke_I$ be the 
subalgebra of $\Hecke$ with basis $h_w\:(w \in W_I)$.
The longest Kazhdan-Lusztig
basis element of $\Hecke_I$
is
\begin{equation}
b_{w_I} = 
\sum_{w \in W_I} q^{\ell(w_I)-\ell(w)} h_w.
\end{equation}
Note that $h_i b_{w_I} = b_{w_I} h_i = q^{-1} b_{w_I}$
for all $i \in I$.
It follows that
\begin{equation}\label{whydiv}
b_{w_I} b_{w_I} = \pi_I b_{w_I}
\end{equation}
where $\pi_I$ is the renormalized Poincar\'e polynomial from \cref{usualpoincare}.

Given $I, J \finsub N$, we are going to consider the set $W_I \backslash W / W_J$
of $(W_I, W_J)$-double cosets in $W$.
We let 
$(W_I \backslash W / W_J)_{\min}$
denote the set of minimal length coset representatives.
For $d \in (W_I \backslash W / W_J)_{\min}$ we often use the shorthand
$IdJ$ to succinctly encode the triple of $I$, $d$ and $J$, hence, the $(W_I, W_J)$-double coset $W_I d W_J$. On the other hand, given just the set $W_I d W_J$, one can recover $d$ as the minimal length element, but
$I$ and $J$ themselves are not uniquely determined.
For $d \in (W_I \backslash W / W_J)_{\min}$,
let 
\begin{equation}
\LR{I}{d}{J} := \{i \in I\:|\:
s_i d  = d s_j\text{ for some } j \in J\},
\end{equation}
which we call the {\em left redundancy}
 of the double coset $W_I d W_J$.
 We have that $W_I \cap d W_J d^{-1} = W_{\LR{I}{d}{J}}$ (Kilmoyer's theorem);
this is the stabilizer of $d W_J$ 
under the natural left action of $W_I$ on $W / W_J$.
Elements of $W_I d W_J$ can be written uniquely as
$u d v$ for $u \in (W_I / W_{\LR{I}{d}{J}})_{\min}$ and $v \in W_J$, and we have that $\ell(udv) = \ell(u)+\ell(d)+\ell(v)$ (Howlett's theorem).
The unique element of maximal length in the double coset $W_I d W_J$ is
\begin{equation}\label{dclongest}
w_{IdJ} := w_I w_{\LR{I}{d}{J}}^{-1} d w_J
\end{equation}
with
\begin{equation}\label{usefuldifference}
\ell(w_{IdJ}) = \ell(w_I) - \ell(w_{\LR{I}{d}{J}})+\ell(d)+\ell(w_J).
\end{equation}

\begin{rem}
There is also the 
{\em right redundancy} 
$\RR{I}{d}{J}
:= \{j \in J\:|\: d s_j = s_i d\text{ for some }i \in I\}$.
We have that
$d^{-1} W_I d  \cap W_J = W_{\RR{I}{d}{J}}$.
%,hence, $d W_{\RR{I}{d}{J}} =  W_{\LR{I}{d}{J}} d$.
Elements of $W_I d W_J$ can be written uniquely as
$u d v$ for $u \in W_I$ and $v \in (W_{\RR{I}{d}{J}} \backslash W_J)_{\min}$, and $\ell(udv) = \ell(u)+\ell(d)+\ell(v)$.
In our exposition, we only need $\LR{I}{d}{J}$, but everything could easily be reformulated in terms of $\RR{I}{d}{J}$,
e.g.,
$w_{IdJ} = w_I d w_{\RR{I}{d}{J}}^{-1} w_J$
with
$\ell(w_{IdJ}) = \ell(w_I) +\ell(d)+\ell(w_J)- \ell(w_{\RR{I}{d}{J}})$.
\end{rem}

\begin{defin}[{\cite[Sec.~2.3]{Will}}]\label{heckealgebroid}
For $I,J \finsub N$, let ${_{I\hspace{.8pt}}}\Schur_J := 
(b_{w_I} \Hecke) \cap (\Hecke
b_{w_J})$, which is a free $\Z[q,q^{-1}]$-submodule of $\Hecke$.
The (generalized) {\em $q$-Schur algebra}\footnote{In \cite{Will}, this is called the {\em Hecke algebroid}, 
and others refer to it as the {\em Schur algebroid}.} is 
the free $\Z[q,q^{-1}]$-module
\begin{equation}\label{schurd}
\Schur := \bigoplus_{I, J \finsub N} 
{_{I\hspace{.8pt}}}\Schur_J
\end{equation}
viewed as a $\Z[q,q^{-1}]$-algebra with multiplication $\m:\Schur\times \Schur \rightarrow \Schur$ defined by
\begin{equation}\label{schurm}
x \m y := 
\begin{cases}
\frac{1}{\pi_J} xy&\text{if $J=J'$}\\
0&\text{otherwise,}
\end{cases}
\end{equation}
$I,J,J',K \subseteq N$, $x \in {_{I\hspace{.8pt}}}\Schur_J$
and $y \in {_{J'}}\Schur_K$; 
in the first case here, we have used \cref{whydiv} to see that 
$xy$ is divisible by $\pi_J$.
We denote the element $b_{w_I}$ of the summand ${_{I\hspace{.8pt}}}\Schur_I$ by $1_I$ (by itself, $b_{w_I}$ is ambiguous as an element of $\Schur$).
By \cref{whydiv} again, we have that
$1_I \m 1_J = \delta_{I,J} 1_I$.
The elements $1_I\:(I \finsub N)$
are mutually orthogonal idempotents summing to the identity in $\Schur$, and ${_{I\hspace{.8pt}}}\Schur_J=1_I \m \Schur \m 1_J$. 
\end{defin}

The {\em bar involution} 
$\Schur\rightarrow \Schur,\:
h \mapsto \overline{h}$
is the anti-linear algebra involution
defined 
on the summand ${_{I\hspace{.8pt}}}\Schur_J$ by the restriction of the bar involution from $\Hecke$.
There is also a 
linear algebra anti-involution
\begin{align}
\label{tau}
\rho&: \Schur \rightarrow \Schur
\end{align}
defined by the maps
${_{I\hspace{.8pt}}}\Schur_J \rightarrow {_{J\hspace{.2pt}}}\Schur_I$ obtained by the restrictions
of $\rho:\Hecke \rightarrow \Hecke$ for all $I,J \finsub N$.
Again, $\rho$ commutes with the bar involution, so the composition
$\omega:\Schur\rightarrow \Schur$ of $\rho$ and the bar involution (in either order)
is an anti-linear
algebra anti-involution of $\Schur$.

For our purposes, there are two important bases of $\Schur$,
the {\em standard basis} 
$\{h_{IdJ}\}$
and the {\em Kazhdan-Lusztig basis}
$\{b_{IdJ}\}$,
both indexed by the symbols $IdJ$ for $I,J \finsub N$ and $d \in (W_I \backslash W / W_J)_{\min}$.
They are defined by
\begin{align}
h_{IdJ} &:=
\sum_{w \in W_I d W_J}
q^{\ell(w_{IdJ})-\ell(w)}h_w,
&b_{IdJ} &:=b_{w_{IdJ}},
\end{align}
both viewed as elements of the summand 
${_{I\hspace{.8pt}}}\Schur_J$ of $\Schur$.
We have that $b_{I1I}=h_{I1I} = 1_I$.
Also
\begin{align}\label{starbucksagain}
\rho(h_{IdJ}) &= h_{J d^{-1} I},&
\rho(b_{IdJ}) &= b_{J d^{-1} I}.
\end{align}
In general, the Kazhdan-Lusztig basis element
$b_{IdJ}$ is the unique bar-invariant element 
of the set $$
h_{IdJ} +\sum_{\substack{d' \in (W_I \backslash W / W_J)_{\min}\\ d' < d}} 
q\Z[q] h_{Id'J}.
$$
Since $1$ is minimal, we 
have that $b_{I1J} = h_{I1J}$.
We denote this instead by $b_{I,J}$; this parallels the notation 
for the bimodules $B_{I,J}$ in \cref{appearedalready} and their generalizations introduced just after \cref{inductionrestrictionbimodule}.
These special elements 
generate $\Schur$ as an algebra.
In fact, 
$\Schur$ is already generated 
by $b_{I,Ii}$ and $b_{Ii,I}$ for $I \subset N$ and $i \in N - I$ with $Ii$ finitary.
We also let
\begin{equation}
b^K_{I,J} := b_{I,K} \m b_{K,J}\label{bselements}
\end{equation}
for $I, J \finsub N$ and $K \subseteq I\cap J$; this parallels \cref{inductionrestrictionbimodule}.

\begin{lem}[{\cite[Prop.~2.8]{Will}}]\label{multrule}
The following hold for any $I, J, K \finsub N$ 
with $I \subseteq J$ or $I \supseteq J$ 
and $d \in (W_J \backslash W / W_K)_{\min}$:
\begin{enumerate}
\item
If $I \subseteq J$ then
$$
b_{I,J} \m h_{JdK} = \sum_{d' \in 
W_J d W_K \cap (W_I \backslash W / W_K)_{\min}}
q^{\ell(w_{JdK})-\ell(w_{Id'K})} h_{Id'K}.
$$
\item
If $I \supseteq J$ then
$$
b_{I,J} \m h_{JdK} = q^{\ell(d')-\ell(d)}
\frac{\pi_{\LR{I}{d'}{K}}}{\pi_{\LR{J}{d}{K}}} h_{Id'K}
$$
where $d'$ is the minimal length element of $W_I d W_K$.
\end{enumerate}
\end{lem}

\begin{cor}\label{bscorollary}
For $I, J \finsub N$ and $K \subseteq I \cap J$, we have that
$b_{I,J}^K = \frac{\pi_{I \cap J}}{\pi_K}
b_{I,J}$.
\end{cor}

The final basic notion is a bilinear form
$(-,-):\Schur \times \Schur \rightarrow \Z[q,q^{-1}]$.
Let 
$\tr:\Hecke \rightarrow \Z[q,q^{-1}]$ be
the usual symmetric Frobenius trace on the Hecke algebra, i.e.,
the $\Z[q,q^{-1}]$-linear map
with $\tr(h_w) = \delta_{w,1}$ for all $w \in W$.
Then, $(-,-)$ is defined so that 
different summands $_{I\hspace{.8pt}}\Schur_J$ are orthogonal to each other, and
\begin{equation}\label{altbilform}
(x,y):= 
\frac{1}{\pi^+_I}
\tr(\rho(x) y)
\end{equation}
for $x,y \in {_{I\hspace{.8pt}}}\Schur_J$.
From the definition of the form, we have
that
$(zx,y) = (x,\rho(z)y)$
for any $x,y,z \in \Schur$. 
The standard basis of $\Schur$ is an orthogonal basis with respect to this form, with
\begin{align}\label{justforjon}
(h_{IdJ}, h_{IdJ})&:= 
\frac{\pi^+_J}{\pi^+_{\LR{I}{d}{J}}} \in 1 + q^2\N[q^2].
\end{align}
In particular, the form is symmetric. 
This all follows from \cite[Lem.~2.13]{Will}, noting that the form $\langle-,-\rangle$ there is 
related to ours by $\langle x,y\rangle 
= q^{\ell(w_I)} (\rho(x), \rho(y))$.

\vspace{2mm}
Now we come to the results
from \cite{Will} relating the $q$-Schur algebra to singular Soergel bimodules.
Suppose that $I,J \subseteq N$.
As well as singular Soergel bimodules,
it is important to consider 
the graded $(R^I, R^J)$-bimodule
$R^{\LR{I}{d}{J}}_d$ for $d \in (W_I\backslash W / W_J)_{\min}$.
This denotes
$R^{\LR{I}{d}{J}}$ viewed as a graded
$(R^I, R^J)$-bimodule,
with action defined by $a \cdot v \cdot b := 
a v d(b)$.
By \cite[Lem.~4.5]{Will}, we have for $d \in (W_I \backslash W / W_J)_{\min}$ that
\begin{equation}\label{howstandard}
R_d \cong 
\pi^+_{I \cap d J d^{-1}} R^{I \cap d J d^{-1}}_d 
\end{equation}
as graded $(R^I,R^J)$-bimodules,
where the $R_d$ on the left hand side denotes $R$ 
viewed as a graded $(R^I, R^J)$-bimodule by
$a \cdot v \cdot b :=  a v d(b)$ again.
By \cite[Lem.~4.2(1)]{Will}, we have that
\begin{align}\label{gatwick}
R_d^{\LR{I}{d}{J}} &\cong 
\begin{cases}
\frac{\pi^+_I}{\pi^+_{\LR{I}{d}{J}}}
R^I&\text{as a graded left $R^I$-module}\\
\frac{\pi^+_J}{\pi^+_{\LR{I}{d}{J}}}
R^J&\text{as a graded right $R^J$-module.}
\end{cases}
\end{align}
Also \cite[Cor.~4.13]{Will} gives that
\begin{align}
\Hom_{R^I\dash R^J}(R_d^{\LR{I}{d}{J}},
R_{d'}^{\LR{I}{d'}{J}})
&\cong 
\begin{cases}
R^{\LR{I}{d}{J}}_d&\text{if $d=d'$}\\
\{0\}&\text{otherwise,}
\end{cases}
\end{align}
for any $d,d' \in (W_I \backslash W / W_J)_{\min}$.
The {\em costandard} and {\em standard bimodules} are
\begin{align}\label{nabla}
\nabla_{IdJ} &:=
q^{\ell(w_{\LR{I}{d}{J}})-\ell(w_I)-\ell(d)}R^{\LR{I}{d}{J}}_d=
q^{\ell(w_J)-\ell(w_{IdJ})}
R^{\LR{I}{d}{J}}_d,\\\label{delta}
\Delta_{IdJ} &:=  q^{2\ell(d)} \nabla_{IdJ},
\end{align}
respectively.
In the special case $d=1$, 
$\nabla_{I1J}$ and $\Delta_{I1J}$ are both equal to the
graded $(R^I, R^J)$-bimodule $B_{I,J}$ introduced just
after \cref{inductionrestrictionbimodule}.

The 2-functor $\tR$ defined in \cref{tR} makes sense on all graded $(R^I, R^J)$-bimodules, not just on Bott-Samelsons. It reverses tensor products, mapping
a graded $(R^I, R^J)$-bimodule $M$
to $q^{\ell(w_I)-\ell(w_J)} M$ viewed as a graded $(R^J, R^I)$-bimodule. It is clear from \cref{delta,nabla} that
\begin{align}\label{orderreversal}
\tR(\Delta_{IdJ}) &= \Delta_{Jd^{-1}I},&
\tR(\nabla_{IdJ}) &= \nabla_{Jd^{-1}I}.
\end{align}
In particular, we have that $\tR(B_{I,J}) = B_{J,I}$.
There is also a duality functor $\tD$
defined on a graded $(R^I, R^J)$-bimodule 
$M$ by
\begin{equation}\label{rightduality}
\tD(M) := q^{2\ell(w_J)-2\ell(w_I)}\Hom_{\dash R^J}(M, R^J),
\end{equation}
viewed as a graded $(R^I, R^J)$-bimodule 
so that $r \in R^I$ and $s \in R^J$
act on $f:M \rightarrow R^J$ by $(rf)(v) = f(rv)$ and
$(fs)(v) = sf(v)$ (this uses commutativity of these algebras).
Arguing as in the proof of
\cite[Prop.~6.17]{Will}, remembering the modified grading shifts in \cref{nabla,delta}, we have that
\begin{align}\label{standardduality}
\tD(\nabla_{IdJ}) &\cong \Delta_{IdJ},
&\tD(\Delta_{IdJ}) &\cong \nabla_{IdJ}.
\end{align}
In particular, this shows that $\tD(B_{I,J})
\cong B_{I,J}$.
Another important point is that 
\begin{align}\label{cwithduality}
\tD (M \otimes_{R^J} M')
\cong \tD(M)  \otimes_{R^J} \tD (M'),
\end{align}
for graded $(R^I, R^J)$- and 
$(R^J,R^K)$-bimodules $M$ and $M'$ such that $M$ is a singular Soergel bimodule.
It suffices to prove this 
in the special case that $M$ is the self-dual
bimodule $B_{I,J}$ with $I \subseteq J$ or $I \supseteq J$. When $I \supseteq J$, the result follows using that $\Res^J_I$
commutes with $\Hom_{\dash J}(-,R^J)$ (tautologically). Then one deduces the result when $I \subseteq J$ using the adjunctions \cref{adj1,adj2};
cf. \cite[Prop.~6.15]{Will}.
This argument also shows that
all of the singular Bott-Samelson bimodules \cref{abs} are self-dual.

\begin{rem}\label{drivingmecrazy}
We have changed the grading shift in the definition of duality compared to \cite{Will}.
Our choice of this grading shift is one of two sensible ways to correct a minor error in Williamson's original setup. 
This is also explained in \cite[Sec.~0.4]{Willerratum},
where it is described as the ``more comprehensive fix''.
It forces several other changes,
including to the degree shifts in \cref{nabla,delta}, which are
not quite the same as in \cite{Will} but
are exactly as in \cite[Sec.~0.4]{Willerratum}.
In view of \cref{cwithduality}, our version of the duality functor $\tD$ is convenient for inductive arguments using ``translation on the left".
If one favors ``translation on the right" as in \cite{Will}, it is more convenient to work from the outset with another duality functor $\tDl$ 
defined on a graded $(R^I, R^J)$-bimodule $M$ by
\begin{equation}\label{leftduality}
\tDl(M) := \Hom_{R^I\dash}(M, R^I).
\end{equation}
This is what is used in \cite[Sec.~0.4]{Willerratum}.
Note that $\tDl = \tR \circ \tD \circ \tR$.
From this, \cref{standardduality,orderreversal},
it follows that $\tDl(\nabla_{IdJ}) \cong \Delta_{IdJ}$.
From \cref{cwithduality},
it follows that
$\tDl (M \otimes_{R^J} M')
\cong \tDl(M)  \otimes_{R^J} \tDl(M')$
for graded $(R^I, R^J)$- and 
$(R^J,R^K)$-bimodules $M$ and $M'$ such that $M'$ is a singular Soergel bimodule.
Since all singular Bott-Samelsons 
are $\tD$-self-dual (explained above)
and $\tR(B_{I,J}) = B_{J,I}$, it follows
that Bott-Samelsons are also $\tDl$-self-dual.
Hence, $\tD(M) \cong \tDl(M)$
for a singular Soergel bimodule $M$.
\end{rem}

Let $\preceq$ be a total order on $(W_I \backslash W / W_J)_{\min}$ refining the Bruhat order $\leq$.
For $d \in (W_I \backslash W / W_J)_{\min}$
and an $(R^I, R^J)$-bimodule $M$,
we define the subquotients
$\Gamma_{\preceq d} M / \Gamma_{\prec d}M, \Gamma_{\succeq d}M / \Gamma_{\succ d}M$
of $M$ in terms of supports
as explained in \cite[Sec.~4.5]{Will}.
\begin{itemize}
\item
Following \cite[Def.~6.1]{Will}, a graded $(R^I, R^J)$-bimodule $M$
has a {\em $\nabla$-flag} if
it is finitely generated both as a left $R^I$-module
and as a right $R^J$-module, and
there exist 
Laurent polynomials $(M:\nabla_{IdJ})_q \in \N[q,q^{-1}]$ such that
$$
\Gamma_{\preceq d} M / \Gamma_{\prec d} M
\cong (M:\nabla_{IdJ})_q\:\nabla_{IdJ}
$$ 
for all $d \in (W_I \backslash W / W_J)_{\min}$, 
with $(M:\nabla_{IdJ})_q = 0$ for all but finitely many $d$.
Assuming $M$ has a $\nabla$-flag, we define its {\em $\nabla$-character}
to be
\begin{equation}
\ch_\nabla(M) := \sum_{d \in (W_I \backslash W / W_J)_{\min}}
(M:\nabla_{IdJ})_q\:
h_{IdJ} \in {_{I\hspace{.8pt}}}\Schur_J.
\end{equation}
\item
Following \cite[Def.~6.12]{Will},
a graded $(R^I, R^J)$-bimodule $M$ has a 
{\em $\Delta$-flag} if it is finitely generated both as a left $R^I$-module
and as a right $R^J$-module, and
there exist Laurent polynomials
$(M:\Delta_{IdJ})_q \in \N[q,q^{-1}]$ such that
$$
\Gamma_{\succeq d} M / \Gamma_{\succ d} M
\cong (M:\Delta_{IdJ})_q\:\Delta_{IdJ}$$ 
for all $d \in (W_I \backslash W / W_J)_{\min}$, with $(M:\Delta_{IdJ})_q = 0$ for all but finitely many $d$.
Assuming $M$ has a $\Delta$-flag, its {\em $\Delta$-character} is
\begin{equation}
\ch_\Delta(M) := \sum_{d \in (W_I \backslash W / W_J)_{\min}}
(M:\Delta_{IdJ})_q\:
\overline{h_{IdJ}} \in {_{I\hspace{.8pt}}}\Schur_J.
\end{equation}
\end{itemize}
The notion of a bimodule $M$ possessing a $\nabla$- or $\Delta$-flag, and its
character $\ch_\nabla(M)$ or $\ch_\Delta(M)$, is 
independent of the choice of the total order $\preceq$ 
thanks to the hin-und-her lemmas in \cite{Will}.

\begin{rem}
Our $\ch_\nabla$ is the same as in \cite[Sec.~0.4]{Willerratum},
but our $\ch_\Delta$ is the one there post-composed with the bar involution.
\end{rem}

By \cite[Prop.~6.16]{Will}, 
the duality $\tD$ interchanges the two sorts of flag. Moreover, using also \cref{standardduality}, we have that
\begin{equation}\label{harmless}
\overline{\ch_\Delta(M)}
= \ch_\nabla(\tD(M))
\end{equation}
on bimodules with $\Delta$-flags.

Now we can state the main results.
These identify the $q$-Schur algebra $\Schur$ with  $K_0(\SBim)$,
the split Grothendieck ring 
consisting of $\cong$-classes of singular Soergel bimodules.

\begin{theo}[Homomorphism formula {\cite[Th.~7.9]{Will}}]\label{homformula}
Let $M$, $M'$ be graded $(R^I, R^J)$-bimodules such that either $M$ is a singular Soergel bimodule and $M'$ has a $\nabla$-flag, or $M$ has a $\Delta$-flag
and $M'$ is a singular Soergel bimodule.
There is an isomorphism
$$
\Hom_{R^I\dash R^J}(M,M')
\cong
\big(\overline{\ch_\Delta(M)}, 
\ch_\nabla (M')\big) R^J
$$
of graded $R^J$-modules.
\end{theo}

\begin{theo}[Classification of indecomposables {\cite[Th.~7.10]{Will}}]\label{class}
For $d \in (W_I \backslash W / W_J)_{\min}$,
there is a unique (up to isomorphism)
indecomposable singular Soergel bimodule
$B_{IdJ}$ in $\one_I \SBim \one_J$
such that $B_{IdJ} / \Gamma_{\prec d} B_{IdJ}
\cong \nabla_{IdJ}$.
We have that $\tD(B_{IdJ}) \cong B_{IdJ}$
and $\tR(B_{IdJ}) \cong B_{J d^{-1} I}$.
Moreover, 
these bimodules for all $I,J \finsub N$
and $d \in (W_I \backslash W / W_J)_{\min}$
give a full set of pairwise inequivalent indecomposable 
singular Soergel bimodules (up to grading shifts).
\end{theo}

\begin{theo}[Categorification theorem {\cite[Th.~7.12]{Will}}]\label{catthm}
If $M$ is a singular Soergel bimodule 
then it has both a $\Delta$-flag and a $\nabla$-flag, and
\begin{equation}
\ch(M) :=  
\ch_\nabla(M)=\ch_\Delta(M).
\end{equation}
The map $\ch$ induces a $\Z[q,q^{-1}]$-algebra isomorphism
$\ch:K_0(\SBim)
\stackrel{\sim}{\rightarrow} \Schur$
taking 
the isomorphism class of 
the graded $(R^I, R^I)$-bimodule $R^I$
to $1_I$ for each $I \subseteq N$. 
It intertwines the anti-linear involution
of $K_0(\SBim)$ arising from the duality $\tD$ with the bar involution on $\Schur$,
and it intertwines the linear anti-involution arising from the symmetry $\tR$ with $\rho$.
\end{theo}

From these results, it follows that 
$\ch(B_{IdJ})$
is a bar-invariant element of 
${_{I\hspace{.8pt}}}\Schur_J$ which is equal to $h_{IdJ}$
plus an $\N[q,q^{-1}]$-linear combination of
$h_{Id'J}$ for $d' < d$.
Also there is $b_{IdJ}$, which is a bar-invariant element of
${_{I\hspace{.8pt}}}\Schur_J$
equal to 
$h_{IdJ}$ plus a $q\Z[q]$-linear combination of $h_{I d' J}$ for $d' < d$.
We say that {\em Soergel's conjecture holds for $W_I\backslash W / W_J$}
if $\ch(B_{IdJ}) = b_{IdJ}$ for all $d \in (W_I\backslash W / W_J)_{\min}$.

\begin{lem}[{\cite[Prop.~7.11]{Will}}]
Suppose that $I, J \finsub N$ and $I' \subseteq I$.
For $d \in 
(W_I \backslash W / W_J)_{\min}$, we have that
\begin{align}
\Ind_{I}^{I'} B_{IdJ} 
&= B_{I',I} \otimes_{R^I} B_{IdJ}
\cong
B_{I'd'J},\\
\Res^{I'}_I B_{I' d' J} &=
q^{\ell(w_I)-\ell(w_{I'})}
B_{I,I'} \otimes_{R^{I'}} B_{I'd'J}
\cong
\frac{\pi^+_I}{\pi^+_{I'}} 
B_{IdJ},\label{hamilton}
\end{align}
where
$d'$ is the longest element of 
$(W_{I'} \backslash W / W_{J})_{\min} \cap W_I d W_J$.
\end{lem}

The following corollary is discussed in the final sentence of \cite{Will}.

\begin{cor}\label{notours}
Suppose that $I, J \finsub N$, $I' \subseteq I$ and $J' \subseteq J$.
If Soergel's conjecture holds for $W_{I'}\backslash W / W_{J'}$
then it holds for $W_I\backslash W / W_J$.
\end{cor}

\begin{proof}
Using the symmetry $\tR$, it suffices to show that
Soergel's conjecture holds 
for $W_I \backslash W / W_J$ if it holds for
$W_{I'}\backslash W / W_J$.
This follows because
$$
\ch(B_{I d J})
= \frac{\pi_{I'}}{\pi_{I}} 
\ch (B_{I,I'} \otimes_{R^{I'}} B_{I' d' J})
= \frac{\pi_{I'}}{\pi_{I}} 
\ch(B_{I,I'}) * \ch(B_{I'd'J})
= \frac{\pi_{I'}}{\pi_{I}} 
b_{I,I'} * b_{I'd'J}
= b_{IdJ}
$$
for $d \in (W_I \backslash W / W_J)_{\min}$
and $d'$ that is
the longest element of 
$(W_{I'} \backslash W / W_{J})_{\min} \cap W_I d W_J$.
Here, we used \cref{hamilton} for the first equality, \cref{catthm}
for the second, and the assumption that Soergel's conjecture holds for $(I',J)$ 
for the third equality.
The final equality
follows because 
$w_{I d J} = w_{I' d' J}$,
hence, we have that $\frac{\pi_{I'}}{\pi_{I}}
b_{I,I'} * b_{I'd'J}=\frac{1}{\pi_{I}}
b_{w_I} b_{w_{I'd'J}} =
\frac{1}{\pi_{I}} b_{w_I} b_{w_{IdJ}}=1_I * b_{IdJ}=b_{IdJ}$.
\end{proof}

\begin{rem}
If $\kk$ is a field of characteristic 0
then Soergel's conjecture holds for $W_I \backslash W / W_J$
for {\em all} $I,J \finsub N$.
When $\kk = \mathbb{R}$, this follows
from \cite{EW} and \cref{notours}.
To deduce it for other fields of characteristic 0, 
the Soergel conjecture is equivalent to the nondegeneracy of certain local intersection pairings (see [EW14]), which is unaffected by base change between fields of characteristic 0; this follows from Soergel's homomorphism formula which 
implies that the Soergel category is flat over $\Q$. 
%For the very special case of this result that is relevant for our purposes, we will give a self-contained proof in \cref{bigtheorem} below, via an argument which is valid also in characteristic $p > 2$.
\end{rem}

%% file: s4-typeD.tex
\setcounter{section}{3}
\section{Extended singular Soergel bimodules in type D\texorpdfstring{$_l$}{}}\label{secD}

Now we specialize to the case of interest in 
this article: singular Soergel bimodules for the Weyl group $W$ of type D$_l$. We will actually work with a variant which we call {\em extended} singular Soergel bimodules, the purpose of which we explain in \cref{ssec-interlude}.
We also introduce an extended version of the $q$-Schur algebra, and compute
explicitly the Kazhdan-Lusztig basis of the part of it that relates to nil-Brauer. Finally, we use the extended diagrammatic calculus to derive some 
difficult relations.

\subsection{Realization of the root system of type \texorpdfstring{$\operatorname{D}_l$}{}}\label{ssec4-1}
Fix $l \geq 2$
and consider the root system of type $\operatorname{D}_l$. We label the nodes of the Dynkin diagram
by the set $N := \{\pm 1, 2, \dots, l-1\}$ as in the introduction.
We adopt the usual reflection realization of the corresponding Coxeter group $W$,
which satisfies all of the hypotheses assumed in \cref{ssec3-1}.
So we let $\mathfrak{h}$ be the vector space with basis $x_1,\dots,x_l$,
and identify $\mathfrak{h}^*$ with $\mathfrak{h}$ via the non-degenerate symmetric bilinear form 
defined by declaring that this basis is orthonormal.
The root system $\Phi$ is
$\{\pm x_j \pm x_i\:|\:1 \leq i, j \leq l, i \neq j\}$,
and we have simply that $\alpha^\vee = \alpha$ for each $\alpha \in \Phi$.
We make the following choice for simple roots:
\begin{align*}
\begin{array}{rccccc}
\alpha_1:=x_2-x_1\\
&
\alpha_2 &
:= x_3 - x_2&
\alpha_3 := x_4-x_3&\quad\cdots\quad
&\alpha_{l-1} := x_{l}-x_{l-1}\\
\alpha_{-1}:=x_2+x_1
\end{array}
\end{align*}
The corresponding set of positive roots
is $\{x_j \pm x_i\:|\:1 \leq i< j \leq l\}$.
For $i \in N$, the simple reflection $s_i$ is the reflection in the hyperplane orthogonal to $\alpha_i$.
Except when $i = -1$, 
 $s_i$ permutes $x_i$ and $x_{i+1}$ fixing all other $x_j$, while $s_{-1}$ permutes $x_1$ and $-x_2$ fixing all other $x_j$.
 
A useful additional feature is the existence of the {\em graph automorphism}
$\gamma:N\rightarrow N$.
This switches $1$ and $-1$ and fixes all other elements.
Let $s_0 : \mathfrak{h} \rightarrow \mathfrak{h}$ be the reflection with $s_0(x_1) = -x_1$, fixing all other $x_j$. This corresponds to $\gamma$ 
in that $s_0(\alpha_i) = \alpha_{\gamma(i)}$ for each $i \in N$. The reflections 
$s_0,s_1,\dots,s_{l-1}$ generate the Weyl group of 
type B$_l$, with the Weyl group $W$ of type D$_l$
being the subgroup generated by $s_{-1} = s_0 s_1 s_0,
s_1,\dots,s_{l-1}$.
All $s_i$
 extend to algebra automorphisms of $R=\kk[x_1,\dots,x_l]$ (defined in the previous section to be the
  symmetric algebra of $\mathfrak{h}^*$).
  We prefer to denote the automorphism of $R$ defined by 
  $s_0$ by
  \begin{align}\label{prettydumbconvention}
  \gamma:R &\rightarrow R&f &\mapsto s_0(f).
  \end{align}

We reserve the unusual letter $\O$ for the special subset
$\{2,\dots,l-1\}\subset N$, since it appears very often. Remember that we write simply $Ii$ for $I \cup \{i\}$
(assuming $i \notin I$) and $I \hat\imath$ for $I - \{i\}$ (assuming $i \in I$).
Mainly 
because they look nicer in diagrams, we will use the following shorthands:
\begin{align}\label{sledgehammer}
I+ &:= I1, &\quad I- &:= I (-1),&I\pm &:= I 1 (-1)
\end{align}
for $I \subseteq \O$. So we are using the symbol $+$ to indicate inclusion of $1$, the symbol $-$ to indicate inclusion of $-1$, and $\pm$ to indicate the inclusion of both. The graph automorphism $\gamma$ interchanges $+$ and $-$.
Recall also for $I \subseteq N$ that $R^I$ denotes the invariant subalgebra 
$R^{W_I}$ of $R$.
The following table gives some useful facts about some of these algebras and the corresponding parabolic subgroups of $W$, for $n \in \O$:
\begin{equation}\label{lengths}
\begin{array}{|l|l|l|l|}
\hline
I&\text{Type of $W_I$}&\ell(w_I)&\pi_I \\\hline
\O\pm&D_l&2\binom{l}{2}&[2]_q^{l-1} [l-1]^!_{q^2} [l]_q\\
\O{+}\text{ or } \O{-}&A_{l-1}&\binom{l}{2}&[l]_q^!\\
\O&A_{l-2}&\binom{l-1}{2}&[l-1]_q^!\\\hline
\O{\hat n\pm} & A_{l-n-1}\times D_{n}&
\binom{l-n}{2}+2\binom{n}{2}&[l-n]_q^![2]_q^{n-1} [n-1]^!_{q^2} [n]_q\\
\O{\hat n+}\text{ or }\O{\hat n-} & A_{l-n-1}\times A_{n-1}&\binom{l-n}{2}+\binom{n}{2}&[l-n]_q^![n]_q^!\\
\O{\hat n}& A_{l-n-1}\times A_{n-2}&\binom{l-n}{2}+\binom{n-1}{2}&[l-n]_q^![n-1]^!_q\\
\hline
\end{array}
\end{equation}
Also, here are explicit descriptions of some of the algebras of invariants:
\begin{itemize}
\item
$R^{N}$ is free with generators given by
the elementary symmetric polynomials $e_r(x_1^2,\dots,x_l^2)$ for $r=1,\dots,l-1$
together with $e_l(x_1,\dots,x_l) = x_1 \cdots x_l$.
\item $R^{\O+}$ is the 
algebra $\Sym$ of symmetric polynomials which appeared already in the introduction and \cref{simonriche}.
It is freely generated by 
$e_r(x_1,\dots,x_l)$ for $r=1,\dots,l$.
\item
$R^\O = \kk[x_1,x_2,\dots,x_{l}]^{S_1 \times S_{l-1}}$, which is generated by $x_1$ and
$e_r(x_2,\dots,x_l)$ for $r=1, \dots, l-1$.
This was seen in the introduction when the $(A,A)$-bimodule $B$ was defined.
\end{itemize}

\subsection{Motivational interlude} \label{ssec-interlude}
We will be interested primarily in double cosets
$W_I \backslash W / W_J$ where $I, J \in \{O+, O-\}$, and in the corresponding categories of singular Soergel bimodules, that is, we wish to study in depth the sub-bicategory of $\SBim$ with two objects $\{R^{O+}, R^{O-}\}$, and the corresponding piece of the $q$-Schur algebra 
$\Schur$. Our goal is to relate these bimodules to the cyclotomic nilBrauer category $\CNB_l$,
but there is some awkward book-keeping involved in this since we are comparing a monoidal category with a bicategory with two objects.

The approach we found to be most convenient is to transform any $R^{O-}$-module into an $R^{O+}$-module using the graph automorphism. In this way, we are able to replace the aforementioned bicategory-with-two-objects simply with a monoidal category of $(R^{O+}, R^{O+})$-bimodules. See \cref{ssec4-3} for details. The Grothendieck group of this formal construction is a piece of the extended $q$-Schur algebra described in the next section. The general approach taken also fits well with the iSchur-Weyl duality from \cite[Ch.~5]{BW18KL}, although we will not discuss this further here.

An alternate approach would be to instead replace $\cNB_t$ with a 2-category with two objects. This is analogous to the use of the two-colored Temperley-Lieb algebra in \cite{Esatake}. There are advantages to this approach, e.g., one can eliminate the pesky sign in the first relations from \cref{rels4,rels4b} by rescaling the dot based on the ambient colors. 
However, the disadvantage of having to keep track of two versions of each bimodule seemed to us to be more burdensome.
For example, in the two-colored approach, one has to work with two generating bimodules
$B_{O-,O} \otimes_{R^{O}} 
B_{O,O+}$ and 
$B_{O+,O} \otimes_{R^{O}} 
B_{O,O-}$. In the extended category, due to the added twists, they both yield the same bimodule over $(R^{O+},R^{O+})$.

\subsection{The extended \texorpdfstring{$q$}{}-Schur algebra}\label{ssec4-2}

Let $\Schur$ be the $q$-Schur algebra of type D$_l$
as in \cref{heckealgebroid}.
The graph automorphism defines a $\Z[q,q^{-1}]$-algebra involution $\gamma:\Hecke \rightarrow \Hecke,
h_i \mapsto h_{\gamma(i)}$ of the Hecke algebra,
which extends to an automorphism $\gamma:\Schur \rightarrow \Schur$ 
taking $1_I$ to $1_{\gamma(I)}$.
The {\em extended $q$-Schur algebra} is the $\Z[q,q^{-1}]$-module $\Schur \oplus s_0 \Schur$, 
where $s_0 \Schur$ is a copy of $\Schur$ with elements denoted $s_0 x\:(x \in \Schur)$, viewed as a $\Z[q,q^{-1}]$-algebra with the extended multiplication $\m$ defined by
\begin{equation}\label{likeme}
(x + s_0 y) \m (u+s_0 v) := 
\left(x\m u + \gamma(y) \m v\right)
+s_0 \left(y\m u + \gamma(x) \m v\right)
\end{equation}
for $x,y,u,v \in \Schur$.
From the standard and Kazhdan-Lusztig bases for $\Schur$,
we obtain standard and Kazhdan-Lusztig bases also for $\Schur \oplus s_0 \Schur$. These consist of the elements
$\{h_{IdJ}, s_0 h_{IdJ}\}$ or the elements
$\{b_{IdJ}, s_0 b_{IdJ}\}$ for $I,J \subseteq N$
and $d \in (W_I \backslash W / W_J)_{\min}$.
The Kazhdan-Lusztig basis is invariant under the bar involution on $\Schur \oplus s_0 \Schur$, which is
defined simply by 
\begin{equation}
\overline{x+s_0 y} := \overline{x} + s_0 \overline{y}.
\end{equation}
The linear anti-involution $\rho$ from \cref{tau}
extends to $\Schur \oplus s_0 \Schur$ by setting
\begin{equation}
\rho(x+s_0 y) := \rho(x) + s_0 \gamma(\rho(y)).
\end{equation}
Finally, we extend the symmetric bilinear form on $\Schur$ from \cref{altbilform} to $\Schur \oplus s_0\Schur$ by 
\begin{equation}
(x+s_0 y, u + s_0 v) := (x,u)+(y,v).
\end{equation}

\begin{rem}
The Weyl group $W$, its Hecke algebra $\Hecke$ and the $q$-Schur algebra $\Schur$ here are of type D$_l$.
The Weyl group of type B$_l$ is 
$W \sqcup s_0 W$,
and $\Hecke \oplus s_0 \Hecke$ with twisted multiplication defined like in \cref{likeme} is the Hecke algebra of
type B$_l$ at unequal parameters, with the short root parameter being 1.
Similarly, our extended $q$-Schur algebra $\Schur \oplus s_0 \Schur$ can be related to the $(q,1)$-Schur algebra
of type B$_l$.
We will not use this observation here so leave the formulation of a more precise statement to the reader; the construction in \cite{Bao} is also somewhat relevant.
\end{rem}

\subsection{The subalgebra \texorpdfstring{$\V$}{} and its Kazhdan-Lusztig basis}
In the remainder of the paper, the focus will be on
the subalgebra
\begin{equation}\label{thisisbetter}
\V := 1_{\O+} \m (\Schur \oplus s_0 \Schur) \m 1_{\O+}=
{_{\O+}}\Schur_{\O+} \oplus s_0 {_{\O-}} \Schur_{\O+}
\end{equation}
of $\Schur\oplus s_0 \Schur$, 
which has the idempotent $1_{\O+}$ as its identity element.
We are now going to calculate its Kazhdan-Lusztig basis explicitly.
In fact, we will show that it is combinatorially the same as the icanonical basis of
$\V(l)$ seen already in \cref{crazy}.
We denote the special element $s_0 b_{\O-,\O+} \in \V$
simply by $b$ (it will turn out to be ``the same'' as the $b$ in \cref{secnilbrauer}):
\begin{equation}\label{killingme}
b := s_0 b_{\O-,\O+} \stackrel{\cref{bselements}}{=}
s_0 b_{\O-,\O} \m b_{\O,\O+}.
\end{equation}
Since $\gamma(b_{\O-,\O+}) = b_{\O+,\O-}$,
we have that
\begin{equation}
b^{\m n} = 
\label{wrongspot}
\begin{cases}
1_{\O+}&\text{if $n = 0$}\\
(b_{\O+,\O-}\m b_{\O-,\O+})^{\m \frac{n}{2}}
&\text{if $n$ is even}\\
s_0 b_{\O-,\O+}\m
(b_{\O+,\O-} \m b_{\O-,\O+})^{\m \frac{n-1}{2}}&\text{if $n$ is odd}.
\end{cases}
\end{equation}
By \cref{cruisier} below, the elements 
$b^{\m n}\:(0 \leq n \leq l)$
give a rational basis for $\V$, i.e., they are a basis for the $\Q(q)$-vector space $\Q(q) \otimes_{\Z[q,q^{-1}]} \V$.

We will start now to draw elements of $W \sqcup s_0 W$ by permutation-type diagrams with symmetry about the middle axis, labeling the boundary in order from left to right by $-x_{l},\dots,-x_2,-x_1;x_1,x_2,\dots,x_{l}$.
For example, the following are the pictures for 
$s_0, s_1$ and $s_2$ when $l=3$:
\begin{align*}
s_0 &=
\begin{tikzpicture}[anchorbase,scale=2.5]
\draw[dashed,gray] (0.5,.6) to (0.5,-.1);
\draw[semithick] (0,.5) to (0,0);
\draw[semithick] (0.2,.5) to (.2,0);
\draw[semithick] (0.4,.5) to (.6,0);
\draw[semithick] (0.6,.5) to (.4,0);
\draw[semithick] (0.8,.5) to (.8,0);
\draw[semithick] (1,.5) to (1,0);
\node at (-.1,.55) {$\stringlabel{-x_3}$};
\node at (.13,.55) {$\stringlabel{-x_2}$};
\node at (.34,.55) {$\stringlabel{-x_1}$};
\node at (.6,.55) {$\stringlabel{x_1}$};
\node at (.8,.55) {$\stringlabel{x_2}$};
\node at (1,.55) {$\stringlabel{x_3}$};
\node at (-.1,-.05) {$\stringlabel{-x_3}$};
\node at (.13,-.05) {$\stringlabel{-x_2}$};
\node at (.34,-.05) {$\stringlabel{-x_1}$};
\node at (.6,-.05) {$\stringlabel{x_1}$};
\node at (.8,-.05) {$\stringlabel{x_2}$};
\node at (1,-.05) {$\stringlabel{x_3}$};
\end{tikzpicture}
,&
s_1 &=
\begin{tikzpicture}[anchorbase,scale=2.5]
\draw[semithick] (0,.5) to (0,0);
\draw[dashed,gray] (0.5,.6) to (0.5,-.1);
\draw[semithick] (0.2,.5) to (.4,0);
\draw[semithick] (0.4,.5) to (.2,0);
\draw[semithick] (0.6,.5) to (.8,0);
\draw[semithick] (0.8,.5) to (.6,0);
\draw[semithick] (1,.5) to (1,0);
\node at (-.1,.55) {$\stringlabel{-x_3}$};
\node at (.13,.55) {$\stringlabel{-x_2}$};
\node at (.34,.55) {$\stringlabel{-x_1}$};
\node at (.6,.55) {$\stringlabel{x_1}$};
\node at (.8,.55) {$\stringlabel{x_2}$};
\node at (1,.55) {$\stringlabel{x_3}$};
\node at (-.1,-.05) {$\stringlabel{-x_3}$};
\node at (.13,-.05) {$\stringlabel{-x_2}$};
\node at (.34,-.05) {$\stringlabel{-x_1}$};
\node at (.6,-.05) {$\stringlabel{x_1}$};
\node at (.8,-.05) {$\stringlabel{x_2}$};
\node at (1,-.05) {$\stringlabel{x_3}$};
\end{tikzpicture},
&
s_2 &=
\begin{tikzpicture}[anchorbase,scale=2.5]
\draw[dashed,gray] (0.5,.6) to (0.5,-.1);
\draw[semithick] (0,.5) to (0.2,0);
\draw[semithick] (0.2,.5) to (0,0);
\draw[semithick] (0.4,.5) to (.4,0);
\draw[semithick] (0.6,.5) to (.6,0);
\draw[semithick] (0.8,.5) to (1,0);
\draw[semithick] (1,.5) to (.8,0);
\node at (-.1,.55) {$\stringlabel{-x_3}$};
\node at (.13,.55) {$\stringlabel{-x_2}$};
\node at (.34,.55) {$\stringlabel{-x_1}$};
\node at (.6,.55) {$\stringlabel{x_1}$};
\node at (.8,.55) {$\stringlabel{x_2}$};
\node at (1,.55) {$\stringlabel{x_3}$};
\node at (-.1,-.05) {$\stringlabel{-x_3}$};
\node at (.13,-.05) {$\stringlabel{-x_2}$};
\node at (.34,-.05) {$\stringlabel{-x_1}$};
\node at (.6,-.05) {$\stringlabel{x_1}$};
\node at (.8,-.05) {$\stringlabel{x_2}$};
\node at (1,-.05) {$\stringlabel{x_3}$};
\end{tikzpicture}.
\end{align*}
For $0 \leq n \leq l$, let 
\begin{align*}
d_n &:=
\begin{tikzpicture}[baseline=0,scale=1.5]
\draw[dashed,gray] (0,.6) to (0,-.6);
\draw[semithick] (-1.2,.5) to (-1.2,-.5);
\draw[semithick] (-1.4,.5) to (-1.4,-.5);
\draw[semithick] (-1,.5) to (-1,-.5);
\draw[semithick] (-.8,.5) to (.2,-.5);
\draw[semithick] (-.6,.5) to (.4,-.5);
\draw[semithick] (-.4,.5) to (.6,-.5);
\draw[semithick] (-.2,.5) to (.8,-.5);
\draw[semithick] (1.2,.5) to (1.2,-.5);
\draw[semithick] (1.4,.5) to (1.4,-.5);
\draw[semithick] (1,.5) to (1,-.5);
\draw[semithick] (.8,.5) to (-.2,-.5);
\draw[semithick] (.6,.5) to (-.4,-.5);
\draw[semithick] (.4,.5) to (-.6,-.5);
\draw[semithick] (.2,.5) to (-.8,-.5);
\draw [blue,thick,decoration={brace,mirror},decorate] (1,-.55) -- (1.4,-.55) node [pos=0.5,anchor=north,yshift=-0.5] {$\color{blue}\scriptstyle l-n$}; 
\draw [blue,thick,decoration={brace,mirror},decorate] (.2,-.55) -- (.8,-.55) node [pos=0.5,anchor=north,yshift=-0.5] {$\color{blue}\scriptstyle n$}; 
\draw [blue,thick,decoration={brace,mirror},decorate] (-1.4,-.55) -- (-1,-.55) node [pos=0.5,anchor=north,yshift=-0.5] {$\color{blue}\scriptstyle l-n$}; 
\draw [thick,blue,decoration={brace,mirror},decorate] (-.8,-.55) -- (-.2,-.55) node [pos=0.5,anchor=north,yshift=-0.5] {$\color{blue}\scriptstyle n$}; 
\end{tikzpicture}
\end{align*}
Assuming $1 \leq n \leq l$, we have that
$$
s_0 d_n = \begin{tikzpicture}[baseline=0,scale=1.5]
\draw[dashed,gray] (0,.6) to (0,-.6);
\draw[semithick] (-1.4,.5) to (-1.4,-.5);
\draw[semithick] (-1.6,.5) to (-1.6,-.5);
\draw[semithick] (-1.2,.5) to (-1.2,-.5);
\draw[semithick] (-1,.5) to (.2,-.5);
\draw[semithick] (-.8,.5) to (.4,-.5);
\draw[semithick] (-.6,.5) to (.6,-.5);
\draw[semithick] (-.4,.5) to (.8,-.5);
\draw[semithick] (.2,.5) to [out=-90,in=135](1,-.5);
\draw[semithick] (1.4,.5) to (1.4,-.5);
\draw[semithick] (1.6,.5) to (1.6,-.5);
\draw[semithick] (1.2,.5) to (1.2,-.5);
\draw[semithick] (1,.5) to (-.2,-.5);
\draw[semithick] (.8,.5) to (-.4,-.5);
\draw[semithick] (.6,.5) to (-.6,-.5);
\draw[semithick] (.4,.5) to (-.8,-.5);
\draw[semithick,-] (-.2,.5) to [out=-90,in=45] (-1,-.5);
\draw [thick,blue,decoration={brace,mirror},decorate] (1.2,-.55) -- (1.6,-.55) node [pos=0.5,anchor=north,yshift=-0.5] {$\color{blue}\scriptstyle l-n$}; 
\draw [thick,blue,decoration={brace,mirror},decorate] (.2,-.55) -- (.8,-.55) node [pos=0.5,anchor=north,yshift=-0.5] {$\color{blue}\scriptstyle n-1$}; 
\draw [thick,blue,decoration={brace,mirror},decorate] (-1.6,-.55) -- (-1.2,-.55) node [pos=0.5,anchor=north,yshift=-0.5] {$\color{blue}\scriptstyle l-n$}; 
\draw [thick,blue,decoration={brace,mirror},decorate] (-.8,-.55) -- (-.2,-.55) node [pos=0.5,anchor=north,yshift=-0.5] {$\color{blue}\scriptstyle n-1$}; 
\end{tikzpicture}
$$
Note that $d_n \in W$ if and only if $n$ is even,
and $s_0 d_n \in W$ if and only if $n$ is odd.
Hence, $s_0^n d_n \in W$ for all $0 \leq n \leq l$.
It is worthwhile to verify that
$\ell\big(s_0^n d_n\big) = \binom{n}{2}$, using the length function in type D$_l$.

\begin{lem}\label{pure}
The set $(W_{\O} \backslash W / W_{\O+})_{\min}$
is equal to $\{s_0^n d_n\:|\:1 \leq n \leq l\}$.
Also $$
\LR{\O}{s_0^n d_n}{(\O+)}
=
\begin{cases}\O&\text{if $n=1$ or $n=l$}\\
\O\, \hat n&\text{for $2\leq n\leq l-1$.}
\end{cases}
$$
\end{lem}

\begin{proof}
Consider the left action of $W$
on the set $X$ of sign sequences
$(\sigma_1,\dots,\sigma_{l}) \in \{+,-\}^l$ 
defined so that $s_i$ switches $\sigma_i$ and $\sigma_{i+1}$
for $i=1,\dots,l-1$ and
$$
s_{-1} \cdot (\sigma_1,\sigma_2,\dots,\sigma_{l})
= (-\sigma_1,-\sigma_2,\sigma_3,\dots,\sigma_{l}).
$$
The stabilizer of $(+^l) = (+,\dots,+)$ is the subgroup $W_{\O+}$, and the $W$-orbit of this point is the set $X$ of all sign sequences with an even number of minus signs.
The corresponding 
orbit map identifies $W / W_{\O+}$
with $X$. A set of representatives
for the orbits of $W_\O$ on $X$
are given by the sequences
$$
s_0^n d_n \cdot (+^l)
= 
\begin{cases}
(-^{n},+,^{l-n})&\text{if $n$ is even}\\
(+,-^{n-1},+,^{l-n})&\text{if $n$ is odd,}
\end{cases}
$$ 
for $1 \leq n \leq l$.
This shows that $s_0^n d_n\:(1 \leq n \leq l)$
is a full set of $W_\O \backslash W / W_{\O+}$-double coset representatives.
Each one clearly gets longer if one acts on the left or right by some $s_i\:(2 \leq i \leq l-1)$
or on the right by $s_1$,
hence, these are the minimal length double coset representatives.
Finally, to compute the left redundancy, $W_\O \cap s_0^n d_n W_{\O+} \big(s_0^n d_n\big)^{-1}$ is the stabilizer 
in $W_\O$ of the point
$s_0^n d_n \cdot (+^l)$, which is easily computed.
\end{proof}

\begin{lem}\label{pure1}
The set
$(W_{\O+} \backslash W / W_{\O+})_{\min}$ equals $\{d_n\:|\:0 \leq n \leq l\text{ with }n\text{ even}\}$.
For these values of $n$, we have that
$$
\LR{(\O+)}{d_n}{(\O+)}=
\begin{cases}
\O+&\text{if $n=0$ or $n=l$}\\
\O\hat n+&\text{otherwise.}
\end{cases}
$$
Also
\begin{align*}
W_{\O+} d_0 W_{\O+}&=W_{\O} d_0 W_{\O+},\\
W_{\O+} d_n W_{\O+}&=W_\O d_n W_{\O+} \sqcup W_\O s_0 d_{n} W_{\O+}&&\text{for even $2 \leq n \leq l-1$,}\\
W_{\O+} d_l W_{\O+}&=W_\O d_l W_{\O+}&&\text{assuming $l$ is even.}
\end{align*}
\end{lem}

\begin{proof}
This follows in a similar way to the proof of \cref{pure}, considering 
$W_{\O+}$-orbits on the set $X$ introduced there.
\end{proof}

\begin{lem}\label{pure2}
The set
$(W_{\O-} \backslash W / W_{\O+})_{\min}$ equals $\{s_0 d_n\:|\:1 \leq n \leq l\text{ with }n\text{ odd}\}$.
For these values of $n$, we have that
\begin{align*}
\LR{(\O-)}{s_0 d_n}{(\O+)}&=
\begin{cases}
\O&\text{if $n=1$}\\
\O-&\text{if $n=l$}\\
\O\hat n-&\text{otherwise.}
\end{cases}
\end{align*}
Also
\begin{align*}
W_{\O-} s_0 d_n W_{\O+}&=W_\O s_0 d_n W_{\O+}\sqcup W_\O d_{n+1} W_{\O+}&&
\text{for odd $1 \leq n \leq l-1$,}\\
W_{\O-} s_0 d_l W_{\O+}&=W_\O s_0 d_l W_{\O+}
&&\text{when $l$ is odd.}
\end{align*}
\end{lem}

\begin{proof}
Similar arguments to the previous two lemmas.
\end{proof}

Recall the $\Ui$-module $\V(l)$ from \cref{ssec2-6}
with its monomial basis $f^{(n)} \eta_l\:(0 \leq n \leq l)$
and icanonical basis $b^{(n)} \eta_l\:(0 \leq n \leq l)$,
and the $\Z[q,q^{-1}]$-algebra $\V$ from \cref{thisisbetter}, with its standard and Kazhdan-Lusztig bases.
Both spaces have also been equipped with a bar involution and a symmetric bilinear form.

\begin{theo}\label{asleftmodule}
There is a unique
isomorphism of free $\Z[q,q^{-1}]$-modules
\begin{align}
\phi:\V(l) &\stackrel{\sim}{\rightarrow}
\V,\notag\\
\overline{f^{(n)} \eta_l}&\mapsto
\begin{cases}
h_{(\O+) d_n (\O+)}&\text{for even $0 \leq n \leq l$}\\
s_0 h_{(\O-) s_0 d_n (\O+)}&\text{for odd $1 \leq n \leq l$,}
\end{cases}\label{rain1}\\
b^{(n)} \eta_l&\mapsto
\begin{cases}
b_{(\O+) d_n (\O+)}&\text{for even $0 \leq n \leq l$}\\
s_0 b_{(\O-) s_0 d_n (\O+)}&\text{for odd $1 \leq n \leq l$.}
\end{cases}\label{rain2}
\end{align}
In fact, viewing $\V(l)$ as a $\Z[q,q^{-1}]$-algebra by identifying it with the quotient $\Ui / \I_l$ as explained in \cref{oldpeculiar}, $\phi$ is an algebra isomorphism.
Also,
$\phi(bv) = b \m \phi(v)$,
$\phi(\overline{v}) = \overline{\phi(v)}$,
and $\overline{(\overline{v},\overline{w})_l} = 
\big(\phi(v), \phi(w)\big)$
for all $v,w \in \V(l)$.
\end{theo}

\begin{proof}
We define $\phi$ by \cref{rain1}. It follows immediately that $\phi$ is a $\Z[q,q^{-1}]$-module isomorphism since it takes a basis to a basis.
First we claim that $\phi(bv) = b \phi(v)$ for any $v \in \V(l)$. Note that the claim 
implies that $\phi$ is also an algebra isomorphism since $\Q(q) \otimes_{\Z[q,q^{-1}]} \Ui$ is generated by $b$.
To prove it,
we may assume that $v = \overline{f^{(n)} \eta_l}$
for some $0 \leq n \leq l$.
We have
that $b_{\O-,\O+}=b_{\O-,\O}\m b_{\O,\O+}$ and 
$b_{\O+,\O-}=
b_{\O+,\O}\m b_{\O,\O-}$.
If $n$ is even then, in view also of \cref{rakes}, 
the proof reduces to checking that
$$
b_{\O-,\O} \m b_{\O,\O+} \!\m h_{(\O+)d_n(\O+)}
= 
[n+1]_q 
h_{(\O-)s_0 d_{n+1}(\O+)}
+q^{l-2n+1} [l-n+1]_q
h_{(\O-)s_0 d_{n-1}(\O+)},
$$
omitting the first term if $n=l$ and the last term if $n=0$.
If $n$ is odd, we need to show instead that
$$
b_{\O+,\O} \m b_{\O,\O-} \!\m h_{(\O-)s_0 d_n(\O+)}
= 
[n+1]_q 
h_{(\O+)d_{n+1}(\O+)}
+q^{l-2n+1} [l-n+1]_q
h_{(\O+)d_{n-1}(\O+)},
$$
omitting the first term if $n=l$.
This follows in both cases
by a computation using \cref{multrule}, \cref{lengths}, and the information about double cosets in \cref{pure,pure1,pure2}.
For example, if $n$ is even and $0 < n < l$,
then \cref{multrule}(1) implies that
$$
b_{\O,\O+} \m h_{(\O+)d_n(\O+)}
=
q^{l-n} h_{\O d_n(\O+)}
+
h_{\O s_0 d_{n+1}(\O+)},
$$
there being two terms in the summation ($d' = d_n$ and $d'=d_{n+1}$).
Then we use \cref{multrule}(2) to get
$$
b_{\O-,\O} \m b_{\O,\O+} \m h_{(\O+)d_n(\O+)}
=
q^{1-n} [l-n+1]_q q^{l-n}
h_{(\O-)s_0 d_{n-1}(\O+)}+
[n+1]_q
h_{(\O-)s_0 d_{n+1}(\O+)},
$$
as required. 
The other cases follow by similar calculations.

Next we check that $\phi(\overline{v}) = \overline{\phi(v)}$.
Since $\V(l)$ is generated as a $\Ui$-module by the highest weight vector 
$\eta_l$, and the elements denoted $b$ 
on both sides of the picture are bar-invarant, 
the proof of this reduces using the claim established in the previous paragraph just to checking that
$\phi(\overline{\eta_l}) = \overline{\phi(\eta_l)}$.
This is clear since $\eta_l$
and $\phi(\eta_l) = 1_{\O+}$ are both bar-invariant.

Now we can show that \cref{rain2} holds.
Applying the bar involution to \cref{characterize}
shows that
$b^{(n)} \eta_l$ is bar-invariant and equals $\overline{f^{(n)} \eta_l} +$ (a $q\Z[q]$-linear combination of other $\overline{f^{(i)} \eta_l}$). We deduce from \cref{rain1,rain2}
that $\phi(b^{(n)} \eta_l)$
is bar-invariant and equals
$h_{(\O+)d_n(\O+)}+(*)$ if $n$ is even or $s_0 h_{(\O-) s_0 d_n (\O+)}+(*)$ if $i$ is odd, where $(*)$
is a $q\Z[q]$-linear combination of other standard basis elements.
It follows that $\phi(b^{(n)} \eta_l)$ equals
$b_{(\O+)d_n(\O+)}$ if $n$ is even 
or $s_0 b_{(\O-)s_0 d_n (\O+)}$ if $n$ is odd, since these properties characterize this Kazhdan-Lusztig basis.

Finally, we check that
$\overline{(\overline{v},\overline{w})_l} = \big(\phi(v),\phi(w)\big)$.
We may assume that $v = \overline{f^{(m)} \eta_l}$
and $w = \overline{f^{(n)} \eta_l}$ for $0 \leq m,n \leq l$.
By the orthogonality of the respective standard bases, both sides are 0 if $m \neq n$,
so assume that $m=n$.
Using \cref{bilform}, we are reduced to 
 showing that
$$
q^{n(l-n)} \qbinom{l}{n}_q
= 
\begin{cases}
\big(h_{(\O+)d_n(\O+)}, h_{(\O+)d_n(\O+)}\big)
&\text{if $n$ is even}\\
\big(h_{(\O-)s_0 d_n(\O+)}, h_{(\O-)s_0 d_n(\O+)}\big)
&\text{if $n$ is odd.}
\end{cases}
$$
This follows from \cref{justforjon,lengths},
using also the descriptions of stabilizers in \cref{pure1,pure2}.
\end{proof}

The next corollaries follow from the theorem using
\cref{ohio,cute1,cute2}, respectively.
In their statements, we use $t$ to denote the unique element of
$\{0,1\}$ such that $l \equiv t \pmod{2}$, as it was in \cref{ssec2-6}.

\begin{cor}\label{cruisier}
For any $n \geq 0$, we have in $\Schur \oplus s_0 \Schur$ that
\begin{align*}
b^{\m n}
&= 
\begin{cases}
\displaystyle\sum_{\substack{0\leq i \leq \frac{n}{2}\\ n-2i \leq l}}[n-2i]_q^!
\left(\sum_{\lambda \in \Par_{\not\equiv t}(i\times (n-2i))}
[\lambda_1+1]_q^2\cdots [\lambda_i+1]_q^2\right) b_{(\O+)d_{n-2i}(\O+)}&
\text{if $n$ is even}\\
\displaystyle\sum_{\substack{0\leq i \leq \frac{n-1}{2}\\ n-2i \leq l}}[n-2i]_q^!
\left(\sum_{\lambda \in \Par_{\not\equiv t}(i\times (n-2i))}
[\lambda_1+1]_q^2\cdots [\lambda_i+1]_q^2\right) s_0 b_{(\O-)s_0 d_{n-2i}(\O+)}&
\text{if $n$ is odd.}
\end{cases}
\end{align*}
\end{cor}

\begin{cor}\label{formhecke}
For $0 \leq n \leq \lfloor\frac{l}{2}\rfloor$, we have that
$\big(1_{\O+},b_{(\O+)d_{2n}(\O+)}\big) = 
q^{n(l+t-1)}\qbinom{(l-t)/2}{n}_{q^2}$.
\end{cor}

\begin{cor}\label{significant}
The lowest degree term
of the Laurent polynomial 
$\big(b^{\m n},b^{\m n}\big)$
is $q^{-2\binom{n}{2}}$ when $0 \leq n \leq l$,
and it is $q^{-(2n-l)(l-1)}$ when $n > l$.
\end{cor}

\subsection{Extended singular Soergel bimodules}\label{ssec4-3}
Let $\BSBim$ (resp., $\SBim$) be the graded bicategory of singular Bott-Samelson bimodules (resp., singular Soergel bimodules) from \cref{newssbimdef} for the realization
of D$_l$ fixed in \cref{ssec4-1}. 
We remind again that $\gamma:N \rightarrow N$ 
denotes the graph automorphism. 

\begin{defin} For $I \subseteq N$, let $s_0 R^I$ be the
graded $(R^{\gamma(I)}, R^I)$-bimodule that is a copy 
$\{s_0 f\:|\:f \in R^I\}$ of the graded vector space $R^I$ (here $s_0$ is a formal symbol).
It is equipped with the natural right action of $R^I$, and a twisted left action of $R^{\gamma(I)}$ defined by $x \cdot (s_0 f) := s_0 \gamma(x) f$.
Similarly, there is a graded $(R^I, R^{\gamma(I)})$-bimodule $R^I s_0 = \{f s_0\:|\:f \in R^I\}$ with the natural left action and the twisted right action. 
\end{defin}

The functor $s_0 R^I \otimes_{R^I}-$ 
(resp., $-\otimes_{R^I} R^I s_0$)
amounts to twisting the left (resp., right) 
action of $R^I$ into an action of $R^{\gamma(I)}$
by pull-back along the isomorphism 
$\gamma:R^{\gamma(I)} \stackrel{\sim}{\rightarrow} R^I$.
We refer to $s_0 R^I$ and $R^I s_0$ as {\em twisting bimodules}.
There is an obvious isomorphism of graded $(R^{I}, R^{\gamma(I)})$-bimodules
$s_0 R^{\gamma(I)} \stackrel{\sim}{\rightarrow} R^{I} s_0$ which sends $s_0 f  \mapsto \gamma(f)s_0$. Because of this isomorphism we avoid using $R^I s_0$ entirely in the definitions below, having chosen to prefer $s_0 R^{\gamma(I)}$.

\begin{defin}\label{anywayyouwant}
Let $\EBSBim$, the
graded bicategory of {\em extended singular Bott-Samelson bimodules} (of type D$_l$),
be the full sub-bicategory of $\gBim$
with the same objects $R^I\:(I \subseteq N)$ as $\BSBim$, but with 1-morphisms 
given by tensor products of the bimodules $B_{I,J}^K$
for all $I \supseteq K \subseteq J$ as before plus the new twisting bimodules $s_0 R^I$ 
for $I \subseteq N$.
Then the graded bicategory $\ESBim$ of {\em extended singular Soergel bimodules} is the graded Karoubian closure of
$\EBSBim$ in $\gBim$.
\end{defin}

In string diagrams,
we use a {\em dashed string}
$\begin{tikzpicture}[anchorbase]
\draw[twist] (0,0) to (0,.5);
\node at (.3,.25) {$\regionlabel I$};
\node at (-.4,.25) {$\regionlabel \gamma(I)$};
\end{tikzpicture}$
separating 2-cells labelled $I$ to the right and $\gamma(I)$ to the left to denote the identity endomorphism of the twisting bimodule $s_0 R^I$. In contrast, our previous diagrams were built from {\em undashed strings}.
We have the obvious bubble slides:
\begin{equation}\label{marchesud}
\begin{tikzpicture}[anchorbase]
  \draw[twist] (-.04,-.4) to (-.04,.4);
\node at (-.4,0) {$\regionlabel \gamma(I)$};
\node at (.6,0) {$\regionlabel I$};
\node [draw,rounded corners,inner sep=2pt,black,fill=gray!20!white] at (.25,0) {$\scriptstyle f$};
\end{tikzpicture}=
\begin{tikzpicture}[anchorbase]
  \draw[twist] (0.04,-.4) to (0.04,.4);
\node at (-1.25,0) {$\regionlabel \gamma(I)$};
\node at (.2,0) {$\regionlabel I$};
\node [draw,rounded corners,inner sep=2pt,black,fill=gray!20!white] at (-.5,0) {$\scriptstyle \gamma(f)$};
\end{tikzpicture}
\end{equation}
for $f \in R^I$ (so $\gamma(f) \in R^{\gamma(I)}$).
Left-tensoring with the twisting bimodule $s_0 R^I$ is obviously an equivalence of categories with quasi-inverse given by left-tensoring with $s_0 R^{\gamma(I)}$.
One of the resulting adjunctions defines degree 0 graded bimodule isomorphisms $s_0 R^{\gamma(I)} \otimes_{R^{\gamma(I)}} s_0 R^{I}
\stackrel{\sim}{\rightarrow} R^I$
and $R^I \stackrel{\sim}{\rightarrow}
s_0 R^I \otimes_{R^I} s_0 R^{\gamma(I)}$
represented in string diagrams
by the dashed caps and cups
\begin{align*}
\begin{tikzpicture}[anchorbase,scale=2]
\draw[twist] (-.15,-.15) to [out=90,in=180] (0,.15) to [out=0
,in=90] (.15,-.15);
\node at (0.25,-.05) {$\regionlabel I$};
\node at (0,-.05) {$\regionlabel \gamma(I)$};
\end{tikzpicture}
\hspace{.5mm}
&:s_0 R^{\gamma(I)} \otimes_{R^{\gamma(I)}} s_0 R^I \rightarrow R^I,&
f \otimes g &\mapsto \gamma(f)g,\\
\begin{tikzpicture}[anchorbase,scale=2]
\draw[twist] (-.15,.15) to [out=-90,in=180] (0,-.15) to [out=0,in=-90] (.15,.15);
\node at (0,0.05) {$\regionlabel I$};
\node at (0.35,0.05) {$\regionlabel \gamma(I)$};
\end{tikzpicture}
&:R^{\gamma(I)}\rightarrow s_0 R^{I} \otimes_{R^I} s_0 R^{\gamma(I)},&
1 &\mapsto 1 \otimes 1.
\end{align*}
These satisfy the zig-zag identities, in addition to 
being mutually inverse isomorphisms:
\begin{align}
\begin{tikzpicture}[anchorbase,scale=1.4]
  \draw[twist] (0.3,0) to (0.3,.4);
	\draw[twist] (0.3,0) to[out=-90, in=0] (0.1,-0.2);
	\draw[twist] (0.1,-0.2) to[out = 180, in = -90] (-0.1,0);
	\draw[twist] (-0.1,0) to[out=90, in=0] (-0.3,0.2);
	\draw[twist] (-0.3,0.2) to[out = 180, in =90] (-0.5,0);
  \draw[twist] (-0.5,0) to (-0.5,-.4);
\node at (0.56,0) {$\regionlabel \gamma(I)$};
\node at (-.68,0) {$\regionlabel I$};
\end{tikzpicture}
&=
\begin{tikzpicture}[anchorbase,scale=1.4]
  \draw[twist] (0,-0.4) to (0,.4);
\node at (0.29,0) {$\regionlabel \gamma(I)$};
\node at (-.18,0) {$\regionlabel I$};
\end{tikzpicture}=
\begin{tikzpicture}[anchorbase,scale=1.4]
  \draw[twist] (0.3,0) to (0.3,-.4);
	\draw[twist] (0.3,0) to[out=90, in=0] (0.1,0.2);
	\draw[twist] (0.1,0.2) to[out = 180, in = 90] (-0.1,0);
	\draw[twist] (-0.1,0) to[out=-90, in=0] (-0.3,-0.2);
	\draw[twist] (-0.3,-0.2) to[out = 180, in =-90] (-0.5,0);
  \draw[twist] (-0.5,0) to (-0.5,.4);
\node at (0.56,0) {$\regionlabel \gamma(I)$};
\node at (-.68,0) {$\regionlabel I$};
\end{tikzpicture},&
\begin{tikzpicture}[anchorbase,scale=1.4]
\draw[twist] (0,-.4) to[out=90,in=90,looseness=2] (.4,-.4);
\draw[twist] (0,.4) to[out=-90,in=-90,looseness=2] (.4,.4);
\node at (0.2,-.3) {$\regionlabel I$};
\node at (0.2,.3) {$\regionlabel I$};
\node at (0.2,0) {$\regionlabel \gamma(I)$};
\end{tikzpicture}
&=
\begin{tikzpicture}[anchorbase,scale=1.4]
  \draw[twist] (0,-0.4) to (0,.4);
    \draw[twist] (0.4,-0.4) to (0.4,.4);
\node at (0.2,0) {$\regionlabel I$};
\node at (-.28,0) {$\regionlabel \gamma(I)$};
\node at (.68,0) {$\regionlabel \gamma(I)$};
\end{tikzpicture}
,&
\begin{tikzpicture}[anchorbase]
\draw[twist] (.4,0) arc[start angle=0, end angle=360,radius=.4];
\node at (.75,0) {$\regionlabel I$};
\node at (0,0) {$\regionlabel \gamma(I)$};
\end{tikzpicture}
&= \id_{R^I}.
\end{align}
 Using these relations, one can argue that any two diagrams built entirely from dashed strings and their cups and caps, and having the same boundary (thus representing morphisms between the same tensor products of twisting bimodules), are actually equal.
 (Formally, this augmented string calculus is encoding a categorical action of the cyclic group $C_2$.)

More interestingly there are {\em mixed crossings}, that is, crossings of dashed strings with undashed strings. These represent the following graded bimodule homomorphisms of degree 0:
\begin{align*}
\begin{tikzpicture}[anchorbase,scale=1.8]
	\draw[<-] (0.3,.3) to (-0.3,-.3);
	\draw[twist] (0.3,-.3) to (-0.3,.3);
 \node at (.25,0) {$\regionlabel Ii$};
 \node at (0,.25) {$\regionlabel I$};
 \node at (0,-.25) {$\regionlabel \gamma(I i)$};
 \node at (-.25,0) {$\regionlabel \gamma(I)$};
\end{tikzpicture}
&:B_{\gamma(I), \gamma(Ii)} \otimes_{R^{\gamma(Ii)}}
s_0 R^{Ii}
\rightarrow s_0 R^{I}\otimes_{R^I}  B_{I,Ii},
&f \otimes 1 &\mapsto 1\otimes \gamma(f),\\
\begin{tikzpicture}[anchorbase,scale=1.8]
	\draw[twist] (0.3,.3) to (-0.3,-.3);
	\draw[->] (0.3,-.3) to (-0.3,.3);
 \node at (.25,0) {$\regionlabel Ii$};
 \node at (0,.25) {$\regionlabel \gamma(Ii)$};
 \node at (0,-.25) {$\regionlabel I$};
 \node at (-.25,0) {$\regionlabel \gamma(I)$};
\end{tikzpicture}
&:s_0 R^I \otimes_{R^I} B_{I,Ii}
\rightarrow B_{\gamma(I),\gamma(Ii)} \otimes_{R^{\gamma(Ii)}} s_0 R^{Ii},
&1 \otimes f &\mapsto \gamma(f)\otimes 1,\\
\begin{tikzpicture}[anchorbase,scale=1.8]
	\draw[->] (0.3,.3) to (-0.3,-.3);
	\draw[twist] (0.3,-.3) to (-0.3,.3);
 \node at (.25,0) {$\regionlabel I$};
 \node at (0,.25) {$\regionlabel Ii$};
 \node at (0,-.25) {$\regionlabel \gamma(I)$};
 \node at (-.25,0) {$\regionlabel \gamma(Ii)$};
\end{tikzpicture}
&: B_{\gamma(Ii),\gamma(I)}\otimes_{R^{\gamma(I)}} s_0 R^{I} 
\rightarrow s_0 R^{Ii}\otimes_{R^{Ii}} B_{Ii,I},
&f \otimes 1 &\mapsto 1\otimes \gamma(f),\\
\begin{tikzpicture}[anchorbase,scale=1.8]
	\draw[twist] (0.3,.3) to (-0.3,-.3);
	\draw[<-] (0.3,-.3) to (-0.3,.3);
 \node at (.25,0) {$\regionlabel I$};
 \node at (0,.25) {$\regionlabel \gamma(I)$};
 \node at (0,-.25) {$\regionlabel Ii$};
 \node at (-.25,0) {$\regionlabel \gamma(Ii)$};
\end{tikzpicture}
&:s_0 R^{Ii} \otimes_{R^{Ii}} B_{Ii,I}
\rightarrow B_{\gamma(Ii),\gamma(I)} \otimes_{R^{\gamma(I)}} s_0 R^{I},
&1 \otimes f &\mapsto \gamma(f) \otimes 1.
\end{align*}
Although we are not labelling strings with their colors, in all four of the above diagrams, the undashed strings change color from $i$ to $\gamma(i)$ at the point that they cross the dashed string.

The mixed crossings are mutual inverses, that is, they satisfy Reidemeister II, in the obvious ways:
\begin{align}
\begin{tikzpicture}[anchorbase,scale=1.6]
	\draw[twist] (0.25,-.5) to[looseness=2,out=150,in=-150] (0.25,.5);
	\draw[->] (-.25,-.5) to[looseness=2,out=30,in=-30] (-0.25,.5);
 \node at (.4,0) {$\regionlabel Ii$};
 \node at (0,.54) {$\regionlabel \gamma(Ii)$};
 \node at (0,0) {$\regionlabel I$};
 \node at (0,-.54) {$\regionlabel \gamma(Ii)$};
 \node at (-.45,0) {$\regionlabel \gamma(I)$};
\end{tikzpicture}
&=
\begin{tikzpicture}[anchorbase,scale=1.6]
	\draw[twist] (0.25,-.5) to (0.25,.5);
	\draw[->] (-.2,-.5) to (-0.2,.5);
 \node at (.4,0) {$\regionlabel Ii$};
 \node at (-.4,0) {$\regionlabel \gamma(I)$};
 \node at (0.02,0) {$\regionlabel \gamma(Ii)$};
 \end{tikzpicture},&
 \begin{tikzpicture}[anchorbase,scale=1.6]
	\draw[->] (0.25,-.5) to[looseness=2,out=150,in=-150] (0.25,.5);
	\draw[twist] (-.25,-.5) to[looseness=2,out=30,in=-30] (-0.25,.5);
 \node at (.4,0) {$\regionlabel Ii$};
 \node at (0,.54) {$\regionlabel I$};
 \node at (0,0) {$\regionlabel I$};
 \node at (0,-.54) {$\regionlabel I$};
 \node at (-.45,0) {$\regionlabel \gamma(I)$};
\end{tikzpicture}
&=
\begin{tikzpicture}[anchorbase,scale=1.6]
	\draw[->] (0.25,-.5) to (0.25,.5);
	\draw[twist] (-.2,-.5) to (-0.2,.5);
 \node at (.4,0) {$\regionlabel Ii$};
 \node at (-.4,0) {$\regionlabel \gamma(I)$};
 \node at (0.02,0) {$\regionlabel I$};
 \end{tikzpicture},\\
 \begin{tikzpicture}[anchorbase,scale=1.6]
	\draw[twist] (0.25,-.5) to[looseness=2,out=150,in=-150] (0.25,.5);
	\draw[<-] (-.25,-.5) to[looseness=2,out=30,in=-30] (-0.25,.5);
 \node at (.4,0) {$\regionlabel I$};
 \node at (0,.54) {$\regionlabel \gamma(I)$};
 \node at (0,0) {$\regionlabel Ii$};
 \node at (0,-.54) {$\regionlabel \gamma(I)$};
 \node at (-.5,0) {$\regionlabel \gamma(Ii)$};
\end{tikzpicture}
&=
\begin{tikzpicture}[anchorbase,scale=1.6]
	\draw[twist] (0.25,-.5) to (0.25,.5);
	\draw[<-] (-.2,-.5) to (-0.2,.5);
 \node at (.4,0) {$\regionlabel I$};
 \node at (-.45,0) {$\regionlabel \gamma(Ii)$};
 \node at (0.02,0) {$\regionlabel \gamma(I)$};
 \end{tikzpicture},&
 \begin{tikzpicture}[anchorbase,scale=1.6]
	\draw[<-] (0.25,-.5) to[looseness=2,out=150,in=-150] (0.25,.5);
	\draw[twist] (-.25,-.5) to[looseness=2,out=30,in=-30] (-0.25,.5);
 \node at (-.5,0) {$\regionlabel \gamma(Ii)$};
 \node at (0,.54) {$\regionlabel Ii$};
 \node at (0,0) {$\regionlabel \gamma(I)$};
 \node at (0,-.54) {$\regionlabel Ii$};
 \node at (.4,0) {$\regionlabel I$};
\end{tikzpicture}
&=
\begin{tikzpicture}[anchorbase,scale=1.6]
	\draw[<-] (0.25,-.5) to (0.25,.5);
	\draw[twist] (-.2,-.5) to (-0.2,.5);
 \node at (-.45,0) {$\regionlabel \gamma(Ii)$};
 \node at (.4,0) {$\regionlabel I$};
 \node at (0.02,0) {$\regionlabel Ii$};
 \end{tikzpicture}.\label{alistair}
\end{align}
There are also Reidemeister III relations for diagrams involving dashed strings. These are obvious and we will not spell them out further. Consequently,
when at least one dashed string is involved,
we can even draw triple intersections without any ambiguity.

The upshot of this discussion is that dashed strings are easy to manipulate graphically---they slide across all other features, applying the graph automorphism to bubbles in the process, and they can be cut open and restitched in any reasonable way.

The graded 2-functor $\tR$ from \cref{tR}
extends to a graded 2-functor
\begin{align}\label{tRnew}
\tR&:\ESBim\rightarrow \ESBim^{\rev}
\end{align}
This fixes objects,
takes the graded $(R^I, R^J)$-bimodule $M$ to $q^{\ell(w_I)-\ell(w_J)} M$
viewed as an $(R^J, R^I)$-bimodule, and takes the twisting bimodule $s_0 R^I$ to $s_0 R^{\gamma(I)}$.
To define $\tR$ on a bimodule homomorphism $\theta$, the same underlying function defines
a bimodule homomorphism between the bimodules obtained instead by replacing each $s_0 R^I$ with $R^I s_0$; one then needs to conjugate with the appropriate
isomorphisms 
$s_0 R^{\gamma(I)} \stackrel{\sim}{\rightarrow} R^I s_0$ 
to obtain from this 
the bimodule homomorphism $\tR(\theta)$.
In terms of string diagrams, the functor $\tR$ reflects diagrams in a vertical axis and reverses the orientation on all undashed strings as before; dotted bubbles of the form $\bubble{f}$ in 2-cells are unchanged.
(The duality $\tD$ from \cref{rightduality} obviously extends too, but we will not use this below.)

The Soergel-Williamson \cref{homformula,class,catthm} are easily adapted to the extended setting. Although we do not believe that extended singular Soergel bimodules are yet in the literature, one can find a study of extended (ordinary) Soergel bimodules in \cite[Ch.3, Appendix]{Egaitsgory}, and it is straightforward to adapt the proofs to our setting. We just make a few comments.

 Any graded bimodule that is a 1-morphism in $\EBSBim$ is a tensor product of singular Bott-Samelson bimodules and some number of twisting bimodules. We say that the bimodule is {\em untwisted} or {\em twisted} according to whether the number of twisting bimodules in the tensor product is even or odd. These two types of bimodule can easily be distinguished since left and right multiplication by $e_l(x_1,\dots,x_l) = x_1 \cdots x_l \in R^N$ are equal on untwisted bimodules, and they differ by a sign on twisted ones.
In particular, no non-zero untwisted bimodule is isomorphic to a twisted one. 

Mixed crossings give an isomorphism between any 1-morphism in $\EBSBim$ and a bimodule of the form $M \otimes X$, where $M$ is a tensor product of twisting bimodules and $X$ is a singular Bott-Samelson bimodule. Dashed cups and caps give additional isomorphisms, letting one assume that the tensor factor $M$ is either the identity bimodule (if untwisted) or a single twisting bimodule (if twisted). Finally, we claim that tensoring with $\id_M$ gives an isomorphism
\begin{equation}  \Hom_{\BSBim}(X,X') \cong \Hom_{\EBSBim}(M \otimes X, M \otimes X'). \end{equation}
A priori, the right-hand side also contains diagrams with dashed strings appearing willy-nilly, but the relations allow one to slide all dashed strings to the left and then pull them straight.

The following is easily deduced from this discussion:
\begin{itemize}
\item
The classification of indecomposable untwisted bimodules in $\ESBim$ is exactly the same as in $\SBim$, as described by \cref{class}. The classification of indecomposable twisted bimodules is almost the same, except one needs to twist once so that, up to isomorphism and grading shift, the indecomposable twisted $(R^I, R^J)$-bimodules in $\ESBim$ are the bimodules $s_0 R^I \otimes_{R^I} B_{IdJ}$
for $I, J \subseteq N$ and $d \in (W_I \backslash W / W_J)_{\min}$.
\item
There are no non-zero homomorphisms between a twisted and an untwisted bimodule.
Given two untwisted bimodules,
the homomorphism formula is the same as before.
Given twisted bimodules, there is an analogous formula, but one needs to allow $\Delta$- and $\nabla$-flags which involve the twisted standard and costandard bimodules 
$s_0 R^I \otimes_{R^I} \Delta_{I d J}$ and
$s_0 R^I \otimes_{R^I} \nabla_{I d J}$, leading to a slightly modified definitions of the character maps $\ch_\Delta$ and $\ch_\nabla$.
(Alternatively, one can just tensor one more time with a twisting bimodule to reduce to the usual homomorphism formula in the untwisted setting.)
\item
Finally, the categorification theorem is easily adapted. 
The modified character map $\ch$ takes the indecomposable untwisted bimodule $B_{IdJ}$ to $b_{IdJ} \in \Schur$ as before, and it takes the twisted bimodule
$s_0 R^I \otimes_{R^I} B_{I d J}$ to $s_0 b_{I d J}$.
One recovers the extended $q$-Schur algebra:
\begin{equation}
K_0(\ESBim) \cong \Schur \oplus s_0 \Schur.
\end{equation}
\end{itemize}

\subsection{Diagrammatics for 
extended singular Bott-Samelson bimodules}
We now start to use the diagrammatic calculus explained in \cref{ssec3-3} more seriously
for making calculations in the graded bicategory $\EBSBim$ from \cref{anywayyouwant}.
Traditionally in this subject an unlabelled 2-cell in a string diagram
is used to denote the 2-cell labelled by the empty set,
but from now on:

\vspace{2mm}
{\begin{center}
\fbox{\begin{minipage}{29em}
An unlabelled 2-cell in a string diagram denotes the 2-cell labelled by $\O$.
\end{minipage}}
\end{center}}

\vspace{2mm}
\noindent
Thus $\O$ is our ``ground state''.
Then we decorate the 2-cell with $+$, $-$ or $\pm$
to denote the 2-cell labelled by $\O+, \O-$ or $\O\pm$, and we decorate it with $\hat n, 
\hat n+, \hat n+$ or $\hat n\pm$
to denote the 2-cell labelled by $\O{\hat n},
\O{\hat n+}, \O{\hat n-}$ or $\O{\hat n\pm}$
for $n \in \O$,
using the shorthands \cref{sledgehammer}.
In fact, moving forwards, we will rarely need 2-cells with labels of the form $I-$ for $I \subseteq \O$; it will be enough to consider string diagrams with 2-cells labelled $I, I+$ or $I \pm$ for $I \subseteq \O$. 
Typically $I$ will be $\O$ or $\O \hat n$.

Now we introduce an important notational shorthand to the string calculus for extended singular Soergel bimodules. 
Suppose that $I \subseteq \O$, so that $\gamma(I) = I$.
We will use the {\em special undotted
strings}
\begin{align}\label{fire1}
\begin{tikzpicture}[anchorbase,scale=1.2]
\draw[colp,->] (0,0) to (0,.6);
\node at (-0,.3) {$\regionlabel I\phantom{+}$};
\node at (.22,.3) {$\regionlabel I+$};
\end{tikzpicture}\!\!
&:=\!\!
\begin{tikzpicture}[anchorbase,scale=1.2]
\draw[->] (0,0) to (0,.6);
\node at (-0,.3) {$\regionlabel I\phantom{+}$};
\node at (.22,.3) {$\regionlabel I+$};
\end{tikzpicture},
&\begin{tikzpicture}[anchorbase,scale=1.2]
\draw[colp,->] (0,0) to (0,.6);
\node at (-.2,.3) {$\regionlabel I+$};
\node at (.24,.3) {$\regionlabel I\pm$};
\end{tikzpicture}
&\!\!:=\!\!
\begin{tikzpicture}[anchorbase,scale=1.2]
\draw[->] (0,0) to (0,.6);
\node at (-.2,.3) {$\regionlabel I+$};
\node at (.24,.3) {$\regionlabel I\pm$};
\end{tikzpicture},&
\begin{tikzpicture}[anchorbase,scale=1.2]
\draw[colp,<-] (0,0) to (0,.6);
\node at (-.22,.3) {$\regionlabel I\pm$};
\node at (.22,.3) {$\regionlabel I+$};
\end{tikzpicture}\!\!
&:=\!\!
\begin{tikzpicture}[anchorbase,scale=1.2]
\draw[<-] (0,0) to (0,.6);
\node at (-.22,.3) {$\regionlabel I\pm$};
\node at (.22,.3) {$\regionlabel I+$};
\end{tikzpicture},
&\begin{tikzpicture}[anchorbase,scale=1.2]
\draw[colp,<-] (0,0) to (0,.6);
\node at (-.2,.3) {$\regionlabel I+$};
\node at (.24,.3) {$\regionlabel I\phantom{+}$};
\end{tikzpicture}\!\!
&:=\!\!
\begin{tikzpicture}[anchorbase,scale=1.2]
\draw[<-] (0,0) to (0,.6);
\node at (-.2,.3) {$\regionlabel I+$};
\node at (.24,.3) {$\regionlabel I\phantom{+}$};
\end{tikzpicture}\\\intertext{to denote
the identity endomorphisms of the untwisted bimodules
$B_{I, I+}$, $B_{I+,I\pm}$,
$B_{I\pm,I+}$ and $B_{I+,I}$. 
We use the {\em special dotted strings}}
\label{fire2}
\begin{tikzpicture}[anchorbase,scale=1.2]
\draw[colm,->] (0,0) to (0,.6);
\node at (-.0,.3) {$\regionlabel I\phantom{+}$};
\node at (.2,.3) {$\regionlabel I+$};
\end{tikzpicture}\!\!
&:=\!\!
\begin{tikzpicture}[anchorbase,scale=1.2]
\draw[->] (0,0) to (0,.6);
\draw[twist] (-.2,0) to (-.2,.6);
\node at (0,.3) {$\regionlabel I\phantom{+}$};
\node at (-.2,.3) {$\regionlabel I\phantom{+}$};
\node at (.2,.3) {$\regionlabel I+$};
\end{tikzpicture},
&\begin{tikzpicture}[anchorbase,scale=1.2]
\draw[colm,->] (0,0) to (0,.6);
\node at (-.2,.3) {$\regionlabel I+$};
\node at (.22,.3) {$\regionlabel I\pm$};
\end{tikzpicture}\!\!
&:=\!\!
\begin{tikzpicture}[anchorbase,scale=1.2]
\draw[->] (-0.2,0) to (-0.2,.6);
\node at (-.37,.3) {$\regionlabel I+$};
\draw[twist] (0.2,0) to (0.2,.6);
\node at (0,.3) {$\regionlabel I\pm$};
\node at (.4,.3) {$\regionlabel I\pm$};
\end{tikzpicture},&
\begin{tikzpicture}[anchorbase,scale=1.2]
\draw[colm,<-] (0,0) to (0,.6);
\node at (-.2,.3) {$\regionlabel I\pm$};
\node at (.2,.3) {$\regionlabel I+$};
\end{tikzpicture}\!\!
&:=\!\!
\begin{tikzpicture}[anchorbase,scale=1.2]
\draw[<-] (0.2,0) to (0.2,.6);
\draw[twist] (-0.2,0) to (-0.2,.6);
\node at (0,.3) {$\regionlabel I\pm$};
\node at (-.4,.3) {$\regionlabel I\pm$};
\node at (.4,.3) {$\regionlabel I+$};
\end{tikzpicture},
&\begin{tikzpicture}[anchorbase,scale=1.2]
\draw[colm,<-] (0,0) to (0,.6);
\node at (-.2,.3) {$\regionlabel I+$};
\node at (.25,.3) {$\regionlabel I\phantom{+}$};
\end{tikzpicture}\!\!
&:=\!\!
\begin{tikzpicture}[anchorbase,scale=1.2]
\draw[twist] (0.03,0) to (0.03,.6);
\draw[<-] (-0.2,0) to (-0.2,.6);
\node at (-.1,.3) {$\regionlabel I$};
\node at (-.4,.3) {$\regionlabel I+$};
\node at (.25,.3) {$\regionlabel I\phantom{+}$};
\end{tikzpicture}
\end{align}
to denote the identity endomorphisms of the twisted bimodules 
$s_0 R^I \otimes_{R^I} B_{I,I+}$,
$B_{I_+, I_{\pm}} \otimes_{R^{I \pm}} s_0 R^{I\pm}$,
$s_0 R^{I\pm}\otimes_{R^{I\pm}} B_{I\pm,I+}$
and $s_0 R^{I+} \otimes_{R^{I+}} B_{I+, I}$.
It is helpful to remember that the dashed line, when it appears in these diagrams, passes through the 2-cells whose label is invariant under $\gamma$.
We reiterate that these special dotted strings are not identity morphisms of new objects in the category, they are merely a shorthand for the identity morphisms of existing objects. In this shorthand, the equality \cref{adele} is tautological, and 
\cref{needsabitofthought} is a consequence of \cref{marchesud}.
Actually, the special undotted string is an ordinary string labelled either $+$ or $-$, but we find the addition of color to it makes the diagrams in which it arises later on more readable.

Notice that
$\begin{tikzpicture}[anchorbase,scale=1]
\draw[colm,<-] (0,0.05) to (0,.5);
\node at (-.24,.3) {$\regionlabel I+$};
\node at (.23,.3) {$\regionlabel I$};
\draw[colp,->] (0.45,0.05) to (0.45,.5);
\node at (.71,.3) {$\regionlabel I+$};
\end{tikzpicture} = 
\begin{tikzpicture}[anchorbase,scale=1]
\draw[<-] (-.2,0.05) to (-.2,.5);
\node at (-.44,.3) {$\regionlabel I+$};
\draw[twist] (.2,0.05) to (.2,.5);
\node at (0,.3) {$\regionlabel I$};
\node at (.4,.3) {$\regionlabel I$};
\draw[->] (0.6,0.05) to (0.6,.5);
\node at (.86,.3) {$\regionlabel I+$};
\end{tikzpicture}
= \begin{tikzpicture}[anchorbase,scale=1]
\draw[colp,<-] (0,0.05) to (0,.5);
\node at (-.24,.3) {$\regionlabel I+$};
\node at (.23,.3) {$\regionlabel I$};
\draw[colm,->] (0.45,0.05) to (0.45,.5);
\node at (.71,.3) {$\regionlabel I+$};
\end{tikzpicture}$, indeed, these represent the identity endomorphism of the twisted $(R^{I+},R^{I+})$-bimodule
$B_{I+,I} \otimes_{R^{I}} s_0 R^{I} \otimes_{R^I}
B_{I,I^+}$.
Most important is the case $I = \O$ here, when
\begin{equation}\label{adele}
\begin{tikzpicture}[anchorbase,scale=1.2]
\draw[colm,<-] (0,0.05) to (0,.5);
\node at (-.24,.3) {$\labelp$};
\draw[colp,->] (0.3,0.05) to (0.3,.5);
\node at (.56,.3) {$\labelp$};
\end{tikzpicture}=\begin{tikzpicture}[anchorbase,scale=1.2]
\draw[<-] (0,0.05) to (0,.5);
\draw[twist] (.225,0.05) to (.225,.5);
\node at (-.24,.3) {$\labelp$};
\draw[->] (0.45,0.05) to (0.45,.5);
\node at (.71,.3) {$\labelp$};
\end{tikzpicture}
= \begin{tikzpicture}[anchorbase,scale=1.2]
\draw[colp,<-] (0,0.05) to (0,.5);
\node at (-.24,.3) {$\labelp$};
\draw[colm,->] (0.3,0.05) to (0.3,.5);
\node at (.56,.3) {$\labelp$};
\end{tikzpicture}
\end{equation}
is the identity endomorphism of the twisted $(R^{\O+}, R^{\O+})$-bimodule
\begin{equation}\label{words}
B := B_{\O+,\O} \otimes_{R^{\O}} s_0 R^{\O} \otimes_{R^\O}
B_{\O,\O^+}.
\end{equation}
This is isomorphic via an obvious contraction to the bimodule $B$ defined in the introduction. We will use the present definition of $B$ from now on.
At the level of characters, $B$ corresponds to the element $b \in s_0 \Schur$ from \cref{killingme}: 
\begin{equation}
\ch(B) = b_{\O+,\O} \m (s_0 b_{\O,\O+})
= s_0 b_{\O-,\O} \m b_{\O,\O+} = s_0 b_{\O-,\O+} = b.
\end{equation}
Note also that
\begin{equation}\label{needsabitofthought}
\begin{tikzpicture}[anchorbase,scale=1.4]
\draw[colp,<-] (-.05,0.05) to (-.05,.5);
\node at (-.24,.3) {$\labelp$};
\draw[colm,->] (0.45,0.05) to (0.45,.5);
\node at (.66,.3) {$\labelp$};
\node [draw,rounded corners,inner sep=2pt,black,fill=gray!20!white] at (.2,.275) {$\scriptstyle f$};
\end{tikzpicture}
=
\begin{tikzpicture}[anchorbase,scale=1.4]
\draw[colm,<-] (-.2,0.05) to (-.2,.5);
\node at (-.39,.3) {$\labelp$};
\draw[colp,->] (0.6,0.05) to (0.6,.5);
\node at (.81,.3) {$\labelp$};
\node [draw,rounded corners,inner sep=2pt,black,fill=gray!20!white] at (.2,.275) {$\scriptstyle \gamma(f)$};
\end{tikzpicture}\ .
\end{equation}
The symmetry $\tR$ from \cref{tRnew} acts on string diagrams as before, preserving the type (undotted or dotted) of special strings.
For example:
\begin{align}
\tR\left(
\begin{tikzpicture}[anchorbase,scale=1.4]
\draw[colp,<-] (-0.05,0.05) to (-.05,.5);
\node at (-.24,.3) {$\labelp$};
\draw[colm,->] (0.45,0.05) to (0.45,.5);
\node at (.66,.3) {$\labelp$};
\node [draw,rounded corners,inner sep=2pt,black,fill=gray!20!white] at (.2,.275) {$\scriptstyle f$};
\end{tikzpicture}\right)&=
\begin{tikzpicture}[anchorbase,scale=1.4]
\draw[colm,<-] (-.05,0.05) to (-.05,.5);
\node at (-.24,.3) {$\labelp$};
\draw[colp,->] (0.45,0.05) to (0.45,.5);
\node at (.66,.3) {$\labelp$};
\node [draw,rounded corners,inner sep=2pt,black,fill=gray!20!white] at (.2,.275) {$\scriptstyle f$};
\end{tikzpicture}\ ,&
\tR\left(
\begin{tikzpicture}[anchorbase,scale=1.4]
\draw[colm,<-] (-0.05,0.05) to (-.05,.5);
\node at (-.24,.3) {$\labelp$};
\draw[colp,->] (0.45,0.05) to (0.45,.5);
\node at (.66,.3) {$\labelp$};
\node [draw,rounded corners,inner sep=2pt,black,fill=gray!20!white] at (.2,.275) {$\scriptstyle f$};
\end{tikzpicture}\right)&=
\begin{tikzpicture}[anchorbase,scale=1.4]
\draw[colp,<-] (-.05,0.05) to (-.05,.5);
\node at (-.24,.3) {$\labelp$};
\draw[colm,->] (0.45,0.05) to (0.45,.5);
\node at (.66,.3) {$\labelp$};
\node [draw,rounded corners,inner sep=2pt,black,fill=gray!20!white] at (.2,.275) {$\scriptstyle f$};
\end{tikzpicture}\ .
\end{align}

There are several more diagrammatic shorthands to be introduced involving the two new types of colored strings.
These should all seem unsurpising by this point.
Continuing with any $I \subseteq \O$, 
there are colored caps and cups satisfying zig-zag identities:
\begin{align*}
&\begin{tikzpicture}[anchorbase,scale=2]
\draw[colp,->] (-.15,-.15) to [out=90,in=180] (0,.15) to [out=0
,in=90] (.15,-.15);
\node at (0.25,0) {$\regionlabel I$};
\node at (0,0) {$\regionlabel I+$};
\end{tikzpicture},
&&
\begin{tikzpicture}[anchorbase,scale=2]
\draw[colp,->] (-.15,.15) to [out=-90,in=180] (0,-.15) to [out=0
,in=-90] (.15,.15);
\node at (0.3,0) {$\regionlabel I+$};
\node at (0,0) {$\regionlabel I$};
\end{tikzpicture},&&
\begin{tikzpicture}[anchorbase,scale=2]
\draw[colp,<-] (-.15,-.15) to [out=90,in=180] (0,.15) to [out=0
,in=90] (.15,-.15);
\node at (0.3,0) {$\regionlabel I+$};
\node at (0,0) {$\regionlabel I$};
\end{tikzpicture},
&&
\begin{tikzpicture}[anchorbase,scale=2]
\draw[colp,<-] (-.15,.15) to [out=-90,in=180] (0,-.15) to [out=0
,in=-90] (.15,.15);
\node at (0.25,0) {$\regionlabel I$};
\node at (0,0) {$\regionlabel I+$};
\end{tikzpicture},&
&\begin{tikzpicture}[anchorbase,scale=2]
\draw[colp,->] (-.15,-.15) to [out=90,in=180] (0,.15) to [out=0
,in=90] (.15,-.15);
\node at (0.3,0) {$\regionlabel I+$};
\node at (0,0) {$\regionlabel I\pm$};
\end{tikzpicture},
&&
\begin{tikzpicture}[anchorbase,scale=2]
\draw[colp,->] (-.15,.15) to [out=-90,in=180] (0,-.15) to [out=0
,in=-90] (.15,.15);
\node at (0.3,0) {$\regionlabel I\pm$};
\node at (0,0) {$\regionlabel I+$};
\end{tikzpicture},&&
\begin{tikzpicture}[anchorbase,scale=2]
\draw[colp,<-] (-.15,-.15) to [out=90,in=180] (0,.15) to [out=0
,in=90] (.15,-.15);
\node at (0.3,0) {$\regionlabel I\pm$};
\node at (0,0) {$\regionlabel I+$};
\end{tikzpicture},
&&
\begin{tikzpicture}[anchorbase,scale=2]
\draw[colp,<-] (-.15,.15) to [out=-90,in=180] (0,-.15) to [out=0
,in=-90] (.15,.15);
\node at (0.3,0) {$\regionlabel I+$};
\node at (0,0) {$\regionlabel I\pm$};
\end{tikzpicture},\\
&\begin{tikzpicture}[anchorbase,scale=2]
\draw[colm,->] (-.15,-.15) to [out=90,in=180] (0,.15) to [out=0
,in=90] (.15,-.15);
\node at (0.25,0) {$\regionlabel I$};
\node at (0,0) {$\regionlabel I+$};
\end{tikzpicture},
&&
\begin{tikzpicture}[anchorbase,scale=2]
\draw[colm,->] (-.15,.15) to [out=-90,in=180] (0,-.15) to [out=0
,in=-90] (.15,.15);
\node at (0.3,0) {$\regionlabel I+$};
\node at (0,0) {$\regionlabel I$};
\end{tikzpicture},
&&
\begin{tikzpicture}[anchorbase,scale=2]
\draw[colm,<-] (-.15,-.15) to [out=90,in=180] (0,.15) to [out=0
,in=90] (.15,-.15);
\node at (0.3,0) {$\regionlabel I+$};
\node at (0,0) {$\regionlabel I$};
\end{tikzpicture},
&&
\begin{tikzpicture}[anchorbase,scale=2]
\draw[colm,<-] (-.15,.15) to [out=-90,in=180] (0,-.15) to [out=0
,in=-90] (.15,.15);
\node at (0.25,0) {$\regionlabel I$};
\node at (0,0) {$\regionlabel I+$};
\end{tikzpicture},&&
\begin{tikzpicture}[anchorbase,scale=2]
\draw[colm,->] (-.15,-.15) to [out=90,in=180] (0,.15) to [out=0
,in=90] (.15,-.15);
\node at (0.3,0) {$\regionlabel I+$};
\node at (0,0) {$\regionlabel I\pm$};
\end{tikzpicture},
&&
\begin{tikzpicture}[anchorbase,scale=2]
\draw[colm,->] (-.15,.15) to [out=-90,in=180] (0,-.15) to [out=0
,in=-90] (.15,.15);
\node at (0.3,0) {$\regionlabel I\pm$};
\node at (0,0) {$\regionlabel I+$};
\end{tikzpicture},&&
\begin{tikzpicture}[anchorbase,scale=2]
\draw[colm,<-] (-.15,-.15) to [out=90,in=180] (0,.15) to [out=0
,in=90] (.15,-.15);
\node at (0.3,0) {$\regionlabel I\pm$};
\node at (0,0) {$\regionlabel I+$};
\end{tikzpicture},
&&
\begin{tikzpicture}[anchorbase,scale=2]
\draw[colm,<-] (-.15,.15) to [out=-90,in=180] (0,-.15) to [out=0
,in=-90] (.15,.15);
\node at (0.3,0) {$\regionlabel I+$};
\node at (0,0) {$\regionlabel I\pm$};
\end{tikzpicture}.
\end{align*}
For example, the last two
are defined by
\begin{align*}
\begin{tikzpicture}[anchorbase,scale=2]
\draw[colm,<-] (-.15,-.15) to [out=90,in=180] (0,.15) to [out=0
,in=90] (.15,-.15);
\node at (0.3,0) {$\regionlabel I\pm$};
\node at (0,0) {$\regionlabel I+$};
\end{tikzpicture}
&:=\begin{tikzpicture}[anchorbase,scale=2]
\draw[<-] (-.15,-.15) to [out=90,in=180] (0,.15) to [out=0
,in=90] (.15,-.15);
\draw[twist] (-.4,-.15) to [out=90,in=180] (0,.4) to [out=0
,in=90] (.4,-.15);
\node at (0.25,0) {$\regionlabel I\pm$};
\node at (0.5,0) {$\regionlabel I\pm$};
\node at (0,0) {$\regionlabel I+$};
\end{tikzpicture}
,&
\begin{tikzpicture}[anchorbase,scale=2]
\draw[colm,<-] (-.15,.15) to [out=-90,in=180] (0,-.15) to [out=0
,in=-90] (.15,.15);
\node at (0.3,0) {$\regionlabel I+$};
\node at (0,0) {$\regionlabel I\pm$};
\end{tikzpicture}&:=
\begin{tikzpicture}[anchorbase,scale=2]
\draw[twist] (-.15,.15) to [out=-90,in=180] (0,-.15) to [out=0
,in=-90] (.15,.15);
\draw[<-] (-.4,.15) to [out=-90,in=180] (0,-.4) to [out=0
,in=-90] (.4,.15);
\node at (0.5,0) {$\regionlabel I+$};
\node at (0.25,0) {$\regionlabel I\pm$};
\node at (0,0) {$\regionlabel I\pm$};
\end{tikzpicture}
\ ,
\end{align*}
and they satisfy the zig-zag identities by combining the ones for ordinary and dashed strings seen before.
The other cases are similar.

Now we discuss crossings. 
Undotted special strings are ordinary strings (with added color for clarity), so we have already defined crossings between them and other ordinary strings. There are also crossings of dotted special strings with others strings which we call
{\em special crossings}.
For $i \in I \subseteq \O$, 
the upward special crossings are defined as follows:
\begin{align*}
\begin{tikzpicture}[anchorbase,scale=2.2]
	\draw[<-] (0.3,.3) to (-0.3,-.3);
	\draw[colm,->] (0.3,-.3) to (-0.3,.3);
 \node at (.25,0) {$\regionlabel Ii+$};
 \node at (0,.25) {$\regionlabel I+$};
 \node at (0,-.25) {$\regionlabel Ii$};
 \node at (-.25,0) {$\regionlabel I$};
\end{tikzpicture} &:= 
\begin{tikzpicture}[anchorbase,scale=2.2]
	\draw[<-] (0.6,.3) to (0,-.3);
	\draw[twist] (0.6,-.3) to (0,.3);
 \draw[->] (0.3,-.3) to (0.3,.3);
 \node at (.55,0) {$\regionlabel Ii+$};
  \node at (.4,-.23) {$\regionlabel Ii-$};
 \node at (.2,.23) {$\regionlabel I$};
  \node at (.4,.23) {$\regionlabel I+$};
 \node at (0.21,-.23) {$\regionlabel Ii$};
 \node at (0.05,0) {$\regionlabel I$};
\end{tikzpicture}\ ,&
\begin{tikzpicture}[anchorbase,scale=2.2]
	\draw[<-] (-0.3,.3) to (0.3,-.3);
	\draw[colm,->] (-0.3,-.3) to (0.3,.3);
 \node at (.25,0) {$\regionlabel Ii+$};
 \node at (0,.25) {$\regionlabel Ii$};
 \node at (0,-.25) {$\regionlabel I+$};
 \node at (-.25,0) {$\regionlabel I$};
\end{tikzpicture} &:= 
\begin{tikzpicture}[anchorbase,scale=2.2]
	\draw[twist] (0.6,.3) to (0,-.3);
	\draw[->] (0.6,-.3) to (0,.3);
 \draw[->] (0.3,-.3) to (0.3,.3);
 \node at (.55,0) {$\regionlabel Ii +$};
  \node at (.4,-.23) {$\regionlabel I+$};
 \node at (.2,.23) {$\regionlabel Ii$};
  \node at (.4,.23) {$\regionlabel Ii-$};
 \node at (0.21,-.23) {$\regionlabel I$};
 \node at (0.05,0) {$\regionlabel I$};
\end{tikzpicture}\ ,\\
\begin{tikzpicture}[anchorbase,scale=2.2]
	\draw[colp,<-] (0.3,.3) to (-0.3,-.3);
	\draw[colm,->] (0.3,-.3) to (-0.3,.3);
 \node at (.25,0) {$\regionlabel I\pm$};
 \node at (0,.25) {$\regionlabel I+$};
 \node at (0,-.25) {$\regionlabel I+$};
 \node at (-.25,0) {$\regionlabel I$};
\end{tikzpicture} &:= 
\begin{tikzpicture}[anchorbase,scale=2.2]
	\draw[<-] (0.6,.3) to (0,-.3);
	\draw[twist] (0.6,-.3) to (0,.3);
 \draw[->] (0.3,-.3) to (0.3,.3);
 \node at (.55,0) {$\regionlabel I\pm$};
  \node at (.4,-.23) {$\regionlabel I\pm$};
 \node at (.2,.23) {$\regionlabel I$};
  \node at (.4,.23) {$\regionlabel I+$};
 \node at (0.21,-.23) {$\regionlabel I+$};
 \node at (0.05,0) {$\regionlabel I$};
\end{tikzpicture}\ ,&
\begin{tikzpicture}[anchorbase,scale=2.2]
	\draw[colp,<-] (-0.3,.3) to (0.3,-.3);
	\draw[colm,->] (-0.3,-.3) to (0.3,.3);
 \node at (.25,0) {$\regionlabel I\pm$};
 \node at (0,.25) {$\regionlabel I+$};
 \node at (0,-.25) {$\regionlabel I+$};
 \node at (-.25,0) {$\regionlabel I$};
\end{tikzpicture} &:= 
\begin{tikzpicture}[anchorbase,scale=2.2]
	\draw[twist] (0.6,.3) to (0,-.3);
	\draw[->] (0.6,-.3) to (0,.3);
 \draw[->] (0.3,-.3) to (0.3,.3);
 \node at (.55,0) {$\regionlabel I\pm$};
  \node at (.4,-.23) {$\regionlabel I+$};
 \node at (.2,.23) {$\regionlabel I+$};
  \node at (.4,.23) {$\regionlabel I\pm$};
 \node at (0.21,-.23) {$\regionlabel I$};
 \node at (0.05,0) {$\regionlabel I$};
\end{tikzpicture}\ ,
\end{align*}
referring to 
the discussion after \cref{alistair} for interpretation of the triple crossings. One defines downward and sideways special crossings by rotating the ones above, noting the cyclicity of the string calculus.

The special crossings satisfy the easy Reidemeister II relations, as follows directly from \cref{easyRII} and the definitions, using the simple behavior of dashed strings.
Other relations involving the new colored strings can be worked out in just the same way.
Here are a few more examples, mainly 
to draw attention to the importance of the type of special strings (i.e., undotted vs. dotted):
\begin{align}\label{howareyou}
\begin{tikzpicture}[anchorbase]
  \draw[colp,->] (-.04,-.4) to (-.04,.4);
\node at (.7,0) {$\labelp$};
\node [draw,rounded corners,inner sep=2pt,black,fill=gray!20!white] at (.25,0) {$\scriptstyle f$};
\end{tikzpicture}&=
\begin{tikzpicture}[anchorbase]
  \draw[colp,->] (0.04,-.4) to (0.04,.4);
\node at (.3,0) {$\labelp$};
\node [draw,rounded corners,inner sep=2pt,black,fill=gray!20!white] at (-.25,0) {$\scriptstyle f$};
\end{tikzpicture},&
\begin{tikzpicture}[anchorbase]
  \draw[colm,->] (-.04,-.4) to (-.04,.4);
\node at (.7,0) {$\labelp$};
\node [draw,rounded corners,inner sep=2pt,black,fill=gray!20!white] at (.25,0) {$\scriptstyle f$};
\end{tikzpicture}&=
\begin{tikzpicture}[anchorbase]
  \draw[colm,->] (0.04,-.4) to (0.04,.4);
\node at (.3,0) {$\labelp$};
\node [draw,rounded corners,inner sep=2pt,black,fill=gray!20!white] at (-.54,0) {$\scriptstyle \gamma(f)$};
\end{tikzpicture}
&&\text{for $f \in R^{\O+}$,}\\
\begin{tikzpicture}[anchorbase]
\draw[colp,->] (.5,0) arc[start angle=0, end angle=360,radius=.5];
\node [draw,rounded corners,inner sep=2pt,black,fill=gray!20!white] at (0,0) {$\scriptstyle f$};
\node at (.75,0) {$\labelp$};
\end{tikzpicture}&=\begin{tikzpicture}[anchorbase]
\node [draw,rounded corners,inner sep=2pt,black,fill=gray!20!white] at (0,0) {$\scriptstyle \tr_{\O+}^{\O}(f)$};
\node at (.75,0) {$\labelp$};
\end{tikzpicture},&
\begin{tikzpicture}[anchorbase]
\draw[colm,->] (.5,0) arc[start angle=0, end angle=360,radius=.5];
\node [draw,rounded corners,inner sep=2pt,black,fill=gray!20!white] at (0,0) {$\scriptstyle f$};
\node at (.75,0) {$\labelp$};
\end{tikzpicture}&=\begin{tikzpicture}[anchorbase]
\node [draw,rounded corners,inner sep=2pt,black,fill=gray!20!white] at (0,0) {$\scriptstyle \tr_{\O+}^{\O}(\gamma(f))$};
\node at (1,0) {$\labelp$};
\end{tikzpicture}=
\begin{tikzpicture}[anchorbase]
\node [draw,rounded corners,inner sep=2pt,black,fill=gray!20!white] at (0,0) {$\scriptstyle \gamma(\tr_{\O-}^{\O}(f))$};
\node at (1.05,0) {$\labelp$};
\end{tikzpicture}
&&\text{for $f \in R^{\O}$,}\\
\begin{tikzpicture}[anchorbase]
\draw[colp,<-] (.5,0) arc[start angle=0, end angle=360,radius=.3];
\node at (.2,0) {$\labelp$};
\end{tikzpicture}&=
\begin{tikzpicture}[anchorbase]
\node [draw,rounded corners,inner sep=2pt,black,fill=gray!20!white] at (0,0) {$\scriptstyle \eta^\O_{\O+}$};
\end{tikzpicture}\ ,
& \begin{tikzpicture}[anchorbase]
\draw[colm,<-] (.5,0) arc[start angle=0, end angle=360,radius=.3];
\node at (.2,0) {$\labelp$};
\end{tikzpicture}
&=
\begin{tikzpicture}[anchorbase]
\node [draw,rounded corners,inner sep=2pt,black,fill=gray!20!white] at (0,0) {$\scriptstyle \gamma(\eta^\O_{\O+})$};
\end{tikzpicture}=
\begin{tikzpicture}[anchorbase]
\node [draw,rounded corners,inner sep=2pt,black,fill=gray!20!white] at (0,0) {$\scriptstyle \eta^\O_{\O-}$};
\end{tikzpicture}\ .\label{iamwell}
\end{align}
This calculus will be used extensively in the rest of the paper.

\begin{rem}\label{benhomework}
Since thinking in a foreign language is hard, we sometimes find it easier to replace the 
special strings with ordinary and dashed strings according to \cref{fire1,fire2}, slide the resulting dashed strings to the edges, then work with 
the ordinary (i.e., non-extended) string calculus after that. 
A roughly equivalent technique is to use the following mutually inverse isomorphisms to relate $B^{\otimes 2}$ with the 
untwisted singular Soergel bimodule
$B_{O+,O-}\otimes_{R^{O-}}B_{O-,O+}$: 
\begin{align}
\begin{tikzpicture}[anchorbase,scale=.9]
\node at (1.4,.3) {$\labelp$};
\node at (-.2,.3) {$\labelp$};
\node at (.6,.15) {$\labelp$};
\node at (.6,.45) {$\labelm$};
\draw[->] (0,0) to (0,.6);
\draw[<-] (.4,0) to (.4,.6);
\draw[->] (.8,0) to (.8,.6);
\draw[<-] (1.2,0) to (1.2,.6);
\draw[twist] (.2,0) to[out=90,in=90,looseness=1.4] (1,0);
\end{tikzpicture}\!&:
B^{\otimes 2}\stackrel{\sim}{\rightarrow} B_{O+,O-}\otimes_{R^{O-}}B_{O-,O+},&
\begin{tikzpicture}[anchorbase,scale=.9]
\node at (1.4,-.3) {$\labelp$};
\node at (-.2,-.3) {$\labelp$};
\node at (.6,-.15) {$\labelp$};
\node at (.6,-.45) {$\labelm$};
\draw[->] (0,0) to (0,-.6);
\draw[<-] (.4,0) to (.4,-.6);
\draw[->] (.8,0) to (.8,-.6);
\draw[<-] (1.2,0) to (1.2,-.6);
\draw[twist] (.2,0) to[out=-90,in=-90,looseness=1.4] (1,0);
\end{tikzpicture}\!&:
B_{O+,O-}\otimes_{R^{O-}}B_{O-,O+}
\stackrel{\sim}{\rightarrow}B^{\otimes 2}.\label{sanitycheck}
\end{align}
There are similar isomorphisms removing extra twisting bimodules in tensor products $B^{\otimes n}$ for all larger values of $n$, e.g., $B^{\otimes 3} \cong s_0 R_{O-} \otimes_{R^{O-}} B_{O-,O_+}\otimes_{R^{O+}}
B_{O+,O-}\otimes_{R^{O-}} B_{O-,O+}$.
Using these, relations involving extended string calculus become equivalent to relations involving ordinary string calculus. For an example, see \cref{chatGPTdoesmyhomework}.
\end{rem}

\subsection{Calculations}
Many of our interesting calculations require some additional understanding of the underlying Frobenius extensions. As a first example, 
we look again at 
the first equality in \cref{iamwell} .
The circle $\begin{tikzpicture}[anchorbase]
\draw[<-,colp] (.6,0) arc[start angle=0, end angle=360,radius=.2];
\node at (.4,0) {$\labelp$};
\end{tikzpicture}$ is, by definition, equal to multiplication by the product-coproduct element for the Frobenius extension $R^{O+} \hookrightarrow R^\O$. For practice, let us re-check that this product-coproduct element is $\eta^{\O}_{\O+}$, as we already known by the general theory. The longest element of $(W_{\O+} / W_\O)_{\min}$
is $s_{l-1} \cdots s_2 s_1$,
so the trace of the Frobenius
extension $R^{\O +} \hookrightarrow R^\O$ is
$\partial_{l-1} \cdots \partial_2 \partial_1$. For $i=1,\dots,l-1$, $\partial_i(f) = \frac{f-s_i(f)}{x_{i+1}-x_{i}}$ by our choice of simple roots in \cref{ssec4-1}. 
The reader can verify (a standard exercise in Schubert calculus) that an explicit pair of dual bases for
$R^\O$ as a free $R^{\O+}$-module are given by
$\big((-1)^r x_1^r\big)_{0 \leq r \leq l-1}$
and $\big(e_{l-1-r}(x_2,\dots,x_l)\big)_{0 \leq r \leq l-1}$.
Consequently
$\begin{tikzpicture}[anchorbase]
\draw[<-,colp] (.6,0) arc[start angle=0, end angle=360,radius=.2];
\node at (.4,0) {$\labelp$};
\end{tikzpicture}$
is multiplication by 
\[ \sum_{r=0}^{l-1} (-1)^r x_1^r e_{l-1-r}(x_2,\dots,x_l)= (x_2-x_1) \cdots (x_{l}-x_1) = \eta^{\O}_{\O+}. \]

Now recall that $\Sym = R^{\O+} = \kk[x_1,\dots,x_l]^{S_l}$.
The elementary symmetric polynomial 
$e_r(x_1,\dots,x_l)$ 
is the image of the elementary symmetric function $e_r \in \Lambda$ under the 
the natural evaluation map
$\ev_l:\Lambda \twoheadrightarrow \Sym$ from \cref{maggiesmith}.
We also let $h_r(x_1,\dots,x_l)$ and $q_r(x_1,\dots,x_l)$
be the evaluations of the $r$th complete symmetric function and the $q$-function $q_r$ from \cref{qfunc}, respectively.

\begin{lem}\label{missedflight}
For $r \geq 0$, we have that
$\partial_{l-1} \cdots \partial_2 \partial_1\big((-x_1)^r\big) = (-1)^{r-l+1} h_{r-l+1}(x_1,\dots,x_l)$ (interpreted as $0$ in case
$r < l-1$).
\end{lem}

\begin{proof}
Given that $x_1^r = h_r(x_1)$, this follows from the well-known fact that \[ \partial_i h_k(x_1, \ldots, x_i) = - h_{k-1}(x_1, \ldots, x_{i+1}).\]
\end{proof}

\begin{lem}
We have that 
\begin{align*}
\tr^\O_{\O+}\left(x_1^n\eta^{\O}_{\O-}\right) =
\begin{cases}
\frac{(-1)^{l-1}+1}{2} \id_{\Sym}&\text{if $n=0$}\\
-(-1)^{l}\frac{1}{2}q_n(x_1,\dots,x_l)\id_{\Sym}
&\text{if $n > 0$.}
\end{cases}
\end{align*}
(Later on, the number 
$t = \frac{(-1)^{l-1}+1}{2}$, i.e., 
the element of $\{0,1\}$ 
which is $\equiv l\pmod{2}$, appearing here and in the next corollary will match
the parameter in the nilBrauer category $\cNB_t$.)
\end{lem}

\begin{proof}
We use the identity $e_r(x_2,\dots,x_l) + x_1 e_{r-1}(x_2,\dots,x_l)= e_r(x_1,\dots,x_l)$,
which is valid for all $r \geq 0$ if one interprets $e_{-1}(x_2,\dots,x_l)$ as 0, to see that
\begin{align*}
\tr^\O_{\O+}&\left(x_1^n\eta^{\O}_{\O-}\right) =\partial_{l-1} \cdots \partial_1\left(
x_1^{n} (x_2+x_1)\cdots (x_{l}+x_1)\right)\\
&=\sum_{r=0}^{l-1}
\partial_{l-1} \cdots \partial_1\left(
x_1^{n+l-1-r} e_r(x_2,\dots,x_{l})\right)\\
&={\textstyle\frac{1}{2}}\sum_{r=0}^{l-1}
\partial_{l-1} \cdots \partial_1\left(
x_1^{n+l-1-r} e_r(x_2,\dots,x_l)\right)
+
{\textstyle\frac{1}{2}}
\sum_{r=1}^{l}
\partial_{l-1} \cdots \partial_1\left(
x_1^{n+l-r} e_{r-1}(x_2,\dots,x_l)\right)\\
&= {\textstyle\frac{1}{2}}\sum_{r=0}^{l-1}
\partial_{l-1} \cdots \partial_1\left(
x_1^{n+l-1-r} \big(e_r(x_2,\dots,x_l)+ x_1 e_{r-1}(x_2,\dots,x_l)\big)\right)+\\&
\hspace{80mm}{\textstyle\frac{1}{2}}\partial_{l-1} \cdots \partial_1\left(
x_1^{n} e_{l-1}(x_2,\dots,x_l)\right)\\
&= {\textstyle\frac{1}{2}}\sum_{r=0}^{l-1}
\partial_{l-1} \cdots \partial_1\left(x_1^{n+l-1-r} e_r(x_1,\dots,x_l)\right)+
{\textstyle\frac{1}{2}}\partial_{l-1} \cdots \partial_1\left(
x_1^{n} e_{l-1}(x_2,\dots,x_l)\right)\\
&= {\textstyle\frac{1}{2}}\sum_{r=0}^{l-1}
e_r(x_1,\dots,x_l)\partial_{l-1} \cdots \partial_1\left(x_1^{n+l-1-r} \right)+
{\textstyle\frac{1}{2}}\partial_{l-1} \cdots \partial_1\left(
x_1^{n} e_{l-1}(x_2,\dots,x_l)\right).
\end{align*}
Now we use \cref{missedflight} to simplify further.
When $n=0$, 
using also the easy identity
\begin{equation}
\partial_{l-1} \cdots \partial_1\big(x_2 x_3 \cdots x_l\big) = 1,
\end{equation}
it reduces to
$\frac{(-1)^{l-1} + 1}{2}$. When $n > 0$,
we obviously have that $x_1^n e_{l-1}(x_2,\dots,x_l) = x_1^{n-1} e_l(x_1,\dots,x_l)$,
and the expression simplifies to
\begin{multline*}
{\textstyle\frac{1}{2}}\sum_{r=0}^{l}
e_r(x_1,\dots,x_l)\partial_{l-1} \cdots \partial_1\left(x_1^{n+l-1-r} \right)
= 
(-1)^{l-1} {\textstyle\frac{1}{2}}\sum_{r=0}^{\min(l,n)}
e_r(x_1,\dots,x_l) h_{n-r}(x_1,\dots,x_l)\\ = (-1)^{l-1} {\textstyle\frac{1}{2}}\sum_{r=0}^{n}
e_r(x_1,\dots,x_l) h_{n-r}(x_1,\dots,x_l)= -(-1)^l {\textstyle\frac{1}{2}}q_n(x_1,\dots,x_l).
\end{multline*}
\end{proof}

\begin{cor}\label{coffeeisgood}
We have that
\begin{align*}
\begin{tikzpicture}[anchorbase]
\draw[colp,->] (.7,0) arc[start angle=0, end angle=360,radius=.7];
\draw[<-,colm] (-.12,0) arc[start angle=0, end angle=360,radius=.2];
\node [draw,rounded corners,inner sep=2pt,black,fill=gray!20!white] at (.3,0) {$\scriptstyle x_1^n$};
\node at (1,0) {$\labelp$};
\node at (-0.33,0) {$\labelp$};
\end{tikzpicture}
&=
\begin{cases}
\frac{(-1)^{l-1}+1}{2} \id_{\Sym}&\text{if $n=0$}\\
-(-1)^{l}\frac{1}{2}q_n(x_1,\dots,x_l)\id_{\Sym}
&\text{if $n > 0$.}
\end{cases}\end{align*}
\end{cor}

\begin{lem}\label{marriage}
Assuming $l \geq 3$, we have that
\begin{align}\label{marry1}
\begin{tikzpicture}[anchorbase,scale=1.2]
\draw[<-,colp] (.3,-.6) to (.3,.6);
\draw[->,colp] (-.3,-.6) to (-.3,.6);
\draw[-] (0,0) ellipse (.8 and .4);
\draw[->] (.8,0)
arc[start angle=-.5, end angle=.5,radius=.8];
\node at (.5,0) {$\labeltwo$};
\node at (-.5,0) {$\labeltwo$};
\node at (0,0) {$\labeltwo\labelp$}; 
\node at (0,-0.54) {$\labelp$};
\node at (0,0.54) {$\labelp$};
\end{tikzpicture}
&=\begin{tikzpicture}[anchorbase,scale=1.2]
\draw[<-,colp] (.25,-.6) to (.25,.6);
\draw[->,colp] (-.25,-.6) to (-.25,.6);
\node at (0,0) {$\labelp$};
\end{tikzpicture}.
\end{align}
\end{lem}

\begin{proof}
This is an instance of \cref{switchback}. 
One needs to check the switchback relation $\partial^{\O}_{\O+} = 
\partial^{\O\hat 2}_{\O}
\circ \partial^{\O\hat 2}_{\O\hat 2+}$, which follows
because both sides equal
$\partial_{l-1} \cdots \partial_1$.
\end{proof}

\begin{lem}\label{puppies}
Assuming $l \geq 3$,
the following relations hold:
\begin{align}\label{puppies1}
\begin{tikzpicture}[anchorbase,scale=1.2]
\draw[->,colm] (.9,-.6) to (.9,.6);
\draw[<-,colm] (-.9,-.6) to (-.9,.6);
\draw[<-,colp] (.3,-.6) to (.3,.6);
\draw[->,colp] (-.3,-.6) to (-.3,.6);
\node at (0,0) {$\labelp$};
\node at (-1.1,0) {$\labelp$};
\node at (1.1,0) {$\labelp$};
\end{tikzpicture}
&=
\frac{1}{2}
\left(
\begin{tikzpicture}[anchorbase,scale=1.7]
\draw[colp,<-] (.15,-.8) to (.15,.8);
\draw[colp,->] (-.15,-.8) to (-.15,.8);
	\draw[colm,<-] (.8,.8) to[out=-90, in=0,looseness=1.4] (0,.15) to[out = 180, in = -90,looseness=1.4] (-0.8,.8);
	\draw[colm,->] (.8,-.8) to[looseness=1.4,out=90, in=0] (0,-.15) to[looseness=1.4,out = 180, in = 90] (-0.8,-.8);
\draw (0,0) ellipse (.85 and .6);
\draw[->] (.85,0)
arc[start angle=-.5, end angle=.5,radius=.9];
\node at (1.05,0) {$\labelp$};
\node at (.68,0) {$\labelp$};\node at (.5,0.04) {$\labeltwo$};
\node at (-.5,0) {$\labelp$};\node at (-.68,0.04) {$\labeltwo$};
\node at (0.06,0) {$\labelpm$};
\node at (-.07,0) {$\labeltwo$};
    \node at (-1,0) {$\labelp$};
\node at (0,-0.7) {$\labelp$};
\node at (0,0.7) {$\labelp$};
\node at (.38,0.4) {$\labeltwo$};
\node at (-.21,0.4) {$\labeltwo$};
\node at (.38,-0.4) {$\labeltwo$};
\node at (-.38,-0.4) {$\labeltwo$};
\node at (.06,0.4) {$\labelp$};\node at (-.07,0.4) {$\labeltwo$};
\node at (.06,-0.4) {$\labelp$};\node at (-.07,-0.4) {$\labeltwo$};
\node [draw,rounded corners,inner sep=1.5pt,black,fill=gray!20!white] at (-.43,.35) {$\scriptstyle \alpha_{-1}$};
\end{tikzpicture}
+\begin{tikzpicture}[anchorbase,scale=1.7]
\draw[colp,<-] (.15,-.8) to (.15,.8);
\draw[colp,->] (-.15,-.8) to (-.15,.8);
	\draw[colm,<-] (.8,.8) to[out=-90, in=0,looseness=1.4] (0,.15) to[out = 180, in = -90,looseness=1.4] (-0.8,.8);
	\draw[colm,->] (.8,-.8) to[looseness=1.4,out=90, in=0] (0,-.15) to[looseness=1.4,out = 180, in = 90] (-0.8,-.8);
\draw (0,0) ellipse (.85 and .6);
\draw[->] (.85,0)
arc[start angle=-.5, end angle=.5,radius=.9];
\node at (1.05,0) {$\labelp$};
\node at (.68,0) {$\labelp$};\node at (.5,0.04) {$\labeltwo$};
\node at (-.5,0) {$\labelp$};\node at (-.68,0.04) {$\labeltwo$};
\node at (0.06,0) {$\labelpm$};
\node at (-.07,0) {$\labeltwo$};
    \node at (-1,0) {$\labelp$};
\node at (0,-0.7) {$\labelp$};
\node at (0,0.7) {$\labelp$};
\node at (.38,0.4) {$\labeltwo$};
\node at (-.38,0.4) {$\labeltwo$};
\node at (.21,-0.4) {$\labeltwo$};
\node at (-.38,-0.4) {$\labeltwo$};
\node at (.06,0.4) {$\labelp$};\node at (-.07,0.4) {$\labeltwo$};
\node at (.06,-0.4) {$\labelp$};\node at (-.07,-0.4) {$\labeltwo$};
\node [draw,rounded corners,inner sep=1.5pt,black,fill=gray!20!white] at (.43,-.35) {$\scriptstyle \alpha_{-1}$};
\end{tikzpicture}
\right),\\
\begin{tikzpicture}[anchorbase,scale=1.2]
\draw[->,colm] (.7,-.9) to[looseness=1.8,out=90,in=90] (-.7,-.9);
\draw[->,colm] (-.7,.9) to[looseness=1.7,out=-90,in=-90] (.7,.9);
\draw[<-,colp] (.3,-.9) to[out=90,in=90,looseness=1.8] (-.3,-.9);
\draw[<-,colp] (-.3,.9) to[looseness=1.8,out=-90,in=-90] (.3,.9);
\node at (0,.8) {$\labelp$};
\node at (0,-.8) {$\labelp$};
\node at (0,0) {$\labelp$};
\end{tikzpicture}\:\;\:
&=
\frac{1}{2}
\left(
\begin{tikzpicture}[anchorbase,scale=1.7]
\draw[colp,<-] (.15,-.8) to (.15,.8);
\draw[colp,->] (-.15,-.8) to (-.15,.8);
	\draw[colm,<-] (.8,.8) to[out=-90, in=0,looseness=1.4] (0,.15) to[out = 180, in = -90,looseness=1.4] (-0.8,.8);
	\draw[colm,->] (.8,-.8) to[looseness=1.4,out=90, in=0] (0,-.15) to[looseness=1.4,out = 180, in = 90] (-0.8,-.8);
\draw (0,0) ellipse (.85 and .6);
\draw[->] (.85,0)
arc[start angle=-.5, end angle=.5,radius=.9];
\node at (1.05,0) {$\labelp$};
\node at (.68,0) {$\labelp$};\node at (.5,0) {$\labeltwo$};
\node at (-.5,0) {$\labelp$};\node at (-.68,0) {$\labeltwo$};
\node at (0.06,0) {$\labelpm$};
\node at (-.07,0) {$\labeltwo$};
    \node at (-1,0) {$\labelp$};
\node at (0,-0.7) {$\labelp$};
\node at (0,0.7) {$\labelp$};
\node at (.38,0.42) {$\labeltwo$};
\node at (-.21,0.42) {$\labeltwo$};
\node at (.38,-0.42) {$\labeltwo$};
\node at (-.38,-0.42) {$\labeltwo$};
\node at (.06,0.42) {$\labelp$};\node at (-.07,0.42) {$\labeltwo$};
\node at (.06,-0.42) {$\labelp$};\node at (-.07,-0.42) {$\labeltwo$};
\node [draw,rounded corners,inner sep=1.5pt,black,fill=gray!20!white] at (-.43,.34) {$\scriptstyle \alpha_{1}$};
\end{tikzpicture}
+\begin{tikzpicture}[anchorbase,scale=1.7]
\draw[colp,<-] (.15,-.8) to (.15,.8);
\draw[colp,->] (-.15,-.8) to (-.15,.8);
	\draw[colm,<-] (.8,.8) to[out=-90, in=0,looseness=1.4] (0,.15) to[out = 180, in = -90,looseness=1.4] (-0.8,.8);
	\draw[colm,->] (.8,-.8) to[looseness=1.4,out=90, in=0] (0,-.15) to[looseness=1.4,out = 180, in = 90] (-0.8,-.8);
\draw (0,0) ellipse (.85 and .6);
\draw[->] (.85,0)
arc[start angle=-.5, end angle=.5,radius=.9];
\node at (1.05,0) {$\labelp$};
\node at (.68,0) {$\labelp$};\node at (.5,0) {$\labeltwo$};
\node at (-.5,0) {$\labelp$};\node at (-.68,0) {$\labeltwo$};
\node at (0.06,0) {$\labelpm$};
\node at (-.07,0) {$\labeltwo$};
    \node at (-1,0) {$\labelp$};
\node at (0,-0.7) {$\labelp$};
\node at (0,0.7) {$\labelp$};
\node at (.38,0.42) {$\labeltwo$};
\node at (-.38,0.42) {$\labeltwo$};
\node at (.21,-0.42) {$\labeltwo$};
\node at (-.38,-0.42) {$\labeltwo$};
\node at (.06,0.42) {$\labelp$};\node at (-.07,0.42) {$\labeltwo$};
\node at (.06,-0.42) {$\labelp$};\node at (-.07,-0.42) {$\labeltwo$};
\node [draw,rounded corners,inner sep=1.5pt,black,fill=gray!20!white] at (.43,-.34) {$\scriptstyle \alpha_{1}$};
\end{tikzpicture}
\right).\label{puppies2}
\end{align}
There are analogous relations in the case $l=2$: one just has to omit $2$-colored circles and all labels $\hat 2$.
\end{lem}

\begin{proof}
Suppose first that $l \geq 3$.
Using that
$\eta^{\O\hat 2}_{\O\hat 2 \pm} = 1$ 
and taking the basis $\B_{\O\hat 2 \pm}^{\O\hat 2+}$ to be $1,\alpha_{-1}$
with dual basis $\frac{1}{2}\alpha_{-1}, \frac{1}{2}$, we have that
\begin{multline*}
\begin{tikzpicture}[anchorbase,scale=1.2]
\draw[->,colm] (.9,-.6) to (.9,.6);
\draw[<-,colm] (-.9,-.6) to (-.9,.6);
\draw[<-,colp] (.3,-.6) to (.3,.6);
\draw[->,colp] (-.3,-.6) to (-.3,.6);
\node at (0,0) {$\labelp$};
\node at (-1.1,0) {$\labelp$};
\node at (1.1,0) {$\labelp$};
\end{tikzpicture}
\stackrel{\cref{marry1}}{=}\begin{tikzpicture}[anchorbase,scale=1.2]
\node at (-1.1,0) {$\labelp$};
\node at (1.1,0) {$\labelp$};
\draw[->,colm] (.9,-.6) to (.9,.6);
\draw[<-,colm] (-.9,-.6) to (-.9,.6);
\draw[<-,colp] (.3,-.6) to (.3,.6);
\draw[->,colp] (-.3,-.6) to (-.3,.6);
\draw[-] (0,0) ellipse (.75 and .4);
\draw[->] (.75,0)
arc[start angle=-.5, end angle=.5,radius=.75];
\node at (.5,0) {$\labeltwo$};
\node at (-.5,0) {$\labeltwo$};
\node at (0,0) {$\labeltwo\labelp$}; 
\node at (0,-0.54) {$\labelp$};
\node at (0,0.54) {$\labelp$};
\end{tikzpicture}
\stackrel{\cref{easyRII}}{=}
\begin{tikzpicture}[anchorbase,scale=1.2]
\node at (-1.2,0) {$\labelp$};
\node at (1.2,0) {$\labelp$};
\draw[->,colm] (.6,-.6) to (.6,.6);
\draw[<-,colm] (-.6,-.6) to (-.6,.6);
\draw[<-,colp] (.25,-.6) to (.25,.6);
\draw[->,colp] (-.25,-.6) to (-.25,.6);
\draw[-] (0,0) ellipse (1 and .4);
\draw[->] (1,0)
arc[start angle=-.5, end angle=.5,radius=1];
\node at (.42,0) {$\labeltwo$};
\node at (-.42,0) {$\labeltwo$};
\node at (0,0) {$\labeltwo\labelp$}; 
\node at (.77,0) {$\labeltwo\labelp$}; 
\node at (-.78,0) {$\labeltwo\labelp$}; 
\node at (0,-0.54) {$\labelp$};
\node at (0,0.54) {$\labelp$};
\end{tikzpicture}\\
\stackrel{\cref{distantRII}}{=}\begin{tikzpicture}[anchorbase,scale=1.6]
\node at (-1.2,0) {$\labelp$};
\node at (1.2,0) {$\labelp$};
\draw[->,colm] (.6,-.8) to[looseness=2,out=90,in=-90] (.2,0) to[looseness=2,out=90,in=-90] (.6,.8);
\draw[<-,colm] (-.6,-.8) to[looseness=2,out=90,in=-90] (-.2,0) to[out=90,in=-90,looseness=2] (-.6,.8);
\draw[<-,colp] (.2,-.8) to[out=90,in=-90,looseness=2] (.6,0) to[out=90,in=-90,looseness=2] (.2,.8);
\draw[->,colp] (-.2,-.8) to[out=90,in=-90,looseness=2] (-.6,0) to[out=90,in=-90,looseness=2] (-.2,.8);
\draw[-] (0,0) ellipse (1 and .66);
\draw[->] (1,0)
arc[start angle=-.5, end angle=.5,radius=1];
\node at (.32,0) {$\labeltwo$};\node at (.47,0) {$\labelpm$};
\node at (-.35,0) {$\labelpm$};\node at (-.5,0) {$\labeltwo$};
\node at (.4,.51) {$\labeltwo$};
\node at (.4,-.5) {$\labeltwo$};
\node at (-.4,.51) {$\labeltwo$};
\node at (-.4,-.5) {$\labeltwo$};
\node at (0,0) {$\labeltwo\labelp$}; 
\node at (.78,0) {$\labeltwo\labelp$}; 
\node at (-.79,0) {$\labeltwo\labelp$}; 
\node at (0,-0.75) {$\labelp$};
\node at (0,0.75) {$\labelp$};
\end{tikzpicture}\stackrel{\cref{cupcap}}{=}
\frac{1}{2}\left(
\begin{tikzpicture}[anchorbase,scale=1.6]
\draw[colp,<-] (.25,-.8) to (.25,.8);
\draw[colp,->] (-.25,-.8) to (-.25,.8);
	\draw[colm,<-] (.8,.8) to[out=-90, in=0,looseness=1.4] (0,.15) to[out = 180, in = -90,looseness=1.4] (-0.8,.8);
	\draw[colm,->] (.8,-.8) to[looseness=1.4,out=90, in=0] (0,-.15) to[looseness=1.4,out = 180, in = 90] (-0.8,-.8);
\draw (0,0) ellipse (.85 and .6);
\draw[->] (.85,0)
arc[start angle=-.5, end angle=.5,radius=.9];
\node at (1.05,0) {$\labelp$};
\node at (.68,0) {$\labelp$};\node at (.5,0) {$\labeltwo$};
\node at (-.5,0) {$\labelp$};\node at (-.68,0) {$\labeltwo$};
\node at (0.06,0) {$\labelpm$};
\node at (-.07,0) {$\labeltwo$};
    \node at (-1,0) {$\labelp$};
\node at (0,-0.7) {$\labelp$};
\node at (0,0.7) {$\labelp$};
\node at (.38,0.4) {$\labeltwo$};
\node at (-.38,0.4) {$\labeltwo$};
\node at (.38,-0.4) {$\labeltwo$};
\node at (-.38,-0.4) {$\labeltwo$};
\node at (.06,0.48) {$\labelp$};\node at (-.07,0.5) {$\labeltwo$};
\node at (.06,-0.4) {$\labelp$};\node at (-.07,-0.4) {$\labeltwo$};
\node [draw,rounded corners,inner sep=1.5pt,black,fill=gray!20!white] at (0,.29) {$\scriptstyle \alpha_{-1}$};
\end{tikzpicture}
+\begin{tikzpicture}[anchorbase,scale=1.6]
\draw[colp,<-] (.25,-.8) to (.25,.8);
\draw[colp,->] (-.25,-.8) to (-.25,.8);
	\draw[colm,<-] (.8,.8) to[out=-90, in=0,looseness=1.4] (0,.15) to[out = 180, in = -90,looseness=1.4] (-0.8,.8);
	\draw[colm,->] (.8,-.8) to[looseness=1.4,out=90, in=0] (0,-.15) to[looseness=1.4,out = 180, in = 90] (-0.8,-.8);
\draw (0,0) ellipse (.85 and .6);
\draw[->] (.85,0)
arc[start angle=-.5, end angle=.5,radius=.9];
\node at (1.05,0) {$\labelp$};
\node at (.68,0) {$\labelp$};\node at (.5,0) {$\labeltwo$};
\node at (-.5,0) {$\labelp$};\node at (-.68,0) {$\labeltwo$};
\node at (0.06,0) {$\labelpm$};
\node at (-.07,0) {$\labeltwo$};
    \node at (-1,0) {$\labelp$};
\node at (0,-0.7) {$\labelp$};
\node at (0,0.7) {$\labelp$};
\node at (.38,0.4) {$\labeltwo$};
\node at (-.38,0.4) {$\labeltwo$};
\node at (.38,-0.4) {$\labeltwo$};
\node at (-.38,-0.4) {$\labeltwo$};
\node at (.06,-0.5) {$\labelp$};\node at (-.07,-0.5) {$\labeltwo$};
\node at (.06,0.4) {$\labelp$};\node at (-.07,0.4) {$\labeltwo$};
\node [draw,rounded corners,inner sep=1.5pt,black,fill=gray!20!white] at (0,-.29) {$\scriptstyle \alpha_{-1}$};
\end{tikzpicture}\right).
\end{multline*}
Now we 
slide the $\alpha_{-1}$ bubbles across the special undotted strings using \cref{howareyou}.
To deduce the second identity, 
recall the equality \cref{adele}.
By expanding the special undotted and dotted strings in terms of ordinary and dashed strings,
then
using the easy relations satisfied by dashed strings,
the identity just displayed 
can be written in an equivalent form as
\begin{equation} \label{streamofconsciousness}
\begin{tikzpicture}[anchorbase,scale=1.2]
\draw[->,colp] (.9,-.6) to (.9,.6);
\draw[<-,colp] (-.9,-.6) to (-.9,.6);
\draw[<-,colm] (.3,-.6) to (.3,.6);
\draw[->,colm] (-.3,-.6) to (-.3,.6);
\node at (0,0) {$\labelp$};
\node at (-1.1,0) {$\labelp$};
\node at (1.1,0) {$\labelp$};
\end{tikzpicture}
=
\frac{1}{2}
\left(
\begin{tikzpicture}[anchorbase,scale=1.6]
\draw[colm,<-] (.25,-.8) to (.25,.8);
\draw[colm,->] (-.25,-.8) to (-.25,.8);
	\draw[colp,<-] (.8,.8) to[out=-90, in=0,looseness=1.4] (0,.15) to[out = 180, in = -90,looseness=1.4] (-0.8,.8);
	\draw[colp,->] (.8,-.8) to[looseness=1.4,out=90, in=0] (0,-.15) to[looseness=1.4,out = 180, in = 90] (-0.8,-.8);
\draw (0,0) ellipse (.85 and .6);
\draw[->] (.85,0)
arc[start angle=-.5, end angle=.5,radius=.9];
\node at (1.05,0) {$\labelp$};
\node at (.68,0) {$\labelp$};\node at (.5,0) {$\labeltwo$};
\node at (-.5,0) {$\labelp$};\node at (-.68,0) {$\labeltwo$};
\node at (0.06,0) {$\labelpm$};
\node at (-.07,0) {$\labeltwo$};
    \node at (-1,0) {$\labelp$};
\node at (0,-0.7) {$\labelp$};
\node at (0,0.7) {$\labelp$};
\node at (.38,0.4) {$\labeltwo$};
\node at (-.38,0.4) {$\labeltwo$};
\node at (.38,-0.4) {$\labeltwo$};
\node at (-.38,-0.4) {$\labeltwo$};
\node at (.06,0.48) {$\labelp$};\node at (-.07,0.5) {$\labeltwo$};
\node at (.06,-0.4) {$\labelp$};\node at (-.07,-0.4) {$\labeltwo$};
\node [draw,rounded corners,inner sep=1.5pt,black,fill=gray!20!white] at (0,.29) {$\scriptstyle \alpha_{-1}$};
\end{tikzpicture}
+\begin{tikzpicture}[anchorbase,scale=1.6]
\draw[colm,<-] (.25,-.8) to (.25,.8);
\draw[colm,->] (-.25,-.8) to (-.25,.8);
	\draw[colp,<-] (.8,.8) to[out=-90, in=0,looseness=1.4] (0,.15) to[out = 180, in = -90,looseness=1.4] (-0.8,.8);
	\draw[colp,->] (.8,-.8) to[looseness=1.4,out=90, in=0] (0,-.15) to[looseness=1.4,out = 180, in = 90] (-0.8,-.8);
\draw (0,0) ellipse (.85 and .6);
\draw[->] (.85,0)
arc[start angle=-.5, end angle=.5,radius=.9];
\node at (1.05,0) {$\labelp$};
\node at (.68,0) {$\labelp$};\node at (.5,0) {$\labeltwo$};
\node at (-.5,0) {$\labelp$};\node at (-.68,0) {$\labeltwo$};
\node at (0.06,0) {$\labelpm$};
\node at (-.07,0) {$\labeltwo$};
    \node at (-1,0) {$\labelp$};
\node at (0,-0.7) {$\labelp$};
\node at (0,0.7) {$\labelp$};
\node at (.38,0.4) {$\labeltwo$};
\node at (-.38,0.4) {$\labeltwo$};
\node at (.38,-0.4) {$\labeltwo$};
\node at (-.38,-0.4) {$\labeltwo$};
\node at (.06,-0.5) {$\labelp$};\node at (-.07,-0.5) {$\labeltwo$};
\node at (.06,0.4) {$\labelp$};\node at (-.07,0.4) {$\labeltwo$};
\node [draw,rounded corners,inner sep=1.5pt,black,fill=gray!20!white] at (0,-.29) {$\scriptstyle \alpha_{-1}$};
\end{tikzpicture}\right).
\end{equation}
Now we rotate counterclockwise through 90$^\circ$ using the cyclic structure, then slide the $\alpha_{-1}$ bubbles across the special dotted strings, which flips them to $\alpha_1$ bubbles.

This argument is easily modified to treat the case $l=2$, when there is no color 2.
One can skip the first line altogether since
\cref{distantRII} can be applied directly (without needing to remove any simple reflections).
\end{proof}

\begin{rem} \label{chatGPTdoesmyhomework} Pre- and post-composing \cref{puppies1} with the isomorphisms from \cref{benhomework}, we get the following equivalent relation in the ordinary string calculus:
\begin{equation*}
\begin{tikzpicture}[anchorbase,scale=1.2]
\draw[->] (.9,-.6) to (.9,.6);
\draw[<-] (-.9,-.6) to (-.9,.6);
\draw[<-] (.3,-.6) to (.3,.6);
\draw[->] (-.3,-.6) to (-.3,.6);
\node at (0,0) {$\labelm$};
\node at (-1.1,0) {$\labelp$};
\node at (1.1,0) {$\labelp$};
\end{tikzpicture}
=
\frac{1}{2}
\left(
\begin{tikzpicture}[anchorbase,scale=1.7]
\draw[<-] (.15,-.8) to (.15,.8);
\draw[->] (-.15,-.8) to (-.15,.8);
	\draw[<-] (.8,.8) to[out=-90, in=0,looseness=1.4] (0,.15) to[out = 180, in = -90,looseness=1.4] (-0.8,.8);
	\draw[->] (.8,-.8) to[looseness=1.4,out=90, in=0] (0,-.15) to[looseness=1.4,out = 180, in = 90] (-0.8,-.8);
\draw (0,0) ellipse (.85 and .6);
\draw[->] (.85,0)
arc[start angle=-.5, end angle=.5,radius=.9];
\node at (1.05,0) {$\labelp$};
\node at (.68,0) {$\labelp$};\node at (.5,0.04) {$\labeltwo$};
\node at (-.5,0) {$\labelp$};\node at (-.68,0.04) {$\labeltwo$};
\node at (0.06,0) {$\labelpm$};
\node at (-.07,0) {$\labeltwo$};
    \node at (-1,0) {$\labelp$};
\node at (0,-0.7) {$\labelm$};
\node at (0,0.7) {$\labelm$};
\node at (.38,0.4) {$\labeltwo$};
\node at (-.21,0.4) {$\labeltwo$};
\node at (.38,-0.4) {$\labeltwo$};
\node at (-.38,-0.4) {$\labeltwo$};
\node at (.06,0.4) {$\labelm$};\node at (-.07,0.4) {$\labeltwo$};
\node at (.06,-0.4) {$\labelm$};\node at (-.07,-0.4) {$\labeltwo$};
\node [draw,rounded corners,inner sep=1.5pt,black,fill=gray!20!white] at (-.43,.35) {$\scriptstyle \alpha_{1}$};
\end{tikzpicture}
+\begin{tikzpicture}[anchorbase,scale=1.7]
\draw[<-] (.15,-.8) to (.15,.8);
\draw[->] (-.15,-.8) to (-.15,.8);
	\draw[<-] (.8,.8) to[out=-90, in=0,looseness=1.4] (0,.15) to[out = 180, in = -90,looseness=1.4] (-0.8,.8);
	\draw[->] (.8,-.8) to[looseness=1.4,out=90, in=0] (0,-.15) to[looseness=1.4,out = 180, in = 90] (-0.8,-.8);
\draw (0,0) ellipse (.85 and .6);
\draw[->] (.85,0)
arc[start angle=-.5, end angle=.5,radius=.9];
\node at (1.05,0) {$\labelp$};
\node at (.68,0) {$\labelp$};\node at (.5,0.04) {$\labeltwo$};
\node at (-.5,0) {$\labelp$};\node at (-.68,0.04) {$\labeltwo$};
\node at (0.06,0) {$\labelpm$};
\node at (-.07,0) {$\labeltwo$};
    \node at (-1,0) {$\labelp$};
\node at (0,-0.7) {$\labelm$};
\node at (0,0.7) {$\labelm$};
\node at (.38,0.4) {$\labeltwo$};
\node at (-.38,0.4) {$\labeltwo$};
\node at (.21,-0.4) {$\labeltwo$};
\node at (-.38,-0.4) {$\labeltwo$};
\node at (.06,0.4) {$\labelm$};\node at (-.07,0.4) {$\labeltwo$};
\node at (.06,-0.4) {$\labelm$};\node at (-.07,-0.4) {$\labeltwo$};
\node [draw,rounded corners,inner sep=1.5pt,black,fill=gray!20!white] at (.43,-.35) {$\scriptstyle \alpha_{1}$};
\end{tikzpicture}
\right).
\end{equation*}
This illustrates the alternate approach mentioned in \cref{ssec-interlude}.
\end{rem}

Now we consider the $n$th tensor power
$B^{\otimes n}$ of the $(\Sym, \Sym)$-bimodule $B$ (the tensor product being over $\Sym=R^{\O+}$).
Its character is the element $b^{\m n} \in \V$ from \cref{wrongspot}.
We are going to investigate certain endomorphisms which will be shown later to be primitive homogeneous idempotents.

\begin{lem}\label{handy3}
If $0 \leq n \leq l$ then the endomorphism algebra of 
$B^{\otimes n}$
is 1-dimensional in degree $-2\binom{n}{2}$
and 0 in all smaller degrees.
If $n > l$ then this endomorphism algebra 
is 0 in all degrees $\leq -2\binom{n}{2}$.
\end{lem}

\begin{proof}
If $n$ is even, this
follows from
the homomorphism formula (\cref{homformula})
and \cref{significant}. If $n$ is odd, it follows from the modified homomorphism formula discussed in \cref{ssec4-3} (or the original one can be applied 
after left-tensoring with $s_0 \Sym$). When $n > l$ the degree bound from the homomorphism formula gives that the lowest degree is $-(2n-l)(l-1)$, so one also needs to observe that $-2\binom{n}{2} < -(2n-l)(l-1)$. This may be seen by replacing $n$ by $l+k$ for $k > 0$, then checking that
$$
2 \binom{l+k}{2} - (2(l+k)-l)(l-1) = 
(l+k)(l+k-1) - (l+2k)(l-1) = k(k+1)> 0.
$$
\end{proof}

For $1 \leq n \leq l$, we define an endomorphism 
$\u_n \in \End_{\Sym\dash \Sym}(B^{\otimes n})$ by setting
\begin{align}\label{aplus}
\u_n &:= 
\begin{cases}
\begin{tikzpicture}[anchorbase,scale=1.2]
\node at (-.45,0) {$\labelp$};
\node at (.45,0) {$\labelp$};
\node at (0,0) {$\labelpm$};
	\draw[colp] (0.28,0) to[out=90,in=-90] (-0.28,.6);
	\draw[colm,->] (-0.28,0) to[out=90,in=-90] (0.28,.6);
	\draw[colm] (0.28,-.6) to[out=90,in=-90] (-0.28,0);
	\draw[colp,<-](-0.28,-.6) to[out=90,in=-90] (0.28,0);
\end{tikzpicture}
&\text{if $n=1$}\\
\begin{tikzpicture}[anchorbase,scale=1.6]
\node at (1,.6) {$\phantom{space!}$};
\node at (1,-.6) {$\phantom{space!}$};
\node at (1.5,.4) {$\labelp$};
\node at (1.1,.4) {$\labelp$};
\node at (-.9,.4) {$\labelp$};
\node at (-.5,.4) {$\labelp$};
\node at (-0.1,.4) {$\labelp$};
\node at (1.9,0) {$\labelp$};
\node at (.5,0) {$\labelpm$};
\node at (-1.3,0) {$\labelp$};
	\draw[colm,->] (-1.2,.6) to[out=-90, in=180] (-.9,.2) to[out = 0, in = -90] (-0.6,.6);
	\draw[colp,->] (-.8,.6) to[out=-90, in=180] (-.5,.2) to[out = 0, in = -90] (-0.2,.6);
	\draw[colm,->] (-.4,.6) to[out=-90, in=180] (-.1,.2) to[out = 0, in = -90] (0.2,.6);
 \node at (.5,.5) {$\cdots$};
\node at (1.1,-.4) {$\labelp$};
\node at (1.5,-.4) {$\labelp$};
\node at (-.9,-.4) {$\labelp$};
\node at (-.5,-.4) {$\labelp$};
\node at (-0.1,-.4) {$\labelp$};
 \draw[colp,->] (.8,.6) to[out=-90, in=180] (1.1,.2) to[out = 0, in = -90] (1.4,.6);
	\draw[colm,->] (1.2,.6) to[out=-90, in=180]
 (1.5,.2) to[out = 0, in = -90] (1.8,.6);
	\draw[colm,<-] (-1.2,-.6) to[out=90, in=180] (-.9,-.2) to[out = 0, in = 90] (-0.6,-.6);
	\draw[colp,<-] (-.8,-.6) to[out=90, in=180] (-.5,-.2) to[out = 0, in = 90] (-0.2,-.6);
	\draw[colm,<-] (-.4,-.6) to[out=90, in=180] (-.1,-.2) to[out = 0, in = 90] (0.2,-.6);
 \node at (.5,-.5) {$\cdots$};
	\draw[colp,<-] (.8,-.6) to[out=90, in=180] (1.1,-.2) to[out = 0, in = 90] (1.4,-.6);
	\draw[colm,<-] (1.2,-.6) to[out=90, in=180] (1.5,-.2) to[out = 0, in = 90] (1.8,-.6);
 \draw[colp,<-] (1.6,-.6) to (1.6,.6);
 \draw[colp,->] (-1,-.6) to (-1,.6);
\end{tikzpicture}&\text{if $n$ is even}\\
\begin{tikzpicture}[anchorbase,scale=1.6]
\node at (1.5,.4) {$\labelp$};
\node at (1.1,.4) {$\labelp$};
\node at (-.5,.4) {$\labelp$};
\node at (-0.1,.4) {$\labelp$};
\node at (-.9,0) {$\labelp$};
\node at (.5,0) {$\labelpm$};
\node at (1.9,0) {$\labelp$};
\node at (1.5,-.4) {$\labelp$};
\node at (1.1,-.4) {$\labelp$};
\node at (-.5,-.4) {$\labelp$};
\node at (-0.1,-.4) {$\labelp$};
	\draw[colp,->] (-.8,.6) to[out=-90, in=180] (-.5,.2) to[out = 0, in = -90] (-0.2,.6);
	\draw[colm,->] (-.4,.6) to[out=-90, in=180] (-.1,.2) to[out = 0, in = -90] (0.2,.6);
 \node at (.5,.5) {$\cdots$};
	\draw[colp,->] (.8,.6) to[out=-90, in=180] (1.1,.2) to[out = 0, in = -90] (1.4,.6);
	\draw[colm,->] (1.2,.6) to[out=-90, in=180]
 (1.5,.2) to[out = 0, in = -90] (1.8,.6);
	\draw[colp,<-] (-.8,-.6) to[out=90, in=180] (-.5,-.2) to[out = 0, in = 90] (-0.2,-.6);
	\draw[colm,<-] (-.4,-.6) to[out=90, in=180] (-.1,-.2) to[out = 0, in = 90] (0.2,-.6);
 \node at (.5,-.5) {$\cdots$};
	\draw[colp,<-] (.8,-.6) to[out=90, in=180] (1.1,-.2) to[out = 0, in = 90] (1.4,-.6);
	\draw[colm,<-] (1.2,-.6) to[out=90, in=180] (1.5,-.2) to[out = 0, in = 90] (1.8,-.6);
 \draw[colp,<-] (1.6,-.6) to (1.6,.6);
 \draw[colm,->] (-.6,-.6) to (-.6,.6);
 \end{tikzpicture}
 &\text{if $n\neq 1$ is odd,}
\end{cases}
\end{align}
where there are a total of 
$n$ $\down\up$-pairs of strings at the top and bottom, $(n-1)$ cups at the top, and $(n-1)$ caps at the bottom.
We have that $\u_n = \check\u_n \circ \hat\u_n$
where $\check\u_n$ denotes the top half of the above diagram
and $\hat\u_n$ is its bottom half.

\begin{lem}\label{play1}
For $1 \leq n \leq l$, both of the endomorphisms $\check\u_{n}$
and $\hat\u_n$
are of degree $\binom{l}{2}-n(l-1)$.
This equals $-\binom{n}{2}$ if $n=l-1$ or $n=l$.
\end{lem}

\begin{proof}
It is clear that $\check\u_n$ and $\hat\u_n$
are of the same degree. 
To compute the degree of $\check\u_n$, 
we use 
\cref{lengths} to see:
\begin{itemize}
\item
Each of the $(n-1)$ cups is
of degree $\ell(w_{\O+}) - \ell(w_{\O\pm})
=\ell(w_{\O-}) - \ell(w_{\O\pm})=-\binom{l}{2}$.
\item
Each of the $n$ crossings is of degree
$\ell(w_{\O\pm}) + \ell(w_\O) - \ell(w_{\O+}) -\ell(w_{\O-})
= \binom{l-1}{2}$.
\end{itemize}
This gives the total degree
$-(n-1)\binom{l}{2}+n\binom{l-1}{2}=
\binom{l}{2} - n(l-1)$.
\end{proof}

Let $\v_0 := \id_{\Sym}$,
$\v_1 := \begin{tikzpicture}[anchorbase,scale=1]
\draw[colm,<-] (0,0.05) to (0,.5);
\node at (-.24,.3) {$\labelp$};
\draw[colp,->] (0.3,0.05) to (0.3,.5);
\node at (.56,.3) {$\labelp$};
\end{tikzpicture}$,
and $\v_l := \u_l$.
Then for $2 \leq n \leq l-1$,
we add an $n$-colored circle to the middle of the diagram defining $\u_{n}$ to obtain
\begin{align}\label{bplus}
\v_n := 
\begin{cases}
\begin{tikzpicture}[anchorbase,scale=2.4]
\node at (1.5,.7) {$\labelp$};
\node at (1.1,.7) {$\labelp$};
\node at (-.9,.7) {$\labelp$};
\node at (-.5,.7) {$\labelp$};
\node at (-0.1,.7) {$\labelp$};
\node at (1.05,.4) {$\labeln$};\node at (1.13,.4) {$\labelp$};
\node at (1.45,.4) {$\labeln$};\node at (1.53,.4) {$\labelp$};
\node at (-.95,.4) {$\labeln$};\node at (-.87,.4) {$\labelp$};
\node at (-.55,.4) {$\labeln$};\node at (-.47,.4) {$\labelp$};
\node at (-.15,.4) {$\labeln$};\node at (-0.07,.4) {$\labelp$};
\node at (1.3,.4) {$\labeln$};
\node at (1.68,.4) {$\labeln$};
\node at (-1.08,.4) {$\labeln$};
\node at (-.7,.4) {$\labeln$};
\node at (-.3,.4) {$\labeln$};
\node at (1.75,0) {$\labeln$};\node at (1.83,0) {$\labelp$};
\node at (.5,0) {$\labeln$};\node at (.6,0) {$\labelpm$};
\node at (2.1,0) {$\labelp$};
\node at (-1.25,0) {$\labeln$};\node at (-1.17,0) {$\labelp$};
\node at (-1.5,0) {$\labelp$};
	\draw[colm,->] (-1.2,.8) to[out=-90, in=180] (-.9,.2) to[out = 0, in = -90] (-0.6,.8);
	\draw[colp,->] (-.8,.8) to[out=-90, in=180] (-.5,.2) to[out = 0, in = -90] (-0.2,.8);
	\draw[colm,->] (-.4,.8) to[out=-90, in=180] (-.1,.2) to[out = 0, in = -90] (0.2,.8);
 \node at (.5,.7) {$\cdots$};
\node at (1.1,-.7) {$\labelp$};
\node at (1.5,-.7) {$\labelp$};
\node at (-.9,-.7) {$\labelp$};
\node at (-.5,-.7) {$\labelp$};
\node at (-0.1,-.7) {$\labelp$};
\node at (1.05,-.4) {$\labeln$};\node at (1.13,-.4) {$\labelp$};
\node at (1.45,-.4) {$\labeln$};\node at (1.53,-.4) {$\labelp$};
\node at (-.95,-.4) {$\labeln$};\node at (-.87,-.4) {$\labelp$};
\node at (-.55,-.4) {$\labeln$};\node at (-.47,-.4) {$\labelp$};
\node at (-.15,-.4) {$\labeln$};\node at (-0.07,-.4) {$\labelp$};
\node at (1.3,-.4) {$\labeln$};
\node at (1.68,-.4) {$\labeln$};
\node at (-1.08,-.4) {$\labeln$};
\node at (-.7,-.4) {$\labeln$};
\node at (-.3,-.4) {$\labeln$};
 \draw[colp,->] (.8,.8) to[out=-90, in=180] (1.1,.2) to[out = 0, in = -90] (1.4,.8);
	\draw[colm,->] (1.2,.8) to[out=-90, in=180]
 (1.5,.2) to[out = 0, in = -90] (1.8,.8);
	\draw[colm,<-] (-1.2,-.8) to[out=90, in=180] (-.9,-.2) to[out = 0, in = 90] (-0.6,-.8);
	\draw[colp,<-] (-.8,-.8) to[out=90, in=180] (-.5,-.2) to[out = 0, in = 90] (-0.2,-.8);
	\draw[colm,<-] (-.4,-.8) to[out=90, in=180] (-.1,-.2) to[out = 0, in = 90] (0.2,-.8);
 \node at (.5,-.7) {$\cdots$};
	\draw[colp,<-] (.8,-.8) to[out=90, in=180] (1.1,-.2) to[out = 0, in = 90] (1.4,-.8);
	\draw[colm,<-] (1.2,-.8) to[out=90, in=180] (1.5,-.2) to[out = 0, in = 90] (1.8,-.8);
 \draw[colp,<-] (1.6,-.8) to (1.6,.8);
 \draw[colp,->] (-1,-.8) to (-1,.8);
 \draw[-] (-1.4,0) to[out=90,looseness=1,in=180] (0.3,.6);
 \draw[dotted] (0.3,.6) to (.7,.6);
 \draw[-] (.7,.6) to[out=0,in=90,looseness=1] (2,0);
\draw[-] (-1.4,0) to[out=-90,in=180,looseness=1] (0.3,-.6);
 \draw[dotted] (0.3,-.6) to (.7,-.6);
\draw[->] (.7,-.6) to[out=0,in=-90,looseness=1] (2,0);
\end{tikzpicture}
&\text{if $n$ is even}\\
\begin{tikzpicture}[anchorbase,scale=2.4]
\node at (1,.8) {$\phantom{space!}$};
\node at (1.5,.7) {$\labelp$};
\node at (1.1,.7) {$\labelp$};
\node at (-.5,.7) {$\labelp$};
\node at (-0.1,.7) {$\labelp$};
\node at (1.05,.4) {$\labeln$};\node at (1.13,.4) {$\labelp$};
\node at (1.45,.4) {$\labeln$};\node at (1.53,.4) {$\labelp$};
\node at (-.55,.4) {$\labeln$};\node at (-.47,.4) {$\labelp$};
\node at (-.15,.4) {$\labeln$};\node at (-0.07,.4) {$\labelp$};
\node at (1.3,.4) {$\labeln$};
\node at (1.68,.4) {$\labeln$};
\node at (-.7,.4) {$\labeln$};
\node at (-.3,.4) {$\labeln$};
\node at (1.75,0) {$\labeln$};\node at (1.83,0) {$\labelp$};
\node at (.5,0) {$\labeln$};\node at (.6,0) {$\labelpm$};
\node at (2.1,0) {$\labelp$};
\node at (1.5,-.7) {$\labelp$};
\node at (1.1,-.7) {$\labelp$};
\node at (-.5,-.7) {$\labelp$};
\node at (-0.1,-.7) {$\labelp$};
\node at (1.05,-.4) {$\labeln$};\node at (1.13,-.4) {$\labelp$};
\node at (1.45,-.4) {$\labeln$};\node at (1.53,-.4) {$\labelp$};
\node at (-.55,-.4) {$\labeln$};\node at (-.47,-.4) {$\labelp$};
\node at (-.15,-.4) {$\labeln$};\node at (-0.07,-.4) {$\labelp$};
\node at (1.3,-.4) {$\labeln$};
\node at (1.68,-.4) {$\labeln$};
\node at (-.7,-.4) {$\labeln$};
\node at (-.3,-.4) {$\labeln$};
\node at (-.85,0) {$\labeln$};\node at (-.77,0) {$\labelp$};
\node at (-1.1,0) {$\labelp$};
	\draw[colp,->] (-.8,.8) to[out=-90, in=180] (-.5,.2) to[out = 0, in = -90] (-0.2,.8);
	\draw[colm,->] (-.4,.8) to[out=-90, in=180] (-.1,.2) to[out = 0, in = -90] (0.2,.8);
 \node at (.5,.7) {$\cdots$};
	\draw[colp,->] (.8,.8) to[out=-90, in=180] (1.1,.2) to[out = 0, in = -90] (1.4,.8);
	\draw[colm,->] (1.2,.8) to[out=-90, in=180]
 (1.5,.2) to[out = 0, in = -90] (1.8,.8);
	\draw[colp,<-] (-.8,-.8) to[out=90, in=180] (-.5,-.2) to[out = 0, in = 90] (-0.2,-.8);
	\draw[colm,<-] (-.4,-.8) to[out=90, in=180] (-.1,-.2) to[out = 0, in = 90] (0.2,-.8);
 \node at (.5,-.7) {$\cdots$};
	\draw[colp,<-] (.8,-.8) to[out=90, in=180] (1.1,-.2) to[out = 0, in = 90] (1.4,-.8);
	\draw[colm,<-] (1.2,-.8) to[out=90, in=180] (1.5,-.2) to[out = 0, in = 90] (1.8,-.8);
 \draw[colp,<-] (1.6,-.8) to (1.6,.8);
 \draw[colm,->] (-.6,-.8) to (-.6,.8);
 \draw[-] (-1,0) to[out=90,looseness=1,in=180] (.3,.6);
 \draw[dotted] (.3,.6) to (.7,.6);
 \draw[-] (.7,.6) to[out=0,in=90,looseness=1] (2,0);
\draw[-] (-1,0) to[out=-90,in=180,looseness=1] (.3,-.6);
\draw[dotted] (.3,-.6) to (.7,-.6);
\draw[->] (.7,-.6) to[out=0,in=-90,looseness=1] (2,0);
\end{tikzpicture}&\text{if $n$ is odd.}
\end{cases}
\end{align}
We reiterate that there are a total of 
$n$ $\down\up$-pairs and $(n-1)$ cups and caps at the top and bottom of these diagrams,
so $\v_n$ is an endomorphism of
$B^{\otimes n}$.
We have that $\v_n = \check\v_n \circ \hat\v_n$
where $\check\v_n$ denotes the top half of the above diagram
and $\hat\v_n$ is its bottom half.

\begin{lem}\label{play2}
For $0 \leq n \leq l$, both of the endomorphisms $\check\v_{n}$ and $\hat\v_n$ are of degree $-\binom{n}{2}$.
\end{lem}

\begin{proof}
It is clear that $\check\v_n$ and $\hat\v_n$ are of the same degree.
They have degree 0
if $n=0$ or $n=1$, and the result follows from \cref{play1} when $n=l$.
So now we assume that $2 \leq n \leq l-1$, and proceed to compute the degree of $\check\v_n$.
Using \cref{lengths} again, we see:
\begin{itemize}
\item Each of the $n-1$ special cups 
is of degree
$\ell(w_{\O\hat n+}) - \ell(w_{\O\hat n\pm})=\ell(w_{\O\hat n-}) - \ell(w_{\O\hat n\pm})= -\binom{n}{2}$.
\item
The $n$ crossings of special strings
are of degree $\ell(w_{\O\hat n\pm})
+\ell(w_{\O\hat n})- \ell(w_{\O\hat n+})
- \ell(w_{\O\hat n-})=\binom{n-1}{2}$;
\item The ordinary cap of color $n$
is of degree $\ell(w_{\O\hat n})-\ell(w_\O)=-(n-1)(l-n)$ if $n$ is odd
or $\ell(w_{\O\hat n+})-\ell(w_{\O+})=\ell(w_{\O\hat n-})-\ell(w_{\O-})= -n(l-n)$ if $n$ is even.
\item Each of the crossings involving the ordinary string and a special string has degree $\ell(w_{\O+}) + \ell(w_{\O\hat n})-\ell(w_{\O\hat n +}) - \ell(w_\O)=\ell(w_{\O-})+ \ell(w_{\O\hat n})-\ell(w_{\O\hat n -}) - \ell(w_\O) =l-n$,
and there are $(n-1)$ such if $n$ is odd or
$n$ such if $n$ is even.
\end{itemize}
The total degree coming from the last two points is 0. We are left with
$-(n-1)\binom{n}{2}+n\binom{n}{2} = 
-\binom{n}{2}$.
\end{proof}

We need two more 
families of endomorphisms.
For $n \geq 1$, let
\begin{align}\label{wn}
\w_n &:=
\begin{cases}
\begin{tikzpicture}[anchorbase,scale=2.5]
\node at (.7,.1) {$\labelp$};
\node at (-1.1,.1) {$\labelp$};
\node at (-.1,.07) {$\labelp$};
\node at (-.7,.07) {$\labelp$};
\node at (.3,.07) {$\labelp$};
\draw[->,colm] (-1,.25) to (-1,0);
\draw[->,colp] (-.8,0) to[looseness=3,out=90,in=90] (-.6,0);
\draw[->,colm] (-.2,0) to[looseness=3,out=90,in=90] (0,0);
\node at (-0.38,.06) {$\cdots$};
\draw[->,colp] (0.2,0) to[looseness=3,out=90,in=90] (.4,0);
\draw[<-,colm] (.6,.25) to (.6,0);
\node at (1,0.02) {$\phantom{space!}$};
\node at (1,.23) {$\phantom{space!}$};
\end{tikzpicture}
&\text{if $n$ is even}\\
\begin{tikzpicture}[anchorbase,scale=2.5]
\node at (.7,.1) {$\labelp$};
\node at (-1.1,.1) {$\labelp$};
\node at (-.1,.07) {$\labelp$};
\node at (-.7,.07) {$\labelp$};
\node at (.3,.07) {$\labelp$};
\draw[->,colp] (-1,.25) to (-1,0);
\draw[->,colm] (-.8,0) to[looseness=3,out=90,in=90] (-.6,0);
\draw[->,colm] (-.2,0) to[looseness=3,out=90,in=90] (0,0);
\node at (-0.38,.06) {$\cdots$};
\draw[->,colp] (0.2,0) to[looseness=3,out=90,in=90] (.4,0);
\draw[<-,colm] (.6,.25) to (.6,0);
\node at (1,.23) {$\phantom{spaceee!!}$};
\node at (1,0.02) {$\phantom{spaceee!!}$};
\end{tikzpicture}
&\text{if $n$ is odd,}
\end{cases}
\\\intertext{where there are $n-1$ caps in total. 
For $n \geq 0$, let}\label{rn}
\r_n &:= 
\begin{cases}
\begin{tikzpicture}[anchorbase,scale=2.4]
\node at (1.06,.88) {$\phantom{space!}$};
\node at (1.06,.52) {$\phantom{space!}$};
\node at (1.56,.7) {$\labelp$};
\draw[->,colm] (1.46,.5) to (1.46,.9);
\draw[<-,colp] (1.26,.5) to (1.26,.9);
\node at (1.16,.7) {$\labelp$};
\node [draw,rounded corners,inner sep=2pt,black,fill=gray!20!white] at (.93,0.7) {$\scriptstyle x_1$};
\draw[<-,colm] (.8,.5) to (.8,.9);
\node at (.7,.7) {$\labelp$};
\draw[->,colp] (1.06,.5) to (1.06,.9);
\node [draw,rounded corners,inner sep=2pt,black,fill=gray!20!white] at (.45,0.7) {$\scriptstyle x_1^2$};
\draw[->,colm] (.6,.5) to (.6,.9);
\draw[<-,colp] (.3,.5) to (.3,.9);
\node at (.2,.7) {$\labelp$};
\node at (-.04,.7) {$\cdots$};
\node at (-.3,.7) {$\labelp$};
\draw[->,colp] (-.4,.5) to (-.4,.9);
\node [draw,rounded corners,inner sep=2pt,black,fill=gray!20!white] at (-.6,0.7) {$\scriptstyle x_1^{n-1}$};
\draw[<-,colm] (-.8,.5) to (-.8,.9);
\node at (-.9,.7) {$\labelp$};
\end{tikzpicture}
&\text{if $n$ is even}\\
\begin{tikzpicture}[anchorbase,scale=2.4]
\node at (1.06,.88) {$\phantom{space!}$};
\node at (1.06,.52) {$\phantom{space!}$};
\node at (1.56,.7) {$\labelp$};
\draw[->,colm] (1.46,.5) to (1.46,.9);
\draw[<-,colp] (1.26,.5) to (1.26,.9);
\node at (1.16,.7) {$\labelp$};
\node [draw,rounded corners,inner sep=2pt,black,fill=gray!20!white] at (.93,0.7) {$\scriptstyle x_1$};
\draw[<-,colm] (.8,.5) to (.8,.9);
\node at (.7,.7) {$\labelp$};
\draw[->,colp] (1.06,.5) to (1.06,.9);
\node [draw,rounded corners,inner sep=2pt,black,fill=gray!20!white] at (.45,0.7) {$\scriptstyle x_1^2$};
\draw[->,colm] (.6,.5) to (.6,.9);
\draw[<-,colp] (.3,.5) to (.3,.9);
\node at (.2,.7) {$\labelp$};
\node at (-.04,.7) {$\cdots$};
\node at (-.3,.7) {$\labelp$};
\draw[->,colm] (-.4,.5) to (-.4,.9);
\node [draw,rounded corners,inner sep=2pt,black,fill=gray!20!white] at (-.6,0.7) {$\scriptstyle x_1^{n-1}$};
\draw[<-,colp] (-.8,.5) to (-.8,.9);
\node at (-.9,.7) {$\labelp$};
\end{tikzpicture}&\text{if $n$ is odd.}
\end{cases}
\end{align}
In particular, $\r_0 = \id_{\Sym}$.
Note that the endomorphism $\r_n$ is of degree $2\binom{n}{2}$.

\begin{lem}\label{first_impossible_lemma}
Suppose that $n=l$ or $n=l-1$.
Then 
$$
\sum_{b \in \B^{\O+}_{\O\pm}}
\tr^\O_{\O+}
\left(x_1^{n-1}\tr^\O_{\O-}\left(\cdots \left(x_1^3\tr^\O_{\O-}\left(
x_1^2 \tr^\O_{\O+}\left(x_1 \tr^\O_{\O-}\left(b\right)\right)\right)\cdots\right)\right)\right)
\otimes b^\vee =  (-1)^{\binom{n}{2}} 1 \otimes 1
$$
if $n$ is even, and
$$
\sum_{b \in \B^{\O+}_{\O\pm}}
\tr^\O_{\O-}
\left(x_1^{n-1}\tr^\O_{\O+}\left(\cdots \left(x_1^3\tr^\O_{\O-}\left(
x_1^2 \tr^\O_{\O+}\left(x_1 \tr^\O_{\O-}\left(b\right)\right)\right)\cdots\right)\right)\right)
\otimes b^\vee =  (-1)^{\binom{n}{2}} 1 \otimes 1
$$
if $n$ is odd.
\end{lem}

\begin{proof}
Let
$\partial_{\pm 1}^{(i)} := \partial_{i}\cdots\partial_2 \partial_{\pm 1}$, interpreted as $\partial_{\pm 1}$ if $i=1$
and as $1$ if $i=0$.
A reduced expression for the longest element of $\left(W_{\O\pm}/W_{\O+}\right)_{\min}$ is
$s_{(-1)^{l-1}}(s_2 s_{(-1)^{l-2}}) \cdots (s_{l-2}\cdots s_2 s_1)
(s_{l-1} \cdots s_2 s_{- 1})$.
Hence, recalling \cref{tracedef}, $
\tr^{\O+}_{\O\pm}=
\partial_{(-1)^{l-1}}^{(1)}
\partial_{(-1)^{l-2}}^{(2)}
\cdots \partial_1^{(l-2)} \partial_{-1}^{(l-1)}.$
Also $\tr^{\O}_{\O+} 
= \partial_1^{(l-1)}$
and $\tr^\O_{\O-} = \partial_{-1}^{(l-1)}$.

By degree considerations, all of the terms in the summation to be computed are 0 except for the one in which $b$ is of maximal degree.
We choose this top degree $b$ so that the corresponding $b^\vee$ equals 1.
Thus, $b$ is an element of $R^{\O+}$
such that \begin{equation}\label{nap}
\partial_{(-1)^{l-1}}^{(1)}
\cdots \partial_{1}^{(l-2)} \partial_{-1}^{(l-1)}(b) = 1.
\end{equation}
Now the result we are trying to prove is reduced to checking that
$$
\partial_{(-1)^n}^{(l-1)} 
\left(x_1^{n-1} \partial_{(-1)^{n-1}}^{(l-1)}\left(\cdots\left(x_1^2 \partial_1^{(l-1)} \left(x_1 \partial_{-1}^{(l-1)}\left(b\right)\right)\right) \cdots\right)\right)
= (-1)^{\binom{n}{2}}
$$
for $n=l-1$ or $n=l$.
In view of \cref{nap}, this follows from the $n=l-1$ or $n=l$ cases of the following more general statement:
for $1 \leq n \leq l$ we have that
\begin{equation}\label{hahaha}
\partial_{(-1)^n}^{(l-1)} 
\left(x_1^{n-1} \partial_{(-1)^{n-1}}^{(l-1)}\left(\cdots\left(x_1^2 \partial_1^{(l-1)} \left(x_1 \partial_{-1}^{(l-1)}\left(b\right)\right)\right) \cdots\right)\right)
= (-1)^{\binom{n}{2}}
\partial_{(-1)^n}^{(l-n)} 
\cdots \partial^{(l-2)}_1
\partial^{(l-1)}_{-1}(b).
\end{equation}
Now we proceed to prove \cref{hahaha} by induction on $n=1,\dots,l$.

The base case $n=1$ is vacuous---it is just asserting that $\partial^{(l-1)}_{-1}(b)=\partial^{(l-1)}_{-1}(b)$.
For the induction step, we assume that \cref{hahaha} is true for some $1 \leq n \leq l-1$, and must show that
$$
\partial_{(-1)^{n+1}}^{(l-1)}\left(x_1^n
\cdot (-1)^{\binom{n}{2}}
\partial_{(-1)^n}^{(l-n)} 
\cdots \partial^{(l-2)}_1
\partial^{(l-1)}_{-1}(b)\right)
=
(-1)^{\binom{n+1}{2}}
\partial_{(-1)^{n+1}}^{(l-n-1)} \partial_{(-1)^n}^{(l-n)} \cdots \partial^{(l-2)}_1
\partial^{(l-1)}_{-1}(b).
$$
Equivalently,
\begin{equation}
\partial_{l-1} \cdots \partial_2 \partial_{(-1)^{n+1}}
\left((-x_1)^n
\cdot 
\partial_{(-1)^n}^{(l-n)} 
\cdots \partial^{(l-2)}_1
\partial^{(l-1)}_{-1}(b)\right)
=
\partial_{(-1)^{n+1}}^{(l-n-1)} \partial_{(-1)^n}^{(l-n)} \cdots \partial^{(l-2)}_1
\partial^{(l-1)}_{-1}(b).\label{needy}
\end{equation}
To prove this, we will need the identity
\begin{align}
\partial_i\partial_{(-1)^m}^{(l-m)} \cdots \partial_1^{(l-2)} \partial_{-1}^{(l-1)}(b)=0,\label{needy1}
\end{align}
for $2 \leq i \leq l-1$ and $1 \leq m \leq l-1$. To prove \eqref{needy1}, we observe that
$$
s_i (s_{l-m} \cdots s_2 s_{(-1)^m})
\cdots (s_{l-2} \cdots s_2 s_1) 
(s_{l-1} \cdots s_2 s_{-1})
$$
is either not a reduced expression 
or it is braid-equivalent to $w s_j$ for some $j \in \{1,2,\dots,l-1\}$. This can be seen
by drawing the permutation diagram, or one can argue algebraically 
by induction on $m$.
In the former case, \cref{needy1} follows as the unreduced product of Demazures is 0.
In the latter case, \cref{needy1} follows because $\partial_j(b) = 0$.

To prove \cref{needy}, we use
the product rule
$\partial_i(f g)
= \partial_i(f) g + s_i(f) \partial_i(g)$ 
to apply $\partial_{(-1)^{n+1}},\partial_2,\dots,\partial_{l-1}$ in order to the product
$(-x_1)^n
\cdot 
\partial_{(-1)^n}^{(l-n)} 
\cdots \partial^{(l-2)}_1
\partial^{(l-1)}_{-1}(b)$, starting with
$\partial_{(-1)^{n+1}}$. In this way, we expand
the left hand side of \cref{needy}
as a sum of $2^{l-1}$ terms.
As soon as one of the Demazure operators
acts on the term $(-x_1)^n$ on the left, all remaining Demazure operators must also act on this term, else we get 0 by \cref{needy1}.
Thus, the only potentially non-zero terms are
$$
\partial_{l-1} \cdots \partial_{l-m}\left((-x_{l-m})^n\right) \cdot
\partial_{l-m-1} \cdots \partial_2 \partial_{(-1)^{n+1}}
\partial^{(l-n)}_{(-1)^n}
\cdots \partial_1^{(l-2)}\partial^{(l-1)}_{-1}(b).
$$
for $0 \leq m \leq l$.
If $m > n$ then $\partial_{l-1} \cdots \partial_{l-m}\left((-x_{l-m})^n\right) =0$ by degree considerations. So we may assume that $m \leq n$.
If $m < n=l-1$
the term is 0 since
$$
\partial_{l-m-1} \cdots \partial_2 \partial_{(-1)^l}
\partial^{(1)}_{(-1)^{l-1}}
\cdots \partial_1^{(l-2)}\partial^{(l-1)}_{-1}(b)
\stackrel{\cref{nap}}{=}
\partial_{l-m-1} \cdots \partial_2 \partial_{(-1)^l}(1) = 0.
$$
If $m < n < l-1$ the term is 0 since
$$
\partial_{l-m-1} \cdots \partial_2 \partial_{(-1)^{n+1}}
\partial^{(l-n)}_{(-1)^n}
\cdots \partial_1^{(l-2)}\partial^{(l-1)}_{-1}(b)
=
\partial_{l-m-1} \cdots \partial_{l-n} \partial_{(-1)^{n+1}}^{(l-n-1)}
\partial^{(l-n)}_{(-1)^n}
\cdots \partial_1^{(l-2)}\partial^{(l-1)}_{-1}(b),
$$
which is 0 by \cref{needy1} (taking $i$ there to be $l-n$).
Thus we are left just with the $m=n$ term
$$
\partial_{l-1} \cdots \partial_{l-n} \left((-x_{l-n-1})^n\right)\cdot
\partial^{(l-n-1)}
_{(-1)^{n+1}}
\partial^{(l-n)}
_{(-1)^n} \cdots \partial_1^{(l-2)} \partial_{-1}^{(l-1)}(b),
$$
which is equal to the right hand side
of \cref{needy} by a variant of \cref{missedflight}.
\end{proof}

\begin{cor}\label{liberty}
For $n=l-1$ or $n=l$, we have that
$$
\hat\u_n \circ \r_n \circ \check \u_n = 
\begin{cases}
\begin{tikzpicture}[anchorbase,scale=1.2]
\node at (.4,0) {$\labelp$};
\node at (0,0.) {$\labelpm$};
\node at (-.4,0) {$\labelp$};
\draw[->,colp] (-.2,-.2) to (-.2,.2);
\draw[->,colp] (.2,.2) to (.2,-.2);
\node at (1,.18) {$\phantom{space!}$};
\node at (1,-.18) {$\phantom{space!}$};
\end{tikzpicture}
&\text{if $n$ is even}\\
\begin{tikzpicture}[anchorbase,scale=1.2]
\node at (.4,0) {$\labelp$};
\node at (0,0) {$\labelpm$};
\node at (-.4,-.02) {$\labelp$};
\draw[->,colm] (-.2,-.2) to (-.2,.2);
\draw[->,colp] (.2,.2) to (.2,-.2);
\node at (1,.18) {$\phantom{space!}$};
\node at (1,-.18) {$\phantom{space!}$};
\end{tikzpicture}&\text{if $n$ is odd.}
\end{cases}
$$
\end{cor}

\begin{proof}
This follows from \cref{first_impossible_lemma} by a calculation which is similar to the proof of \cref{switchback}. We illustrate with an example, taking $n=5$. We recommend thinking also about an example in which $n$ is even, at which point the pattern is clear.
We need to show that
$$
\begin{tikzpicture}[anchorbase,scale=1.65]
\node at (.7,1) {$\labelp$};
\node at (-.6,1) {$\labelp$};
\node at (-1.83,1) {$\labelp$};
\node at (0.04,1) {$\labelp$};
\node at (-1.2,1) {$\labelp$};
\node at (-2.5,1) {$\labelp$};
\node at (-.9,0.4) {$\labelpm$};
\node at (-.9,1.6) {$\labelpm$};
\draw[->,colm] (-1.95,.2) to (-1.95,1.8);
\draw[->,colp] (.2,1.8) to (.2,0.2);
   \draw[colm,decoration={markings, mark=at position 0.0 with {\arrow{>}}},
    postaction={decorate}] (0,1) circle(0.5);
  \draw[colp,decoration={markings, mark=at position 0.0 with {\arrow{>}}},
    postaction={decorate}] (-.6,1) circle(0.5);  \draw[colm,decoration={markings, mark=at position 0.0 with {\arrow{>}}},
    postaction={decorate}] (-1.2,1) circle(0.5);  \draw[colp,decoration={markings, mark=at position 0.0 with {\arrow{>}}},
    postaction={decorate}] (-1.8,1) circle(0.5);
    \node [draw,rounded corners,inner sep=2pt,black,fill=gray!20!white] at (-.3,1) {$\scriptstyle x_1$};
       \node [draw,rounded corners,inner sep=2pt,black,fill=gray!20!white] at (-.92,1) {$\scriptstyle x_1^2$};
          \node [draw,rounded corners,inner sep=2pt,black,fill=gray!20!white] at (-1.51,1) {$\scriptstyle x_1^3$};
             \node [draw,rounded corners,inner sep=2pt,black,fill=gray!20!white] at (-2.13,1) {$\scriptstyle x_1^4$};
    \end{tikzpicture}=\begin{tikzpicture}[anchorbase,scale=1.2]
\node at (.4,0) {$\labelp$};
\node at (0,0) {$\labelpm$};
\node at (-.4,-.02) {$\labelp$};
\draw[->,colm] (-.2,-.2) to (-.2,.2);
\draw[->,colp] (.2,.2) to (.2,-.2);
\end{tikzpicture}.
$$
This is a situation in which we find it easiest to  expand the special undotted and dotted strings as ordinary and dashed strings using \cref{fire1,fire2}.
Remembering signs which arise because $\gamma(x_1)=-x_1$,
this reduces the problem to proving that
$$
\begin{tikzpicture}[anchorbase,scale=1.65]
\node at (.7,1) {$\labelp$};
\node at (-.6,1) {$\labelp$};
\node at (-1.83,1) {$\labelp$};
\node at (0.04,1) {$\labelm$};
\node at (-1.2,1) {$\labelm$};
\node at (-2.5,1) {$\labelm$};
\node at (-.9,0.4) {$\labelpm$};
\node at (-.9,1.6) {$\labelpm$};
\draw[->] (-1.95,.2) to (-1.95,1.8);
\draw[->] (.2,1.8) to (.2,0.2);
\draw[decoration={markings, mark=at position 0.0 with {\arrow{>}}}, postaction={decorate}] (0,1) circle(0.5);  
\draw[decoration={markings, mark=at position 0.0 with {\arrow{>}}},   postaction={decorate}] (-.6,1) circle(0.5);  
\draw[decoration={markings, mark=at position 0.0 with {\arrow{>}}},  postaction={decorate}] (-1.2,1) circle(0.5);
\draw[decoration={markings, mark=at position 0.0 with {\arrow{>}}}, postaction={decorate}] (-1.8,1) circle(0.5);
\node [draw,rounded corners,inner sep=2pt,black,fill=gray!20!white] at (-.3,1) {$\scriptstyle x_1$};
\node [draw,rounded corners,inner sep=2pt,black,fill=gray!20!white] at (-.92,1) {$\scriptstyle x_1^2$};
\node [draw,rounded corners,inner sep=2pt,black,fill=gray!20!white] at (-1.51,1) {$\scriptstyle x_1^3$};
\node [draw,rounded corners,inner sep=2pt,black,fill=gray!20!white] at (-2.13,1) {$\scriptstyle x_1^4$};
\end{tikzpicture}=(-1)^{\binom{n}{2}}\big(\begin{tikzpicture}[anchorbase,scale=1.2]
\node at (.4,0) {$\labelp$};
\node at (0,0) {$\labelpm$};
\node at (-.4,-.02) {$\labelm$};
\draw[->] (-.2,-.2) to (-.2,.2);
\draw[->] (.2,.2) to (.2,-.2);
\end{tikzpicture}\big)
$$
in the ordinary string calculus for $\BSBim$.
Now we simplify the left hand side of this working from right to left, starting with
an application of
hard RII \cref{hardRII}, to get
$$
\sum_{b \in \B^{\O+}_{\O\pm}}
\begin{tikzpicture}[anchorbase,scale=1.65]
\node at (1.3,1) {$\labelp$};
\node at (-.59,1) {$\labelp$};
\node at (-1.83,1) {$\labelp$};
\node at (0.1,.65) {$\labelm$};
\node at (-1.2,1) {$\labelm$};
\node at (-2.5,1) {$\labelm$};
\node at (-.9,0.4) {$\labelpm$};
\node at (-.9,1.6) {$\labelpm$};
\draw[->] (-1.95,.2) to (-1.95,1.8);
\draw[->] (.8,1.8) to (.8,0.2);
   \draw[decoration={markings, mark=at position 0.0 with {\arrow{>}}},
    postaction={decorate}] (0.12,1) ellipse(0.6 and 0.5);
  \draw[decoration={markings, mark=at position 0.0 with {\arrow{>}}},
    postaction={decorate}] (-.6,1) circle(0.5);  \draw[decoration={markings, mark=at position 0.0 with {\arrow{>}}},
    postaction={decorate}] (-1.2,1) circle(0.5);  \draw[decoration={markings, mark=at position 0.0 with {\arrow{>}}},
    postaction={decorate}] (-1.8,1) circle(0.5);
 \node [draw,rounded corners,inner sep=2pt,black,fill=gray!20!white] at (.3,1) {$\scriptstyle \tr^\O_{\O-}(b)$};
 \node [draw,rounded corners,inner sep=2pt,black,fill=gray!20!white] at (1,1) {$\scriptstyle b^\vee$};
    \node [draw,rounded corners,inner sep=2pt,black,fill=gray!20!white] at (-.3,1) {$\scriptstyle x_1$};
       \node [draw,rounded corners,inner sep=2pt,black,fill=gray!20!white] at (-.92,1) {$\scriptstyle x_1^2$};
          \node [draw,rounded corners,inner sep=2pt,black,fill=gray!20!white] at (-1.51,1) {$\scriptstyle x_1^3$};
             \node [draw,rounded corners,inner sep=2pt,black,fill=gray!20!white] at (-2.13,1) {$\scriptstyle x_1^4$};
    \end{tikzpicture}.
$$
Then we use \cref{bubbleslides} and easy RII \cref{easyRII} to move the rightmost bubble inside the second bubble from the right, and evaluating the internal bubble using \cref{evaluations} to reduce to
$$
\sum_{b \in \B^{\O+}_{\O\pm}}
\begin{tikzpicture}[anchorbase,scale=1.65]
\node at (1.5,1) {$\labelp$};
\node at (-1.83,1) {$\labelp$};
\node at (0.1,.65) {$\labelp$};
\node at (-1.2,1) {$\labelm$};
\node at (-2.5,1) {$\labelm$};
\node at (-.82,0.4) {$\labelpm$};
\node at (-.82,1.6) {$\labelpm$};
\draw[->] (-1.95,.2) to (-1.95,1.8);
\draw[->] (.97,1.8) to (.97,0.2);
  \draw[decoration={markings, mark=at position 0.0 with {\arrow{>}}},
    postaction={decorate}] (-.11,1) ellipse(1 and 0.5);
    \draw[decoration={markings, mark=at position 0.0 with {\arrow{>}}},
    postaction={decorate}] (-1.2,1) circle(0.5);  \draw[decoration={markings, mark=at position 0.0 with {\arrow{>}}},
    postaction={decorate}] (-1.8,1) circle(0.5);
 \node [draw,rounded corners,inner sep=2pt,black,fill=gray!20!white] at (0.1,1) {$\scriptstyle \tr^{\O}_{\O+}\left(x_1\tr^\O_{\O-}(b)\right)$};
 \node [draw,rounded corners,inner sep=2pt,black,fill=gray!20!white] at (1.2,1) {$\scriptstyle b^\vee$};
       \node [draw,rounded corners,inner sep=2pt,black,fill=gray!20!white] at (-.92,1) {$\scriptstyle x_1^2$};
          \node [draw,rounded corners,inner sep=2pt,black,fill=gray!20!white] at (-1.51,1) {$\scriptstyle x_1^3$};
             \node [draw,rounded corners,inner sep=2pt,black,fill=gray!20!white] at (-2.13,1) {$\scriptstyle x_1^4$};
    \end{tikzpicture}.
$$
Repeating in this way eventually yields
$$
\sum_{b \in \B^{\O+}_{\O\pm}}
\begin{tikzpicture}[anchorbase,scale=1.65]
\node at (1.5,1) {$\labelp$};
\node at (-3.05,1) {$\labelm$};
\node at (.78,1) {$\labelpm$};
\draw[->] (.65,.8) to (.65,1.2);
\draw[->] (.9,1.2) to (.9,0.8);
 \node [draw,rounded corners,inner sep=2pt,black,fill=gray!20!white] at (-1.2,1) {$\scriptstyle \tr^\O_{\O-}\left(x_1^4\tr^\O_{\O+}\left(x_1^3\tr^\O_{\O-}\left(x_1^2 \tr^{\O}_{\O+}\left(x_1\tr^\O_{\O-}(b)\right)\right)\right)\right)$};
 \node [draw,rounded corners,inner sep=2pt,black,fill=gray!20!white] at (1.2,1) {$\scriptstyle b^\vee$};
    \end{tikzpicture}.
$$
It remains to apply \cref{first_impossible_lemma} to obtain the result.
\end{proof}

\begin{cor}\label{mutual}
For $2 \leq n \leq l-1$, we have that
$$
\hat \v_n \circ \r_n \circ \check\u_n = 
\begin{cases}
\begin{tikzpicture}[anchorbase,scale=1.4]
\node at (.29,.3) {$\labeln$};
\node at (.43,.3) {$\labelp$};
\node at (-.47,.3) {$\labeln$};
\node at (-.33,.3) {$\labelp$};
\node at (.6,0) {$\labelp$};
\node at (0,-0.2) {$\labelpm$};
\node at (-.6,0) {$\labelp$};
\node at (-.09,.3) {$\labeln$};
\node at (.07,0.3) {$\labelpm$};
\draw[->,colp] (-.2,-.4) to (-.2,.5);
\draw[->,colp] (.2,.5) to (.2,-.4);
\draw[->] (-.6,.5) to[out=-90,in=-90,looseness=1.5] (.6,.5);
\node at (1,-.38) {$\phantom{space!}$};
\node at (1,.48) {$\phantom{space!}$};
\end{tikzpicture}&\text{if $n$ is even}\\
\begin{tikzpicture}[anchorbase,scale=1.4]
\node at (1,-.38) {$\phantom{space!}$};
\node at (1,.48) {$\phantom{space!}$};
\node at (.29,.3) {$\labeln$};
\node at (.43,.3) {$\labelp$};
\node at (-.47,.3) {$\labeln$};
\node at (-.33,.3) {$\labelp$};
\node at (.6,0) {$\labelp$};
\node at (0,-0.2) {$\labelpm$};
\node at (-.6,0) {$\labelp$};
\node at (-.09,.3) {$\labeln$};
\node at (.07,0.3) {$\labelpm$};
\draw[->,colm] (-.2,-.4) to (-.2,.5);
\draw[->,colp] (.2,.5) to (.2,-.4);
\draw[->] (-.6,.5) to[out=-90,in=-90,looseness=1.5] (.6,.5);
\end{tikzpicture}
&\text{if $n$ is odd.}
\end{cases}
$$\end{cor}

\begin{proof}
Again we explain the calculation with an example, taking $n = 4$. By the definitions, we have that
$$
\hat\v_4 \circ \r_4 \circ \check \u_4
= 
\begin{tikzpicture}[anchorbase,scale=1.6]
\node at (-.3,1.3) {$\labeln$};
\node at (-.9,1.3) {$\labeln$};
\node at (-1.38,1.38) {$\labeln$};
\node at (.2,1.38) {$\labeln$};
\node at (.7,1) {$\labelp$};
\node at (-.6,.75) {$\labelp$};
\node at (-1.83,1) {$\labelp$};
\node at (0.04,1) {$\labelp$};
\node at (-1.2,1) {$\labelp$};
\node at (-.6,0.35) {$\labelpm$};
\node at (-.65,1.7) {$\labeln$};
\node at (-.55,1.7) {$\labelpm$};
\node at (-.65,1.3) {$\labeln$};
\node at (-.54,1.3) {$\labelp$};
\node at (-1.21,1.35) {$\labeln$};
\node at (-1.11,1.35) {$\labelp$};
\node at (-.09,1.35) {$\labeln$};
\node at (0.01,1.35) {$\labelp$};
\node at (.3,1.7) {$\labeln$};
\node at (.4,1.7) {$\labelp$};
\node at (-1.65,1.7) {$\labeln$};
\node at (-1.55,1.7) {$\labelp$};
\draw[->,colp] (-1.3,.2) to (-1.3,1.8);
\draw[->,colp] (.14,1.8) to (.14,0.2);
   \draw[colm,decoration={markings, mark=at position 0.0 with {\arrow{>}}},
    postaction={decorate}] (0,1) circle(0.5);
  \draw[colp,decoration={markings, mark=at position 0.0 with {\arrow{>}}},
    postaction={decorate}] (-.6,1) circle(0.5);  \draw[colm,decoration={markings, mark=at position 0.0 with {\arrow{>}}},
    postaction={decorate}] (-1.2,1) circle(0.5);
    \node [draw,rounded corners,inner sep=2pt,black,fill=gray!20!white] at (-.3,1) {$\scriptstyle x_1$};
       \node [draw,rounded corners,inner sep=2pt,black,fill=gray!20!white] at (-.92,1) {$\scriptstyle x_1^2$};
          \node [draw,rounded corners,inner sep=2pt,black,fill=gray!20!white] at (-1.51,1) {$\scriptstyle x_1^3$};
           \draw[->] (-1.9,1.8) to[out=-90,looseness=.8,in=-90] (.7,1.8);
    \end{tikzpicture}.
$$
Now we use easy RII, easy RIII and \cref{bubbleslides} to slide the big rightward cap downward to obtain
$$
\begin{tikzpicture}[anchorbase,scale=1.6]
\node at (-.3,1.5) {$\labeln$};
\node at (-.9,1.5) {$\labeln$};
\node at (-1.45,1.5) {$\labeln$};
\node at (.25,1.5) {$\labeln$};
\node at (.63,.8) {$\labelp$};
\node at (-1.83,.8) {$\labelp$};
\node at (-.7,.45) {$\labeln$};
\node at (-.6,0.45) {$\labelpm$};
\node at (-.6,0.05) {$\labelpm$};
\node at (-.6,1.83) {$\labeln$};
\node at (-.5,1.83) {$\labelpm$};
\node at (-.65,1.5) {$\labeln$};
\node at (-.54,1.5) {$\labelp$};
\node at (-1.21,1.5) {$\labeln$};
\node at (-1.11,1.5) {$\labelp$};
\node at (-.09,1.5) {$\labeln$};
\node at (0.01,1.5) {$\labelp$};
\node at (.39,1.8) {$\labeln$};
\node at (.49,1.8) {$\labelp$};
\node at (-1.65,1.8) {$\labeln$};
\node at (-1.55,1.8) {$\labelp$};
\draw[->,colp] (-1.3,0) to (-1.3,1.9);
\draw[->,colp] (.14,1.9) to (.14,0);
   \draw[colm,decoration={markings, mark=at position 0.0 with {\arrow{>}}},
    postaction={decorate}] (0,1.2) circle(0.5);
  \draw[colp,decoration={markings, mark=at position 0.0 with {\arrow{>}}},
    postaction={decorate}] (-.6,1.2) circle(0.5);  \draw[colm,decoration={markings, mark=at position 0.0 with {\arrow{>}}},
    postaction={decorate}] (-1.2,1.2) circle(0.5);
    \node [draw,rounded corners,inner sep=2pt,black,fill=gray!20!white] at (-.3,1.2) {$\scriptstyle x_1$};
       \node [draw,rounded corners,inner sep=2pt,black,fill=gray!20!white] at (-.92,1.2) {$\scriptstyle x_1^2$};
          \node [draw,rounded corners,inner sep=2pt,black,fill=gray!20!white] at (-1.51,1.2) {$\scriptstyle x_1^3$};
           \draw[->] (-1.95,1.9) to[out=-90,looseness=2.15,in=-90] (.75,1.9);
    \end{tikzpicture}.
$$
The circles can now be removed by applying \cref{liberty}---the presence of the $\hat n$ means that we are effectively working inside a D$_n$ subsystem, so the hypothesis of \cref{liberty} is valid.
\end{proof}

\begin{lem}\label{xmasday}
For $2 \leq n \leq l-1$,
we have that
$$
\sum_{b \in \B^\O_{\O+}}
b^\vee\tr^{\O\hat n}_{\O+}
\left(
b \left(\eta^{\O\hat n}_{\O-}\right)^{\lceil\frac{n-1}{2}\rceil}
\left(\eta^{\O\hat n}_{\O+}\right)^{\lfloor\frac{n-1}{2}\rfloor}\right)
= 1.
$$
\end{lem}

\begin{proof}
A reduced expression for the 
longest element of $\left(W_{\O+} / W_{\O\hat n}\right)_{\min}$ is
$$
(s_{l-1} \cdots s_2 s_1)
(s_{l-n+1} \cdots s_{l-1}) \cdots (s_3 \cdots s_{n+1}) (s_2 \cdots s_n).
$$
The corresponding product of Demazure operators computes $\tr^{\O\hat n}_{\O+}$.
We also have that $\eta^{\O\hat n}_{\O+} = (x_{n+1}-x_1)(x_{n+2}-x_1)\cdots (x_l-x_1)$
and
$\eta^{\O\hat n}_{\O-} = (x_{n+1}+x_1)(x_{n+2}+x_1)\cdots (x_{l}+x_1)$.

By degree considerations, 
all terms in the sum vanish except for
the one in which $b$ is of maximal degree.
For this, we take $b = (-x_1)^{l-1}$ and $b^\vee = 1$. So the expression we are trying to calculate is
$$
(\partial_{l-1} \cdots \partial_1)
(\partial_{l-n+1} \cdots \partial_{l-1}) \cdots (\partial_3 \cdots \partial_{n+1}) (\partial_2 \cdots \partial_n)\!
\left(\!\!(-x_1)^{l-1}\!\!\!
\prod_{i=n+1}^{l} \!\left[(x_{i}+x_1)^{\lceil\frac{n-1}{2}\rceil}
(x_i-x_1)^{\lfloor\frac{n-1}{2}\rfloor}\right]\!\right)\!.
$$
First, note that $(-x_1)^{l-1}$ is invariant under $s_i$ for $i=2, \dots, l-1$, and therefore commutes past the product of Demazure operators $(\partial_{l-n+1} \cdots\partial_{l-1}) \cdots (\partial_3 \cdots \partial_{n+1}) (\partial_2 \cdots \partial_n)$. Now we observe that $\partial_2\cdots\partial_n\left(
(x_{n+1}+x_1)^{\lceil\frac{n-1}{2}\rceil}
(x_{n+1}+x_1)^{\lfloor\frac{n-1}{2}\rfloor}
\right)=1$.
All that matters for this is that the expression inside the brackets is a monic polynomial of degree $(n-1)$ in $x_{n+1}$ with coefficients in $\kk[x_1]$,
and $\partial_2\cdots\partial_n(x_{n+1}^{n-1}) = 1$.
Similarly, $\partial_3\cdots\partial_{n+2}
\left(
(x_{n+2}+x_1)^{\lceil\frac{n-1}{2}\rceil}
(x_{n+2}+x_1)^{\lfloor\frac{n-1}{2}\rfloor}\right)
= 1$, and so on. In this way, we reduce until we are left with 
$\partial_{l-1} \cdots\partial_1\left((-x_1)^{l-1}\right)$.
This is 1 by \cref{missedflight}.
\end{proof}

\begin{cor}\label{xmaseve1}
For $2 \leq n \leq l-1$, we have that
$$
\w_n \circ \v_n \circ \r_n \circ \check\u_n = 
\begin{cases}
\:\:\:\:\:\;\begin{tikzpicture}[anchorbase,scale=1.3]
\node at (0,0) {$\labelp$};
\node at (0,-0.3) {$\labelpm$};
\draw[->,colp] (-.2,-.45) to[out=90,in=90,looseness=2.8] (.2,-.45);
\draw[->,colm] (-.2,.45) to[out=-90,in=-90,looseness=2.8] (.2,.45);
\node at (1,.43) {$\phantom{space!}$};
\node at (1,-.43) {$\phantom{space!}$};
\end{tikzpicture}&\text{if $n$ is even}\\
\begin{tikzpicture}[anchorbase,scale=1.3]
\node at (.4,-.1) {$\labelp$};
\node at (0,-0.3) {$\labelpm$};
\node at (-.4,-.1) {$\labelp$};
\draw[->,colm] (-.3,-.36) to (.3,.2);
\draw[->,colp] (-.3,.2) to (.3,-.36);
\node at (1,.18) {$\phantom{space!}$};
\node at (1,-.34) {$\phantom{space!}$};
\end{tikzpicture}
&\text{if $n$ is odd.}
\end{cases}
$$
\end{cor}

\begin{proof}
We first write $\v_n=\check \v_n\circ \hat\v_n$ and 
use \cref{mutual} to simplify the bottom $\hat \v_n \circ \r_n \circ \check\u_n$
part of 
$\w_n\circ \check \v_n \circ \hat \v_n \circ \r_n \circ \check\u_n$. For example, if $n=6$, it shows that
$$
\w_6 \circ \v_6 \circ \r_6 \circ \check \u_6
=
\begin{tikzpicture}[anchorbase,scale=2.6]
\draw[-,decoration={markings, mark=at position 0.0 with {\arrow{>}}},
    postaction={decorate}] (-1.1,1.15) ellipse (1.75 and .3);
    \draw[->,colp] (-2.6,.7) to (-2.6,1) to [out=90,in=180] (-2.3,1.5) to [out=0,in=180] (-1.7,1) to [out=0,in=180] (-1.1,1.5) to [out=0,in=180] (-.5,1) to [out=0,in=180] (.1,1.5) to [out=0,in=90] (.4,1) to (.4,.7);
    \draw[->,colm] (-2.8,1.6) to [out=-90,in=180] (-2.3,1) to [out=0,in=180] (-1.7,1.5) to [out=0,in=180] (-1.1,1) to [out=0,in=180] (-.5,1.5) to [out=0,in=180] (.1,1) to [out=0,in=-90] (.6,1.6);
\node at (.46,1.1) {$\labeln$};\node at (.53,1.1) {$\labelp$};
\node at (-2.73,1.1) {$\labeln$};\node at (-2.67,1.1) {$\labelp$};
\node at (.77,1) {$\labelp$};
\node at (-2.97,1) {$\labelp$};
\node at (-1.1,0.77) {$\labelpm$};
\node at (-1.13,.92) {$\labeln$};
\node at (-1.05,.92) {$\labelpm$};
\node at (.07,1.25) {$\labeln$};\node at (.15,1.25) {$\labelm$};
\node at (-.53,1.25) {$\labeln$};\node at (-.45,1.25) {$\labelp$};
\node at (-1.13,1.25) {$\labeln$};\node at (-1.05,1.25) {$\labelm$};
\node at (-1.73,1.25) {$\labeln$};\node at (-1.65,1.25) {$\labelp$};
\node at (-2.33,1.25) {$\labeln$};\node at (-2.25,1.25) {$\labelm$};
\node at (-2.63,1.23) {$\labeln$};
\node at (.42,1.23) {$\labeln$};
\node at (-.2,1.35) {$\labeln$};
\node at (-2,1.35) {$\labeln$};
\node at (-1.4,1.4) {$\labeln$};
\node at (-.8,1.4) {$\labeln$};
\end{tikzpicture}
$$
Now it seems simplest to rewrite the special strings as ordinary strings using \cref{fire1,fire2}, sliding the dashed string to the top, to reduce the proof to showing using the ordinary string calculus that
$$
\begin{tikzpicture}[anchorbase,scale=2.6]
\draw[-,decoration={markings, mark=at position 0.0 with {\arrow{>}}},
    postaction={decorate}] (-1.1,1.15) ellipse (1.75 and .3);
    \draw[->] (-2.6,.7) to (-2.6,1) to [out=90,in=180] (-2.3,1.5) to [out=0,in=180] (-1.7,1) to [out=0,in=180] (-1.1,1.5) to [out=0,in=180] (-.5,1) to [out=0,in=180] (.1,1.5) to [out=0,in=90] (.4,1) to (.4,.7);
    \draw[->] (-2.8,1.6) to [out=-90,in=180] (-2.3,1) to [out=0,in=180] (-1.7,1.5) to [out=0,in=180] (-1.1,1) to [out=0,in=180] (-.5,1.5) to [out=0,in=180] (.1,1) to [out=0,in=-90] (.6,1.6);
\node at (.46,1.1) {$\labeln$};\node at (.53,1.1) {$\labelp$};
\node at (-2.73,1.1) {$\labeln$};\node at (-2.67,1.1) {$\labelp$};
\node at (.77,1) {$\labelp$};
\node at (-2.97,1) {$\labelp$};
\node at (-1.1,0.77) {$\labelpm$};
\node at (-1.13,.92) {$\labeln$};
\node at (-1.05,.92) {$\labelpm$};
\node at (.07,1.25) {$\labeln$};\node at (.15,1.25) {$\labelm$};
\node at (-.53,1.25) {$\labeln$};\node at (-.45,1.25) {$\labelp$};
\node at (-1.13,1.25) {$\labeln$};\node at (-1.05,1.25) {$\labelm$};
\node at (-1.73,1.25) {$\labeln$};\node at (-1.65,1.25) {$\labelp$};
\node at (-2.33,1.25) {$\labeln$};\node at (-2.25,1.25) {$\labelm$};
\node at (-2.63,1.23) {$\labeln$};
\node at (.42,1.23) {$\labeln$};
\node at (-.2,1.35) {$\labeln$};
\node at (-2,1.35) {$\labeln$};
\node at (-1.4,1.4) {$\labeln$};
\node at (-.8,1.4) {$\labeln$};
\end{tikzpicture}= 
\begin{tikzpicture}[anchorbase,scale=1.3]
\node at (0,0) {$\labelp$};
\node at (0,-0.3) {$\labelpm$};
\draw[->] (-.2,-.45) to[out=90,in=90,looseness=2.8] (.2,-.45);
\draw[->] (-.2,.45) to[out=-90,in=-90,looseness=2.8] (.2,.45);
\end{tikzpicture}.
$$
To prove this, 
we use hard RII \cref{hardRII}
to commute the five caps at the top past the big cap to obtain
$$
\begin{tikzpicture}[anchorbase,scale=2.6]
    \draw[-,decoration={markings, mark=at position 0.0 with {\arrow{>}}},
    postaction={decorate}] (-1.13,1.1) ellipse (1.73 and .55);
    \draw[->] (-2.6,.5) to (-2.6,.8) to [out=90,in=180] (-2.3,1.3) to [out=0,in=180] (-1.7,.8) to [out=0,in=180] (-1.1,1.3) to [out=0,in=180] (-.5,.8) to [out=0,in=180] (.1,1.3) to [out=0,in=90] (.4,.8) to (.4,.5);
    \draw[->] (-2.8,1.6) to [out=-90,in=180] (-2.3,.8) to [out=0,in=180] (-1.7,1.3) to [out=0,in=180] (-1.1,.8) to [out=0,in=180] (-.5,1.3) to [out=0,in=180] (.1,.88) to [out=0,in=-90] (.6,1.6);
\node at (.46,1) {$\labeln$};\node at (.52,1) {$\labelp$};
\node at (-2.76,1) {$\labeln$};\node at (-2.71,1) {$\labelp$};
\node at (.7,1) {$\labelp$};
\node at (-3,1) {$\labelp$};
\node at (-1.1,0.5) {$\labelpm$};
\node at (-1.13,.7) {$\labeln$};
\node at (-1.05,.7) {$\labelpm$};
\node at (.07,1.08) {$\labeln$};\node at (.15,1.07) {$\labelm$};
\node at (-.53,1.08) {$\labeln$};\node at (-.45,1.07) {$\labelp$};
\node at (-1.13,1.08) {$\labeln$};\node at (-1.05,1.07) {$\labelm$};
\node at (-1.73,1.08) {$\labeln$};\node at (-1.65,1.07) {$\labelp$};
\node at (-2.33,1.08) {$\labeln$};\node at (-2.25,1.07) {$\labelm$};
\node[draw,rounded corners,inner sep=2pt,black,fill=gray!20!white] at (-.195,1.38) {$\scriptstyle\eta^{\O\hat n}_{\O-}$};
\node[draw,rounded corners,inner sep=2pt,black,fill=gray!20!white] at (-.65,1.45) {$\scriptstyle\eta^{\O\hat n}_{\O+}$};
\node[draw,rounded corners,inner sep=2pt,black,fill=gray!20!white] at (-2,1.38) {$\scriptstyle\eta^{\O\hat n}_{\O-}$};
\node[draw,rounded corners,inner sep=2pt,black,fill=gray!20!white] at (-1.55,1.45) {$\scriptstyle\eta^{\O\hat n}_{\O+}$};
\node[draw,rounded corners,inner sep=2pt,black,fill=gray!20!white] at (-1.09,1.47) {$\scriptstyle\eta^{\O\hat n}_{\O-}$};
\node at (-2.61,1.1) {$\labeln$};
\end{tikzpicture}.
$$
The internal caps and cups can now be  separated and straightened using easy RII and \cref{zigzag} to give
$$
\begin{tikzpicture}[anchorbase,scale=2.5]
    \draw[-,decoration={markings, mark=at position 0.0 with {\arrow{>}}},
    postaction={decorate}] (-1.1,1) ellipse (.7 and .38);
    \node[draw,rounded corners,inner sep=1.5pt,black,fill=gray!20!white] at (-1.09,1.15) {$\left(\scriptstyle\eta^{\O\hat n}_{\O-}\right)^3\left(\scriptstyle\eta^{\O\hat n}_{\O+}\right)^2$};
    \draw[->] (-1.8,1.5) to [out=-90,in=-90,looseness=1.52] (-.4,1.5);
    \draw[->] (-1.8,0.5) to [out=90,in=90,looseness=.57] (-.4,0.5);
\node at (-.28,1) {$\labelp$};
\node at (-1.9,1) {$\labelp$};
\node at (-1.1,.95) {$\labeln$};
\node at (-1.11,.68) {$\labeln$};\node at (-1.04,.7) {$\labelp$};\node at (-1.04,.65) {$\labelm$};
\node at (-1.11,.8) {$\labeln$};\node at (-1.04,.8) {$\labelp$};
\node at (-1.07,.55) {$\labelpm$};
\end{tikzpicture}\stackrel{\cref{easyRII}}{=}
\begin{tikzpicture}[anchorbase,scale=2.4]
    \draw[-,decoration={markings, mark=at position 0.0 with {\arrow{>}}},
    postaction={decorate}] (-1.1,1.1) ellipse (.6 and .38);
    \node[draw,rounded corners,inner sep=1.5pt,black,fill=gray!20!white] at (-1.09,1.22) {$\left(\scriptstyle\eta^{\O\hat n}_{\O-}\right)^3\left(\scriptstyle\eta^{\O\hat n}_{\O+}\right)^2$};
    \draw[->] (-1.7,1.5) to [out=-90,in=-90,looseness=1.8] (-.5,1.5);
    \draw[->] (-1.6,0.55) to [out=90,in=90,looseness=.49] (-.6,0.55);
\node at (-.38,1) {$\labelp$};
\node at (-1.8,1) {$\labelp$};
\node at (-1.1,.95) {$\labeln$};
\node at (-1.11,.8) {$\labeln$};\node at (-1.04,.8) {$\labelp$};
\node at (-1.07,.58) {$\labelpm$};
\end{tikzpicture}.
    $$
    Now we apply hard RII one more time to slide the cup to the top. Then the argument concludes by applying \cref{evaluations} and \cref{xmasday}.
    All other cases with $n$ even proceed in the same way as this. We leave it to the reader to check that a similar argument works when $n$ is odd.
\end{proof}

\begin{cor}\label{mincies}
For $1 \leq n \leq l$, we have that
$\v_n\left(\r_n\left(1 \otimes\cdots\otimes 1\right)\right)
=
1 \otimes\cdots\otimes 1$
in $B^{\otimes n}$.
\end{cor}

\begin{proof}
Since $\v_n$ is of degree $-2\binom{n}{2}$ by \cref{play2} and $\r_n$ is of degree $2 \binom{n}{2}$,
the composition $\v_n \circ \r_n$ is of degree 0, so
we have that
$\v_n\left(\r_n\left(1 \otimes\cdots\otimes 1\right)\right)
=
 c 1 \otimes\cdots\otimes 1$ for some $c \in \kk$.
 Applying $\w_n$, which preserves the 1-tensor since it is built out of rightward caps,
 we deduce that
$\w_n\left(\v_n\left(\r_n\left(1\otimes\cdots\otimes 1\right)\right)\right)
 = c 1 \otimes 1$.
Since $\check\u_n$ is built out of rightward cups, we have that $\check\u_n(1 \otimes 1)=1 \otimes\cdots\otimes 1$,
so this is equivalent to the equation
$\w_n\left(\v_n\left(\r_n\left(\check\u_n\left(1\otimes 1\right)\right)\right)\right)
 = c 1 \otimes 1$.
The left hand side is $1 \otimes 1$ by \cref{xmaseve1},
so $c=1$.
\end{proof}

\begin{lem}\label{xmaseve2}
If $n=l-1$ or $n=l$ then 
$$
\w_n \circ \u_n \circ \r_n \circ \check\u_n = 
\begin{cases}
\:\:\:\:\begin{tikzpicture}[anchorbase,scale=1.4]
\node at (0,0) {$\labelp$};
\node at (0,-0.3) {$\labelpm$};
\draw[->,colp] (-.2,-.45) to[out=90,in=90,looseness=2.8] (.2,-.45);
\draw[->,colm] (-.2,.45) to[out=-90,in=-90,looseness=2.8] (.2,.45);
\node at (1,.43) {$\phantom{space!}$};
\node at (1,-.43) {$\phantom{space!}$};
\end{tikzpicture}&\text{if $n$ is even}\\
\begin{tikzpicture}[anchorbase,scale=1.4]
\node at (.25,-.1) {$\labelp$};
\node at (0,-0.3) {$\labelpm$};
\node at (-.25,-.1) {$\labelp$};
\draw[->,colm] (-.3,-.36) to (.3,.2);
\draw[->,colp] (-.3,.2) to (.3,-.36);
\node at (1,.18) {$\phantom{space!}$};
\node at (1,-.34) {$\phantom{space!}$};
\end{tikzpicture}
&\text{if $n$ is odd.}
\end{cases}
$$
\end{lem}

\begin{proof}
We first 
use \cref{liberty} to simplify the bottom $\hat \u_n \circ \r_n \circ \check\u_n$
part of 
$\w_n\circ \check \u_n \circ \hat \u_n \circ \r_n \circ \check\u_n$, then the result follows in an obvious way using easy RII.
If $n=2$ then the proof is already complete after the application of \cref{liberty}. For a typical example, suppose that $n=3$. Then
the application of \cref{liberty} reduces to showing that
$$
\begin{tikzpicture}[anchorbase,scale=2]
\node at (.65,0) {$\labelp$};
\node at (-.35,0) {$\labelp$};
\node at (0.15,-0.13) {$\labelpm$};
\node at (0,0.1) {$\labelp$};
\node at (0.3,0.1) {$\labelp$};
\draw[->,colm] (-.25,-.3) to[out=90,in=-120] (-.15,0.1) to [out=60,in=120,looseness=1.4] (.15,0.1) to [out=-60,in=-120,looseness=1.4] (.45,0.1) to[out=60,in=-90] (.5,.3);
\draw[->,colp] (-.2,.3) to[out=-90,in=120] (-.15,.1) to [out=-60,in=-120,looseness=1.4] (.15,.1) to [out=60,in=120,looseness=1.4] (.45,.1) to[out=-60,in=90] (.55,-.3);
\end{tikzpicture}
=\begin{tikzpicture}[baseline=-6.4,scale=2]
\node at (.22,-.05) {$\labelp$};
\node at (0,-0.26) {$\labelpm$};
\node at (-.22,-.05) {$\labelp$};
\draw[->,colm] (-.3,-.36) to (.3,.2);
\draw[->,colp] (-.3,.2) to (.3,-.36);
\end{tikzpicture}\:,
$$
which follows by easy RII.
\end{proof}

\begin{theo}\label{popping}
$\v_{l-1} = \u_{l-1}$.
\end{theo}

\begin{proof}
If $l=2$ (when there is no bubble to pop) this follows from the definitions and an application of distant RII \cref{distantRII}.
Now assume that $l \geq 3$.
\cref{xmaseve1,xmaseve2} imply that
\begin{equation} \label{jonlikessillynames}
\w_{l-1} \circ \v_{l-1} \circ \r_{l-1} \circ \check\u_{l-1}
=\w_{l-1} \circ \u_{l-1} \circ \r_{l-1} \circ \check\u_{l-1}.
\end{equation}
Since $\w_{l-1}$ and $\check\u_{l-1}$ are built from rightward cups and caps, they preserve the 1-tensor. So \cref{mincies} implies that the left-hand side of \cref{jonlikessillynames} takes $1 \otimes 1 \mapsto 1\otimes 1$. In particular, it is non-zero, hence, so are
$\v_{l-1}$ and $\u_{l-1}$.
The maps $\v_{l-1}$ and $\u_{l-1}$ are both of degree
$-2\binom{l-1}{2}$ by \cref{play1,play2}.
By \cref{handy3}, the
morphism space $\End_{R^{\O+}\dash R^{\O+}}(B^{\otimes(l-1)})$ in this (lowest) degree is one-dimensional.
Hence, $\u_{l-1}$ and $\v_{l-1}$ are equal up to multiplication by a scalar. 
They become equal on postcomposing with $\w_{l-1}$
and precomposing with $\r_{l-1} \circ \check\u_{l-1}$, so the scalar is 1.
\end{proof}

%% file: s5-cyclotomic.tex
\setcounter{section}{4}
\section{Main results}\label{seccyclotomic}

Now we prove the main results,
which explain the combinatorial coincidences observed in \cref{asleftmodule}.
Throughout, $l \geq 2$ is an integer with
$l \equiv t\pmod{2}$.
Nil-Brauer notation is as in \cref{secnilbrauer}, and
all other notation is for the
type D$_l$ root system as in \cref{secD}.
In particular, $\Sym = \kk[x_1,\dots,x_l]^{S_l}=R^{\O+}$.

\subsection{Construction of the monoidal functor \texorpdfstring{$\Theta$}{}}\label{ssec5-2}
Recall the subalgebra $\Lambda^{[2]}$ of the algebra of symmetric functions defined just before \cref{clanky}. 
As discussed in \cref{simonriche}, 
the image of $\Lambda^{[2]}$ 
under the evaluation homomorphism 
$\ev_l:\Lambda \rightarrow \Sym$
is the subalgebra 
generated by $e_r(x_1^2,\dots,x_l^2)\:(1 \leq r \leq l)$, that is,
the algebra of invariants for the Weyl group of type B$_l$. 
As in the introduction, we denote this subalgebra by $\Center$.
The action of $\Lambda^{[2]}$ factors through this quotient (via $\ev_l$) to make the cyclotomic nil-Brauer category $\CNB_l$ from \cref{letters}
into a strict $\Center$-linear graded monoidal category. 

The endomorphism category $\END_{\ESBim}(\Sym) = \one_{\O+} \ESBim \one_{\O+}$
is a graded Karoubian 
$\Center$-linear monoidal category with objects that are certain graded $(\Sym,\Sym)$-bimodules, and whose tensor product 
is $-\otimes_{\Sym}-$.
This will be our main focus from now on: we will show it categorifies the $\Z[q,q^{-1}]$-algebra $\V$ 
appearing in \cref{asleftmodule}.
We remark that $\END_{\ESBim}(\Sym)$
is generated as a graded Karoubian category by the special bimodule $B$ from \cref{words}. This follows from \cref{cruisier}.

\begin{theo}\label{okinawa}
There is a strict $\Center$-linear 
graded monoidal functor
$\Theta:\CNB_l \rightarrow \END_{\ESBim}(\Sym)$
mapping the generating object $B$ of 
$\CNB_l$ to the bimodule $B$, and defined
on generating 
morphisms by
\begin{align}\label{easier1}
\begin{tikzpicture}[anchorbase]
	\draw[semithick] (0.08,-.3) to (0.08,.4);
    \closeddot{0.08,.05};
\end{tikzpicture}
&\mapsto
\begin{tikzpicture}[anchorbase]
    \node at (-.6,0.05) {$\labelp$};
	\draw[colm,<-] (-0.35,-.3) to (-0.35,.4);
	\draw[colp,->] (0.25,-.3) to (0.25,.4);
\node [draw,rounded corners,inner sep=2pt,black,fill=gray!20!white] at (-.05,0.05) {$\scriptstyle x_1$};
    \node at (.5,0.05) {$\labelp$};
\end{tikzpicture}=-\left( \begin{tikzpicture}[anchorbase]
    \node at (-.6,0.05) {$\labelp$};
	\draw[colp,<-] (-0.35,-.3) to (-0.35,.4);
	\draw[colm,->] (0.25,-.3) to (0.25,.4);
\node [draw,rounded corners,inner sep=2pt,black,fill=gray!20!white] at (-.05,0.05) {$\scriptstyle x_1$};
    \node at (.5,0.05) {$\labelp$};
\end{tikzpicture}\right),
\\\label{easier2}
\begin{tikzpicture}[anchorbase]
	\draw[semithick] (0.4,0.4) to[out=-90, in=0] (0.1,0);
	\draw[semithick] (0.1,0) to[out = 180, in = -90] (-0.2,0.4);
\end{tikzpicture}
&\mapsto\begin{tikzpicture}[anchorbase]
\node at (.2,0.02) {$\phantom{space!}$};
\node at (.2,-.58) {$\phantom{space!}$};
	\draw[colm,->] (0.4,0) to[out=-90, in=0] (0.1,-0.4) to[out = 180, in = -90] (-0.2,0);
	\draw[colp,<-] (0.7,0) to[out=-90, in=0] (0.1,-0.65) to[out = 180, in = -90] (-0.5,0);
    \node at (.9,-0.25) {$\labelp$};
    \node at (.1,-0.15) {$\labelp$};
\end{tikzpicture}=\begin{tikzpicture}[anchorbase]
\node at (.2,0.02) {$\phantom{space!}$};
\node at (.2,-.58) {$\phantom{space!}$};
\draw[colp,->] (0.4,0) to[out=-90, in=0] (0.1,-0.4) to[out = 180, in = -90] (-0.2,0);
	\draw[colm,<-] (0.7,0) to[out=-90, in=0] (0.1,-0.65) to[out = 180, in = -90] (-0.5,0);
    \node at (.9,-0.25) {$\labelp$};
    \node at (.1,-0.15) {$\labelp$};
\end{tikzpicture},\\\label{easier3}
\begin{tikzpicture}[anchorbase]
	\draw[semithick] (0.4,0) to[out=90, in=0] (0.1,0.4);
	\draw[semithick] (0.1,0.4) to[out = 180, in = 90] (-0.2,0);
\end{tikzpicture}
&\mapsto
\begin{tikzpicture}[anchorbase]
	\draw[colm,<-] (0.4,0) to[out=90, in=0] (0.1,0.4) to[out = 180, in = 90] (-0.2,0);
	\draw[colp,->] (0.7,0) to[out=90, in=0] (0.1,0.65) to[out = 180, in = 90] (-0.5,0);
    \node at (.9,0.25) {$\labelp$};
    \node at (.1,0.15) {$\labelp$};
\end{tikzpicture}=
\begin{tikzpicture}[anchorbase]
	\draw[colp,<-] (0.4,0) to[out=90, in=0] (0.1,0.4) to[out = 180, in = 90] (-0.2,0);
	\draw[colm,->] (0.7,0) to[out=90, in=0] (0.1,0.65) to[out = 180, in = 90] (-0.5,0);
    \node at (.9,0.25) {$\labelp$};
    \node at (.1,0.15) {$\labelp$};
\end{tikzpicture}
,\\\label{easier4}
\begin{tikzpicture}[anchorbase]
\draw[semithick] (0.28,-.3) to (-0.28,.4);
	\draw[semithick] (-0.28,-.3) to (0.28,.4);
 \end{tikzpicture}
&\mapsto
\begin{tikzpicture}[anchorbase,scale=1.3]
\draw[colm,<-] (.3,-.8) to (.3,.8);
\draw[colm,->] (-.3,-.8) to (-.3,.8);
	\draw[colp,<-] (.8,.8) to[out=-90, in=0] (0,.2) to[out = 180, in = -90] (-0.8,.8);
	\draw[colp,->] (.8,-.8) to[out=90, in=0] (0,-.2) to[out = 180, in = 90] (-0.8,-.8);
\draw (0,0) ellipse (.85 and .6);
\draw[->] (.85,0)
arc[start angle=-0.5, end angle=.5,radius=.9];
\node at (1.05,0) {$\labelp$};
\node at (.68,0) {$\labelp$};\node at (.5,0.04) {$\labeltwo$};
\node at (-.5,0) {$\labelp$};\node at (-.68,0.04) {$\labeltwo$};
\node at (.09,0) {$\labelpm$};\node at (-.09,0) {$\labeltwo$}; 
\node at (-1,0) {$\labelp$};
\node at (0,-0.7) {$\labelp$};
\node at (0,0.7) {$\labelp$};
\node at (.38,0.4) {$\labeltwo$};
\node at (-.38,0.4) {$\labeltwo$};
\node at (.38,-0.4) {$\labeltwo$};
\node at (-.38,-0.4) {$\labeltwo$};
\node at (.09,0.4) {$\labelp$};\node at (-.09,0.4) {$\labeltwo$};
\node at (.09,-0.4) {$\labelp$};\node at (-.09,-0.4) {$\labeltwo$};
\end{tikzpicture}=\begin{tikzpicture}[anchorbase,scale=1.3]
\draw[colp,<-] (.3,-.8) to (.3,.8);
\draw[colp,->] (-.3,-.8) to (-.3,.8);
	\draw[colm,<-] (.8,.8) to[out=-90, in=0] (0,.2) to[out = 180, in = -90] (-0.8,.8);
	\draw[colm,->] (.8,-.8) to[out=90, in=0] (0,-.2) to[out = 180, in = 90] (-0.8,-.8);
\draw (0,0) ellipse (.85 and .6);
\draw[->] (.85,0)
arc[start angle=-0.5, end angle=.5,radius=.9];
\node at (1.05,0) {$\labelp$};
\node at (.68,0) {$\labelp$};\node at (.5,0.04) {$\labeltwo$};
\node at (-.5,0) {$\labelp$};\node at (-.68,0.04) {$\labeltwo$};
\node at (.09,0) {$\labelpm$};\node at (-.09,0) {$\labeltwo$};    
\node at (-1,0) {$\labelp$};
\node at (0,-0.7) {$\labelp$};
\node at (0,0.7) {$\labelp$};
\node at (.38,0.4) {$\labeltwo$};
\node at (-.38,0.4) {$\labeltwo$};
\node at (.38,-0.4) {$\labeltwo$};
\node at (-.38,-0.4) {$\labeltwo$};
\node at (.09,0.4) {$\labelp$};\node at (-.09,0.4) {$\labeltwo$};
\node at (.09,-0.4) {$\labelp$};\node at (-.09,-0.4) {$\labeltwo$};
\end{tikzpicture}.
\end{align}
In the case $l=2$,
the circle of color 2 and all labels $\hat 2$ here should simply be omitted.
Furthermore:
\begin{enumerate}
\item
For any symmetric polynomial $f \in \Sym$,
$\Theta\left(\:\Bubble{f}\:\right) = \bubble{f}\begin{tikzpicture}[anchorbase]\node at (0,0) {$\labelp$};\end{tikzpicture}$, the endomorphism of $\Sym$ defined by multiplication by $f$.
\item
The functor $\Theta$ intertwines
the symmetry $\tR$ of $\CNB_l$
arising from \cref{R} with the symmetry $\tR$ of $\END_{\ESBim}(\Sym)$ arising from \cref{tRnew}.
\end{enumerate}
\end{theo}

\begin{proof}
We first justify the equalities in the formulae for the images of the generators under $\Theta$ in the statement of the theorem.
This follows in each case by expanding the definitions of the undotted and dotted special strings like in \cref{adele}.
For the dot, it is a special case of \cref{needsabitofthought}
since $\gamma(x_1)=-x_1$. The other three are obvious once the definitions have been expanded, using the easy properties of the dashed string; for the crossing, we already used a similar strategy in the proof of \cref{puppies}.
These equalities show that
$\tR \circ \Theta = \Theta \circ \tR$ on all four of the generating morphisms of $\CNB_l$. Once the well-definedness of $\Theta$ has been established, this is all that is needed to prove the property (2).

Next we check that the degrees are consistent.
The dot has degree 2, which is the same as the degree of $x_1$.
For the cap, the clockwise cap is of degree
$\ell(w_{\O-})-\ell(w_\O) = \ell(w_{\O+})-\ell(w_\O)$
and the counterclockwise cap is of degree $\ell(w_\O)-\ell(w_{\O+})=\ell(w_\O)-\ell(w_{\O-})$, giving the required total degree 0.
The cup is similar. For the crossing,
note from the definition 
that 
$\Theta\left(\ \begin{tikzpicture}[anchorbase]
\draw[semithick] (0.28,-.3) to (-0.28,.4);
	\draw[semithick] (-0.28,-.3) to (0.28,.4);
\end{tikzpicture}\ \right)$ is $\v_2$ from \cref{bplus},
which is of the desired degree $-2$ thanks to
\cref{play2}.

To prove the existence of $\Theta$,
it suffices to check the defining
relations for $\cNB_t$, that is, the eight
relations
\cref{rels2,rels3,rels4,rels1} plus the cyclotomic relations that are the generators of the ideal $\cI_l$ from \cref{bensdumbideatocallitfoursided}.
We are first going to 
check \cref{rels2,rels3,rels4,rels1}.
This establishes the existence of a 
graded $\Lambda^{[2]}$-linear monoidal functor
$\widehat{\Theta}:\cNB_t\rightarrow \END_{\ESBim}(\Sym)$,
viewing $\END_{\ESBim}(\Sym)$ as a $\Lambda^{[2]}$-linear monoidal category via the evaluation homomorphism $\Lambda^{[2]}\twoheadrightarrow \Center$.
After that, we will discuss property (1), leaving the cyclotomic relation to the end.
\begin{itemize}
\item
The zig-zag relations from \cref{rels2}
follow easily using \cref{zigzag} and the zig-zag relations for the dashed string.
\item
The bubble relation from \cref{rels2} follows
by \cref{coffeeisgood}, taking $n=0$.
\item
The first relation from \cref{rels4} follows by the Interchange Law. The sign appears due to \cref{needsabitofthought} 
since $\gamma(x_1) = -x_1$.
\item
The first relation from \cref{rels3} follows
by cyclicity of the diagrammatic calculus for $\EBSBim$:
$$
\begin{tikzpicture}[anchorbase,scale=1.2]
\draw[colm,<-] (.3,-.8) to (.3,.8)to[out=90,in=90,looseness=1.5] (1.6,.8) to (1.6,-.8);
\draw[colm,->] (-.3,-.8) to (-.3,1.4);
	\draw[colp,<-] (1.25,-.8) to (1.25,.8) to[out=90,in=90,looseness=2] (.8,.8) to[out=-90, in=0] (0,.2) to[out = 180, in = -90] (-0.8,.8) to (-.8,1.4);
	\draw[colp,->] (.8,-.8) to[out=90, in=0] (0,-.2) to[out = 180, in = 90] (-0.8,-.8);
\draw (0,0) ellipse (.85 and .6);
\draw[->] (.85,0)
arc[start angle=-0.5, end angle=.5,radius=.9];
\node at (1.05,0) {$\labelp$};
\node at (.68,0) {$\labelp$};\node at (.5,0.04) {$\labeltwo$};
\node at (-.5,0) {$\labelp$};\node at (-.68,0.04) {$\labeltwo$};
\node at (0.09,0) {$\labelpm$};
\node at (-.09,0) {$\labeltwo$};
    \node at (-1,0) {$\labelp$};
\node at (0,-0.7) {$\labelp$};
\node at (1.8,0) {$\labelp$};
\node at (.38,0.4) {$\labeltwo$};
\node at (-.38,0.4) {$\labeltwo$};
\node at (.38,-0.4) {$\labeltwo$};
\node at (-.38,-0.4) {$\labeltwo$};
\node at (.09,0.4) {$\labelp$};\node at (-.09,0.4) {$\labeltwo$};
\node at (.09,-0.4) {$\labelp$};\node at (-.09,-0.4) {$\labeltwo$};
\end{tikzpicture}
=
\begin{tikzpicture}[anchorbase,scale=1.2]
\draw (0,0) ellipse (.85 and .6);
\draw[->] (.85,0)
arc[start angle=-0.5, end angle=.5,radius=.9];
\draw[colp,->] (-.4,1) to[looseness=1.5,out=-70,in=90] (.6,-1.2);
\draw[colm,<-] (.4,1) to[looseness=1.5,out=-110,in=90] (-.6,-1.2);
\draw[colp,->] (.2,-1.2) to (.2,-.9) to[out=90,in=90,looseness=2] (-1,-.7) to (-1,-1.2);
\draw[colm,<-] (-.2,-1.2) to (-.2,-.9) to[out=90,in=90,looseness=2] (1,-.7) to (1,-1.2);
\node at (0.6,-.25) {$\labeltwo$};
\node at (-0.6,-.25) {$\labeltwo$};
\node at (-1,0) {$\labelp$};
\node at (1.05,0) {$\labelp$};
\node at (.08,-0.05) {$\labelpm$};
\node at (-.08,-0.05) {$\labeltwo$};
\node at (0,0.47) {$\labeltwo$};
\node at (0,-0.47) {$\labeltwo$};
\node at (.37,-.42) {$\labelp$};\node at (0.22,-0.38) {$\labeltwo$};
\node at (-.23,-.42) {$\labelp$};\node at (-0.38,-0.38) {$\labeltwo$};
\node at (.5,0.24) {$\labelp$};\node at (.3,0.26) {$\labeltwo$};
\node at (-.3,0.24) {$\labelp$};\node at (-.5,0.26) {$\labeltwo$};
\node at (.38,-0.9) {$\labelp$};
\node at (-.38,-0.9) {$\labelp$};
\end{tikzpicture}
=
\begin{tikzpicture}[anchorbase,scale=1.2]
\draw[colp,<-] (.3,-.8) to (.3,1.4);
\draw[colp,->] (-.3,-.8) to (-.3,.8)
to[out=90,in=90,looseness=1.5] (-1.6,.8) to (-1.6,-.8);
	\draw[colm,<-] (.8,1.4) to (.8,.8) to[out=-90, in=0] (0,.2) to[out = 180, in = -90] (-0.8,.8) to [out=90,in=90,looseness=2] (-1.25,.8) to (-1.25,-.8);
	\draw[colm,->] (.8,-.8) to[out=90, in=0] (0,-.2) to[out = 180, in = 90] (-0.8,-.8);
\draw (0,0) ellipse (.85 and .6);
\draw[->] (.85,0)
arc[start angle=-0.5, end angle=.5,radius=.9];
\node at (1.05,0) {$\labelp$};
\node at (.68,0) {$\labelp$};\node at (.5,0.04) {$\labeltwo$};
\node at (-.5,0) {$\labelp$};\node at (-.68,0.04) {$\labeltwo$};
\node at (0.09,0) {$\labelpm$};
\node at (-.09,0) {$\labeltwo$};
    \node at (-1,0) {$\labelp$};
\node at (0,-0.7) {$\labelp$};
\node at (-1.8,0) {$\labelp$};
\node at (.38,0.4) {$\labeltwo$};
\node at (-.38,0.4) {$\labeltwo$};
\node at (.38,-0.4) {$\labeltwo$};
\node at (-.38,-0.4) {$\labeltwo$};
\node at (.09,0.4) {$\labelp$};\node at (-.09,0.4) {$\labeltwo$};
\node at (.09,-0.4) {$\labelp$};\node at (-.09,-0.4) {$\labeltwo$};
\end{tikzpicture}.
$$
In view of the alternative forms for the cap and the crossing in \cref{easier3,easier4}, which was already established in the opening paragraph of the proof, this is what we want. 
\item
The dot sliding relation from \cref{rels4}
follows because
\begin{align*}
\Theta\Bigg(
\begin{tikzpicture}[anchorbase,scale=1.4]
 	\draw[semithick] (0,-.3) to (0,.3);
	\draw[semithick] (-.3,-.3) to (-0.3,.3);
\end{tikzpicture}
\:-\:
&\begin{tikzpicture}[anchorbase,scale=1.4]
 	\draw[semithick] (-0.15,-.3) to[out=90,in=180] (0,-.1) to[out=0,in=90] (0.15,-.3);
 	\draw[semithick] (-0.15,.3) to[out=-90,in=180] (0,.1) to[out=0,in=-90] (0.15,.3);
\end{tikzpicture}\:\:\Bigg)=
\begin{tikzpicture}[anchorbase,scale=1.2]
\draw[->,colm] (.9,-.6) to (.9,.6);
\draw[<-,colm] (-.9,-.6) to (-.9,.6);
\draw[<-,colp] (.3,-.6) to (.3,.6);
\draw[->,colp] (-.3,-.6) to (-.3,.6);
\node at (0,0) {$\labelp$};
\node at (-1.1,0) {$\labelp$};
\node at (1.1,0) {$\labelp$};
\end{tikzpicture}\:-\!\begin{tikzpicture}[anchorbase,scale=1.2]
\draw[->,colm] (.7,-.9) to[looseness=1.8,out=90,in=90] (-.7,-.9);
\draw[->,colm] (-.7,.9) to[looseness=1.7,out=-90,in=-90] (.7,.9);
\draw[<-,colp] (.3,-.9) to[out=90,in=90,looseness=1.8] (-.3,-.9);
\draw[<-,colp] (-.3,.9) to[looseness=1.8,out=-90,in=-90] (.3,.9);
\node at (0,.8) {$\labelp$};
\node at (0,-.8) {$\labelp$};
\node at (0,0) {$\labelp$};
\end{tikzpicture}\\
&\stackrel{(*)}{=}
\begin{tikzpicture}[anchorbase,scale=1.4]
\draw[colp,<-] (.15,-.8) to (.15,.8);
\draw[colp,->] (-.15,-.8) to (-.15,.8);
	\draw[colm,<-] (.8,.8) to[out=-90, in=0,looseness=1.4] (0,.15) to[out = 180, in = -90,looseness=1.4] (-0.8,.8);
	\draw[colm,->] (.8,-.8) to[looseness=1.4,out=90, in=0] (0,-.15) to[looseness=1.4,out = 180, in = 90] (-0.8,-.8);
\draw (0,0) ellipse (.85 and .6);
\draw[->] (.85,0)
arc[start angle=-0.5, end angle=.5,radius=.9];
\node at (1.05,0) {$\labelp$};
\node at (.68,0) {$\labelp$};\node at (.5,0.04) {$\labeltwo$};
\node at (-.5,0) {$\labelp$};\node at (-.68,0.04) {$\labeltwo$};
\node at (0.06,0) {$\labelpm$};
\node at (-.07,0) {$\labeltwo$};
    \node at (-1,0) {$\labelp$};
\node at (0,-0.7) {$\labelp$};
\node at (0,0.7) {$\labelp$};
\node at (.38,0.4) {$\labeltwo$};
\node at (-.38,0.4) {$\labeltwo$};
\node at (.38,-0.4) {$\labeltwo$};
\node at (-.38,-0.4) {$\labeltwo$};
\node at (.06,0.4) {$\labelp$};\node at (-.07,0.4) {$\labeltwo$};
\node at (.06,-0.4) {$\labelp$};\node at (-.07,-0.4) {$\labeltwo$};
\node [draw,rounded corners,inner sep=2pt,black,fill=gray!20!white] at (-.55,.65) {$\scriptstyle x_1$};
\end{tikzpicture}
+
\begin{tikzpicture}[anchorbase,scale=1.4]
\draw[colp,<-] (.15,-.8) to (.15,.8);
\draw[colp,->] (-.15,-.8) to (-.15,.8);
	\draw[colm,<-] (.8,.8) to[out=-90, in=0,looseness=1.4] (0,.15) to[out = 180, in = -90,looseness=1.4] (-0.8,.8);
	\draw[colm,->] (.8,-.8) to[looseness=1.4,out=90, in=0] (0,-.15) to[looseness=1.4,out = 180, in = 90] (-0.8,-.8);
\draw (0,0) ellipse (.85 and .6);
\draw[->] (.85,0)
arc[start angle=-0.5, end angle=.5,radius=.9];
\node at (1.05,0) {$\labelp$};
\node at (.68,0) {$\labelp$};\node at (.5,0.04) {$\labeltwo$};
\node at (-.5,0) {$\labelp$};\node at (-.68,0.04) {$\labeltwo$};
\node at (0.06,0) {$\labelpm$};
\node at (-.07,0) {$\labeltwo$};
    \node at (-1,0) {$\labelp$};
\node at (0,-0.7) {$\labelp$};
\node at (0,0.7) {$\labelp$};
\node at (.38,0.4) {$\labeltwo$};
\node at (-.38,0.4) {$\labeltwo$};
\node at (.38,-0.4) {$\labeltwo$};
\node at (-.38,-0.4) {$\labeltwo$};
\node at (.06,0.4) {$\labelp$};\node at (-.07,0.4) {$\labeltwo$};
\node at (.06,-0.4) {$\labelp$};\node at (-.07,-0.4) {$\labeltwo$};
\node [draw,rounded corners,inner sep=2pt,black,fill=gray!20!white] at (.55,-.65) {$\scriptstyle x_1$};
\end{tikzpicture}
=\Theta\Bigg(\begin{tikzpicture}[anchorbase,scale=1.4]
	\draw[semithick] (0.25,.3) to (-0.25,-.3);
	\draw[semithick] (0.25,-.3) to (-0.25,.3);
 \closeddot{-0.12,0.145};
\end{tikzpicture}
-
\begin{tikzpicture}[anchorbase,scale=1.4]
	\draw[semithick] (0.25,.3) to (-0.25,-.3);
	\draw[semithick] (0.25,-.3) to (-0.25,.3);
     \closeddot{.12,-0.145};
\end{tikzpicture}\Bigg).
\end{align*}
The difficult step here is the equality $(*)$. Subtracting \cref{puppies2} from \cref{puppies1}, one obtains a result similar to the above except with $\frac{1}{2}(\alpha_{-1}-\alpha_1)$ in certain regions labeled $\labeltwo$. But
$x_1 = \frac{1}{2}(\alpha_{-1}-\alpha_1)$. Since $x_1$ is invariant under $s_2$ (unlike $\alpha_{\pm 1}$), it  slides across the string of color $2$ to the desired region with label $\O$.
\item
For the Reidemeister I relation from \cref{rels3}, we give an easy degree argument. The curl is of degree $-2$, so 
it maps to a graded bimodule homomorphism 
$B\otimes_{\Sym}B
\rightarrow \Sym$ of degree $-2$.
This homomorphism must be 0
because by adjunction we have that
$$\Hom_{\Sym\dash \Sym}(B\otimes_{\Sym} B,\Sym)
\cong
\Hom_{R^{\O+}\dash R^{\O+}}(B_{\O+,\O-}\otimes_{R^{\O-}} B_{\O-,\O+},R^{\O+})\cong
\End_{R^{\O-}\dash R^{\O+}}(B_{\O-,\O+}).$$
Since $B_{\O-,\O+}$ is isomorphic to a shift of $R^{\O}$ as a graded $(R^{\O-}, R^{\O+})$-bimodule, they have the same endomorphism ring. As a bimodule $R^{\O}$ is generated by $1$ (a simple argument on linear generators) and is non-negatively graded, so its endomorphism ring is zero in negative degrees.
\item
For the Reidemeister II relation from \cref{rels1}, we must show that 
the degree $-4$ endomorphism
\begin{align}\label{kidding1}
\Theta\left(\ 
\begin{tikzpicture}[anchorbase]
\draw[semithick] (0.28,0) to[out=90,in=-90] (-0.28,.6);
\draw[semithick] (-0.28,0) to[out=90,in=-90] (0.28,.6);
\draw[semithick] (0.28,-.6) to[out=90,in=-90] (-0.28,0);
\draw[semithick] (-0.28,-.6) to[out=90,in=-90] (0.28,0);
\end{tikzpicture}\ \right)
&=
\begin{tikzpicture}[anchorbase,scale=1.15]
\draw[colp,-] (.3,0) to (.3,1.6);
\draw[colp,->] (-.3,0) to (-.3,1.6);
\draw[colm,<-] (.8,1.6) to[out=-90, in=0] (0,1) to[out = 180, in = -90] (-0.8,1.6);
\draw[colm,-] (.8,0) to[out=90, in=0] (0,.6) to[out = 180, in = 90] (-0.8,0);
\draw (0,0.8) ellipse (.85 and .6);
\draw[->] (.85,0.8) arc[start angle=-0.5, end angle=.5,radius=.9];
\node at (1.05,0) {$\labelp$};
\node at (.68,0.8) {$\labelp$};
\node at (.5,0.8) {$\labeltwo$};
\node at (-.5,0.8) {$\labelp$};
\node at (-.68,.8) {$\labeltwo$};
\node at (0.09,0.8) {$\labelpm$};
\node at (-.09,0.8) {$\labeltwo$};    
\node at (0,0) {$\labelp$};
\node at (0,1.5) {$\labelp$};
\node at (.38,1.2) {$\labeltwo$};
\node at (-.38,1.2) {$\labeltwo$};
\node at (.38,.4) {$\labeltwo$};
\node at (-.38,0.4) {$\labeltwo$};
\node at (.09,1.2) {$\labelp$};
\node at (-.09,1.2) {$\labeltwo$};
\node at (.09,0.4) {$\labelp$};
\node at (-.09,0.4) {$\labeltwo$};
\draw[colp,<-] (.3,-1.6) to (.3,0);
\draw[colp,-] (-.3,-1.6) to (-.3,0);
\draw[colm,<-] (.8,0) to[out=-90, in=0] (0,-.6) to[out = 180, in = -90] (-0.8,0);
\draw[colm,->] (.8,-1.6) to[out=90, in=0] (0,-1) to[out = 180, in = 90] (-0.8,-1.6);
\draw (0,-.8) ellipse (.85 and .6);
\draw[->] (.85,-.8) arc[start angle=-.5, end angle=.5,radius=.9];
\node at (.68,-.8) {$\labelp$};\node at (.5,-0.8) {$\labeltwo$};
\node at (-.5,-.8) {$\labelp$};\node at (-.68,-0.8) {$\labeltwo$};
\node at (0.09,-0.8) {$\labelpm$};\node at (-.09,-0.8) {$\labeltwo$};    
\node at (-1,0) {$\labelp$};
\node at (0,-1.5) {$\labelp$};
\node at (.38,-0.4) {$\labeltwo$};
\node at (-.38,-0.4) {$\labeltwo$};
\node at (.38,-1.2) {$\labeltwo$};
\node at (-.38,-1.2) {$\labeltwo$};
\node at (.09,-0.4) {$\labelp$};\node at (-.09,-0.4) {$\labeltwo$};
\node at (.09,-1.2) {$\labelp$};\node at (-.09,-1.2) {$\labeltwo$};
\end{tikzpicture}
\end{align}
 of the graded $(\Sym,\Sym)$-bimodule
$B\otimes_{\Sym} B$ is 0.
This follows immediately from \cref{handy3}, 
which shows that
the lowest degree of any non-zero endomorphism of this bimodule is $-2$.
We regard this as a hard degree argument since it is applying some heavy machinery. Here is a more elementary proof, which also
prepares for the proof of Reidemeister III in the next point.
Identifying $B$
simply with $q^{l-1} s_0 R^O$ in the obvious way,
$B\otimes_{\Sym} B$
is generated by
the vectors $x_1^r \otimes x_1^s$
for $r,s \geq 0$, so it suffices to show that 
the endomorphism \cref{kidding1} is 0 on $x_1^r \otimes x_1^s$ for all $r,s \geq 0$. 
This follows by easy degree considerations in the case $r=s=0$.
To deduce the vanishing in general,
note that the image of
$x_1^r \otimes x_1^s$
under \cref{kidding1} is the same as the image of $1\otimes 1$ under the following:
\begin{align}\label{kidding2}
\Theta\left(\ 
(-1)^s
\begin{tikzpicture}[anchorbase]
	\draw[semithick] (0.28,0) to[out=90,in=-90] (-0.28,.6);
	\draw[semithick] (-0.28,0) to[out=90,in=-90] (0.28,.6);
	\draw[semithick] (0.28,-.6) to[out=90,in=-90] (-0.28,0);
	\draw[semithick] (-0.28,-.6) to[out=90,in=-90] (0.28,0);
 \closeddot{-.24,-.46};\node at (-.45,-.45) {$\scriptstyle r$};
 \closeddot{.24,-.46};\node at (.45,-.45) {$\scriptstyle s$};
\end{tikzpicture}\ \right)&=
\begin{tikzpicture}[anchorbase,scale=1.25]
\draw[colp,-] (.3,0) to (.3,1.6);
\draw[colp,->] (-.3,0) to (-.3,1.6);
	\draw[colm,<-] (.8,1.6) to[out=-90, in=0] (0,1) to[out = 180, in = -90] (-0.8,1.6);
	\draw[colm,-] (.8,0) to[out=90, in=0] (0,.6) to[out = 180, in = 90] (-0.8,0);
\draw (0,0.8) ellipse (.85 and .6);
\draw[->] (.85,0.8)
arc[start angle=-.5, end angle=.5,radius=.9];
\node at (1.05,0) {$\labelp$};
\node at (.68,0.8) {$\labelp$};\node at (.5,0.8) {$\labeltwo$};
\node at (-.5,0.8) {$\labelp$};\node at (-.68,.8) {$\labeltwo$};
\node at (0.09,0.8) {$\labelpm$};\node at (-.09,0.8) {$\labeltwo$};    
\node at (0,0) {$\labelp$};
\node at (0,1.5) {$\labelp$};
\node at (.38,1.2) {$\labeltwo$};
\node at (-.38,1.2) {$\labeltwo$};
\node at (.38,.4) {$\labeltwo$};
\node at (-.38,0.4) {$\labeltwo$};
\node at (.09,1.2) {$\labelp$};\node at (-.09,1.2) {$\labeltwo$};
\node at (.09,0.4) {$\labelp$};\node at (-.09,0.4) {$\labeltwo$};
\draw[colp,<-] (.3,-1.7) to (.3,0);
\draw[colp,-] (-.3,-1.7) to (-.3,0);
	\draw[colm,<-] (.8,0) to[out=-90, in=0] (0,-.6) to[out = 180, in = -90] (-0.8,0);
	\draw[colm,->] (.8,-1.7) to (.8,-1.6) to[out=90, in=0] (0,-1) to[out = 180, in = 90] (-0.8,-1.6) to (-0.8,-1.7);
\draw (0,-.8) ellipse (.85 and .6);
\draw[->] (.85,-.8)
arc[start angle=-.5, end angle=.5,radius=.9];
\node at (.68,-.8) {$\labelp$};\node at (.5,-0.8) {$\labeltwo$};
\node at (-.5,-.8) {$\labelp$};\node at (-.68,-0.8) {$\labeltwo$};
\node at (0.09,-0.8) {$\labelpm$};\node at (-.09,-0.8) {$\labeltwo$};    
\node at (-1,0) {$\labelp$};
\node at (0,-1.5) {$\labelp$};
\node at (.38,-0.4) {$\labeltwo$};
\node at (-.38,-0.4) {$\labeltwo$};
\node at (.38,-1.2) {$\labeltwo$};
\node at (-.38,-1.2) {$\labeltwo$};
\node at (.09,-0.4) {$\labelp$};\node at (-.09,-0.4) {$\labeltwo$};
\node at (.09,-1.2) {$\labelp$};\node at (-.09,-1.2) {$\labeltwo$};
\node [draw,rounded corners,inner sep=1.5pt,black,fill=gray!20!white] at (-.55,-1.5) {$\scriptstyle x_1^r$};
\node [draw,rounded corners,inner sep=1.5pt,black,fill=gray!20!white] at (.55,-1.5) {$\scriptstyle x_1^s$};\end{tikzpicture}\:.
\end{align}
We have checked that the relations 
\cref{rels2,rels3,rels4} are satisfied in
$\END_{\ESBim}(\Sym)$, so we already have in our hands a monoidal functor from the monoidal category defined in the same way as $\cNB_t$ but omitting the relations \cref{rels1} to $\END_{\ESBim}(\Sym)$.
From this and \cref{puppy1}, we deduce that
the endomorphism
\cref{kidding2} is equal to
\begin{align}\label{kidding3}
\Theta\left(\ (-1)^s
\begin{tikzpicture}[anchorbase]
	\draw[semithick] (0.28,0) to[out=90,in=-90] (-0.28,.6);
	\draw[semithick] (-0.28,0) to[out=90,in=-90] (0.28,.6);
	\draw[semithick] (0.28,-.6) to[out=90,in=-90] (-0.28,0);
	\draw[semithick] (-0.28,-.6) to[out=90,in=-90] (0.28,0);
 \closeddot{-.24,.46};\node at (-.45,.45) {$\scriptstyle r$};
 \closeddot{.24,.46};\node at (.45,.45) {$\scriptstyle s$};
\end{tikzpicture}\ \right)&=
\begin{tikzpicture}[anchorbase,scale=1.25]
\draw[colp,-] (.3,0) to (.3,1.7);
\draw[colp,->] (-.3,0) to (-.3,1.7);
	\draw[colm,<-] (.8,1.7) to (.8,1.6) to[out=-90, in=0] (0,1) to[out = 180, in = -90] (-0.8,1.6) to (-.8,1.7);
	\draw[colm,-] (.8,0) to[out=90, in=0] (0,.6) to[out = 180, in = 90] (-0.8,0);
\draw (0,0.8) ellipse (.85 and .6);
\draw[->] (.85,0.8)
arc[start angle=-.5, end angle=.5,radius=.9];
\node at (1.05,0) {$\labelp$};
\node at (.68,0.8) {$\labelp$};\node at (.5,0.8) {$\labeltwo$};
\node at (-.5,0.8) {$\labelp$};\node at (-.68,.8) {$\labeltwo$};
\node at (0.09,0.8) {$\labelpm$};\node at (-.09,0.8) {$\labeltwo$};    
\node at (0,0) {$\labelp$};
\node at (0,1.5) {$\labelp$};
\node at (.38,1.2) {$\labeltwo$};
\node at (-.38,1.2) {$\labeltwo$};
\node at (.38,.4) {$\labeltwo$};
\node at (-.38,0.4) {$\labeltwo$};
\node at (.09,1.2) {$\labelp$};\node at (-.09,1.2) {$\labeltwo$};
\node at (.09,0.4) {$\labelp$};\node at (-.09,0.4) {$\labeltwo$};
\draw[colp,<-] (.3,-1.6) to (.3,0);
\draw[colp,-] (-.3,-1.6) to (-.3,0);
	\draw[colm,<-] (.8,0) to[out=-90, in=0] (0,-.6) to[out = 180, in = -90] (-0.8,0);
	\draw[colm,->] (.8,-1.6) to[out=90, in=0] (0,-1) to[out = 180, in = 90] (-0.8,-1.6);
\draw (0,-.8) ellipse (.85 and .6);
\draw[->] (.85,-.8)
arc[start angle=-.5, end angle=.5,radius=.9];
\node at (.68,-.8) {$\labelp$};\node at (.5,-0.8) {$\labeltwo$};
\node at (-.5,-.8) {$\labelp$};\node at (-.68,-0.8) {$\labeltwo$};
\node at (0.09,-0.8) {$\labelpm$};\node at (-.09,-0.8) {$\labeltwo$};    
\node at (-1,0) {$\labelp$};
\node at (0,-1.5) {$\labelp$};
\node at (.38,-0.4) {$\labeltwo$};
\node at (-.38,-0.4) {$\labeltwo$};
\node at (.38,-1.2) {$\labeltwo$};
\node at (-.38,-1.2) {$\labeltwo$};
\node at (.09,-0.4) {$\labelp$};\node at (-.09,-0.4) {$\labeltwo$};
\node at (.09,-1.2) {$\labelp$};\node at (-.09,-1.2) {$\labeltwo$};
\node [draw,rounded corners,inner sep=1.5pt,black,fill=gray!20!white] at (-.55,1.5) {$\scriptstyle x_1^r$};
\node [draw,rounded corners,inner sep=1.5pt,black,fill=gray!20!white] at (.55,1.5) {$\scriptstyle x_1^s$};\end{tikzpicture}\:.
\end{align}
Finally, the endomorphism \cref{kidding3} maps $1 \otimes 1$ to 0 by the $r=s=0$ case treated earlier.
\item
For the Reidemeister III relation from \cref{rels3}, we must show that
the degree $-6$ endomorphism 
\begin{align*}
\Theta\left(\:\begin{tikzpicture}[anchorbase]
	\draw[semithick] (0.45,.6) to (-0.45,-.6);
	\draw[semithick] (0.45,-.6) to (-0.45,.6);
        \draw[semithick] (0,-.6) to[out=90,in=-90] (-.45,0);
        \draw[semithick] (-0.45,0) to[out=90,in=-90] (0,0.6);
\end{tikzpicture}-
\begin{tikzpicture}[anchorbase]
	\draw[semithick] (0.45,.6) to (-0.45,-.6);
	\draw[semithick] (0.45,-.6) to (-0.45,.6);
        \draw[semithick] (0,-.6) to[out=90,in=-90] (.45,0);
        \draw[semithick] (0.45,0) to[out=90,in=-90] (0,0.6);
\end{tikzpicture}\:\right)&=
\begin{tikzpicture}[anchorbase,scale=1.15]
\draw[colp,<-] (.3,-2.4) to (.3,2.4);
\draw[colp,->] (-.3,0.8) to (-.3,2.4);
	\draw[colm,<-] (.8,2.4) to[out=-90, in=0] (0,1.8) to[out = 180, in = -90] (-0.8,2.4);
	\draw[colm,-] (.8,.8) to[out=90, in=0] (0,1.4) to[out = 180, in = 90] (-0.8,0.8);
\draw (0,1.6) ellipse (.85 and .6);
\draw[->] (.85,1.6)
arc[start angle=-.5, end angle=.5,radius=.9];
\node at (1.05,0) {$\labelp$};
\node at (.68,1.6) {$\labelp$};\node at (.5,1.6) {$\labeltwo$};
\node at (-.5,1.6) {$\labelp$};\node at (-.68,1.6) {$\labeltwo$};
\node at (0.09,1.6) {$\labelpm$};\node at (-.09,1.6) {$\labeltwo$};    
\node at (0,0) {$\labelp$};
\node at (0,2.3) {$\labelp$};
\node at (.38,2) {$\labeltwo$};
\node at (-.38,2) {$\labeltwo$};
\node at (.38,1.2) {$\labeltwo$};
\node at (-.38,1.2) {$\labeltwo$};
\node at (.09,2) {$\labelp$};\node at (-.09,2) {$\labeltwo$};
\node at (.09,1.2) {$\labelp$};\node at (-.09,1.2) {$\labeltwo$};
\draw[colm,-] (-.8,-.8) to (-.8,.8);
\draw[colm,->] (-1.4,-2.4) to (-1.4,2.4);
	\draw[colp,-] (-.3,.8) to[out=-90, in=0] (-1.1,.2) to[out = 180, in = -90] (-1.9,.8) to (-1.9,2.4);
	\draw[colp,->] (-.3,-.8) to[out=90, in=0] (-1.1,-.2) to[out = 180, in = 90] (-1.9,-.8) to (-1.9,-2.4);
\draw (-1.1,0) ellipse (.85 and .6);
\draw[->] (-.25,0)
arc[start angle=-.5, end angle=.5,radius=.9];
\node at (-.42,0) {$\labelp$};\node at (-.6,0) {$\labeltwo$};
\node at (-1.6,0) {$\labelp$};\node at (-1.78,0) {$\labeltwo$};
\node at (-1.01,0) {$\labelpm$};
\node at (-1.19,0) {$\labeltwo$};
    \node at (-2.2,0) {$\labelp$};
\node at (-1.1,-1.6) {$\labelp$};
\node at (-1.1,1.6) {$\labelp$};
\node at (-.72,0.4) {$\labeltwo$};
\node at (-.72,0.4) {$\labeltwo$};
\node at (-.72,-0.4) {$\labeltwo$};
\node at (-.72,-0.4) {$\labeltwo$};
\node at (-1.01,0.4) {$\labelp$};\node at (-1.19,0.4) {$\labeltwo$};
\node at (-1.01,-0.4) {$\labelp$};\node at (-1.19,-0.4) {$\labeltwo$};
\draw[colm,->] (.8,-.8) to (.8,0);
\draw[colm,-] (.8,0) to (.8,.8);
\draw[colp,-] (-.3,-2.4) to (-.3,-.8);
	\draw[colm,-] (.8,-.8) to[out=-90, in=0] (0,-1.4) to[out = 180, in = -90] (-0.8,-.8);
	\draw[colm,->] (.8,-2.4) to[out=90, in=0] (0,-1.8) to[out = 180, in = 90] (-0.8,-2.4);
\draw (0,-1.6) ellipse (.85 and .6);
\draw[->] (.85,-1.6)
arc[start angle=-.5, end angle=.5,radius=.9];
\node at (.68,-1.6) {$\labelp$};\node at (.5,-1.6) {$\labeltwo$};
\node at (-.5,-1.6) {$\labelp$};\node at (-.68,-1.6) {$\labeltwo$};
\node at (0.09,-1.6) {$\labelpm$};\node at (-.09,-1.6) {$\labeltwo$};    
\node at (0,-2.3) {$\labelp$};
\node at (.38,-1.2) {$\labeltwo$};
\node at (-.38,-1.2) {$\labeltwo$};
\node at (.38,-2) {$\labeltwo$};
\node at (-.38,-2) {$\labeltwo$};
\node at (.09,-1.2) {$\labelp$};\node at (-.09,-1.2) {$\labeltwo$};
\node at (.09,-2) {$\labelp$};\node at (-.09,-2) {$\labeltwo$};
\end{tikzpicture}-
\begin{tikzpicture}[anchorbase,scale=1.15]
\draw[colm,->] (-.3,-2.4) to (-.3,2.4);
\draw[colm,-] (.3,0.8) to (.3,2.4);
	\draw[colp,->] (-.8,2.4) to[out=-90, in=180] (0,1.8) to[out = 0, in = -90] (0.8,2.4);
	\draw[colp,-] (-.8,.8) to[out=90, in=180] (0,1.4) to[out = 0, in = 90] (0.8,0.8);
\draw (0,1.6) ellipse (.85 and .6);
\draw[<-] (.85,1.6)
arc[start angle=.5, end angle=-.5,radius=.9];
\node at (-1.05,0) {$\labelp$};
\node at (-.5,1.6) {$\labelp$};\node at (-.68,1.6) {$\labeltwo$};
\node at (.68,1.6) {$\labelp$};\node at (.5,1.6) {$\labeltwo$};
\node at (0.09,1.6) {$\labelpm$};\node at (-.09,1.6) {$\labeltwo$};    
\node at (0,0) {$\labelp$};
\node at (0,2.3) {$\labelp$};
\node at (-.38,2) {$\labeltwo$};
\node at (.38,2) {$\labeltwo$};
\node at (-.38,1.2) {$\labeltwo$};
\node at (.38,1.2) {$\labeltwo$};
\node at (.09,2) {$\labelp$};\node at (-.09,2) {$\labeltwo$};
\node at (.09,1.2) {$\labelp$};\node at (-.09,1.2) {$\labeltwo$};
\draw[colp,->] (.8,-.8) to (.8,0);
\draw[colp,-] (.8,0) to (.8,.8);
\draw[colp,<-] (1.4,-2.4) to (1.4,2.4);
	\draw[colm,->] (.3,.8) to[out=-90, in=180] (1.1,.2) to[out =0, in = -90] (1.9,.8) to (1.9,2.4);
	\draw[colm,-] (.3,-.8) to[out=90, in=180] (1.1,-.2) to[out = 0, in = 90] (1.9,-.8) to (1.9,-2.4);
\draw (1.1,0) ellipse (.85 and .6);
\draw[<-] (1.95,0)
arc[start angle=.5, end angle=-.5,radius=.9];
\node at (.6,0) {$\labelp$};\node at (.42,0) {$\labeltwo$};
\node at (1.78,0) {$\labelp$};\node at (1.6,0) {$\labeltwo$};
\node at (1.19,0) {$\labelpm$};
\node at (1.01,0) {$\labeltwo$};
    \node at (2.2,0) {$\labelp$};
\node at (1.1,-1.6) {$\labelp$};
\node at (1.1,1.6) {$\labelp$};
\node at (1.47,0.4) {$\labeltwo$};
\node at (1.47,-0.4) {$\labeltwo$};
\node at (.72,0.4) {$\labeltwo$};
\node at (.72,-0.4) {$\labeltwo$};
\node at (1.19,0.4) {$\labelp$};\node at (1.01,0.4) {$\labeltwo$};
\node at (1.19,-0.4) {$\labelm$};\node at (1.01,-0.4) {$\labeltwo$};
\draw[colp,-] (-.8,-.8) to (-.8,0);
\draw[colp,-] (-.8,0) to (-.8,.8);
\draw[colm,<-] (.3,-2.4) to (.3,-.8);
	\draw[colp,-] (-.8,-.8) to[out=-90, in=180] (0,-1.4) to[out =0, in = -90] (0.8,-.8);
	\draw[colp,<-] (-.8,-2.4) to[out=90, in=180] (0,-1.8) to[out = 0, in = 90] (0.8,-2.4);
\draw (0,-1.6) ellipse (.85 and .6);
\draw[<-] (.85,-1.6)
arc[start angle=0.5, end angle=-.5,radius=.9];
\node at (-.5,-1.6) {$\labelp$};\node at (-.68,-1.6) {$\labeltwo$};
\node at (.68,-1.6) {$\labelp$};\node at (.5,-1.6) {$\labeltwo$};
\node at (0.09,-1.6) {$\labelpm$};\node at (-.09,-1.6) {$\labeltwo$};    
\node at (0,-2.3) {$\labelp$};
\node at (-.38,-1.2) {$\labeltwo$};
\node at (.38,-1.2) {$\labeltwo$};
\node at (-.38,-2) {$\labeltwo$};
\node at (.38,-2) {$\labeltwo$};
\node at (.09,-1.2) {$\labelp$};\node at (-.09,-1.2) {$\labeltwo$};
\node at (.09,-2) {$\labelp$};\node at (-.09,-2) {$\labeltwo$};
\end{tikzpicture}
\end{align*}
of the graded
$(\Sym,\Sym)$-bimodule
$B \otimes_{\Sym} B \otimes_{\Sym} B$ is 0.
Identifying each $B$ with $q^{l-1} s_0 R^O$,
this bimodule
is generated by
$x_1^r \otimes x_1^s \otimes x_1^t$ for $r,s,t\geq 0$, so it suffices to show that the above 
endomorphism sends each of these vectors to 0.
This follows by same strategy as was just used to prove
Reidemeister II, using \cref{puppy2} in place of \cref{puppy1}.
\end{itemize}

At this point, we have in our hands a $\Lambda^{[2]}$-linear  monoidal functor $\widehat{\Theta}:\cNB_t \rightarrow \END_{\ESBim}(\Sym)$.
Next, we show for any $f \in \Lambda$ 
that $\widehat{\Theta}(\Bubble{f})$ is the endomorphism 
$\bubble{\ev_l(f)}\ \labelp$
of $\Sym$ defined by multiplication by the image of $f$
under the evaluation homomorphism $\Lambda \twoheadrightarrow \Sym$. This implies the statement (1).
Since $\widehat{\Theta}$ is $\Lambda^{[2]}$-linear by construction, and $\Lambda$ is generated by $\Lambda^{[2]}$ together with the symmetric functions  $q_n$ for all $n > 0$, it suffices to show that
$\widehat{\Theta}\big(\:\Bubble{q_n}\:) = \bubble{\ev_l(q_n)}\  \labelp$.
By \cref{gammat}, $q_n$ corresponds (up to scalar) to a dotted bubble $\begin{tikzpicture}[anchorbase]
\draw[semithick] (0,0) circle (.2);
\closeddot{.2,0};
\node at (.35,0) {$\scriptstyle n$};
\end{tikzpicture}$ in the nilBrauer category, and  \cref{coffeeisgood} shows that $\widehat{\Theta}$ sends the dotted bubble to $q_n$ (up to the inverse scalar) in $\END_{\ESBim}(\Sym)$.

Now we complete the proof by checking 
the cyclotomic relations, so that $\widehat{\Theta}$ descends to the cyclotomic quotient $\CNB_l$.
Since $e_r(x_1,\dots,x_l) = 0$ for $r > l$, we have that
$\widehat{\Theta}\left(\ \Bubble{e_r}\ \right) = 0$ for $r > l$
by the statement (1) checked in the previous paragraph. 
To show that $\widehat{\Theta}$ maps
\cref{purplesquare} to 0, it suffices to show that the left action of $e_l(x_1,\dots,x_l) = x_1 \cdots x_l$ on $B$ differs from the right action by a sign. This holds because $B$ is a twisted bimodule.
\end{proof}

\begin{rem}
Some of the relation checks in the proof of 
\cref{okinawa} were performed directly using diagrammatic transformations. 
However, we gave indirect arguments for the Reidemeister I, II and III relations. 
In fact, it is not hard to give a direct proof of the Reidemeister I relation, using \cref{evaluations,easyRII,hardRII}. To prove the Reidemeister II relation diagrammatically, one needs additional relations not recorded in \cref{ssec3-3}, namely, the idempotent decompositions which categorify the type A MOY relations. 
Using these, one can
simplify the inside of the middle circle in \cref{kidding3}. However, we have {\em not} been able to work out a purely diagrammatical proof of Reidemeister III, and expect that such an argument would require the use of idempotent decompositions specific to type D.
\end{rem}

\subsection{Characters of indecomposables}
Restricting from the graded bicategory $\ESBim$ to 
the graded monoidal category $\END_{\ESBim}(\Sym)$, and recalling the definition \cref{thisisbetter},
the Soergel-Williamson character map (modified slightly to fit the extended setup)
gives an isomorphism of $\Z[q,q^{-1}]$-algebras
\begin{equation}
\ch:K_0(\END_{\ESBim}(\Sym)) \stackrel{\sim}{\rightarrow} \V.
\end{equation}
Hence, by 
\cref{asleftmodule}, $K_0(\END_{\ESBim}(\Sym))$ is isomorphic as a module to
the integral form
$\V(l)$ for the $(l+1)$-dimensional 
$U_q(\sl_2)$-module, 
and as an algebra to the quotient $\Ui / \I_l$ from \cref{oldpeculiar}.
By \cref{class}
and the classification of double cosets in \cref{pure1,pure2}, we know that 
a full set of self-dual
indecomposable objects
in $\END_{\ESBim}(\Sym)$
are given by the $(\Sym,\Sym)$-bimodules
$B_{(\O+)d_n(\O+)}$ for $0 \leq n \leq l$ with $n$ even
and $s_0 \Sym \otimes_{\Sym} B_{(\O-)s_0 d_n(\O+)}$
for $1 \leq n \leq l$ with $n$ odd.
The next fundamental theorem is an application of the nil-Brauer theory from \cite{BWWiquantum}
to deduce character formulae for these special
singular Soergel bimodules.

From now on, we use $\e_n$ to denote the 
image of the idempotent \cref{endef} 
in the quotient $\CNB_l$ of $\cNB_t$.
So $\e_n$ is a homogeneous idempotent (possibly 0)
in $\End_{\CNB_l}\big(B^{\star n}\big)$.
Let
\begin{align}\label{taipei}
\f_n &:= \Theta(\e_n),
\end{align}
which is a
homogeneous idempotent (possibly 0)
in $\End_{\Sym\dash \Sym}\left(B^{\otimes n}\right)$
defining a summand $\left[B^{\otimes n}, \f_n\right]$ of this bimodule.
We set
\begin{equation}\label{newBB}
B^{[n]} := q^{-\binom{n}{2}} \left[B^{\otimes n}, \f_n\right].
\end{equation}
The similarity of this notation with \cref{poorpolly} is deliberate.

\begin{theo}\label{bigtheorem}
In $\END_{\ESBim}(\Sym)$, we have that
\begin{equation}\label{step}
B^{\otimes n}
\cong 
\bigoplus_{i=0}^{\lfloor\frac{n}{2}\rfloor}
[n-2i]_q^!
\left(\sum_{\lambda \in \Par_{\not\equiv t}(i\times (n-2i))}
[\lambda_1+1]_q^2\cdots [\lambda_i+1]_q^2\right) B^{[n-2i]}.
\end{equation}
Hence, for $0 \leq n \leq l$, we have that
\begin{align}\label{droving}
B^{[n]}
\cong
\begin{cases}
B_{(\O+)d_n(\O+)}&\text{if $n$ is even}\\
s_0 R^{\O-} \otimes_{R^{\O-}} B_{(\O-)s_0 d_n(\O+)}&\text{if $n$ is odd}.
\end{cases}
\end{align}
The characters of these bimodules are equal to the corresponding 
Kazhdan-Lusztig basis elements:
\begin{align}\label{driving}
\ch\big( B^{[n]}\big)
=\begin{cases}
b_{(\O+)d_n (\O+)}&\text{if $n$ is even}\\
s_0 b_{(\O-)s_0 d_n(\O+)}&\text{if $n$ is odd.}
\end{cases}
\end{align}
In particular, these characters are independent of the characteristic of the ground field.
\end{theo}

\begin{proof}
The identity \cref{step} follows by applying $\Theta$
(actually, its extension to the graded Karoubi envelope)
to the decomposition \cref{theworld}
and using the definitions \cref{taipei,,newBB}.
Since $\ch \big(B^{\otimes n}\big) = b^{\m n}$
and the
(triangular) transition matrix in the decomposition \cref{step} is the same as in \cref{cruisier}, we deduce that
$B^{[n]}$
is an extended singular Soergel bimodule with character
$b_{(\O+)d_n(\O+)}$ if $n$ is even or
$s_0 b_{(\O-)s_0 d_n(\O+)}$ if $n$ is odd, proving \cref{driving}.
By the homomorphism formula (\cref{homformula}) and
\cref{asleftmodule}, 
$$
\End_{A\dash A}(B^{[n]})
\cong A^{\oplus \overline{(b^{(n)}\eta_l, b^{(n)} \eta_l)_l}}.
$$
By \cref{bilform,characterize}, the graded dimension of this vector space is an element of $1+q\N\llbracket q\rrbracket$,
hence, its degree 0 component is 1-dimensional.
This implies that $B^{[n]}$
is indecomposable. Consequently, it must be isomorphic to
the indecomposable extended singular Soergel bimodule
$B_{(\O+)d_n(\O+)}$ or $s_0 R^{\O-} \otimes_{R^{\O-}}
B_{(\O-)s_0 d_n(\O+)}$ as in \cref{droving},
in view of the definition of these bimodules  from \cref{class}.
\end{proof}

\begin{cor}\label{handy1}
For $0 \leq n \leq \lfloor\frac{l}{2}\rfloor$, 
$\Hom_{\Sym\dash\Sym}\big(\Sym, B^{[2n]}
\big)$
is a free graded $\Sym$-module of graded rank
$q^{n(l+t-1)}\qbinom{(l-t)/2}{n}_{q^2}$.
\end{cor}

\begin{proof}
This follows from 
the homomorphism formula (\cref{homformula}),
using the identification of the character of  $B^{[2n]}$
obtained in \cref{bigtheorem} and
 the combinatorics from \cref{formhecke}.
\end{proof}

\subsection{Identification of the Grothendieck ring of
\texorpdfstring{$\CNB_l$}{}}
\cref{bigtheorem} implies that the homogeneous
idempotents $\f_n\:(0 \leq n \leq l)$ 
are primitive, and moreover that
any primitive homogeneous idempotent in $\END_{\ESBim}(\Sym)$ is equivalent to $\f_n$
for a unique $0 \leq n \leq l$.
Now we return to the cyclotomic nil-Brauer category
$\CNB_l$. 

\begin{theo}\label{windows}
For $0 \leq n \leq l$, $\e_n$ is a 
primitive homogeneous idempotent in $\CNB_l$,
and any 
primitive homogeneous idempotent in $\CNB_l$
is equivalent to $\e_n$ for a unique $0 \leq n \leq l$. 
\end{theo}

\begin{proof}
By \cref{aboutidempotents}, we know already 
that the idempotents $\e_n\:(n \geq 0)$
give a complete set of primitive homogeneous idempotents in $\cNB_t$.
It follows 
from general principles that the 
{\em non-zero} ones amongst their images give a full set of primitive homogeneous idempotents in $\CNB_l$.
For $n > l$, we know that $\bar\e_n=0$ by \cref{lemmabelow}.
To complete the proof of the theorem, 
it remains to show that
$\e_n$ is non-zero in $\CNB_l$ for $0 \leq n \leq l$.
This follows from \cref{bigtheorem} and the definitions \cref{taipei}.
\end{proof}

\begin{cor}\label{doors}
The monoidal functor $\Theta:\CNB_l \rightarrow \END_{\ESBim}(\Sym)$
induces a $\Z[q,q^{-1}]$-linear algebra isomorphism
$K_0(\gKar(\CNB_l))\stackrel{\sim}{\rightarrow}
K_0(\END_{\ESBim}(\Sym))$
taking $\big[B^{(n)}\big]$ to $\big[B^{[n]}\big]$ for  $0 \leq n \leq l$.
\end{cor}

\subsection{Diagrams for primitive idempotents} 
Recall the bimodule homomorphisms $\v_n$ and $\r_n$
from \cref{bplus,rn}.
The following theorem gives an explicit diagrammatic description of the homogeneous 
primitive idempotents $\f_n\:(0 \leq n \leq l)$:

\begin{theo}\label{bingley}
$\displaystyle
\f_n = \begin{cases}
\r_n \circ \v_n&\text{if $0 \leq n \leq l$}\\
0&\text{if $n > l$.}
\end{cases}
$
\end{theo}

\begin{proof}
By \cref{taipei} and \cref{endef}, we have that
$\f_n =
\r_n \circ 
\Theta\left(\begin{tikzpicture}[anchorbase,scale=1.2]
  \cross{(0,0)};
\draw[ultra thick] (0,.2) to (0,-.2);
\node at (0,-.3) {$\stringlabel{n}$};
\end{tikzpicture}\right).$
Thus, $\f_n$ is the composition of a degree $2\binom{n}{2}$ endomorphism after a
degree $-2\binom{n}{2}$ endomorphism of $B^{\otimes n}$.
It follows that $\f_n = 0$ for $n > l$ since,
by \cref{handy3},
any degree $-2 \binom{n}{2}$ endomorphism of 
$B^{\otimes n}$ is 0 when $n> l$. It is also obvious that
$\f_0 = \r_0 \circ \v_0$.
Finally, suppose that $1 \leq n \leq l$.
Since $\f_n \neq 0$, the degree $-2\binom{n}{2}$ endomorphism 
$\Theta\left(\begin{tikzpicture}[anchorbase,scale=1.2]
  \cross{(0,0)};
\draw[ultra thick] (0,.2) to (0,-.2);
\node at (0,-.3) {$\stringlabel{n}$};
\end{tikzpicture}\right)$ is non-zero.
Also $\v_n$ is of
degree $-2\binom{n}{2}$. This is the lowest non-zero degree of 
$\End_{\Sym\dash\Sym}\big(B^{\otimes n}\big)$, and it is  1-dimensional in this degree 
by \cref{handy3} again. Hence, $\v_n=c
\Theta\left(\begin{tikzpicture}[anchorbase,scale=1.2]
  \cross{(0,0)};
\draw[ultra thick] (0,.2) to (0,-.2);
\node at (0,-.3) {$\stringlabel{n}$};
\end{tikzpicture}\right)$ and
$\r_n \circ \v_n = c \f_n$ for some $c \in \kk$.
To show that $c=1$, we use that $\f_n$ is an idempotent so its only non-zero eigenvalue is 1. Thus, it suffices to observe that $\r_n \circ \v_n$ has an eigenvector of eigenvalue 1:
the vector
$\r_n(1\otimes\cdots\otimes 1)$ is fixed by $\r_n \circ \v_n$ by
\cref{mincies}.
\end{proof}

\begin{rem} \label{rmk:factoredidempotentfn}
Here is some additional commentary on the theory underlying \Cref{bingley}. Recall that $\v_n = \check\v_n \circ \hat\v_n$. Thus the idempotent $\f_n$ factors as a composition of $(\r_n \circ \check\v_n)$ with $\hat\v_n$. When $n$ is even, the object being factored through is the singular Bott--Samelson bimodule $M$ associated to the sequence of parabolic subgroups
$[\O + \supset \O \hat n + \subset \O \hat n \pm \supset \O \hat n + \subset \O +]$.
This sequence is a reduced expression for the double coset $(\O+) d_n (\O+)$ in the sense of \cite{EK} (see also \cite[Ch~1.3]{WillThesis}). Consequently, by \cite[Proof of Theorem 5.4.2]{WillThesis}, the bimodule $M$ has $q^{\binom{n}{2}}B_{(\O+)d_n(\O+)}$ as a direct summand. Thus it is reasonable to expect that the idempotent $\f_n$, whose image is isomorphic to $q^{\binom{n}{2}}B_{(\O+)d_n(\O+)}$ by \Cref{bigtheorem}, should factor through $M$. Indeed, the opposite composition $\hat\v_n \circ (\r_n \circ \check\v_n)$ should yield the idempotent endomorphism of $M$ whose image is $q^{\binom{n}{2}}B_{(\O+)d_n(\O+)}$.
Similar statements can be made when $n$ is odd.
\end{rem}